%% file: __G_Gamma_GmodGamma.tex
\documentclass[12pt,oneside]{book}
\usepackage{geometry}                		
\usepackage{graphicx}									
\usepackage{amssymb}\vspace{0.25cm}
\usepackage{mathtools}
\usepackage{setspace}
\usepackage{tikz-cd}
\usepackage[mathscr]{euscript}
\newcommand{\abs}[1]{\left\vert#1\right\vert}
\usepackage[normalem]{ulem}
\usepackage{manfnt}
\usepackage{pgfplots}
\usepackage{pifont}
\usepackage[safe]{tipa}
\usepackage{marvosym}
\usepackage{scalerel,stackengine}
\usepackage{titlesec}
\usepackage{makeidx}
\usepackage{centernot}
\usepackage{xspace}
\usepackage[noauto]{chappg}
\usepackage[OT2,T1]{fontenc}
\usepackage{hyperref}
\usepackage{amsmath}
\usepackage[amsmath,thmmarks]{ntheorem}
\usepackage{tabularx}
\usepackage{yfonts}
\newcolumntype{Y}{>{\raggedright\arraybackslash}X}

\newcommand\mA{%
$A$\xspace
}
\newcommand\mB{%
$B$\xspace
}
\newcommand\mC{%
$C$\xspace
}

\newcommand\mE{%
$E$\xspace
}
\newcommand\mF{%
$F$\xspace
}
\newcommand\mG{%
$G$\xspace
}
\newcommand\mH{%
$H$\xspace
}

\newcommand\mK{%
$K$\xspace
}
\newcommand\mL{%
$L$\xspace
}
\newcommand\mM{%
$M$\xspace
}
\newcommand\mN{%
$N$\xspace
}

\newcommand\mP{%
$P$\xspace
}
\newcommand\mQ{%
$Q$\xspace
}
\newcommand\mR{%
$R$\xspace
}
\newcommand\mS{%
$S$\xspace
}
\newcommand\mT{%
$T$\xspace
}
\newcommand\mU{%
$U$\xspace
}
\newcommand\mV{%
$V$\xspace
}
\newcommand\mW{%
$W$\xspace
}
\newcommand\mX{%
$X$\xspace
}
\newcommand\mY{%
$Y$\xspace
}

\newcommand\Gsym{%
\abs{-\hspace{-.12cm}\circ \hspace{-.12cm}-}^2
}

\newcommand\fa{%
\mathfrak{a}
}
\newcommand\fF{%
\mathfrak{F}
}
\newcommand\fG{%
\mathfrak{G}
}
\newcommand\fg{%
\mathfrak{g}
}
\newcommand\fh{%
\mathfrak{h}
}

\newcommand\fk{%
\mathfrak{k}
}
\newcommand\fl{%
\mathfrak{l}
}
\newcommand\fn{%
\mathfrak{n}
}
\newcommand\fP{%
\mathfrak{P}
}
\newcommand\fp{%
\mathfrak{p}
}

\newcommand\fs{%
\mathfrak{s}
}

\newcommand\ggS{%
\textgoth{S}
}
\newcommand\gT{%
\textgoth{T}
}

\newcommand\ft{%
\mathfrak{t}
}
\newcommand\fU{%
\mathfrak{U}
}
\newcommand\fu{%
\mathfrak{u}
}

\newcommand\fz{%
\mathfrak{z}
}

\newcommand\gche{%
\Upsilon
}

\newcommand\A{%
\mathbb{A}
}
\newcommand\Cx{%
\mathbb{C}
}
\newcommand\E{%
\mathbb{E}
}

\newcommand\I{%
\mathbb{I}
}

\newcommand\N{%
\mathbb{N}
}

\newcommand\Q{%
\mathbb{Q}
}
\newcommand\R{%
\mathbb{R}
}
\newcommand\T{%
\mathbb{T}
}

\newcommand\Z{%
\mathbb{Z}
}

\newcommand\bS{%
\textbf{S}\xspace
}
\newcommand\bT{%
\textbf{T}\xspace
}

\newcommand\Aut{%
\text{Aut}\hspace{0.05 cm}
}

\newcommand\Hom{%
\text{Hom}
}
\newcommand\MOD{%
\text{MOD}\xspace
}
\newcommand\sA{%
\mathcal{A}
}

\newcommand\sC{%
\mathcal{C}
}
\newcommand\sE{%
\mathcal{E}
}

\newcommand\sH{%
\mathcal{H}
}
\newcommand\sL{%
\mathcal{L}
}
\newcommand\sO{%
\mathcal{O}
}
\newcommand\sP{%
\mathcal{P}
}
\newcommand\sR{%
\mathcal{R}
}
\newcommand\sS{%
\mathcal{S}
}

\newcommand\sV{%
\mathcal{V}
}
\newcommand\sW{%
\mathcal{W}
}

\newcommand\reg{%
\text{reg}
}
\newcommand\sreg{%
\text{{\hspace{0.04cm}}reg}
}

\newcommand\Int{%
\text{Int}
}

\newcommand\Inv{%
\text{Inv}\hspace{0.05cm}
}
\newcommand\CoInv{%
\text{CoInv}\hspace{0.05cm}
}
\newcommand\Com{%
\text{Com}\hspace{0.03cm}
}
\newcommand\CON{%
\text{CON}\hspace{0.03cm}
}
\newcommand\Ker{%
\text{Ker}\hspace{0.03cm}
}
\newcommand\Ind{%
\text{Ind}
}
\newcommand\Inf{%
\text{Inf}
}
\newcommand\Ad{%
\text{Ad}
}

\newcommand\res{%
\text{res}
}
\newcommand\Res{%
\text{Res}
}

\newcommand\id{%
\text{id}
}

\newcommand\td{%
\text{d}
}
\newcommand\ad{%
\text{ad}
}
\newcommand\tdet{%
\text{det}\hspace{0.05cm}
}
\newcommand\tB{%
\text{B}
}
\newcommand\tc{%
\text{c}
}
\newcommand\tC{%
\text{C}
}
\newcommand\tE{%
\text{E}
}
\newcommand\tI{%
\text{I}
}
\newcommand\tJ{%
\text{J}
}
\newcommand\tL{%
\text{L}
}

\newcommand\tss{%
\text{ss}
}
\newcommand\tT{%
\text{T}
}

\newcommand\tZ{%
\text{Z}
}

\newcommand\GL{%
\textbf{GL}
}
\newcommand\SL{%
\textbf{SL}
}
\newcommand\SO{%
\textbf{SO}
}
\newcommand\SU{%
\textbf{SU}
}
\newcommand\ssx{%
\text{ss}
}

\newcommand\UN{%
\text{UN}
}

\newcommand\Tee{%
\mathsf{T}
}

\newcommand\wt{%
\text{wt}
}
\newcommand\rt{%
\text{rt}
}
\newcommand\spanx{%
\text{span}
}
\newcommand\Sphere{%
\textbf{S}
}
\newcommand\Torus{%
\text{$T$}
}
\newcommand\vspx{%
\vspace{0.2cm}
}

\newcommand\xdot{%
\dot{x}
}

\newcommand\ra{%
\rightarrow
}
\newcommand\lra{%
\longrightarrow
}

\newcommand\la{%
\leftarrow
}

\newcommand\ds{%
\displaystyle
}

\newcommand\End{%
\text{End}\hspace{0.05cm}
}
\newcommand\fin{%
\text{fin}\hspace{0.05cm}
}

\newcommand\modx{%
\text{mod}\hspace{0.05cm}
}
\newcommand\spt{%
\text{spt}
}
\newcommand\vol{%
\text{vol}
}
\newcommand\rank{%
\text{rank}
}

\newcommand\tr{%
\text{tr}\hspace{0.05cm}
}

\newcommand\un[1]{%
\underline{#1}\xspace
}
\newcommand\uun[1]{%
\underline{\underline{#1}}\xspace
}

\newcommand\restr[2]{%
{#1}|{#2}
}

\newcommand\hsx{%
\hspace{0.05cm}
}
\newcommand\hsy{%
\hspace{0.03cm}
}
\newcommand\hsz{%
\hspace{0.1cm}
}

\newcommand{\norm}[1]{\left\lVert #1 \right\rVert}

\newcommand{\chisub}[1]{\chi_{_{#1}}}

\newcommand\chis[1]{%
\chi_{_{
\scaleto{
#1
}
{4pt}
}}
}

\newcommand\chiss[1]{%
\chi_{_{
\scaleto{
#1
}
{6pt}
}}
}

\newcommand\stickfigure{%
\raisebox{-.0cm}{{\rule{.5pt}{1ex}}}
\hspace{-.25cm}{
\raisebox{0.15cm}{
\text{$\circ$}
}
}
}

\newcommand\ov[1]{%
\overline{#1}
}

\newcommand\ovha[1]{
\mkern 1.5mu
\overline{\mkern-1.5mu\mbox{$#1$}\raisebox{2.2mm}{}\mkern-1.5mu}
\mkern 1.5mu
}
\newcommand\ovhb[1]{
\mkern 1.5mu
\overline{\mkern-1.5mu\mbox{$#1$}\raisebox{3.2mm}{}\mkern-1.5mu}
\mkern 1.5mu
}

\newcommand\frct[2]{
\frac
{#1}
{\raisebox{-.1cm}{$#2$}}
}

\newcommand\chisubpi{\chi\raisebox{-.1cm}{$\scaleto{_\pi}{4.5pt}$}}
\newcommand\chisubpiStar{\chi\raisebox{-.1cm}{$\scaleto{_{\pi^*}}{6.7pt}$}}
\newcommand\chisubpiLR{\chi\raisebox{-.1cm}{$\scaleto{_{\pi_{L,R}}}{6.5pt}$}}
\newcommand\chisubpiOne{\chi\raisebox{-.1cm}{$\scaleto{_{\pi_{1}}}{6.5pt}$}}
\newcommand\chisubpiTwo{\chi\raisebox{-.1cm}{$\scaleto{_{\pi_{2}}}{6.5pt}$}}

\newcommand\chisubDelta{\chi\raisebox{-.1cm}{$\scaleto{_\Delta}{5.5pt}$}}
\newcommand\chisubtheta{\chi\raisebox{-.1cm}{$\scaleto{_\theta}{5.5pt}$}}
\newcommand\chisubthetaOne{\chi\raisebox{-.1cm}{$\scaleto{_{\theta_{1}}}{6.5pt}$}}
\newcommand\chisubthetaTwo{\chi\raisebox{-.1cm}{$\scaleto{_{\theta_{2}}}{6.5pt}$}}

\newcommand\chisubOne{\chi\raisebox{-.1cm}{$\scaleto{_1}{4.5pt}$}}
\newcommand\chisubTwo{\chi\raisebox{-.1cm}{$\scaleto{_2}{4.5pt}$}}

\newcommand\chisubPi{\chi\raisebox{-.1cm}{$\scaleto{_\Pi}{5.5pt}$}}
\newcommand\chisubPiStar{\chi\raisebox{-.1cm}{$\scaleto{_{\Pi^*}}{5.5pt}$}}
\newcommand\chisubPiZero{\chi\raisebox{-.1cm}{$\scaleto{_{\Pi_0}}{6.5pt}$}}
\newcommand\chisubPiOne{\chi\raisebox{-.1cm}{$\scaleto{_{\Pi_1}}{6.5pt}$}}
\newcommand\chisubPiTwo{\chi\raisebox{-.1cm}{$\scaleto{_{\Pi_2}}{6.5pt}$}}

\newcommand\chisubPii{\chi\raisebox{-.1cm}{$\scaleto{_{\Pi_i}}{6.5pt}$}}
\newcommand\chisubPin{\chi\raisebox{-.1cm}{$\scaleto{_{\Pi_n}}{6.5pt}$}}

\newcommand\chisubPiOneuXTwo{\chi\raisebox{-.1cm}{$\scaleto{_{\Pi_1 \hspace{0.05cm} \un{\otimes} \hspace{0.05cm} \Pi_2}}{6.5pt}$}}

\newcommand\chisubLR{\chi\raisebox{-.1cm}{$\scaleto{_{L,R}}{6.5pt}$}}
\newcommand\chisubC{\chi\raisebox{-.1cm}{$\scaleto{_C}{4.5pt}$}}
\newcommand\chisubCi{\chi\raisebox{-.1cm}{$\scaleto{_{C_i}}{6.5pt}$}}
\newcommand\chisubCj{\chi\raisebox{-.1cm}{$\scaleto{_{C_j}}{6.5pt}$}}
\newcommand\chisubCk{\chi\raisebox{-.1cm}{$\scaleto{_{C_k}}{6.5pt}$}}

\newcommand\chisubL{\chi\raisebox{-.1cm}{$\scaleto{_L}{4.5pt}$}}
\newcommand\chisubLsubGmodGamma{\chi\raisebox{-.1cm}{$\scaleto{_{L_{G/\Gamma}}}{6.5pt}$}}

\newcommand\chisubLambda{\chi\raisebox{-.1cm}{$\scaleto{_\Lambda}{4.5pt}$}}
\newcommand\chisubGmodGamma{\chi\raisebox{-.1cm}{$\scaleto{_{G/\Gamma}}{6.5pt}$}}

\newcommand\chisubfU{\chi\raisebox{-.1cm}{$\scaleto{_\fU}{5.8pt}$}}

\newcommand\chisubk{\chi\raisebox{-.1cm}{$\scaleto{_k}{4.5pt}$}}
\newcommand\chisubi{\chi\raisebox{-.1cm}{$\scaleto{_i}{4.5pt}$}}
\newcommand\chisubs{\chi\raisebox{-.1cm}{$\scaleto{_s}{4.5pt}$}}
\newcommand\ranglesubPi{\rangle\raisebox{-.1cm}{$\scaleto{_\Pi}{5.5pt}$}}
\newcommand\ranglesubGamma{\rangle\raisebox{-.1cm}{$\scaleto{_\Gamma}{5.5pt}$}}
\newcommand\ranglesubGammaOne{\rangle\raisebox{-.1cm}{$\scaleto{_{\Gamma_1}}{6.5pt}$}}
\newcommand\ranglesubGammaTwo{\rangle\raisebox{-.1cm}{$\scaleto{_{\Gamma_2}}{6.5pt}$}}
\newcommand\ranglesubGammaTwoOfs{\rangle\raisebox{-.1cm}{$\scaleto{_{\Gamma_2(s)}}{6.9pt}$}}

\newcommand\ranglesubG{\rangle\raisebox{-.1cm}{$\scaleto{_G}{5.5pt}$}}
\newcommand\ranglesubGxG{\rangle\raisebox{-.1cm}{$\scaleto{_{G \times G}}{5.5pt}$}}
\newcommand\ranglesubGOne{\rangle\raisebox{-.1cm}{$\scaleto{_{G_1}}{5.5pt}$}}
\newcommand\ranglesubGTwo{\rangle\raisebox{-.1cm}{$\scaleto{_{G_2}}{5.5pt}$}}

\newcommand\ranglesubL{\rangle\raisebox{-.1cm}{$\scaleto{_L}{5.5pt}$}}
\newcommand\ranglesubN{\rangle\raisebox{-.1cm}{$\scaleto{_N}{5.5pt}$}}
\newcommand\ranglesubP{\rangle\raisebox{-.1cm}{$\scaleto{_P}{5.5pt}$}}

\newcommand\ranglesubTorus{\rangle\raisebox{-.1cm}{$\scaleto{_\Torus}{5.5pt}$}}
\newcommand\ranglesubHn{\rangle\raisebox{-.1cm}{$\scaleto{_{H_n}}{6.5pt}$}}

\newcommand\ranglesubEE{\rangle\raisebox{-.1cm}{$\scaleto{_\E}{5.5pt}$}}
\newcommand\ranglesubV{\rangle\raisebox{-.1cm}{$\scaleto{_V}{5.5pt}$}}
\newcommand\ranglesubtheta{\rangle\raisebox{-.1cm}{$\scaleto{_\theta}{5.5pt}$}}
\newcommand\ranglesubpi{\rangle\raisebox{-.1cm}{$\scaleto{_\pi}{4.0pt}$}}

\newcommand\esubPi{e\raisebox{-.1cm}{$\scaleto{_\Pi}{5.5pt}$}}
\newcommand\esubPiOne{e\raisebox{-.1cm}{$\scaleto{_{\Pi_1}}{6.5pt}$}}
\newcommand\esubPiTwo{e\raisebox{-.1cm}{$\scaleto{_{\Pi_2}}{6.5pt}$}}

\newcommand\pisubLR{\pi\raisebox{-.1cm}{$\scaleto{_{L, R}}{6.5pt}$}}

\usepackage{scalerel,stackengine}
\stackMath
\newcommand\reallywidehat[1]{%
\savestack{\tmpbox}{\stretchto{%
  \scaleto{%
    \scalerel*[\widthof{\ensuremath{#1}}]{\kern-.6pt\bigwedge\kern-.6pt}%
    {\rule[-\textheight/2]{1ex}{\textheight}}
  }{\textheight}%
}{0.5ex}}%
\stackon[1pt]{#1}{\tmpbox}%
}
\makeatletter
\DeclareRobustCommand\widecheck[1]{{\mathpalette\@widecheck{#1}}}
\def\@widecheck#1#2{%
    \setbox\z@\hbox{\m@th$#1#2$}%
    \setbox\tw@\hbox{\m@th$#1%
       \widehat{%
          \vrule\@width\z@\@height\ht\z@
          \vrule\@height\z@\@width\wd\z@}$}%
    \dp\tw@-\ht\z@
    \@tempdima\ht\z@ \advance\@tempdima2\ht\tw@ \divide\@tempdima\thr@@
    \setbox\tw@\hbox{%
       \raise\@tempdima\hbox{\scalebox{1}[-1]{\lower\@tempdima\box
\tw@}}}%
    {\ooalign{\box\tw@ \cr \box\z@}}}
\makeatother

\newtheoremstyle{xx}
  {4pt}
  {0pt}
  {\upshape}
  {\bfseries}
  {}
  { }
  {}
  
\makeatletter 
 \newtheoremstyle{myu}%
  {\upshape\item[ \indent\indent\bf\underline{\theorem@headerfont ##2:}]}%
\makeatother
\makeatletter 
 \newtheoremstyle{myn}%
  {\item[\hskip\labelsep \ \bf ##1 \theorem@headerfont ##2.]}%
\makeatother
\theoremstyle{myn}
\newtheorem{theoremn}{Theorem} 
\theoremstyle{myu}
{\upshape}
\newtheorem{x}[theoremn]{}

 \newtheoremstyle{mr}%
  {\upshape\item[ \indent{\theorem@headerfont ##2. \hspace{.2cm}}]}%
\makeatother
\theoremstyle{mr}
{\upshape}
\newtheorem{rf}[theoremn]{}

\title{\textbf{\mG, $\Gamma$, $G/\Gamma$}}
\author{Garth Warner\\
Department of Mathematics\\
University of Washington}
\date{}									

\titleformat{\chapter}[display]
{\normalfont\filcenter\huge\bfseries}{}{0pt}{\large}

\titleformat{\chapter}[display]
{\normalfont\filcenter\huge\bfseries}{}{0pt}{\large}

\usepackage[OT2,OT1]{fontenc} 
\newcommand\cyr
{
\renewcommand\rmdefault{wncyr} 
\renewcommand\sfdefault{wncyss} 
\renewcommand\encodingdefault{OT2} 
\normalfont
\selectfont
}

\DeclareTextFontCommand{\textcyr}{\cyr}

\makeindex 
\begin{document}

\maketitle                              

\titlespacing*{\chapter}{0pt}{-50pt}{40pt}
\setlength{\parskip}{0.1em}
\include{__abstract}

\include{_tocX}
\pagenumbering{bychapter}
\setcounter{chapter}{0}
\include{_A}

\include{_B}

\include{_C}

\include{__refs}

\end{document}

%% file: __abstract.tex
\chapter{
ABSTRACT}
\setlength\parindent{2em}
\setcounter{theoremn}{0}

\ \indent 

The purpose of this book is to provide an introduction to one of the fundamental tools of abstract harmonic analysis, namely the Selberg trace formula.


\chapter{
ACKNOWLEDGEMENT}
\setlength\parindent{2em}
\setcounter{theoremn}{0}

\ \indent 

Many thanks to David Clark for his rendering the original transcript into AMS tech.  
Both of us also thank Judith Clare for her meticulous proofreading.


%% file: _tocX.tex
\begingroup
\fontsize{11pt}{11pt}\selectfont

$\text{ }$
\\[-0.75cm]
\[
\textbf{CONTENTS}
\]

A: \ FINITE GROUPS\\

\qquad  \un{I:}\\

\qquad\qquad $\S1.\ $ \qquad ASSOCIATIVE ALGEBRAS\\

\qquad\qquad $\S2.\ $ \qquad REPRESENTATION THEORY\\

\qquad\qquad $\S3.\ $ \qquad CHARACTERS\\

\qquad\qquad $\S4.\ $ \qquad SIMPLE AND SEMISIMPLE ALGEBRAS\\

\qquad  \un{II:}\\

\qquad\qquad $\S1.\ $ \qquad GROUP ALGEBRAS \\

\qquad\qquad $\S2.\ $ \qquad CONTRAGREDIENTS AND TENSOR PRODUCTS \\

\qquad\qquad $\S3.\ $ \qquad FOURIER TRANSFORMS\\

\qquad\qquad $\S4.\ $ \qquad CLASS FUNCTIONS\\

\qquad\qquad $\S5.\ $ \qquad DECOMPOSITION THEORY\\

\qquad\qquad $\S6.\ $ \qquad INTEGRALITY\\

\qquad\qquad $\S7.\ $ \qquad INDUCED CLASS FUNCTIONS\\

\qquad\qquad $\S8.\ $ \qquad MACKEY THEORY\\

\qquad\qquad $\S9.\ $ \qquad INDUCED REPRESENTATIONS\\

\qquad\qquad $\S10.$ \qquad IRREDUCIBILITY OF $\Ind_{\Gamma, \Theta}^G$\\

\qquad\qquad $\S11.$ \qquad BURNSIDE RINGS\\

\qquad\qquad $\S12.$ \qquad BRAUER THEORY\\

\qquad\qquad $\S13.$ \qquad GROUPS OF LIE TYPE\\

\qquad\qquad $\S14.$ \qquad HARISH-CHANDRA THEORY\\

\qquad\qquad $\S15.$ \qquad HOWLETT-LEHRER THEORY\\

\qquad\qquad $\S16.$ \qquad MODULE LANGUAGE\\

\qquad  \un{III:}\\

\qquad\qquad $\S1.\ $ \qquad ORBITAL SUMS \\

\qquad\qquad $\S2.\ $ \qquad THE LOCAL TRACE FORMULA \\

\qquad\qquad $\S3.\ $ \qquad THE GLOBAL PRE-TRACE FORMULA\\

\qquad\qquad $\S4.\ $ \qquad THE GLOBAL TRACE FORMULA\\

B: \ COMPACT GROUPS\\

\qquad  \un{I:}\\

\qquad\qquad $\S1.\ $ \qquad UNITARY REPRESENTATIONS \\

\qquad\qquad $\S2.\ $ \qquad EXPANSION THEORY \\

\qquad\qquad $\S3.\ $ \qquad STRUCTURE THEORY\\

\qquad\qquad $\S4.\ $ \qquad MAXIMAL TORI\\

\qquad\qquad $\S5.\ $ \qquad REGULARITY\\

\qquad\qquad $\S6.\ $ \qquad WEIGHTS AND ROOTS\\

\qquad\qquad $\S7.\ $ \qquad LATTICES\\

\qquad\qquad $\S8.\ $ \qquad WEYL CHAMBERS AND WEYL GROUPS\\

\qquad\qquad $\S9.\ $ \qquad DESCENT\\

\qquad\qquad $\S10.$ \qquad CHARACTER THEORY\\

\qquad\qquad $\S11.$ \qquad THE INVARIANT INTEGRAL\\

\qquad\qquad $\S12.$ \qquad PLANCHEREL\\

\qquad\qquad $\S13.$ \qquad DETECTION\\

\qquad\qquad $\S14.$ \qquad INDUCTION\\

\qquad  \un{II:}\\

\qquad\qquad $\S1.\ $ \qquad ORBITAL INTEGRALS \\

\qquad\qquad $\S2.\ $ \qquad KERNELS \\

\qquad\qquad $\S3.\ $ \qquad THE LOCAL TRACE FORMULA\\

C: \ LOCALLY COMPACT GROUPS\\

\qquad  \un{I:}\\

\qquad\qquad $\S1.\ $ \qquad TOPOLOGICAL TERMINOLOGY \\

\qquad\qquad $\S2.\ $ \qquad INTEGRATION THEORY \\

\qquad\qquad $\S3.\ $ \qquad UNIMODULARITY\\

\qquad\qquad $\S4.\ $ \qquad INTEGRATION ON HOMOGENEOUS SPACES\\

\qquad\qquad $\S5.\ $ \qquad INTEGRATION ON LIE GROUPS\\

\qquad  \un{II:}\\

\qquad\qquad $\S1.\ $ \qquad TRANSVERSALS \\

\qquad\qquad $\S2.\ $ \qquad LATTICES \\

\qquad\qquad $\S3.\ $ \qquad UNIFORMLY INTEGRABLE FUNCTIONS\\

\qquad\qquad $\S4.\ $ \qquad THE SELBERG TRACE FORMULA\\

\qquad\qquad $\S5.\ $ \qquad FUNCTIONS OF REGULAR GROWTH\\

\qquad\qquad $\S6.\ $ \qquad DISCRETE SERIES\\
\\

REFERENCES\\
\endgroup 

%% file: _A.tex
\chapter{
$\boldsymbol{\S}$\textbf{1}.\quad  ASSOCIATIVE ALGEBRAS}
\setlength\parindent{2em}
\setcounter{theoremn}{0}
\renewcommand{\thepage}{A I \S1-\arabic{page}}


\begin{x}{\small\bf DEFINITION} \ 
An \un{associative algebra} over $\Cx$ is a finite dimensional vector space $\sA$ over $\Cx$ equipped with a bilinear map 
\begin{align*}
\mu : \sA \times \sA &\ra \sA \\
(x,y) &\ra \mu(x,y) \ \equiv \ x y
\end{align*}
such that 
$(x y) z = x (y z)$.
\end{x}
\vspace{0.1cm}

\begin{x}{\small\bf DEFINITION} \ 
An associative algebra $\sA$ is said to be \un{unital} if there exists an element $e \in \sA$ with the property that 
$x e = e x  = x$ for all $x \in \sA$.
\vspace{0.1cm}

[Note: \ Such an $e$ is called an \un{identity element} and is denoted by $1_\sA$.]
\end{x}
\vspace{0.1cm}

\begin{x}{\small\bf \un{N.B.}} \ 
Identity elements are unique.
\end{x}
\vspace{0.1cm}

\begin{x}{\small\bf EXAMPLE} \ 
Let \mV be a finite dimensional vector space over $\Cx$ $-$then $\Hom(V)$ 
(the set of all $\Cx$-linear maps of \mV) is a unital associative algebra over $\Cx$ 
(multiplication being composition of linear transformations and the identity element $\id_V$).
\end{x}
\vspace{0.1cm}

Let $\sA$ be an associative algebra over $\Cx$.

\begin{x}{\small\bf DEFINITION} \ 
A 
\un{representation} 
\index{representation} 
of $\sA$ is a pair $(\rho,V)$, where \mV is a finite dimensional vector space over $\Cx$ and 
$\rho:\sA \ra \Hom(V)$ is a morphism of associative algebras.
\vspace{0.1cm}

[Note: \ If $\sA$ is unital, then it will be assumed that $\rho(1_\sA) = \id_V$, thus is a morphism of unital associative algebras.]

\end{x}

\vspace{0.1cm}

\begin{x}{\small\bf DEFINITION} \ 
Let $(\rho,V)$ be a representation of $\sA$ $-$then a linear subspace $U \subset V$ is said to be 
\un{$\rho$-invariant} if $\forall \ x \in \sA$, $\rho(x) U \subset U$.
\end{x}
\vspace{0.1cm}


\begin{x}{\small\bf \un{N.B.}} \ 
A $\rho$-invariant subspace $U \subset V$ gives rise to two representations of $\sA$, viz. by restricting to \mU and passing to the quotient $V/U$.
\end{x}
\vspace{0.1cm}

\begin{x}{\small\bf DEFINITION} \ 
A representation $(\rho,V)$ of \mA is \un{irreducible} if $V \neq \{0\}$ and if the only $\rho$-invariant subspaces are $\{0\}$ and \mV.
\end{x}
\vspace{0.1cm}

\begin{x}{\small\bf NOTATION} \ 
Given a representation $(\rho,V)$ of $\sA$, put
\[
\Ker(\rho) \ = \ \{x \in \sA: \rho(x) = 0\}.
\]
\end{x}
\vspace{0.1cm}

\begin{x}{\small\bf \un{N.B.}} \ 
$\Ker(\rho)$ is a two-sided ideal in $\sA$.
\end{x}
\vspace{0.1cm}

\begin{x}{\small\bf DEFINITION} \ 
A representation $(\rho,V)$ of $\sA$ is 
\un{faithful}
\index{faithful} 
if $\Ker(\rho) = \{0\}$.
\end{x}
\vspace{0.1cm}

\begin{x}{\small\bf DEFINITION} \ 
Let $(\rho,V)$ and $(\sigma,W)$ be representations of $\sA$ $-$then an
\un{intertwining operator}
\index{intertwining operator} 
is a $\Cx$-linear map $T:V \ra W$ such that $T \rho(x) = \sigma(x) T$ for all $x \in \sA$.
\end{x}
\vspace{0.1cm}

\begin{x}{\small\bf NOTATION} \ 
$I_\sA(\rho,\sigma)$ is the set of intertwining operators between $(\rho,V)$ and $(\sigma,W)$.
\end{x}
\vspace{0.1cm}

\begin{x}{\small\bf EXAMPLE} \ 
Let $(\rho,V)$ be a representation of $\sA$ and suppose that $U \subset V$ is a $\rho$-invariant subspace $-$then the inclusion map 
$U \ra V$ is an intertwining operator, as is the quotient map $V \ra V/U$.
\end{x}
\vspace{0.1cm}

\begin{x}{\small\bf DEFINITION} \ 
Representations $(\rho,V)$ and $(\sigma,W)$ of $\sA$ are 
\un{equivalent}
\index{represenation\\equivalent} 
if there exists an invertible operator in $I_\sA(\rho,\sigma)$, in which case we write 
\[
(\rho,V) \ \approx \ (\sigma,W) \qquad (\text{or} \ \rho \approx \sigma).
\]
\end{x}
\vspace{0.1cm}

\begin{x}{\small\bf NOTATION} \ 
$\widehat{\sA}$ is the set of equivalence classes of irreducible representations of $\sA$.
\end{x}

\vspace{0.1cm}
\begin{x}{\small\bf EXAMPLE} \ 
Take $\sA = \Hom(V)$, where \mV is a finite dimensional complex vector space $-$then up to equivalence, the only irreducible representation of $\Hom(V)$ is the representation $\rho$ given by 
\[
\rho(T) v \ = \ Tv \qquad (T \in \Hom(V)).
\]
\end{x}
\vspace{0.1cm}


\chapter{
$\boldsymbol{\S}$\textbf{2}.\quad  REPRESENTATION THEORY}
\setlength\parindent{2em}
\setcounter{theoremn}{0}
\renewcommand{\thepage}{A I \S2-\arabic{page}}

$\text{ }$\\[-1.5cm]

Let $\sA$ be a unital associative algebra over $\Cx$.
\\[-.2cm]

\begin{x}{\small\bf THEOREM} \ 
Let $(\rho,V)$, $(\sigma,W)$ be irreducible representations of $\sA$ $-$then
\[
\dim I_\sA(\rho,\sigma) \ = \ 
\begin{cases}
\ 1 \ \text{if} \  (\rho,V) \approx (\sigma,W)\\
\ 0 \ \text{if} \  (\rho,V) \not\approx (\sigma,W)
\end{cases}
.
\]
\\
\end{x}

\begin{x}{\small\bf THEOREM} \ 
Let $(\rho,V)$  be an irreducible representation of $\sA$ $-$then $\rho(\sA) = \Hom(V)$.
\\
\end{x}

\begin{x}{\small\bf DEFINITION} \ 
A representation $(\rho,V)$ of $\sA$ is 
\un{completely reducible} 
\index{completely reducible} 
if for every $\rho$-invariant subspace $V_1 \subset V$ there exists a $\rho$-invariant subspace $V_2 \subset V$ such that 
$V = V_1 \oplus V_2$.
\\
\end{x}

\begin{x}{\small\bf LEMMA} \ 
Suppose that  $(\rho,V)$ is a  representation of $\sA$ $-$then $(\rho,V)$ is completely reducible iff there is a decomposition 
\[
V \ = \ V_1 \oplus \cdots \oplus V_s,
\]
where each $V_i$ is $\rho$-invariant and irreducible.
\\
\end{x}

Let $\sA$ be an associative algebra over $\Cx$.
\\

\begin{x}{\small\bf LEMMA} \ 
Suppose that  $(\rho,V)$ is a  representation of $\sA$ $-$then $(\rho,V)$ is completely reducible iff there is a decomposition 
\[
V \ = \ U_1 + \cdots + U_t,
\]
where each $U_j$ is $\rho$-invariant and irreducible.
\\
\end{x}

\begin{x}{\small\bf DEFINITION} \ 
Let \mV be a finite dimensional vector space over $\Cx$.  
Given a subset $\sS$ of $\Hom(V)$, put
\[
\Com(\sS) \ = \ \{T \in \Hom(V) : T S = S T \ \forall \ S \in \sS\},
\]
the 
\underline{commutant}
\index{commutant} 
of 
$\sS$.
\\
\end{x}


\begin{x}{\small\bf \un{N.B.}} \ 
$\Com(\sS)$ is a unital associative algebra over $\Cx$.
\\
\end{x}

\begin{x}{\small\bf THEOREM} \ 
Suppose that \mV is a finite dimensional vector space over $\Cx$ and let $\sV \subset \Hom(V)$ be an associative algebra over $\Cx$ with identity $\id_V$.  
Assume: \mV is completely reducible per the canonical action of $\sV$ $-$then
\[
\Com(\Com(\sV) \ = \sV.
\]

[Note: \ A priori, 
\[
\sV \subset \Com(\Com(\sV)).]
\]
\\[-.5cm]
\end{x}

\begin{x}{\small\bf NOTATION} \ 
Let $(\rho,V)$ be a completely reducible representation of $\sA$.  
Given $\delta \in \widehat{\sA}$, put 
\[
V_\delta \ = \  \sum\limits_{U \subset V: [U] = \delta} U,
\] 
the subspaces \mU being $\rho$-invariant and irreducible, $[U]$ standing for the equivalence class in $\widehat{\sA}$ determined by \mU.
\\
\end{x}

\begin{x}{\small\bf THEOREM} \ 
Let $(\rho,V)$ be a completely reducible representation of $\sA$ and let
\[
V \ = \ V_1 \oplus \cdots \oplus V_s
\]
be a decomposition, where each $V_i$ is $\rho$-invariant and irreducible $-$then $\forall \ \delta \in \widehat{\sA}$, 
\[
V_\delta \ = \  \bigoplus\limits_{[V_i] = \delta} V_i,
\] 
thus
\[
V \ = \ \bigoplus\limits_{\delta \in \widehat{\sA}} V_\delta.
\]
\vspace{0.1cm}

[Note: \ An empty sum is taken to be zero.]
\\
\end{x}

\begin{x}{\small\bf DEFINITION} \ 
The decomposition 
\[
V \ = \ \bigoplus\limits_{\delta \in \widehat{\sA}} V_\delta
\]
is the 
\un{primary decomposition}
\index{primary decomposition} 
of $V$ and $V_\delta$ is the 
\un{$\delta$-isotypic}
\index{$\delta$-isotropic} 
subspace of \mV.
\end{x}
\vspace{0.1cm}

\begin{x}{\small\bf DEFINITION} \ 
The cardinality $m_V(\delta)$ of 
\[
\{i:[V_i] = \delta\}
\]
is the 
\un{multiplicity}
\index{multiplicity} 
of $\delta$ in \mV.
\\
\end{x}

\begin{x}{\small\bf NOTATION} \ 
Given $\delta \in \widehat{\sA}$, let $U(\delta)$ be an element in the class $\delta$.
\\
\end{x}

\begin{x}{\small\bf LEMMA} \ 
\[
m_V(\delta) \ = \ \dim I_\sA(U(\delta),V) \ = \ \dim I_\sA (V,U(\delta)).
\]
\end{x}


\chapter{
$\boldsymbol{\S}$\textbf{3}.\quad  CHARACTERS}
\setlength\parindent{2em}
\setcounter{theoremn}{0}
\renewcommand{\thepage}{A I \S3-\arabic{page}}

$\text{ }$\\[-1.5cm]

Let $\sA$ be a unital associative algebra over $\Cx$.
\\[-.2cm]

\begin{x}{\small\bf DEFINITION} \ 
Let $(\rho,V)$ be a representation of $\sA$ $-$then its 
\un{character}
\index{character} 
is the linear functional
\[
\chi_\rho : \sA \ra \Cx
\]
given by the prescription
\[
\chi_\rho (x) \ = \ \tr(\rho(x)) \qquad (x \in \sA).
\]
\\
\end{x}

\begin{x}{\small\bf LEMMA} \ 
\[
\chi_\rho (1_\sA) \ = \ \dim V.
\]
\end{x}
\vspace{0.1cm}

\begin{x}{\small\bf LEMMA} \ 
$\forall \ x,\ y \in \sA$, 
\[
\chi_\rho (x y) \ = \ \chi_\rho (y x).
\]
\\
\end{x}

\begin{x}{\small\bf DEFINITION} \ 
Let $(\rho,V)$ be a representation of $\sA$ $-$then a 
\un{composition series}
\index{composition series} 
for $\rho$ is a sequence of $\rho$-invariant subspaces 
\[
\{0\} \ = \ V_0 \subset V_1 \subset \cdots \subset V_s \ = \ V
\]
such that 
\[
\{0 \} \ \neq \ V_i / V_{i - 1} \qquad (i = 1, \ldots, s)
\]
is irreducible.
\\
\end{x}

\begin{x}{\small\bf LEMMA} \ 
Composition series exist.
\\
\end{x}

\begin{x}{\small\bf DEFINITION} \ 
The 
\un{semisimplification}
\index{semisimplification} 
of $(\rho,V)$ is the direct sum 
\[
V_{\ssx} \ = \ \bigoplus\limits_{i = 1}^s (V_i / V_{i - 1})
\]
equipped with the canonical operations.
\end{x}
\vspace{0.1cm}

\begin{x}{\small\bf DEFINITION} \ 
The irreducible quotients $V_i / V_{i-1}$ are the 
\un{composition factors} 
\index{composition factors} 
of $(\rho,V)$.
\\
\end{x}

Let $\rho_\ssx$ be the representation of $\sA$ per $V_\ssx$ and let $\rho_i$ be the representation of $\sA$ per $V_i/V_{i-1}$.
\\

\begin{x}{\small\bf LEMMA} \ 
\[
\chisub{\rho_\ssx} \ = \ \sum\limits_{i = 1}^s \chisub{\rho_i} \ = \ \chisub{\rho}.
\]
\\
\end{x}

\begin{x}{\small\bf LEMMA} \ 
Suppose that $(\sigma_1,U_1), \ldots, (\sigma_r,U_r)$ are irreducible representations of $\sA$.  
Assume: \ $(\sigma_k,U_k)$ is not equivalent to $(\sigma_\ell,U_\ell)$  $(k \neq \ell)$ $-$then the set 
\[
\{\chisub{\sigma_1}, \ldots, \chisub{\sigma_r}\}
\]
is linearly independent.
\\
\end{x}

\begin{x}{\small\bf SCHOLIUM} \ 
The composition factors in a composition series for $\rho$ are unique up to isomorphism and order of appearance and 
$(\rho_{\ssx},V_\ssx)$ is uniquely determined by 
$\chisub{\rho}$ 
up to isomorphism.
\end{x}



\chapter{
$\boldsymbol{\S}$\textbf{4}.\quad  SIMPLE AND SEMISIMPLE ALGEBRAS}
\setlength\parindent{2em}
\setcounter{theoremn}{0}
\renewcommand{\thepage}{A I \S4-\arabic{page}}
\vspace{-.5cm}

\ \indent 
Let $\sA$ be a unital associative algebra over $\Cx$.
\\[-.25cm]

\begin{x}{\small\bf DEFINITION} \ 
$\sA$ is 
\un{simple}
\index{simple algebra}
\index{algebra // simple} 
if the only two-sided ideals in $\sA$ are $\{0\}$ and $\sA$.
\\[-.2cm]
\end{x}

\begin{x}{\small\bf LEMMA} \ 
If \mV is a finite dimensional vector space over $\Cx$, then $\Hom(V)$ is simple.
\\[-.2cm]
\end{x}

\begin{x}{\small\bf THEOREM} \ 
If $\sA$ is simple, then there is a finite dimensional vector space \mV over $\Cx$ such that $\sA \hsx \approx \hsx \Hom(V)$.
\\[-.2cm]
\end{x}

\begin{x}{\small\bf DEFINITION} \ 
$\sA$ is 
\un{semisimple}
\index{semisimple algebra}
\index{algebra // semisimple} 
if it is a finite direct sum of simple algebras.
\\[-.2cm]
\end{x}

Accordingly, if $\sA$ is semisimple, then there is a finite set \mL,  finite dimensional complex vector spaces 
$V_\lambda$ $(\lambda \in L)$, and an isomorphism
\[
\phi : \sA \ra \bigoplus\limits_{\lambda \in L} \hsx \Hom(V_\lambda).
\]

Denote by $E_\lambda$ the element
\[
0 \oplus \cdots \oplus \id_{V_\lambda}  \oplus \cdots \oplus 0
\]
and define a representation $(\rho_\lambda, V_\lambda)$ by the prescription 
\[
\rho_\lambda (x) 
\ = \ 
\phi(x) E_\lambda 
\qquad (x \in \sA).
\]
\\[-.75cm]

\begin{x}{\small\bf LEMMA} \ 
The $(\rho_\lambda, V_\lambda)$ are irreducible.
\\[-.2cm]
\end{x}


\begin{x}{\small\bf THEOREM} \ 
Every irreducible representation of $\sA$ is equivalent to some $(\rho_\lambda, V_\lambda)$.
\\[-.2cm]
\end{x}

\begin{x}{\small\bf \un{N.B.}} \ 
Therefore 
\[
\widehat{\sA} \longleftrightarrow L,
\]
so the term ``$\lambda$-isotypic subspace'' makes sense.
\end{x}
\vspace{0.1cm}

Put
\[
e_\lambda 
\ = \ 
\phi^{-1} (E_\lambda).
\]
Then $e_\lambda$ is a central idempotent and 
\[
\sum\limits_{\lambda \in L} \hsx e_\lambda \ = \ 1_\sA.
\]
\\[-.75cm]

\begin{x}{\small\bf THEOREM} \ 
Suppose that $\sA$ is semisimple and let $(\rho, V)$ be a representation of $\sA$ $-$then 
its $\lambda$-isotypic subspace is $\rho(e_\lambda) V$ and 
\[
V 
\ = \ 
\bigoplus\limits_{\lambda \in L} \hsx 
\rho(e_\lambda) V
\]
is the primary decomposition of \mV.
\\
\end{x}

\begin{x}{\small\bf LEMMA} \ 
Let $\sA$ be a unital associative algebra over $\Cx$ and let $(\rho,V)$ be a completely reducible representation of 
$\sA$ $-$then $\rho(\sA)$ is semisimple.
\\
\end{x}

\begin{x}{\small\bf THEOREM} \ 
Let $\sA$ be a unital associative algebra over $\Cx$ $-$then the following conditions are equivalent:
\\[-.2cm]

1. \ 
The left regular representation $(L,\sA)$ of $\sA$ is completely reducible $(L(x)y = xy)$.
\\[-.3cm]

2. \ 
Every representation of $\sA$ is completely reducible.
\\[-.3cm]

3. \ 
$\sA$ is a semisimple algebra.
\\[-.3cm]


[$1 \implies 3$: \ $L(\sA)$ is semisimple (cf. \#9).  
On the other hand, $\sA \hsx \approx \hsx L(\sA)$, \mL being faithful.
\vspace{0.2cm}

$3 \implies 2$: \ 
Quote \#8 and \S2, \#4.
\vspace{0.1cm}

$2 \implies 1$: \ Obvious.]
\end{x}
\vspace{0.1cm}

\begin{x}{\small\bf THEOREM} \ 
Every representation of a semisimple algebra is uniquely determined by its character up to isomorphism.
\end{x}
\vspace{0.1cm}


\chapter{
$\boldsymbol{\S}$\textbf{1}.\quad  GROUP ALGEBRAS}
\setlength\parindent{2em}
\setcounter{theoremn}{0}
\renewcommand{\thepage}{A II \S1-\arabic{page}}


\begin{x}{\small\bf NOTATION} \ 
If \mX is a finite set, then $\abs{X}$ is the cardinality of \mX and $C(X)$ is the vector space of complex valued functions on \mX.
\end{x}
\vspace{0.1cm}

\begin{x}{\small\bf \un{N.B.}} \ 
The functions $\{\delta_x : x \in X\}$, where
\[
\delta_x(y) \ = \ 
\begin{cases}
\ 1 \qquad (x = y)\\
\ 0 \qquad (x \neq y)
\end{cases}
,
\]
constitute a basis for $C(X)$.  
Therefore
\[
\dim C(X) \ = \ \abs{X}
\]
and every $f \in C(X)$ admits a decomposition
\[
f \ = \ \sum\limits_{x \in X} \hsz f(x) \delta_x.
\]
In particular: If $1_X$ is the function on \mX which is $\ \equiv \hsx 1$, then 
\[
1_X \ = \ \sum\limits_{x \in X} \hsz  \delta_x.
\]
\end{x}
\vspace{0.1cm}

Let \mG be a finite group.
\\[-.2cm]

\begin{x}{\small\bf DEFINITION} \ 
Given $f, \ g \in C(G)$, their 
\un{convolution}
\index{convolution}\index{f * g} 
$f * g$ is the element of $C(G)$ defined by the rule 
\allowdisplaybreaks
\begin{align*}
(f * g)(x) \ 
&= \ \sum\limits_{y \in G} \hsz f(xy^{-1}) g(y) 
\\[11pt]
&=\ \sum\limits_{y \in G} \hsz f(y) g(y^{-1}x).
\end{align*}

[Note: \ $\forall \ x, \ y \in G$, 
\[
\delta_x * \delta_y \ = \ \delta_{x y}.]
\]
\end{x}
\vspace{0.1cm}

\begin{x}{\small\bf LEMMA} \ 
$C(G)$ is an associative algebra over $\Cx$.
\end{x}
\vspace{0.1cm}

\begin{x}{\small\bf \un{N.B.}} \ 
If $e$ is the identity in \mG, then $\delta_e$ is the identity in $C(G)$, which is therefore unital.
\end{x}
\vspace{0.1cm}

\begin{x}{\small\bf LEMMA} \ 
The center of $C(G)$ consists of those $f$ such that 
\[
f(x) \ = \ f(y x y^{-1}) \qquad (x, y \in G).
\]

[Note: \ In other words, the center of $C(G)$ consists of those $f$ that are constant on conjugacy classes, the so-called 
\un{class functions}.]
\index{class functions} 
\\[-.2cm]

E.g.: \ $\forall \ x \in G$, the function
\[
\sum\limits_{y \in G} \delta_{y x y^{-1}}
\]
is a class function.

[Given $z$ in \mG,
\allowdisplaybreaks
\begin{align*}
\bigg( \hsx\sum\limits_{y \in G} \hsz \delta_{y x y^{-1}} \bigg) * \delta_z \ 
&= \ 
\bigg( \hsx\sum\limits_{y \in G} \hsz \delta_{(zy) x (zy)^{-1}} \bigg) * \delta_z
\\[11pt]
&=\ \sum\limits_{y \in G} \hsz \delta_{z y x (z y)^{-1}z}
\\[11pt]
&=\ \sum\limits_{y \in G} \hsz \delta_{z y x y^{-1}}
\\[11pt]
&=\ \delta_z * \bigg( \hsx \sum\limits_{y \in G} \hsz \delta_{y x y^{-1}} \bigg).]
\end{align*}
\\[11pt]
\end{x}
\vspace{0.1cm}

\begin{x}{\small\bf DEFINITION} \ 
A 
\un{representation}
\index{representation} 
of \mG is a pair $(\pi,V)$, where \mV is a finite dimensional vector space over $\Cx$ and 
$\pi:G \ra \GL(V)$ is a morphism of groups.
\end{x}
\vspace{0.1cm}

\begin{x}{\small\bf SCHOLIUM} \ 
Let \mV be a finite dimensional vector space over $\Cx$.
\\[-.2cm]

\qquad \textbullet \quad Every representation $\pi:G \ra \GL(V)$ extends to a representation $\rho$ of $C(G)$ on \mV, viz.
\[
\rho(f) \ = \ \sum\limits_{x \in G} f(x) \hsz \pi(x).
\]
\vspace{0.1cm}

\qquad \textbullet \quad Every representation $\rho:C(G) \ra \Hom(V)$ restricts to a representation $\pi$ of \mG on \mV, viz.
\[
\pi(x) \ = \ \rho(\delta_x).
\]

[Note: \ If $\pi$ is given, it is customary to denote its extension ``$\rho$'' by $\pi$ as well.]
\end{x}

\vspace{0.1cm}

\begin{x}{\small\bf LEMMA} \ 
Let $W \subset V$ be a linear subspace $-$then \mW is invariant under \mG iff \mW is invariant under $C(G)$.
\end{x}
\vspace{0.1cm}

\begin{x}{\small\bf LEMMA} \ 
An operator $T \in \Hom(V)$ commutes with the action of \mG iff it commutes with the action of $C(G)$.
\end{x}
\vspace{0.1cm}

\begin{x}{\small\bf THEOREM} \ 
$C(G)$ is semisimple.
\\[-.2cm]

PROOF \ Let $(\rho,V)$ be a representation of $C(G)$ and suppose that $V_1 \subset V$ is a $\rho$-invariant subspace.  
Fix a linear complement \mU per $V_1: V = V_1 \oplus U$.  
Let $P:V \ra V_1$ be the corresponding projection and put
\[
Q \ = \ \frac{1}{\abs{G}} \ \sum\limits_{x \in G} \hsz \pi(x) P \pi(x)^{-1}.
\]
Then \mQ is a projection with range $V_1$.  
In addition, $\forall \ y \in G$,
\allowdisplaybreaks
\begin{align*}
\pi(y) Q \ 
&=\ \frac{1}{\abs{G}} \ \sum\limits_{x \in G} \hsz \pi(yx) P \pi(x)^{-1}
\\[11pt]
&=\ \frac{1}{\abs{G}} \ \sum\limits_{x \in G} \hsz \pi(x) P \pi(y^{-1}x)^{-1}
\\[11pt]
&=\ \frac{1}{\abs{G}} \ \sum\limits_{x \in G} \hsz  \pi(x) P \pi(x)^{-1} \pi(y)
\\[11pt]
&=\ Q \pi(y).
\end{align*}
Consequently, $\forall \ v \in V$, 
\allowdisplaybreaks
\begin{align*}
\pi(y) (\id_V - Q) v \ 
&=\ \pi(y) (v - Q v)
\\[11pt]
&=\ \pi(y) v - \pi(y) Q v
\\[11pt]
&=\ \pi(y) v - Q \pi(y) v
\\[11pt]
&=\ (\id_V - Q) \pi(y) v,
\end{align*}
thus the range $V_2$ of $\id_V - Q$ is a $\rho$-invariant complement per $V_1$.  
It therefore follows that every representation of $C(G)$ is completely reducible, hence $C(G)$ is semisimple 
(cf. I, \S4, \#10).
\end{x}
\vspace{0.1cm}

\begin{x}{\small\bf DEFINITION} \ 
\\[-.2cm]

\qquad \textbullet \quad The \un{left translation representation} \mL of \mG on $C(G)$ is the prescription 
\[
L(x) f(y) \ = \ f(x^{-1}y) \qquad (\implies L(x) f = \delta_x * f).
\]

\qquad\textbullet \quad The \un{right translation representation} \mR of \mG on $C(G)$ is the prescription 
\[
R(x) f(y) \ = \ f(yx) \qquad (\implies R(x) f = f * \delta_{x^{-1}}).
\]
\end{x}
\vspace{0.1cm}

\begin{x}{\small\bf \un{N.B.}} \ 
Since $C(G)$ is semisimple, both \mL and \mR are completely reducible.
\end{x}
\vspace{0.1cm}

\begin{x}{\small\bf REMARK} \ 
There is also a representation $\pi_{L,R}$ of $G \times G$ on $C(G)$, namely
\[
(\pi_{L,R}(x_1,x_2) f) (x) \ = \ f(x_1^{-1} x x_2).
\]
And it too is completely reducible ($C(G \times G)$ is semisimple).
\end{x}
\vspace{0.1cm}

\begin{x}{\small\bf DEFINITION} \ 
\ Let $(\pi_1,V_1)$, \ $(\pi_2,V_2)$ be representations of \ \mG \ $-$then an 
\un{intertwining operator}
\index{intertwining operator} 
is a $\Cx$-linear map $T:V_1 \ra V_2$ such that $T \pi_1(x) = \pi_2(x) T$ for all $x \in G$.
\end{x}
\vspace{0.1cm}

\begin{x}{\small\bf NOTATION} \ 
$I_G(\pi_1,\pi_2)$ is the set of intertwining operators between $(\pi_1,V_1)$ and $(\pi_2,V_2)$.
\end{x}
\vspace{0.1cm}

\begin{x}{\small\bf \un{N.B.}} \ 
On the basis of the definitions,
\[
I_G(\pi_1,\pi_2) \ = \ I_{C(G)}(\rho_1,\rho_2).
\]
\end{x}
\vspace{0.1cm}

\begin{x}{\small\bf LEMMA} \ 
Let $(\pi_1,V_1)$, $(\pi_2,V_2)$ be irreducible representations of \mG and let 
$T \in I_G(\pi_1,\pi_2)$ $-$then \mT is zero or it is an isomorphism.
\end{x}
\vspace{0.1cm}

\begin{x}{\small\bf LEMMA} \ 
Suppose that $(\pi,V)$ is an irreducible representation of \mG and suppose that $T \in I_G(\pi,\pi)$ $-$then 
\mT is a scalar multiple of $\id_V$.
\end{x}
\vspace{0.1cm}

\begin{x}{\small\bf DEFINITION} \ 
Representations $(\pi_1,V_1)$ and $(\pi_2,V_2)$ of \mG are 
\un{equivalent}
\index{representations \\equivalent} 
if there exists an invertible operator in $I_G(\pi_1,\pi_2)$, in which case we write
\[
(\pi_1, V_1) \ \approx \ (\pi_2, V_2) \qquad (\text{or} \ \pi_1 \approx \pi_2).
\]
\end{x}
\vspace{0.1cm}

\begin{x}{\small\bf NOTATION} \ 
$\widehat{G}$ is the set of equivalence classes of irreducible representations of \mG.

[Note: \ By convention, the zero representation of \mG on $V = \{0\}$ is not to be viewed as irreducible.]
\end{x}
\vspace{0.1cm}

\begin{x}{\small\bf \un{N.B.}} \ 
There is a one-to-one correspondence
\[
\widehat{G} \ \approx \   \widehat{C(G)}.
\]
\end{x}
\vspace{0.1cm}

In the sequel, $\Pi$ stands for an element of $\widehat{G}$ with representation space $V(\Pi)$ of dimension $d_\Pi$.  
Without loss of generality, it can be assumed moreover that $\Pi$ is unitary with respect to a \mG-invariant inner product 
$\langle \ , \ \ranglesubPi$ on $V(\Pi)$.
\\[-.2cm]

[Recall the argument.  Start with an inner product $\langle \ , \ \rangle$ on $V(\Pi)$ and put
\[
\langle v_1  ,v_2 \ \ranglesubPi 
\ = \ 
\frac{1}{\abs{G}} \ \sum\limits_{x \in G} \hsz \langle \Pi(x)v_1, \Pi(x) v_2 \rangle.]
\]
\vspace{0.1cm}

\[
\text{APPENDIX}
\]

Let $(\pi_1,V_1)$, $(\pi_2,V_2)$ be unitary representations of \mG.  
Suppose that there exists an invertible 
\[
T \in I_G(\pi_1,\pi_2).
\]
Then there exists a unitary
\[
U \in I_G(\pi_1,\pi_2).
\]

[Let $T = U \abs{T}$ be the polar decomposition of \mT $-$then $\forall \ x \in G$, 
\[
\abs{T} \pi_1(x) \ = \ \pi_1(x) \abs{T}.
\]
Therefore
\allowdisplaybreaks
\begin{align*}
U \pi_1(x) U^{-1} \ 
&= \ T \abs{T}^{-1} \pi_1(x) \abs{T} T^{-1}
\\
&=\ T \pi_1(x) T^{-1} 
\\
&= \ \pi_2(x).]
\end{align*}

\chapter{
$\boldsymbol{\S}$\textbf{2}.\quad  CONTRAGREDIENTS AND TENSOR PRODUCTS}
\setlength\parindent{2em}
\setcounter{theoremn}{0}
\renewcommand{\thepage}{A II \S2-\arabic{page}}


\begin{x}{\small\bf NOTATION} \ 
Given a finite dimensional vector space \mV over $\Cx$, let $V^*$ be its dual and denote by 
\[
\begin{cases}
\ V^* \times V \ra \Cx \\
\ (v^*,v) \ra \langle v^*, v \rangle \qquad (= v^*(v))
\end{cases}
\]
the evaluation pairing.
\end{x}
\vspace{0.1cm}

Let \mG be a finite group.
\\[-.2cm]

\begin{x}{\small\bf DEFINITION} \ 
Suppose that $\pi:G \ra \GL(V)$ is a representation $-$then its 
\un{contragredient}
\index{contragredient} 
is the representation $\pi^*:G \ra \GL(V^*)$ defined by requiring that $\forall \ x \in G$, 
\[
\pi^*(x) v^* 
\ = \ v^* \circ \pi(x^{-1}) \qquad (v^* \in V^*),
\]
thus $\forall \ v \in V$, 
\[
\langle \pi^*(x) v^* , v \rangle 
\ = \ 
\langle v^*, \pi(x^{-1}) v \rangle .
\]
\end{x}
\vspace{0.1cm}

\begin{x}{\small\bf \un{N.B.}} \ 
The identification $(V^*)^* \approx V$ leads to an equivalence $(\pi^*)^* \approx \pi$.
\end{x}
\vspace{0.1cm}

\begin{x}{\small\bf LEMMA} \ 
$(\pi,V)$ is irreducible iff $(\pi^*,V^*)$ is irreducible.
\end{x}
\vspace{0.1cm}

\begin{x}{\small\bf CONVENTION} \ 
Given $(\Pi, (V(\Pi))$ in $\widehat{G}$, take 
\[
V(\Pi^*) 
\ = \ V(\Pi)^*, 
\quad 
\Pi^*(x) 
\ = \ 
\Pi(x^{-1})^\Tee.
\]
\end{x}
\vspace{0.1cm}

\begin{x}{\small\bf NOTATION} \ 
Given finite dimensional vector spaces $V_1$, $V_2$ over $\Cx$, let $V_1 \otimes V_2$ be their tensor product.
\end{x}
\vspace{0.1cm}

Let \mG be a finite group. 
\\[-.2cm]

\begin{x}{\small\bf DEFINITION} \ 
Suppose that 
$\pi_1:G \ra \GL(V_1)$, 
$\pi_2:G \ra \GL(V_2)$ 
are representations $-$then their 
\un{tensor product}
\index{tensor product} 
is the representation $\pi_1 \otimes \pi_2 : G \ra \GL(V_1 \otimes V_2)$ defined by requiring that $\forall \ x \in G$,
\[
(\pi_1 \otimes \pi_2) (x) (v_1 \otimes v_2) 
\ = \ 
\pi_1(x) v_1 
\ \otimes \ 
\pi_2(x) v_2.
\]
\end{x}
\vspace{0.1cm}

Let $(\pi_1, V_1)$, $(\pi_2, V_2)$ be representations of \mG $-$then the prescription
\[
\pi_{1,2} (x) T
\ = \ 
\pi_2(x) T \pi_1(x^{-1})
\qquad (T \in \Hom(V_1,V_2))
\]
defines a representation $\pi_{1,2}$ of \mG on $\Hom(V_1,V_2)$.
\\[-.2cm]

\begin{x}{\small\bf RAPPEL} \ 
There is a canonical isomorphism
\[
\theta : V_2 \hsx \otimes \hsx V_1^*
\ \approx \ 
\Hom(V_1,V_2).
\]

[Send $v_2 \otimes v_1^*$ to the linear transformation
\[
T(v_2,v_1^*) : v_1 \ra v_1^* (v_1) v_2.]
\]
\end{x}
\vspace{0.1cm}

Consider 
\[
\pi_2(x) v_2 \hsx \otimes \hsx \pi_1^* (x) v_1^* \in V_2 \hsx \otimes \hsx V_1^*.
\]
Then the corresponding element of $\Hom(V_1,V_2)$ is the assignment
\allowdisplaybreaks
\begin{align*}
v_1 \ra (\pi_1^*(x) v_1^*) (v_1) \pi_2(x) v_2 \
&=\ 
v_1^* (\pi_1(x^{-1}) v_1) \pi_2(x) v_2 
\\[15pt]
&=\ 
\pi_2(x) (v_1^*(\pi_1(x^{-1})v_1))v_2
\\[15pt]
&=\ 
\pi_2(x) T(v_2, v_1^*) \pi_1(x^{-1}) v_1.
\end{align*}
\vspace{0.1cm}

\begin{x}{\small\bf LEMMA} \ 
$\pi_2 \otimes \pi_1^*$ is equivalent to $\pi_{1,2}$.
\vspace{0.2cm}

[The isomorphism $\theta$ intertwines $\pi_2 \otimes \pi_1^*$ and $\pi_{1,2} : \forall \ x \in G$, 
\[
\theta \circ (\pi_2(x) \otimes \pi_1^*(x)) 
\ = \ 
\pi_{1,2}(x) \circ \theta.]
\]
\end{x}
\vspace{0.1cm}

Let $G_1$, $G_2$ be finite groups.
\\[-.2cm]

\begin{x}{\small\bf DEFINITION} \ 
Suppose that $\pi_1: G_1 \ra \GL(V_1)$, $\pi_2: G_2 \ra \GL(V_2)$ are representations $-$then their 
\un{outer tensor product}
\index{outer tensor product} 
is the representation 
$\pi_1 \hsx \un{\otimes} \hsx \pi_2 : G_1 \times G_2 \ra \GL(V_1 \otimes V_2)$ defined by requiring that 
$\forall \ x_1 \in G_1$, $\forall \ x_2 \in G_2$, 
\[
(\pi_1 \hsx \un{\otimes} \hsx \pi_2)(x_1, x_2) 
\ = \ 
\pi_1(x_1) \hsx \otimes \hsx \pi_2(x_2).
\]
\end{x}
\vspace{0.1cm}

\begin{x}{\small\bf \un{N.B.}} \ 
If $G_1 = G_2 = G$, then the restriction of the outer tensor product $\pi_1 \un{\otimes} \pi_2$ to the diagonal subgroup
\[
\{(x, x) : x \in G\}
\]
of $G \times G$ is the tensor product $\pi_1 \otimes \pi_2$.
\end{x}
\vspace{0.1cm}

\begin{x}{\small\bf REMARK} \ 
Take $G_1 = G_2 = G$ and define a representation $\pi_{1,2}$ of $G \times G$ on $\Hom(V_1,V_2)$ via the prescription
\[
\pi_{1,2}(x,y) T 
\ = \ 
\pi_2(x) T \pi_1(y^{-1})
\qquad (T \in \Hom(V_1, V_2)).
\]
\end{x}
Then $\pi_2 \un{\otimes} \pi_1^*$ is equivalent to $\pi_{1,2}$.
\vspace{0.1cm}

\begin{x}{\small\bf LEMMA} \ 
If $\pi_1$ and $\pi_2$ are irreducible, then $\pi_1 \hsx \un{\otimes} \hsx \pi_2$ is irreducible.
\vspace{0.2cm}

[To begin with,
\[
C(G_1 \times G_2) 
\ \approx \ 
C(G_1) \otimes C(G_2)
\]
and
\[
\Hom(V_1 \otimes V_2) 
\ \approx \ 
\Hom(V_1) \otimes \Hom(V_2).
\]
Now make the passage
\[
\begin{cases}
\ \pi_1 \ra \rho_1\\[4pt]
\ \pi_2 \ra \rho_2
\end{cases}
.
\]
Then
\[
\begin{cases}
\  \rho_1(C(G_1)) \hsx = \hsx \Hom(V_1)\\[4pt]
\  \rho_2(C(G_2)) \hsx = \hsx \Hom(V_2)
\end{cases}
\qquad \text{(cf. I, \S2, \#2)}.]
\]
\end{x}
\vspace{0.1cm}

Conversely:
\\[-.5cm]

\begin{x}{\small\bf THEOREM} \ 
Every irreducible representation of $G_1 \times G_2$ is equivalent to an outer tensor product 
$\pi_1 \hsx \un{\otimes} \hsx \pi_2$.
\end{x}
\vspace{0.1cm}

\begin{x}{\small\bf SCHOLIUM} \ 
\[
\reallywidehat{G_1 \times G_2} 
\ \approx \
\widehat{G_1} \hsx \times \hsx \widehat{G_2}.
\]
\end{x}
\vspace{0.1cm}


\chapter{
$\boldsymbol{\S}$\textbf{3}.\quad  FOURIER TRANSFORMS}
\setlength\parindent{2em}
\setcounter{theoremn}{0}
\renewcommand{\thepage}{A II \S3-\arabic{page}}

\qquad Let \mG be a finite group.
\\[-.2cm]

\begin{x}{\small\bf DEFINITION} \ 
Given $f \in C(G)$, its 
\un{Fourier transform}\ 
\index{Fourier transform} 
$\widehat{f}$
\index{$\widehat{f}$} 
is that element of 
\[
\bigoplus\limits_{\Pi \in \widehat{G}} \ \Hom(V(\Pi))
\]
whose $\Pi$-component is
\[
\widehat{f} (\Pi) \ \equiv \ \sum\limits_{x \in G} \  f(x) \hsy \Pi(x) \qquad (= \Pi(f)).
\]

E.g.: \ $\forall \ x \in G$, 
\[
\widehat{\delta}_x (\Pi) \ = \ \Pi(x).
\]
\end{x}
\vspace{0.1cm}

\begin{x}{\small\bf LEMMA} \ 
$\forall \ f_1, \ f_2 \ \in C(G)$, 
\[
\reallywidehat{f_1 \hsx * \hsx f_2} \hsy (\Pi) \ = \ f_1(\Pi) \hsx f_2(\Pi).
\]
\end{x}
\vspace{0.1cm}

\begin{x}{\small\bf EXAMPLE}  \ 
$\forall \ x \in G$, 
\[
\begin{cases}
\ \widehat{L(x) f} \hsy (\Pi) \ = \ \widehat{\delta_x \hsx * \hsx f} \hsy (\Pi) \hspace{0.5cm} = \ \Pi(x) \hsy \widehat{f}(\Pi)\\[0.25cm]
\ \widehat{R(x) f} \hsy (\Pi) \ = \ \widehat{f \hsx * \hsx \delta_{x^{-1}}}\hsy (\Pi) \ = \ \widehat{f} (\Pi) \hsy \Pi(x^{-1})
\end{cases}
.
\]
\end{x}
\vspace{0.1cm}

\begin{x}{\small\bf THEOREM} \ 
The Fourier transform 
\[
\wedge : C(G) \ra \bigoplus\limits_{\Pi \in \widehat{G}} \ \Hom(V(\Pi))
\]
is an algebra isomorphism.
\end{x}
\vspace{0.1cm}


\begin{x}{\small\bf APPLICATION} \ 
\[
\abs{G} \ = \ \sum\limits_{\Pi \in \widehat{G}} \ d_\Pi^2.
\]

[In fact, 
\[
\dim C(G) \ = \ \abs{G} \ \quad \text{and} \quad \dim \Hom(V(\Pi)) \ = \ d_\Pi^2.]
\]
\end{x}
\vspace{0.25cm}

As it stands, $C(G)$ is a unital associative algebra over $\Cx$.  
But more is true: $C(G)$ is a $*$-algebra, i.e., admits a conjugate linear antiautomorphism $f \ra f^*$ given by 
$f^*(x) = \ov{f(x^{-1})} (\ x \in G)$.

Each $T \in \Hom(V(\Pi))$ has an adjoint $T^*$ per 
$\langle \ , \ \ranglesubPi$ : $\forall \ v_1, \ v_2 \ \in V(\Pi)$, 
\[
\langle T v_1, v_2\ranglesubPi \ = \ \langle v_1, T^* v_2\ranglesubPi.
\]
Therefore
\[
\bigoplus\limits_{\Pi \in \widehat{G}} \ \Hom(V(\Pi))
\]
admits a conjugate linear antiautomorphism by using the arrow $T \ra T^*$ on each summand.
\\

\begin{x}{\small\bf \un{N.B.}} \ 
It can and will be assumed that
\[
\begin{cases}
\ V(\Pi^*) \ = \ V(\Pi)\\[0.25cm]
\ \Pi^*(x) \ = \ \Pi(x^{-1})
\end{cases}
\qquad (\text{cf.} \ \S2, \ \#5),
\]
hence in terms of adjoints
\[
\Pi(x)^* \ = \ \Pi(x)^{-1} \ = \ \Pi(x^{-1}) \ = \ \Pi^*(x).
\]
\end{x}
\vspace{0.25cm}


\begin{x}{\small\bf LEMMA} \ 
The Fourier transform 
\[
\wedge : C(G) \ra \bigoplus\limits_{\Pi \in \widehat{G}} \ \Hom(V(\Pi))
\]
preserves the $*$-operations: $\forall \ f \in C(G)$, 
\[
f^* \ = \ (\widehat{f}\hsx)^*.
\]
\end{x}
\vspace{0.25cm}

\begin{x}{\small\bf INVERSION FORMULA} \ 
Given $f \in C(G)$, $\forall \ x \in G$, 
\[
f(x) \ = \ \frac{1}{\abs{G}} \ \sum\limits_{\Pi \in \widehat{G}} \ d_\Pi  \hsx \tr(\Pi(x^{-1})\widehat{f}(\Pi)).
\]
In particular: 
\[
f(e) \ = \ \frac{1}{\abs{G}} \ \sum\limits_{\Pi \in \widehat{G}} \ d_\Pi \hsx \tr(\widehat{f}(\Pi)).
\]
\end{x}
\vspace{0.25cm}

\begin{x}{\small\bf PARSEVAL IDENTITY} \ 
Given $\ f_1, \ f_2 \ \in C(G)$,
\[
\sum\limits_{x \in G} \ f_1(x) f_2(x^{-1})
 \ = \ 
 \frac{1}{\abs{G}} \ \sum\limits_{\Pi \in \widehat{G}} \ d_\Pi \hsx \tr(\widehat{f_1}(\Pi) \hsx \widehat{f_2}(\Pi)).
\]
\vspace{0.1cm}

PROOF \ 
Put $f = f_1 * f_2$ $-$then 
\[
f(e) \ = \ \sum\limits_{x \in G} \ f_1(x) f_2(x^{-1}).
\]
On the other hand, 
\begin{align*}
 \frac{1}{\abs{G}} \ 
\sum\limits_{\Pi \in \widehat{G}} \ d_\Pi \hsx\tr(\widehat{f_1}(\Pi) \hsx \widehat{f_2}(\Pi)) \ 
&=\ 
 \frac{1}{\abs{G}} \ 
\sum\limits_{\Pi \in \widehat{G}} \ d_\Pi \hsx \tr(\reallywidehat{f_1 * f_2} \hsx (\Pi))
\\[11pt]
&=\ 
(f_1 * f_2) (e) 
\\[11pt] 
&=\ f(e).
\end{align*}
\end{x}
\vspace{0.25cm}

\begin{x}{\small\bf COMPLETENESS PRINCIPLE} \ 
If $f \in C(G)$ and if $\widehat{f}(\Pi) = 0$ for all $\Pi$, then $f = 0$.
\end{x}


\chapter{
$\boldsymbol{\S}$\textbf{4}.\quad  CLASS FUNCTIONS}
\setlength\parindent{2em}
\setcounter{theoremn}{0}
\renewcommand{\thepage}{A II \S4-\arabic{page}}

$\text{ }$\\[-1.5cm]

Let \mG be a finite group.
\\[-.2cm]

\begin{x}{\small\bf DEFINITION} \ 
Let $(\pi,V)$ be a representation of \mG $-$then its 
\un{character}
\index{character} 
is the function 
\[
\chis{\pi} : G \ra \Cx
\]
given by the prescription
\[
\chis{\pi}(x) \ = \ \tr(\pi(x)) \qquad (x \in G).
\]
\end{x}
\vspace{0.1cm}

\begin{x}{\small\bf \un{N.B.}} \ 
It is clear that characters are class functions and that equivalent representations have equal characters.
\end{x}
\vspace{0.1cm}

\begin{x}{\small\bf LEMMA} \ 
$\forall \ x \in G$, 
\[
\chisubpi (x^{-1}) \ = \ \ovhb{\chisubpi (x)}.
\]
\end{x}
\vspace{0.1cm}

\begin{x}{\small\bf \un{N.B.}} \ 
\[
\chisubpiStar \ = \ \ovha{\chisubpi}.
\]
\end{x}
\vspace{0.1cm}

\begin{x}{\small\bf LEMMA} \ 
Let
\[
\begin{cases}
\ \pi_1: G \ra \GL(V_1)\\[0.25cm]
\ \pi_2: G \ra \GL(V_2)
\end{cases}
\]
be representations of \mG $-$then the character of 
\[
(\pi_1 \hsx \otimes \hsx \pi_2, V_1 \hsx \otimes \hsx V_2)
\]
is $\chiss{\pi_1} \hsx \chiss{\pi_2}$.
\\[-.25cm]

[For the record, the character of 
\[
(\pi_1 \hsx \oplus \hsx \pi_2, V_1 \hsx \oplus \hsx V_2)
\]
is 
$\chiss{\pi_1} + \chiss{\pi_2}$,
implying thereby that a nonnegative integral linear combination of characters is again a character.]
\end{x}

\begin{x}{\small\bf EXAMPLE} \ 
$\pi_{1,2}$ is equivalent to $\pi_2 \hsx \otimes \hsx \pi_1^*$ (cf. \S2, \#9), hence 
\[
\chi\raisebox{-.15cm}{$\scaleto{\pi_{1,2}}{8pt}$}
\ = \ 
\chi\raisebox{-.15cm}{$\scaleto{\pi_2 \hsx \otimes \hsx \pi_1^*}{10pt}$} 
\ = \ 
\chi\raisebox{-.15cm}{$\scaleto{\pi_2}{6pt}$}
\hsx
\chi\raisebox{-.15cm}{$\scaleto{\pi_1^*}{10pt}$}
\ = \
\chi\raisebox{-.15cm}{$\scaleto{\pi_2}{6pt}$}
\hsx 
\ovha{\chi\raisebox{-.15cm}{$\scaleto{\pi_1}{6pt}$}}.
\]
\end{x}
\vspace{0.1cm}

\begin{x}{\small\bf DEFINITION} \ 
The character of an irreducible representation is called an 
\un{irreducible character}
\index{irreducible character}.
\\[-.25cm]

[Note: \ The zero function (i.e., the additive identity of $C(G)$) is a character but it is not an irreducible character 
(cf. \S1, \#21).]
\end{x}
\vspace{0.1cm}

\begin{x}{\small\bf \un{N.B.}} \ 
The irreducible characters are thus the $\chisubPi$ $(\Pi \in \widehat{G})$.
\end{x}
\vspace{0.1cm}

\begin{x}{\small\bf FIRST ORTHOGONALITY RELATION} \ \  
Let $\Pi_i$, $\Pi_j \in \ \widehat{G}$ $-$then
\[
\frac{1}{\abs{G}} \  
\sum\limits_{x \in G} \  
\chi\raisebox{-.1cm}{$\scaleto{_i}{5.5pt}$}(x) 
\chi\raisebox{-.1cm}{$\scaleto{_j}{5.5pt}$}(x^{-1}) 
\ = \ 
\delta_{i j},
\]
where for short
\[
\chi\raisebox{-.1cm}{$\scaleto{_i}{5.5pt}$}
\ = \ 
\chi\raisebox{-.17cm}{$\scaleto{\Pi_i}{7pt}$}, 
\qquad 
\chi\raisebox{-.1cm}{$\scaleto{_j}{5.5pt}$} 
\ = \ 
\chi\raisebox{-.15cm}{$\scaleto{\Pi_j}{8pt}$}.
\]
\end{x}
\vspace{0.1cm}

\begin{x}{\small\bf NOTATION} \ 
Given $x \in G$, write $C(x)$ for its conjugacy class and $G_x$ for its centralizer.
\end{x}
\vspace{0.1cm}

\begin{x}{\small\bf RAPPEL} \ 
The number of conjugates of $x$ in \mG is $[G:G_x]$, i.e., 
\[
\abs{C(x)}
\ = \ 
[G:G_x].
\]


[Note: \ 
The \un{class equation} for \mG is the relation
\[
\abs{G} 
\ = \ 
\sum\limits_i \ [G : G_{x_i}],
\]
one $x_i$ having been chosen from each conjugacy class.]
\end{x}
\vspace{0.1cm}

\begin{x}{\small\bf SECOND ORTHOGONALITY RELATION} \ \ 
Let $x_1$, $x_2 \in \ G$ $-$then
\[
\sum\limits_{\Pi \in \widehat{G}} \ \chisubPi(x_1) \chisubPi(x_2^{-1}) \ = \ 
\begin{cases}
\ \abs{G_x} \ \text{if} \quad x = x_1 = x_2 \\[.25cm]
\ \ 0 \hspace{.6cm}  \text{if} \quad C(x_1) \neq C(x_2)
\end{cases}
.
\]

[Note: \ 
\[
\abs{G_x} 
\ = \ 
\frac{\abs{G}}{[G:G_x]} 
\ = \ 
\frac{\abs{G}}{\abs{C(x)}}
.]
\]
\end{x}
\vspace{0.25cm}

\begin{x}{\small\bf NOTATION} \ 
Given $f$, $g \in \ C(G)$, put
\[
\langle f, g \ranglesubG 
\ = \ 
\frac{1}{\abs{G}} \ 
\sum\limits_{x \in G} \ 
f(x) \ov{g(x)},
\]
the canonical inner product on $C(G)$.
\end{x}
\vspace{0.1cm}

\begin{x}{\small\bf EXAMPLE} \ 
$\forall \ \Pi_1$, $\Pi_2 \in \widehat{G}$,
\[
\langle \chisubPiOne,  \chisubPiTwo \ranglesubG 
\ = \ 
\begin{cases}
\ 1 \ \text{if} \quad \Pi_1 = \Pi_2 \\[4pt]
\ 0 \ \text{if} \quad \Pi_1 \neq \Pi_2 
\end{cases}
.
\]
\end{x}
\vspace{0.1cm}

\begin{x}{\small\bf SCHOLIUM} \ 
The irreducible characters form an orthonormal set, thus are linearly independent (cf. I, $\S3$, $\#9$).
\end{x}
\vspace{0.1cm}


Recall now that the Fourier transform
\[
\wedge: C(G) \ra 
\bigoplus\limits_{\Pi \in \widehat{G}} \hsx \Hom(V(\Pi))
\]
is an algebra isomorphism.  
Since the center of each $\Hom(V(\Pi))$ consists of scalar multiples of the identity operator, it follows that an $f \in C(G)$ 
is a class function iff $\forall \ \Pi \in \widehat{G}$,
\[
\widehat{f} \hsx (\Pi) 
\ = \ 
C_\Pi \hsx\hsx \id_{V(\Pi)} \qquad (C_\Pi \in \Cx).
\]
\\[-.5cm]

\begin{x}{\small\bf INVERSION FORMULA} \ 
Given a class function $f \in C(G)$, $\forall \ x \in G$, 
\[
f(x) 
\ = \ 
\sum\limits_{\Pi \in \widehat{G}} \ 
\langle f, \ovha{\chisubPi} \ranglesubG  \hsx \ovhb{\chisubPi(x)}.
\]

PROOF \ 
\allowdisplaybreaks
\begin{align*}
f(x) \ 
&=\ 
\frac{1}{\abs{G}}
\
\sum\limits_{\Pi \in \widehat{G}} 
\
d_\Pi \hsx \tr(\Pi(x^{-1}) \widehat{f} \hsx (\Pi))
\\[11pt]
&=\ 
\frac{1}{\abs{G}}
\
\sum\limits_{\Pi \in \widehat{G}} \ 
d_\Pi \hsx C_\Pi \hsx \chisubPi(x^{-1})
\\[11pt]
&=\ 
\frac{1}{\abs{G}}
\
\sum\limits_{\Pi \in \widehat{G}}  \
d_\Pi \hsx C_\Pi \hsx \ovhb{\chisubPi(x)}.
\end{align*}
Fix $\Pi_0 \in \widehat{G}$ $-$then
\allowdisplaybreaks
\begin{align*}
\langle f, \ovha{\chisubPiZero} \ranglesubG \
&=\ 
\frac{1}{\abs{G}}
\
\sum\limits_{\Pi \in \widehat{G}} \ 
d_\Pi \hsx C_\Pi \hsx \langle \ovha{\chisubPi}, \ovha{\chisubPiZero} \ranglesubG
\\[11pt]
&=\ 
\frac{1}{\abs{G}}
\
\sum\limits_{\Pi \in \widehat{G}} \ 
d_\Pi \hsx C_\Pi  \hsx \ovhb{\langle \chisubPi, \chisubPiZero \ranglesubG}
\\[11pt]
&=\ 
\frac{1}{\abs{G}}
\
d_{\Pi_0} \hsx  C_{\Pi_0}.
\end{align*}
\end{x}
\vspace{0.1cm}

\begin{x}{\small\bf \un{N.B.}} \ 
$\forall \ x \in G$, 
\allowdisplaybreaks
\begin{align*}
f(x) \ 
&=\ 
\ovhb{\ovha{f(x)}} 
\\[11pt]
&=\ 
\ovhb
{
\sum\limits_{\Pi \in \widehat{G}} \ 
\langle \ovhb{f}, \ovha{\chisubPi} \ranglesubG  \  \ovhb{\chisubPi(x)} \
}
\\[11pt]
&=\ 
\sum\limits_{\Pi \in \widehat{G}} \ 
\ovhb{\langle \ovhb{f}, \ov{\chisubPi} \ranglesubG}  \ \ovhb{\ovha{\chisubPi(x)}}
\\[11pt]
&=\ 
\sum\limits_{\Pi \in \widehat{G}} \ 
\langle f, \chisubPi \ranglesubG  \ \chisubPi(x).
\end{align*}
\\[-.75cm]
\end{x}

The preceding discussion makes it clear that a class function $f$ is a character iff 
$\langle f, \Pi \ranglesubG$ is a nonnegative integer for all $\Pi \in \ \widehat{G}$.
\\[-.2cm]

\begin{x}{\small\bf NOTATION} \ 
$\CON(G)$ is the set of conjugacy classes of \mG.
\end{x}
\vspace{0.1cm}

\begin{x}{\small\bf SCHOLIUM} \ 
The dimension of the space of class functions is equal to $\abs{\CON(G)}$ or still, is equal to $\widehat{\abs{G}}$.
\end{x}
\vspace{0.1cm}

\begin{x}{\small\bf NOTATION}  
Given $C \in  \CON(G)$, let $\chisubC$ be the characteristic function of \mC:
\[
\chisubC (x) \ = \ 
\begin{cases}
\ 1 \quad \text{if} \quad x \in C\\[4pt]
\ 0 \quad \text{if} \quad x \notin C
\end{cases}
.
\]
\end{x}
\vspace{0.1cm}

\begin{x}{\small\bf LEMMA} \ 
\[
\chisubC
\ = \ 
\sum\limits_{y \in C} \ \delta_y \hsx .
\]
\end{x}
\vspace{0.1cm}

\begin{x}{\small\bf \un{N.B.}} \ 
The $\chisubC$ $(C  \in \CON(G))$ are a basis for the class functions on \mG 
(as are the $\chisubPi$ $(\Pi \in \ \widehat{G})$).
\end{x}
\vspace{0.1cm}


\begin{x}{\small\bf LEMMA} \ 
Let $C_1, C_2, \ldots$ be the elements of $\CON(G)$ $-$then there are nonnegative integers $m_{i, j, k}$ such that 
\[
\chisubCi \hsx \chisubCj
\ = \ 
\sum\limits_k \ 
m_{i, j, k} \hsx \chisubCk.
\]

[Note: 
Fixing an $x_k \ \in C_k$, qualitatively $m_{i, j, k}$ is the number of ordered pairs $(x,y)$ with 
$x \in C_i$, $y \in C_j$ and $x y = x_k$ while quantitatively
\[
m_{i, j, k} 
\ = \ 
\frac
{
\abs{C_i}\abs{C_j}
}
{
\abs{G}
} 
\ 
\sum\limits_{\Pi \in \widehat{G}} \
\frac
{
\chisubPi (x_i) 
\chisubPi(x_j) 
\ovhb
{
\chisubPi(x_k)
}
}
{
d_\Pi
}.]
\]
\end{x}
\vspace{0.1cm}

\begin{x}{\small\bf NOTATION} \ 
Given $\Pi \in \widehat{G}$, put
\[
\esubPi
\ = \ 
\wedge^{-1} (E_\Pi) \qquad (E_\Pi \in \ \Hom(V(\Pi)) \qquad \text{(cf. I, $\S4$)).}
\]
\end{x}
\vspace{0.1cm}

\begin{x}{\small\bf LEMMA} \ 
\[
\esubPi
\ = \ 
\frac{d_\Pi}{\abs{G}}  \
\sum\limits_{y \in G} \ 
\chisubPi(y^{-1}) \delta_y \hsx .
\]

[Note: 
In brief, 
\[
\esubPi
\ = \ 
\frac{d_\Pi}{\abs{G}}  \hsx 
\chisubPiStar \hsx .]
\]
\end{x}
\vspace{0.1cm}

\begin{x}{\small\bf LEMMA} \ 
\[
\esubPiOne * \esubPiTwo 
\ = \ 
\begin{cases}
\ \esubPi \hspace{.5cm} \text{if} \quad \Pi = \Pi_1 = \Pi_2 \\[4pt]
\ 0 \hspace{.72cm} \text{if} \quad \Pi_1 \neq \Pi_2
\end{cases}
.
\]
\end{x}
\vspace{0.1cm}

\begin{x}{\small\bf LEMMA} \ 
\[
\delta_e 
\ = \ 
\sum\limits_{\Pi \in \widehat{G}} \ \esubPi.
\]
\end{x}
\vspace{0.1cm}

\chapter{
$\boldsymbol{\S}$\textbf{5}.\quad  DECOMPOSITION THEORY}
\setlength\parindent{2em}
\setcounter{theoremn}{0}
\renewcommand{\thepage}{A II \S5-\arabic{page}}

\ \indent 
Let \mG be a finite group.
\\[-.2cm]

\begin{x}{\small\bf CONSTRUCTION} \ 
Suppose that \mG operates on a finite set \mS, hence for each $x \in G$ there is given a bijection 
$s \ra x \cdot s$ of \mS satisfying the identities
\[
e \cdot s 
\ = \ 
s, 
\quad 
x \cdot (y \cdot s) 
\ = \ 
(x y) \cdot s.
\]
Let $V = \{f : S \ra \Cx \}$ and define a representation $\pi:G \ra \GL(V)$ by 
\[
\pi(x) f(s) 
\ = \ 
f(x^{-1} \cdot s).
\]
Then 
\[
\chisubpi(x) 
\ = \ 
\abs{\{s \in S: x \cdot s = s\}}.
\]
\\[-1.05cm]
\end{x}

\begin{x}{\small\bf EXAMPLE} \ 
Take $S = G$ and write $x \cdot y = xy$ $-$then the role of \mV is played by $C(G)$ and the role of $\pi$ is played by \mL (the left translation representation of \mG (cf. \S1, \#12)), hence
\[
\chisubL(x) 
\ = \ 
\abs{\{y \in G: xy = y\}},
\]
which is $\abs{G}$ if $x = e$ and is 0 otherwise.
\\[-.25cm]
\end{x}

\begin{x}{\small\bf EXAMPLE} \ 
Take $S = G$ but replace \mG by $G \times G$, the action being 
$(x_1, x_2) \cdot y = x_1 y x_2^{-1}$ $-$then the associated representation is $\pi_{L,R}$ 
(cf. \S1, \#14) and
\allowdisplaybreaks 
\begin{align*}
\chisubpiLR (x_1, x_2) \ 
&=\ 
\abs{\{y \in G : x_1 y x_2^{-1} = y\}}
\\[.25cm]
&=\ 
\abs{\{y \in G : x_1 = y x_2 y^{-1}\}}
\\[.25cm]
&=\ 
\abs{G_{x_1}}
\end{align*}
if $x_1$ and $x_2$ are conjugate and is 0 otherwise.
\end{x}
\vspace{0.1cm}


\begin{x}{\small\bf DEFINITION} \ 
Let $(\pi,V)$ be a representation of \mG $-$then by complete reducibility, there is a direct sum decomposition
\[
\pi
\ = \ 
\bigoplus\limits_{\Pi \in \widehat{G}} \hsx 
m(\Pi,\pi) \Pi, 
\]
the nonnegative integer $m(\Pi,\pi)$ being the 
\un{multiplicity}
\index{multiplicity} 
of $\Pi$ in $\pi$.
\end{x}
\vspace{0.1cm}

\begin{x}{\small\bf LEMMA} \ 
$\forall \ \Pi \in \widehat{G}$, 
\[
m(\Pi,\pi)
\ = \ 
\langle \chisubPi, \chisubpi \ranglesubG.
\]
\end{x}
\vspace{0.1cm}

\begin{x}{\small\bf \un{N.B.}} \ 
\[
\dim I_G(\Pi,\pi) 
\ = \ 
m(\Pi,\pi).
\]
\end{x}
\vspace{0.1cm}

\begin{x}{\small\bf REMARK} \ 
The operator
\[
P_\Pi 
\ = \ 
\frac{d_\Pi}{\abs{G}} \ 
\sum\limits_{x \in G} \ 
\ov{\chisubPi(x)} \hsx\hsx \pi(x)
\]
is the projection onto the $\Pi$-isotypic subspace of \mV.
\end{x}
\vspace{0.1cm}

\begin{x}{\small\bf THEOREM} \ 
Each $\Pi \in \widehat{G}$ is contained in \mL with multiplicity $d_\Pi$.
\\[-.2cm]

PROOF \ 
In fact, 
\allowdisplaybreaks 
\begin{align*}
m(\Pi, L ) \ 
&=\ 
\langle \chisubPi, \chisubL \ranglesubG 
\\[11pt]
&=\ 
\frac{1}{\abs{G}} \ 
\sum\limits_{x \in G} \
\chisubPi(x) \hsx \ov{\chisubL(x)} 
\\[11pt]
&=\ 
\frac{1}{\abs{G}}\  \chisubPi(e) \abs{G}
\\[11pt]
&=\ 
\chisubPi(e) 
\\[11pt]
&=\ 
d_\Pi.
\end{align*}
\end{x}
\vspace{0.1cm}

\begin{x}{\small\bf \un{N.B.}} \ 
It is a corollary that
\[
\abs{G} 
\ = \ 
\sum\limits_{\Pi \in \widehat{G}} \
d_\Pi^{\hsy 2} \qquad (\text{cf. \S3, \#5}).
\]
\end{x}
\vspace{0.1cm}

\begin{x}{\small\bf LEMMA} \ 
Let $(\pi_1, V_1)$, $(\pi_2, V_2)$ be representations of \mG.  
Assume: 
$\chisubpiOne = \chisubpiTwo$ $-$then $(\pi_1, V_1) \hsx \approx \hsx (\pi_2, V_2)$.
\\[-.2cm]

PROOF \ 
$\forall \ \Pi \in \widehat{G}$, 
\[
\langle \chisubPi, \chisubpiOne \ranglesubG 
\ = \ 
\langle \chisubPi, \chisubpiTwo \ranglesubG 
\]
or still, $\forall \ \Pi \in \widehat{G}$, 
\[
m(\Pi, \pi_1) 
\ = \ 
m(\Pi, \pi_2),
\]
from which the assertion.
\end{x}
\vspace{0.1cm}

\begin{x}{\small\bf IRREDUCIBILITY CRITERION} \ 
A representation $\pi:G \ra \GL(V)$ is irreducible iff 
$\langle \chisubpi, \chisubpi \ranglesubG = 1$.
\\[-.2cm]

PROOF \ 
The necessity is implied by the first orthogonality relation and the sufficiency follows upon noting that
\[
\langle \chisubpi, \chisubpi \ranglesubG 
\ = \ 
\sum\limits_{\Pi \in \widehat{G}} \
m(\Pi, \pi)^2.
\]
\\[-.5cm]
\end{x}

Let $G_1$, $G_2$ be finite groups and let
\[
\begin{cases}
\ \Pi_1 : G_1 \ra \GL(V_1) \\
\ \Pi_2 : G_2 \ra \GL(V_2)
\end{cases}
\]
be irreducible representations of $G_1$, $G_2$ $-$then the character 
$\chisubPiOneuXTwo$ of 
\[
(\Pi_1 \hsx \un{\otimes} \hsx\hsx \Pi_2, V_1 \hsx \otimes \hsx V_2)
\]
is the function
\[
(x_1, x_2) 
\ra 
\chisubPiOne (x_1) \chisubPiTwo (x_2) \qquad (x_1 \in G_1, \ x_2 \in G_2).
\]
\vspace{0.1cm}

\begin{x}{\small\bf LEMMA} \ 
$\Pi_1 \hsx \un{\otimes} \hsx \Pi_2$ is irreducible (cf. \S2, \#13).
\\[-.2cm]

PROOF \ 
It is a question of applying the irreducibility criterion.  
Thus
\allowdisplaybreaks 
\begin{align*}
\langle \chisubPiOneuXTwo, \chisubPiOneuXTwo \ranglesubGxG \ 
&=\ 
\frac{1}{\abs{G_1 \times G_2}} \hsx
\sum\limits_{(x_1, x_2) \in G_1 \times G_2} \
\chisubPiOne(x_1) \hsx \chisubPiTwo(x_2)
\hsx
\ov{\chisubPiOne(x_1) \hsx \chisubPiTwo(x_2)}
\\[15pt]
&=\ 
\frac{1}{\abs{G_1}} \hsx
\sum\limits_{x_1 \in G_1} \
\chisubPiOne(x_1)
\ov{\chisubPiOne(x_1)}
\hsx \cdot \hsx
\frac{1}{\abs{G_2}} \hsx
\sum\limits_{x_2 \in G_2} \
\chisubPiTwo(x_2)
\hsx
\ov{\chisubPiTwo(x_2)}
\\[15pt]
&=\ 
\langle \chisubPiOne, \chisubPiTwo \ranglesubGOne 
\hsx \cdot \hsx
\langle \chisubPiOne, \chisubPiTwo \ranglesubGTwo 
\\[15pt]
&=\ 
1.
\end{align*}

\end{x}
\vspace{0.1cm}

\begin{x}{\small\bf REMARK} \ 
The cardinality of $\reallywidehat{G_1 \times G_2}$ is $\abs{\CON(G_1 \times G_2)}$ 
(cf. \S4, \#19).
But 
\[
\abs{\CON(G_1 \times G_2)} 
\ = \ 
\abs{\CON(G_1)}\abs{\CON(G_2)}
\]
and the preceding considerations produce
\[
\abs{\CON(G_1)}\abs{\CON(G_2)}
\]
pairwise distinct irreducible characters of $G_1 \times G_2$.  
Therefore every irreducible 
representation of $G_1 \times G_2$ is equivalent to an outer tensor product 
$\Pi_1 \hsx \un{\otimes} \hsx \Pi_2$, where 
$\Pi_1 \in \widehat{G}_1$, $\Pi_2 \in \widehat{G}_2$ 
(cf. \S2, \#14).
\end{x}
\vspace{0.1cm}


\chapter{
$\boldsymbol{\S}$\textbf{6}.\quad  INTEGRABILITY}
\setlength\parindent{2em}
\setcounter{theoremn}{0}
\renewcommand{\thepage}{A II \S6-\arabic{page}}


\begin{x}{\small\bf DEFINITION} \ 
An 
\un{algebraic integer}
\index{algebraic integer} 
is a complex number $\lambda$ which is a root of a polynomial of the form 
\[
x^n + a_{n-1} x^{n-1} + \cdots + a_0,
\]
where $a_i \in \Z$ $(0 \leq i \leq n - 1)$.
\\[-.5cm]

[Note: \ 
Equivalently, an algebraic integer is a complex number $\lambda$ which is a zero of 
\[
\det (A - X I)
\]
for some square matrix \mA with entries in $\Z$.]
\\[-.25cm]
\end{x}

\begin{x}{\small\bf \un{N.B.}} \ 
The rational algebraic integers are precisely the elements of $\Z$.
\\[-.25cm]
\end{x}

\begin{x}{\small\bf LEMMA} \ 
If $\mu$, $\nu$ are algebraic integers, then $\mu + \nu$ and $\mu \nu$ are also algebraic integers.
\\[-.25cm]
\end{x}

Therefore the set of algebraic integers is a subring of $\Cx$.
\\[-.25cm]

\begin{x}{\small\bf EXAMPLE} \ 
Roots of unity are algebraic integers.
\\[-.25cm]
\end{x}

Let \mG be a finite group.
\\[-.25cm]

\begin{x}{\small\bf LEMMA} \ 
Let $(\pi,V)$ be a representation of \mG, $\chisubpi$ its character $-$then $\forall \ x \in G$, $\chisubpi(x)$ is an algebraic integer.  
\\[-.5cm]

[This is because $\chisubpi(x)$ is a finite sum of roots of unity.]
\\[-.25cm]
\end{x}

The center of $C(G)$ (i.e., the class functions) is a unital commutative associative algebra over $\Cx$, thus its irreducible representations are just 
homomorphisms into $\Cx$ and are indexed by the $\Pi \in \widehat{G}$, say $\omega_\Pi$ with 
\[
\omega_{\Pi_1} (e_{\Pi_2}) 
\ = \ 
\delta_{\Pi_1, \Pi_2}.
\]
\vspace{0.2cm}

[Note: \ 
The $e_\Pi$ $(\Pi \in \widehat{G})$ are a basis for the class functions on \mG.]
\\[-.25cm]

\begin{x}{\small\bf THEOREM} \ 
$\forall \ C \in \CON(G)$, $\omega_\Pi(C)$ is an algebraic integer.
\\[-.5cm]

PROOF \ 
In the notation of \S4, \#23, 
\[
\chisubCi \chisubCj 
\ = \ 
\sum\limits_k \ 
m_{i, j, k} \hsx \chisubCk,
\]
hence 
\allowdisplaybreaks
\begin{align*}
\omega_\Pi(\chisubCi) \hsx \omega_\Pi(\chisubCj) \ 
&=\ 
\omega_\Pi(\chisubCi \chisubCj) 
\\[15pt]
&=\ 
\sum\limits_k \hsx 
m_{i, j, k} \hsx \omega_\Pi(\chisubCk)
\end{align*}
\qquad\qquad $\implies$
\[
\sum\limits_k \ 
(m_{i, j, k} - \delta_{j k} \hsx \omega_\Pi(\chisubCi))\omega_\Pi(\chisubCk) 
\ = \ 
0.
\]
But this means that $\omega_\Pi(\chisubCi)$ is an eigenvalue of the matrix $A_i$ whose $(j,k)^{\text{th}}$ entry is 
$m_{i, j, k}$ or still, is a zero of 
\[
\det (A_i - X I),
\]
thus is an algebraic integer.
\\[-.25cm]
\end{x}

\begin{x}{\small\bf LEMMA} \ 
$\forall \ C \in \CON(G)$, 
\[
\omega_\Pi (\chisubC) 
\ = \ 
\frac{\abs{C}}{d_\Pi} \hsx \chisubPi(x) 
\qquad (x \in C).
\]
\\[-.2cm]

PROOF \ 
Owing to \S4, \#25, 
\[
e_\Pi 
\ = \ 
\frac{d_\Pi}{\abs{G}} \ 
\sum\limits_{y \in G}  \
\chisubPi (y^{-1})  \hsx \delta_y, 
\]
so 
\[
\frac{\abs{G}}{d_\Pi} \hsx \chisubPi(x) \hsx e_\Pi 
\ = \ 
\chisubPi(x)  \
\sum\limits_{y \in G}  \
\chisubPi (y^{-1}) \delta_y
\]
\qquad\qquad $\implies$
\allowdisplaybreaks
\begin{align*}
\sum\limits_{\Pi \in \widehat{G}} \ 
\frac{\abs{G}}{d_\Pi} \hsx \chisubPi(x) \hsx e_\Pi  \ 
&=\ 
\sum\limits_{\Pi \in \widehat{G}} \ 
\chisubPi(x) \
\sum\limits_{y \in G} \ 
\chisubPi(y^{-1})
\delta_y
\\[11pt]
&=\ 
\sum\limits_{y \in G} \ 
\bigg(
\sum\limits_{\Pi \in \widehat{G}} \ 
\chisubPi(x) \hsx \chisubPi(y^{-1})
\bigg)
\delta_y
\\[11pt]
&=\ 
\sum\limits_{y \in C} \
\abs{G_x} \hsx
\delta_y
\qquad \text{(cf. \S4, \#12)}
\\[11pt]
&=\ 
\abs{G_x} \
\sum\limits_{y \in C} \
\delta_y
\\[11pt]
&=\ 
\abs{G_x} \hsx \chisubC \qquad \text{(cf. \S4, \#21)}.
\end{align*}
Now fix $\Pi_0 \in \widehat{G}$ $-$then
\allowdisplaybreaks
\begin{align*}
\omega_{\Pi_0} (\chisubC) \ 
&=\ 
\omega_{\Pi_0}
\bigg(
\frac{1}{\abs{G_x}} \ 
\sum\limits_{\Pi \in \widehat{G}} \ 
\frac{\abs{G}}{d_\Pi} \hsx 
\chisubPi(x)  \hsx e_\Pi
\bigg)
\\[11pt]
&=\ 
\frac{\abs{G}}{\abs{G_x}} \ 
\sum\limits_{\Pi \in \widehat{G}} \ 
\frac{1}{d_\Pi}  \hsx 
\chisubPi(x)  \omega_{\Pi_0} (e_\Pi)
\\[11pt]
&=\ 
\frac{\abs{G}}{\abs{G_x}} \hsx 
\frac{\chisubPiZero(x)}{d_{\Pi_0}} \hsx 
\\[11pt]
&=\ 
\frac{\abs{C}}{d_{\Pi_0}}  \hsx 
\chisubPiZero(x).
\end{align*}
\end{x}


Consequently, $\forall \ C \in \CON(G)$, 
\[
\frac{\abs{C}}{d_\Pi} \hsx \chisubPi(x) \qquad (x \in C)
\]
is an algebraic integer.
\\

\begin{x}{\small\bf THEOREM} \ 
$\forall \ \Pi \in \widehat{G}$, 
\[
\frac{\abs{G}}{d_\Pi} \in \Z.
\]

PROOF \ 
In view of \S4, \#9,
\[
\abs{G}
\ = \ 
\sum\limits_{x \in G} \
\chisubPi(x) \chisubPi(x^{-1}).
\]
Given $C \in \CON(G)$, fix an $x_C \in C$ $-$then 
\[
\abs{G}
\ = \ 
\sum\limits_{C \in \CON(G)} \ 
\abs{C} \chisubPi(x_C) \chisubPi(x_C^{-1})
\]
\qquad $\implies$ 
\[
\frac{\abs{G}}{d_\Pi} 
\ = \ 
\sum\limits_{C \in \CON(G)} \ 
\bigg(
\frac{\abs{C}}{d_\Pi}  \hsx 
\chisubPi(x_C)\bigg) \chisubPi(x_C^{-1}),
\]
hence $\ds\frac{\abs{G}}{d_\Pi}$ is a rational algebraic integer, hence is an integer.
\\[-.25cm]
\end{x}

In other words, the $d_\Pi$ divide $\abs{G}$.
\\[-.2cm]

\begin{x}{\small\bf THEOREM} \ 
If \mA is an abelian normal subgroup of \mG, then the $d_\Pi$ divide $[G:A]$.
\\[-.25cm]
\end{x}

\begin{x}{\small\bf APPLICATION} \ 
Let $Z(G)$ be the center of \mG $-$then the $d_\Pi$ divide $[G:Z(G)]$.
\end{x}
\vspace{0.1cm}


\chapter{
$\boldsymbol{\S}$\textbf{7}.\quad  INDUCED CLASS FUNCTIONS}
\setlength\parindent{2em}
\setcounter{theoremn}{0}
\renewcommand{\thepage}{A II \S7-\arabic{page}}

\qquad Let \mG be a finite group, $\Gamma \subset G$ a subgroup.
\\[-.2cm]

\begin{x}{\small\bf NOTATION} \ 
$CL(G)$ is the subspace of $C(G)$ comprised of the class functions and $CL(\Gamma)$ is the subspace of $C(\Gamma)$ comprised of the class functions.
\end{x}
\vspace{0.1cm}

\begin{x}{\small\bf NOTATION} \ 
Extend a function $\phi \in C(\Gamma)$ to a function $\mathring{\phi} \in C(G)$ by writing
\[
\mathring{\phi} (x) \ = \ 
\begin{cases}
\ \phi(x) \ \ \text{if $x \in \Gamma$}\\
\ \ 0 \qquad \text{if $x \notin \Gamma$}
\end{cases}
.
\]
\end{x}
\vspace{0.1cm}

\begin{x}{\small\bf NOTATION} \ 
Given a class function $\phi \in CL(\Gamma)$, put
\allowdisplaybreaks
\begin{align*}
(i_{_{\Gamma \ra G}} \phi) (x) \ 
&=
\frac{1}{\abs{\Gamma}} \
\sum\limits_{y \in G} \ 
\mathring{\phi} (y x y^{-1}) 
\\[15pt]
&=
\frac{1}{\abs{\Gamma}} \
\sum\limits_{\substack{y \in G, \hsx y x y^{-1} \in \hsy \Gamma}} \ 
\phi (y x y^{-1}) .
\end{align*}
\end{x}
\vspace{0.1cm}

\begin{x}{\small\bf LEMMA} \ 
\[
i_{_{\Gamma \ra G}} \phi \in CL(G),
\]
the 
\un{induced}
\index{class function \\ induced} 
class function.
\end{x}
\vspace{0.1cm}

\begin{x}{\small\bf \un{N.B.}} \ 
Therefore
\[
i_{_{\Gamma \ra G}} : CL(\Gamma) \ra CL(G).
\]

[Note: \ 
\[
i_{_{\Gamma \ra G}} c\phi 
\ = \ 
c \hsx i_{_{\Gamma \ra G}} \phi \quad (c \in \Cx), 
\quad 
i_{_{\Gamma \ra G}} (\phi_1 + \phi_2) 
\ = \ 
i_{_{\Gamma \ra G}} \phi_1 + i_{_{\Gamma \ra G}} \phi_2
\]
but in general, 
\[
i_{_{\Gamma \ra G}} (\phi_1 \phi_2) 
\ \neq \ 
(i_{_{\Gamma \ra G}} \phi_1) (i_{_{\Gamma \ra G}} \phi_2).] 
\]
\end{x}
\vspace{0.1cm}

The arrow of restriction $C(G) \ra C(\Gamma)$ leads to a map
\[
r_{_{G \ra \Gamma}} : CL(G) \ra CL(\Gamma).
\]
And:
\\[-.2cm]

\begin{x}{\small\bf FROBENIUS RECIPROCITY} \ 
Let $\phi \in CL(\Gamma)$, $\psi \in CL(G)$ $-$then
\[
\langle
i_{_{\Gamma \ra G}} \phi, \psi \ranglesubG 
\ = \ 
\langle 
\phi, r_{_{G \ra \Gamma}} \psi
\ranglesubGamma. 
\]

PROOF \ 
\allowdisplaybreaks
\begin{align*}
\langle i_{_{\Gamma \ra G}} \phi, \psi \ranglesubG \ 
&=\ 
\frac{1}{\abs{G}} \ 
\sum\limits_{x \in G} \ 
( i_{_{\Gamma \ra G}} \phi) (x) \ov{\psi(x)}
\\[15pt]
&=\ 
\frac{1}{\abs{G}} \ 
\frac{1}{\abs{\Gamma}} \ 
\sum\limits_{x \in G} \ 
\sum\limits_{y \in G} \ 
\mathring{\phi}(y x y^{-1}) \hsx \ov{\psi(x)}
\\[15pt]
&=\ 
\frac{1}{\abs{\Gamma}} \ 
\frac{1}{\abs{G}} \ 
\sum\limits_{y \in G} \ 
\sum\limits_{x \in G} \ 
\mathring{\phi}(x) \hsx \ov{\psi(y^{-1} x y)}
\\[15pt]
&=\ 
\frac{1}{\abs{G}} \ 
\sum\limits_{y \in G} \ 
\frac{1}{\abs{\Gamma}} \ 
\sum\limits_{\gamma \in \Gamma} \ 
\phi(\gamma) \hsx \ov{\psi(\gamma)}
\\[15pt]
&=\ 
\frac{1}{\abs{G}} \ 
\sum\limits_{y \in G} \ 
\langle
\phi, r_{_{G \ra \Gamma}} \hsx \psi
\ranglesubGamma
\\[15pt]
&=\ 
\langle
\phi, r_{_{G \ra \Gamma}} \phi
\ranglesubGamma.
\end{align*}
\end{x}
\vspace{0.1cm}

\begin{x}{\small\bf APPLICATION} \ 
If $\phi$ is a character of $\Gamma$, then $i_{_{\Gamma \ra G}} \phi$ is a character of \mG.
\vspace{0.2cm}

[If $\chi$ is a character of \mG, then $r_{_{G \ra \Gamma}} \chi$ is a character of $\Gamma$, hence
\[
\langle \phi, r_{_{G \ra \Gamma}} \chisubPi \ranglesubGamma
\]
is a nonnegative integer for all $\Pi \in \widehat{G}$ or still, 
\[
\langle i_{_{\Gamma \ra G}} \phi, \chisubPi \ranglesubG
\]
is a nonnegative integer for all $\Pi \in \widehat{G}$ which implies that $i_{_{\Gamma \ra G}} \phi$ is a character of \mG 
(cf. \S4, \#17 ff.).
\end{x}
\vspace{0.1cm}

\begin{x}{\small\bf LEMMA} \ 
Let $\phi \in CL(\Gamma)$, $\psi \in CL(G)$ $-$then
\[
i_{_{\Gamma \ra G}} ((r_{_{G \ra \Gamma}} \psi)\phi) 
\ = \ 
\psi(i_{_{\Gamma \ra G}} \phi).
\]

PROOF \ 
From the definitions, 
\allowdisplaybreaks
\begin{align*}
i_{_{\Gamma \ra G}} (( r_{_{G \ra \Gamma}} \psi)\phi) (x)\ 
&=\ 
\frac{1}{\abs{\Gamma}} \ 
\sum\limits_{y \in G} \ 
\mathring{r_{_{G \ra \Gamma}}\psi} (y x y^{-1}) 
\mathring{\phi} (y x y^{-1})
\\[15pt]
&=\ 
\frac{1}{\abs{\Gamma}} \ 
\sum\limits_{y \in G} \ 
\psi(y x y^{-1}) 
\mathring{\phi} (y x y^{-1})
\\[15pt]
&=\ 
\frac{1}{\abs{\Gamma}} \ 
\sum\limits_{y \in G} \ 
\psi(x) 
\mathring{\phi} (y x y^{-1})
\\[15pt]
&=\ 
\psi(x) \hsx 
\frac{1}{\abs{\Gamma}} \ 
\sum\limits_{y \in G} \ 
\mathring{\phi} (y x y^{-1})
\\[15pt]
&=\
\psi(x) \hsx 
(i_{_{\Gamma \ra G}} \phi) (x).
\end{align*}
\end{x}
\vspace{0.1cm}

\begin{x}{\small\bf APPLICATION} \ 
The image of $i_{_{\Gamma \ra G}}$ is an ideal in $CL(G)$.
\end{x}
\vspace{0.1cm}

Write
\[
G 
\ = \ 
\coprod\limits_{k=1}^n \ 
x_k \Gamma.
\]

\vspace{0.1cm}


\begin{x}{\small\bf LEMMA} \ 
For any $\phi \in CL(\Gamma)$, 
\[
(i_{_{\Gamma \ra G}} \phi) (x) 
\ = \ 
\sum\limits_{k=1}^n \ 
\mathring{\phi} (x_k^{-1} x x_k).
\]

PROOF \ 
In fact, 
\allowdisplaybreaks
\begin{align*}
(i_{_{\Gamma \ra G}} \phi) (x) \ 
&=\ 
\frac{1}{\abs{\Gamma}} \ 
\sum\limits_{y \in G} \ 
\mathring{\phi} (y x y^{-1})
\\[15pt]
&=\ 
\frac{1}{\abs{\Gamma}} \ 
\sum\limits_{y \in G} \ 
\mathring{\phi} (y^{-1} x y)
\\[15pt]
&=\ 
\frac{1}{\abs{\Gamma}} \ 
\sum\limits_{\gamma \in \Gamma} \ 
\sum\limits_{k=1}^n \ 
\mathring{\phi} (\gamma^{-1} x_k^{-1} x x_k \gamma).
\end{align*}
There are then two possibilities.
\\[-.25cm]

\qquad \textbullet \quad $\gamma^{-1} x_k^{-1} x x_k   \gamma \notin \Gamma$
\\[-.2cm]

\hspace{2cm} $\implies x_k^{-1} x x_k \notin \Gamma$
\\[-.2cm]

\hspace{2cm}  $\implies 
\mathring{\phi} (\gamma^{-1} x_k^{-1} x x_k \gamma) 
\ = \ 0 \ = \ 
\mathring{\phi} (x_k^{-1} x x_k).$
\\[-.2cm]

\qquad \textbullet \quad $\gamma^{-1} x_k^{-1} x x_k   \gamma \in \Gamma$
\\[-.2cm]

\hspace{2cm}  $\implies x_k^{-1} x x_k \in \Gamma$
\\[-.2cm]

\hspace{2cm}  $\implies \mathring{\phi} (\gamma^{-1} x_k^{-1} x x_k \gamma) 
\ = \ 
\phi (\gamma^{-1} x_k^{-1} x x_k \gamma)$
\\[-.2cm]

$
\hspace{6.0cm} \ = \ 
\phi (x_k^{-1} x x_k) $
\\[-.2cm]

$
\hspace{6.0cm} \ = \ 
\mathring{\phi}(x_k^{-1} x x_k) $.
\\[-.2cm]

\noindent Therefore the sum $\ds\frac{1}{\abs{\Gamma}} \ \sum\limits_{\gamma \in \Gamma}$ 
disappears, leaving
\[
\sum\limits_{k=1}^n \ 
\mathring{\phi} (x_k^{-1} x x_k).
\]
\vspace{0.2cm}

[Note: \ 
If instead, 
\[
G 
\ = \ 
\coprod\limits_{k=1}^n \ 
\Gamma x_k,
\]
then for any $\phi \in CL(\Gamma)$, 
\[(i_{_{\Gamma \ra G}} \phi) (x) 
\ = \ 
\sum\limits_{k=1}^n \ 
\mathring{\phi} (x_k x x_k^{-1}).]
\]
\end{x}
\vspace{0.1cm}

\begin{x}{\small\bf EXAMPLE} \ 
Let \mS be a transitive \mG-set, $\pi$ the associated representation (cf. \S5, \#1).  
Fix a point $s \in S$ and let $G_s$ be its stabilizer $-$then
\[
\chisubpi 
\ = \ 
i_{G_s \ra G} 1_{G_s},
\]
where $1_{G_s} \in CL(G_s)$ is $\equiv 1$.
\vspace{0.2cm}

[Take $S \ = \ \{1, \ldots, n\}$ and $s = 1$.  Write
\[
G 
\ = \ 
\coprod\limits_{k=1}^n \ 
x_k G_s
\]
with $x_k \cdot 1 = k$ $-$then
\allowdisplaybreaks
\begin{align*}
(i_{G_s \ra G} \hsx 1_{G_s}) (x) \ 
&=\ 
\sum\limits_{k=1}^n \
\mathring{1}_{G_s} (x_k^{-1} x x_k)
\\[15pt]
&=\ 
\sum\limits_{k, x_k^{-1} x x_k \in G_s} \
1
\\[15pt]
&=\ 
\sum\limits_{k, (x_k^{-1} x x_k) \hsy \cdot \hsy 1 \hsy = \hsy 1} \
1
\\[15pt]
&=\ 
\sum\limits_{k, (x x_k) \hsy \cdot \hsy 1   \hsy = \hsy x_k \hsy \cdot \hsy 1} \
1
\\[15pt]
&=\ 
\sum\limits_{k, x \hsy \cdot \hsy k = k} \
1
\\[15pt]
&=\ 
\abs{\{k \in S : x \hsy \cdot \hsy k = k\}}
\\[15pt]
&=\ 
\chisubpi (x) 
\qquad \text{(cf. \S5, \#1)}.]
\end{align*}
\vspace{0.2cm}

[Note: \ 
Here is a ``for instance''.  
Take $S = G / \Gamma$ and write 
\[
G / \Gamma 
\ = \ 
\coprod\limits_{k=1}^n 
\hsx x_k \Gamma.
\]
Then $G / \Gamma$ is a transitive \mG-set and 
\[
G_{x_k \Gamma} 
\ = \ 
x_k \Gamma x_k^{-1}.
\]
In particular: 
Take $x_k = 1$ to get 
\[
\chisubpi 
 \ = \ 
 i_{_{\Gamma \ra G}} 1_\Gamma, 
\]
thus at a given $x \in G$, $(i_{_{\Gamma \ra G}} 1_\Gamma)(x)$ is the number of left cosets of 
$\Gamma$ in \mG fixed by $x$.]
\end{x}
\vspace{0.1cm}

\begin{x}{\small\bf LEMMA} \ 
Suppose that $\Gamma_1 \subset \Gamma_2 \subset G$.  
Let $\phi_1 \in CL(\Gamma_1)$ $-$then 
\[
i_{_{\Gamma_2 \ra G}} (i_{_{\Gamma_1 \ra \Gamma_2}} \phi_1) 
\ = \ 
i_{_{\Gamma_1 \ra G}} \phi_1.
\]

PROOF \ 
Both sides of the putative equality are class functions, thus it suffices to show that 
\[
\langle
i_{_{\Gamma_2 \ra G}} (i_{_{\Gamma_1 \ra \Gamma_2}} \phi_1), \chisubPi
\ranglesubG 
\ = \ 
\langle
i_{_{\Gamma_1 \ra G}} \phi_1, \chisubPi 
\ranglesubG 
\]
for all $\Pi \in \widehat{G}$.  
But the LHS equals
\allowdisplaybreaks
\begin{align*}
\langle i_{_{\Gamma_1 \ra \Gamma_2}} \phi_1, r_{_{G \ra \Gamma_2}} \chisubPi \ranglesubGammaTwo \ 
&=\ 
\langle
\phi_1, r_{_{\Gamma_2 \ra \Gamma_1}} (r_{_{G \ra \Gamma_2}} \chisubPi)
\ranglesubGammaOne
\\[15pt]
&=\ 
\langle
\phi_1, \res_{_{G \ra \Gamma_1}} \chisubPi
\ranglesubGammaOne
\\[15pt]
&=\ 
\langle
i_{_{\Gamma_1 \ra G}} \phi_1, \chisubPi
\ranglesubG,
\end{align*}
which is the RHS.
\end{x}
\vspace{0.1cm}

\begin{x}{\small\bf NOTATION} \ 
Given $x \in G$, put
\[
\Gamma^x 
\ = \ 
x \Gamma x^{-1} 
\ = \ 
\{x \gamma x^{-1} : \gamma \in \Gamma \}.
\]
The range of 
\[
i_{_{\Gamma \ra G}} : CL(\Gamma) \ra CL(G)
\]
is contained in the subspace $\sS_\Gamma$ of $CL(G)$ consisting of those class functions $f \in CL(G)$ that vanish on 
\[
G 
\ - \ 
\bigcup\limits_{x \in G} \hsx \Gamma^x.
\]
\end{x}
\vspace{0.1cm}

\begin{x}{\small\bf LEMMA} \ 
\[
i_{_{\Gamma \ra G}} CL(\Gamma) 
\ = \ 
\sS_\Gamma.
\]

PROOF \ 
Assume not, thus
\[
i_{_{\Gamma \ra G}} CL(\Gamma) 
\ \neq \ 
\sS_\Gamma.
\]
Then there exists a nonzero $f \in \sS_\Gamma$ which is orthogonal to all functions in 
$i_{_{\Gamma \ra G}} CL(\Gamma)\hsx :$  $\forall \ \phi \in CL(\Gamma)$, 
\[
\langle 
i_{_{\Gamma \ra G}} \phi, f
\ranglesubG 
\ = \ 
0
\]
or still, $\forall \ \phi \in CL(\Gamma)$, 
\[
\langle 
\phi, r_{_{G \ra \Gamma}} f 
\ranglesubGamma 
\ = \ 
0.
\]
Now take $\phi = r_{_{G \ra \Gamma}} f$ to get 
\[
\langle 
r_{_{G \ra \Gamma}} f, r_{_{G \ra \Gamma}} f 
\ranglesubGamma 
\ = \ 
0,
\]
hence $r_{_{G \ra \Gamma}} f = 0$, i.e., $f$ vanishes on $\Gamma$.  
But $f \in CL(G)$, so $\forall \ x \in G$, $f$ vanishes on $\Gamma^x$.  
Since $f \in \sS_\Gamma$, it then follows that $f$ vanishes on \mG: $f \equiv 0$, contradicting the supposition that $f$ is nonzero.
\\
\end{x}

\begin{x}{\small\bf APPLICATION} \ 
The image of $i_{_{\Gamma \ra G}}$ is an ideal in $CL(G)$ (cf. \#9).
\end{x}
\vspace{0.3cm}

Let $\phi \in CL(\Gamma)$.  
Given $x \in G$, define $\phi^x \in C(\Gamma^x)$ by 
\allowdisplaybreaks
\begin{align*}
\phi^x (y) \ 
&=\ 
\phi(x^{-1} y x) 
\qquad (y = x \gamma x^{-1}, \ \gamma \in \Gamma)
\\[15pt]
&=\ 
\phi(x^{-1} x \gamma x^{-1} x)
\\[15pt]
&=\ 
\phi (\gamma).
\end{align*}
\vspace{0.1cm}

\begin{x}{\small\bf LEMMA} \ 
\[
\phi^x \in CL(\Gamma^x).
\]

PROOF \ 
Let 
\[
y_1 
\ = \ 
x \gamma_1 x^{-1}, 
\hspace{0.5cm}
y_2 \ = \ x \gamma_2 x^{-1}.
\]
Then
\allowdisplaybreaks
\begin{align*}
\phi^x (y_1 y_2 y_1^{-1}) \ 
&=\ 
\phi(x^{-1} y_1 y_2 y_1^{-1} x)
\\[15pt]
&=\ 
\phi(x^{-1}  (x \gamma_1 x^{-1}) (x \gamma_2 x^{-1}) (x \gamma_1 x^{-1})^{-1} x)
\\[15pt]
&=\ 
\phi(x^{-1} (x \gamma_1 x^{-1}) (x \gamma_2 x^{-1}) (x \gamma_1^{-1} x^{-1})x)
\\[15pt]
&=\ 
\phi(\gamma_1 \gamma_2 \gamma_1^{-1})
\\[15pt]
&=\ 
\phi(\gamma_2)
\\[15pt]
&=\ 
\phi^x(y_2).
\end{align*}
\\[-.5cm]
\end{x}

\begin{x}{\small\bf LEMMA} \ 
$\forall \ x \in G$ and $\forall \ \phi \in CL(\Gamma)$, 
\[
i_{_{\Gamma^x \ra G}} \phi^x 
\ = \ 
i_{_{\Gamma \ra G}} \phi.
\]
\vspace{0.2cm}

PROOF \ 
Write 
\[
G 
\ = \ 
\coprod\limits_{k=1}^n \ 
x_k \Gamma 
\ = \ 
\coprod\limits_{k=1}^n \ 
x x_k x^{-1} \Gamma^x.
\]
Then (cf. \#10)
\allowdisplaybreaks
\begin{align*}
(i_{_{\Gamma^x \ra G}} \phi^x ) (y) \ 
&=\ 
\sum\limits_{k=1}^n \ 
\mathring{\phi}^x (( x x_k x^{-1})^{-1} \hsx y (x x_k x^{-1}))
\\[15pt]
&=\ 
\sum\limits_{k=1}^n \ 
\mathring{\phi}^x ( x x_k^{-1} x^{-1} y x x_k x^{-1})
\\[15pt]
&=\ 
\sum\limits_{k=1}^n \ 
\mathring{\phi} (x^{-1}( x x_k^{-1} x^{-1} y x x_k x^{-1}) x)
\\[15pt]
&=\ 
\sum\limits_{k=1}^n \ 
\mathring{\phi} (x_k^{-1} x^{-1} y x x_k)
\\[15pt]
&=\ 
(i_{_{\Gamma \ra G}} \phi) (x^{-1} y x)
\\[15pt]
&=\ 
(i_{_{\Gamma \ra G}} \phi) (y)
\qquad \text{(cf. \#4)}.
\end{align*}
\end{x}
\vspace{0.1cm}


\chapter{
$\boldsymbol{\S}$\textbf{8}.\quad  MACKEY THEORY}
\setlength\parindent{2em}
\setcounter{theoremn}{0}
\renewcommand{\thepage}{A II \S8-\arabic{page}}

\qquad Let \mG be a finite group, let $\Gamma_1$, $\Gamma_2 \subset G$ be subgroups, and let
\[
G \ = \ \bigcup\limits_{s \in S} \hsx \Gamma_1 \hsx s \hsx  \Gamma_2
\]
be a double coset decomposition of \mG.  
Given $s \in S$, put
\[
\Gamma_2(s) \ = \ \Gamma_2^s \hsx \cap \hsx \Gamma_1 \qquad (= s \Gamma_2 s^{-1} \cap \Gamma_1).
\]
\\[-.5cm]

\begin{x}{\small\bf LEMMA} \ 
Let 
\[
\Gamma_1 
\ = \ 
\bigcup\limits_{t \hsy \in \hsy T(s)} \ 
t \hsx \Gamma_2 (s)
\]
be a left coset decomposition of $\Gamma_1$ $-$then 
\allowdisplaybreaks
\begin{align*}
\Gamma_1 \hsy s \hsy \Gamma_2 \ 
&=\ 
\bigg(
\bigcup\limits_{t \hsy \in \hsy T(s)} \ 
t \hsy \Gamma_2(s) \bigg) \hsy s \hsy \Gamma_2
\\[15pt]
&=\ 
\bigcup\limits_{t \hsy \in \hsy T(s)} \ 
t \Gamma_2(s) \hsy s \hsy \Gamma_2
\\[15pt]
&=\ 
\bigcup\limits_{t \hsy \in \hsy T(s)} \ 
t \Gamma_2(s) (s \hsy \Gamma_2 \hsy s^{-1}) s
\\[15pt]
&=\ 
\bigcup\limits_{t \hsy \in \hsy T(s)} \ 
t (s \hsy \Gamma_2 \hsy s^{-1}) \hsy s
\\[15pt]
&=\ 
\bigcup\limits_{t \hsy \in \hsy T(s)} \ 
t s \Gamma_2
\end{align*}
is a partition of $\Gamma_1 \hsy s \hsy \Gamma_2$.
\vspace{0.2cm}

PROOF \ 
Suppose that 
\[
t_1 \hsy s \hsy \Gamma_2 \hsx \cap \hsx t_2 \hsy s \hsy \Gamma_2 
\ \neq \ 
\emptyset 
\qquad (t_1 \ \neq \ t_2),
\]
so
\[
t_1 \hsy s 
\ = \ 
t_2 \hsy s \hsy \gamma_2 
\qquad (\gamma_2 \in \Gamma_2).
\]
Then
\[
t_1 
\ = \ 
t_2 \hsy s \hsy \gamma_2\hsy  s^{-1} 
\implies 
t_2^{-1} \hsy t_1 \in \Gamma_2^s.
\]
Meanwhile
\[
t_1, \ t_2 \in \Gamma_1 
\implies 
t_2^{-1} t_1 \in \Gamma_1.
\]
Therefore
\allowdisplaybreaks
\begin{align*}
t_2^{-1} t_1 \ 
&\in \Gamma_2^s \hsx \cap \hsx \Gamma_1 \ = \ \Gamma_2(s) 
\\[12pt]
&\implies 
t_1 = t_2.
\end{align*}
\end{x}
\vspace{0.3cm}

Let $R(s) = \{t s : t \in T(s) \} \hsx \equiv \hsx T(s) \hsy s$ and let
\[
R 
\ = \ 
\bigcup\limits_{s \in S} \
R(s).
\]
\vspace{0.3cm}

\begin{x}{\small\bf LEMMA} \ 
\mR is a set of left coset representatives of $\Gamma_2$ in \mG.  
\vspace{0.3cm}

PROOF
Let $x \in G$ $-$then
\allowdisplaybreaks
\begin{align*}
x \in \Gamma_1 s \Gamma_2 \quad (\exists \ s \in S) \ 
&\implies 
x \ = \ t \hsy s \hsy \gamma_2 \qquad (\exists \ t \in T(s))
\\[12pt]
&\implies 
x \ = \ r \hsy \gamma_2 \qquad (r \in R(s), \ r = t s).
\end{align*}
Therefore
\[
G 
\ = \ 
\bigcup\limits_{r \in R}  \
r \hsy \Gamma_2.
\]
Suppose now that
\[
x \in r \hsy \Gamma_2 \hsx \cap \hsx r^\prime \hsy \Gamma_2.
\]
Then 
\[
x \ = \ r \hsy \gamma_2 \ = \ r^\prime \gamma_2^\prime 
\qquad
(r \in R(s), \ r^\prime \in R(s^\prime))
\]

\[
\implies
\begin{cases}
\ x \ = \ t \hsy s \hsy \gamma_2 \hspace{0.75cm} (t \in T(s))\\
\ x \ = t^\prime \hsy s^\prime \hsy \gamma_2^\prime \hspace{0.7cm} (t^\prime \in T(s^\prime))
\end{cases}
.
\] 
\[
\text{But} \ \  
\begin{cases}
\ t \in T(s) \hspace{0.2cm} \implies t \in \Gamma_1\\
\ t^\prime \in T(s^\prime) \implies t^\prime \in \Gamma_1
\end{cases}
\implies x \in \Gamma_1 \hsy s \hsy \Gamma_2 \hsx \cap \hsx \Gamma_1 \hsy s^\prime \hsx \Gamma_2 \implies s \ = \ s^\prime
\hspace{1.75cm}
\]
\allowdisplaybreaks
\begin{align*}
&\implies t \hsy s \hsy \gamma_2 \ = \ t^\prime \hsy s \hsy \gamma_2^\prime
\\[12pt] 
&\implies 
t s = t^\prime \hsy s \hsy \gamma_2^\prime \hsy \gamma_2^{-1} \ = \ t^\prime \hsy s \hsy \gamma_2^{\prime\prime}
\\[12pt] 
&\implies t \ = \ t^\prime 
\\[12pt]
&\implies r \ = \ r^\prime.
\end{align*}
\end{x}
\vspace{0.1cm}

Given $\phi \in CL(\Gamma_2)$, put
\[
\phi_s 
\ = \ 
r_{_{\Gamma_2^s \ra \Gamma_2(s)}} \phi^s.
\]
Here, by definition (cf. \S7, \#16), $\phi^s \in CL(\Gamma_2^s)$, where
\[
\phi^s(y) 
\ = \ 
\phi(\gamma_2) 
\qquad 
(y = s \hsy \gamma_2 \hsy s^{-1}, \gamma_2 \in \Gamma_2).
\]
\begin{x}{\small\bf THEOREM} \ 
Under the above assumptions, 
\[
r_{_{G  \ra \Gamma_1}} (i_{_{\Gamma_2 \ra G}} \phi) 
\ = \ 
\sum\limits_{s \in S} \ 
i_{_{\Gamma_2(s) \ra \Gamma_1}} \phi_s.
\]
\vspace{0.3cm}

PROOF \ 
Since
\[
G 
\ = \ 
\coprod\limits_{r \in R} \ 
r \Gamma_2,
\]
$\forall \ x \in G$, 
\[
(i_{_{\Gamma_2 \ra G}} \phi) (x) 
\ = \ 
\sum\limits_{r \in R} \ 
\mathring{\phi} (r^{-1} \hsy x \hsy r) 
\qquad \text{(cf. \ \S7, \#10)},
\]
so $\forall \ \gamma_1 \in \Gamma_1$, 
\allowdisplaybreaks
\begin{align*}
(r_{_{G \ra \Gamma_1}} (i_{_{\Gamma_2 \ra G}} \phi)) (\gamma_1) \ 
&=\ 
\sum\limits_{r \in R} \
\mathring{\phi} (r^{-1} \hsy \gamma_1 \hsy r)
\\[15pt]
&=\ 
\sum\limits_{\substack{r \hsy \in \hsy R,\\ r^{-1} \hsy \gamma_1 \hsy r \hsy \in \hsy \Gamma_2}} \ 
\phi (r^{-1} \hsy \gamma_1 \hsy r)
\\[15pt]
&=\ 
\sum\limits_{\substack{s \in S,\\r \hsy \in \hsy R(s), \\ r^{-1} \hsy \gamma_1 \hsy r \hsy \in \hsy \Gamma_2}} \
\phi (r^{-1} \hsy \gamma_1 \hsy r)
\\[15pt]
&=\ 
\sum\limits_{\substack{s \in S,\\t \in T(s), \\s^{-1} \hsy t^{-1} \hsy \gamma_1 \hsy t \hsy s \hsy \in \hsy \Gamma_2}} \
\phi (s^{-1} t^{-1} \gamma_1 \hsy t \hsy s)
\\[15pt]
&=\ 
\sum\limits_{\substack{s \in S,\\t \in T(s), \\t^{-1} \hsy \gamma_1 \hsy t \hsy \in \hsy \Gamma_2^s}} \
\phi^s (t^{-1} \hsy \gamma_1 \hsy t )
\\[15pt]
&=\ 
\sum\limits_{s \in S} \
\sum\limits_{t \hsy \in \hsy T(s)} \
\phi_s (t^{-1} \hsy \gamma_1 \hsy t ) 
\\[15pt]
&=\ 
\sum\limits_{s \in S} \
(i_{_{\Gamma_2(s) \ra \Gamma_1}} \phi_s) (\gamma_1) 
\qquad \text{(cf. \S7, \#10)}.
\end{align*}
\end{x}
\vspace{0.1cm}

\begin{x}{\small\bf LEMMA} \ 
Let $\psi \in CL(\Gamma_1)$, $\phi \in CL(\Gamma_2)$ $-$then
\[
\langle i_{_{\Gamma_1 \ra G}} \psi, i_{_{\Gamma_2 \ra G}} \phi \ranglesubG 
\ = \ 
\sum\limits_{s \in S} \ 
\langle r_{_{\Gamma_1 \ra \Gamma_2(s)}} \psi, \phi_s \ranglesubGammaTwoOfs.
\]
\vspace{0.2cm}

PROOF \ 
Taking into account \S7, \#6, 
\allowdisplaybreaks
\begin{align*}
\langle i_{_{\Gamma_1 \ra G}} \psi, i_{_{\Gamma_2 \ra G}} \phi \ranglesubG \ 
&=\ 
\langle
\psi,
r_{_{G \ra \Gamma_1}} (i_{_{\Gamma_2 \ra G}} \phi)
\ranglesubGammaOne
\\[15pt]
&=\ 
\langle
\psi,\hsy
\sum\limits_{s \in S} \ 
i_{_{\Gamma_2(s) \ra \Gamma_1}} \phi_s
\ranglesubGammaOne
\\[15pt]
&=\ 
\sum\limits_{s \in S} \ 
\langle
\psi, i_{_{\Gamma_2(s) \ra \Gamma_1}} \phi_s
\ranglesubGammaOne
\\[15pt]
&=\ 
\sum\limits_{s \in S} \ 
\raisebox{.4cm}
{$
\ov{
\raisebox{-.4cm}
{$
\langle
i_{_{\Gamma_2(s) \ra \Gamma_1}} \phi_s, \psi
\ranglesubGammaOne
$}
}
$}
\\[15pt]
&=\ 
\sum\limits_{s \in S} \ 
\raisebox{.4cm}
{$
\ov{
\raisebox{-.4cm}
{$
\langle
\phi_s, r_{_{\Gamma_1 \ra \Gamma_2(s)}} \psi
\ranglesubGammaTwoOfs
$}
}
$}
\\[15pt]
&=\ 
\sum\limits_{s \in S} \ 
\langle
r_{_{\Gamma_1 \ra \Gamma_2(s)}} \psi, \phi_s 
\ranglesubGammaTwoOfs.
\end{align*}
\end{x}
\vspace{0.3cm}

\begin{x}{\small\bf NOTATION} \ \ 
Given a subgroup $\Gamma \subset G$, \ 
let $1_\Gamma$ stand for the function $\Gamma \ra \Cx$ which is $\equiv 1$, 
that is, the character of the trivial one-dimensional representation of $\Gamma$.  
\end{x}
\vspace{0.3cm}


\begin{x}{\small\bf EXAMPLE} \ 
Take $\Gamma_1 = \Gamma_2 = \Gamma$ $-$then 
\[
\langle i_{_{\Gamma \ra G}} 1_\Gamma, i_{_{\Gamma \ra G}} 1_\Gamma \ranglesubG 
\ = \ 
\abs{\Gamma \backslash G / \Gamma}.
\]
Therefore $i_{_{\Gamma \ra G}} 1_\Gamma$ is not irreducible if 
$\abs{\Gamma \backslash G / \Gamma} > 1$ (cf. \S5, \#11).
\\[-.2cm]

[Note: \ 
$i_{_{\Gamma \ra G}} 1_\Gamma$ is a character of \mG (cf. \S7, \#7).]
\end{x}
\vspace{0.3cm}


\chapter{
$\boldsymbol{\S}$\textbf{9}.\quad  INDUCED REPRESENTATIONS}
\setlength\parindent{2em}
\setcounter{theoremn}{0}
\renewcommand{\thepage}{A II \S9-\arabic{page}}

\qquad Let \mG be a finite group, $\Gamma \subset G$ a subgroup.
\\[-.2cm]

\begin{x}{\small\bf CONSTRUCTION} \ 
Let $(\theta,\E)$ be a unitary representation of $\Gamma$ and denote by $\E_{\Gamma,\theta}^G$ the space of all $\E$-valued functions $f$ on \mG such that $f(x \gamma) = \theta (\gamma^{-1}) f(x)$ $(x \in G, \gamma \in \Gamma)$ $-$then the prescription 
\[
\big(\Ind_{\Gamma,\theta}^G  (x) f\big) (y)
\ =\ 
f(x^{-1} y)
\]
defines a representation $\Ind_{\Gamma,\theta}^G$ of \mG on $\E_{\Gamma,\theta}^G$, the representation of \mG 
\un{induced}
\index{representation \\ induced} 
by $\theta$.
\end{x}
\vspace{0.3cm}

\begin{x}{\small\bf \un{N.B.}} \ 
The inner product 
\[
\langle f,g \ranglesubtheta
\ = \
\frac{1}{\abs{G}} \hsx
\sum\limits_{x \in G} \hsx \hsx
\langle f(x),g(x) \ranglesubEE
\]
equips $\E_{\Gamma,\theta}^G$ with the structure of a Hilbert space and $\Ind_{\Gamma,\theta}^G$ is a unitary representation.
\end{x}
\vspace{0.3cm}

\begin{x}{\small\bf EXAMPLE} \ 
Take $\Gamma = \{e\}$ and take $\theta$ to be the trivial representation of $\Gamma$ on $\E = \Cx$ $-$then 
$\E_{\Gamma,\theta}^G = \tC(G)$ and 
\[
\Ind_{\Gamma,\theta}^G 
\ = \ 
L,
\]
the left translation representation of \mG (cf. \S1, \#12).
\end{x}
\vspace{0.3cm}

\begin{x}{\small\bf EXAMPLE} \ 
Take $\Gamma = G$ and let $(\pi,V)$ be a unitary representation of \mG.   
Define a linear bijection
\[
T : V_{G,\pi}^G \ra V
\]
by sending $f$ to $f(e)$ $-$then $\forall \ x \in G$, 
\allowdisplaybreaks
\begin{align*}
T(\Ind_{G,\pi}^G ) (x) f) \ 
&=\ 
(\Ind_{G,\pi}^G  (x) f) (e) 
\\[15pt]
&=\ 
f(x^{-1} e)
\\[15pt]
&=\ 
f(x^{-1})
\\[15pt]
&=\ 
f(e x^{-1})
\\[15pt]
&=\ 
\pi(x) f(e)
\\[15pt]
&=\ 
\pi(x) (T f).
\end{align*}
Therefore
\[
T 
\hsx \circ \hsx 
\Ind_{G,\pi}^G 
\ = \ 
\pi \hsx \circ \hsx T.
\]
I.e.: 
\[
T \in I_G(\Ind_{G,\pi}^G ,\pi)
\]
is an invertible interwining operator, thus $\Ind_{G,\pi}^G$ is equivalent to $\pi$.
\vspace{0.2cm}

[Note: \ 
\mT is unitary.  In fact, 
\allowdisplaybreaks
\begin{align*}
\langle f, g \ranglesubpi \ 
&=\ 
\frac{1}{\abs{G}} \hsx 
\sum\limits_{x \in G} \hsx
\langle f(x), g(x) \ranglesubV
\\[15pt]
&=\ 
\frac{1}{\abs{G}} \hsx 
\sum\limits_{x \in G} \hsx
\langle f(ex), g(ex) \ranglesubV
\\[15pt]
&=\ 
\frac{1}{\abs{G}} \hsx 
\sum\limits_{x \in G} \hsx
\langle \pi(x^{-1}) f(e), \pi(x^{-1}) g(e) \ranglesubV
\\[15pt]
&=\ 
\frac{1}{\abs{G}} \hsx 
\sum\limits_{x \in G} \hsx
\langle f(e), g(e) \ranglesubV
\\[15pt]
&=\ 
\langle f(e), g(e) \ranglesubV
\\[15pt]
&=\ 
\langle T(f), T(g) \ranglesubV.]
\end{align*}
\end{x}
\vspace{0.3cm}

\begin{x}{\small\bf LEMMA} \ 
The dimension of $\E_{\Gamma,\theta}^G$ equals
\[
\frac{\abs{G}}{\abs{\Gamma}} \hsx  \dim \E.
\]
\vspace{0.2cm}


PROOF \ 
Write
\[
G 
\ = \ 
\coprod\limits_{k=1}^n \hsx x_k \Gamma, 
\]
where $n = \ds\frac{\abs{G}}{\abs{\Gamma}}$, and define a bijection 
\[
\Lambda : \E_{\Gamma,\theta}^G \ra \bigoplus\limits_{k=1}^n \hsx \E
\]
by the stipulation that 
\[
\Lambda f 
\ = \ 
(f(x_1) , \ldots, f(x_n)),
\]
from which the assertion.
\end{x}
\vspace{0.3cm}

For any character $\chi$ of \mG and for any conjugacy class $C \in \CON(G)$, write $\chi(C)$ for the common value of $\chi(x)$ $(x \in C)$ (and analogously if \mG is replaced by $\Gamma$).

Fixing \mC, the intersection $C \hsx \cap \hsx \Gamma$ is a union of elements of $\CON(\Gamma)$, say
\[
C \hsx \cap \hsx \Gamma
\ = \ 
\bigcup\limits_\ell \hsx C_\ell.
\]
\vspace{0.2cm}

[Note: \ 
If $C \hsx \cap \hsx \Gamma = \emptyset$, then the sum that follows is empty and its value is 0.]
\vspace{0.3cm}

\begin{x}{\small\bf THEOREM} \ 
Set $\pi = \Ind_{\Gamma,\theta}^G$ $-$then 
\[
\chisubpi(C) 
\ =\ 
\frac{\abs{G}}{\abs{\Gamma}} \hsx 
\sum\limits_\ell
\frac{\abs{C_\ell}}{\abs{C}} \hsx
\chisubtheta (C_\ell).
\]

\vspace{0.2cm}

PROOF \ 
If $\chisubC$ is the characteristic function of \mC (cf. \S4, \#20), then 
$\ds\chisubC = \sum\limits_{y \in C} \hsx \delta_y$ (cf. \S4, \#21).  
Denoting by $\rho$ the canonical extension of $\pi$ to $C(G)$, it thus follows that
\[\chisubpi(C) 
\ = \ 
\frac{1}{\abs{C}} \hsx 
\tr(\rho(\chisubC)).
\]
Fix an orthonormal basis $\phi_1, \ldots, \phi_m$ in $\E$ and in $\E_{\Gamma,\theta}^G$, let
\[
f_j(x) \ = \ 
\begin{cases}
\ \ds\bigg( \frac{\abs{G}}{\abs{\Gamma}} \bigg)^{1/2} \hsx \theta(\gamma^{-1}) \phi_j \qquad (x = \gamma \in \Gamma)\\[8pt]
\hspace{0.65cm} 0 \hspace{3.2cm}  (x \notin \Gamma)
\end{cases}
.
\]
\vspace{0.1cm}

\qquad \textbullet \ 
The $f_j \ (1 \leq j \leq m)$ are an orthonormal set in $\E_{\Gamma,\theta}^G$.
\\[-.2cm]

\qquad \textbullet \ 
The $\rho(x_k) f_j$ $(1 \leq k \leq n$, $1 \leq j \leq m)$ are an orthonormal basis for $\E_{\Gamma,\theta}^G$.
\\[-.2cm]

Proceeding
\allowdisplaybreaks
\begin{align*}
\tr(\rho(\chisubC)) \ 
&=\ 
\sum\limits_{k=1}^n \hsx
\sum\limits_{j=1}^m \hsx
\langle \rho(\chisubC) \hsx \rho(x_k) f_j, \rho(x_k) f_j \ranglesubtheta
\\[15pt]
&=\ 
\sum\limits_{k=1}^n \hsx
\sum\limits_{j=1}^m \hsx
\langle \rho(x_k^{-1}) \hsx \rho(\chisubC) \hsx \rho(x_k) f_j,f_j \ranglesubtheta
\\[15pt]
&=\ 
\sum\limits_{k=1}^n \hsx
\sum\limits_{j=1}^m \hsx
\langle \rho(\chisubC) f_j,f_j \ranglesubtheta
\\[15pt]
&=\ 
n \hsx
\sum\limits_{j=1}^m \hsx
\langle \rho(\chisubC) f_j,f_j \ranglesubtheta
\\[15pt]
&=\ 
[G : \Gamma] \hsx
\sum\limits_{j=1}^m \hsx
\langle \rho(\chisubC) f_j,f_j \ranglesubtheta
\\[15pt]
&=\ 
\frac{\abs{G}}{\abs{\Gamma}} \hsx 
\sum\limits_{j=1}^m \hsx
\langle \rho(\chisubC) f_j,f_j \ranglesubtheta.
\end{align*}
But
\allowdisplaybreaks
\begin{align*}
\langle \rho(\chisubC) f_j,f_j \ranglesubtheta \ 
&=\ 
\big\langle \ 
\sum\limits_{y \in C} \hsx 
\rho(\delta_y) f_j,f_j 
\hsx \big\ranglesubtheta
\\[15pt]
&=\ 
\sum\limits_{y \in C} \hsx 
\langle
\rho(\delta_y) f_j,f_j 
\ranglesubtheta
\\[15pt]
&=\ 
\sum\limits_{y \in C \cap \Gamma} \hsx 
\langle
\rho(\delta_y) f_j,f_j 
\ranglesubtheta.
\end{align*}
Therefore
\allowdisplaybreaks
\begin{align*}
\tr(\rho(\chisubC)) \ 
&=\ 
\frac{\abs{G}}{\abs{\Gamma}} \hsx
\sum\limits_{j=1}^m \hsx
\sum\limits_{y \in C \cap \Gamma} \hsx
\langle
\rho(\delta_y) f_j,f_j 
\ranglesubtheta
\\[15pt]
&=\ 
\frac{\abs{G}}{\abs{\Gamma}} \hsx
\sum\limits_{j=1}^m \hsx
\sum\limits_\ell \hsx
\sum\limits_{\gamma \in C_\ell} \hsx
\langle
\rho(\gamma) f_j,f_j 
\ranglesubtheta.
\end{align*}
But 
$\forall \ \gamma_0 \in \Gamma$, 
\allowdisplaybreaks
\begin{align*}
\langle
\rho(\gamma_0) f_j,f_j 
\ranglesubtheta \ 
&=\ 
\frac{1}{\abs{G}} \hsx
\sum\limits_{\gamma \in \Gamma} \hsx
\langle
f_j(\gamma_0^{-1} \gamma), f_j(\gamma)
\ranglesubEE
\\[15pt]
&=\ 
\frac{1}{\abs{G}} \hsx
\sum\limits_{\gamma \in \Gamma} \hsx
\frac{\abs{G}}{\abs{\Gamma}} \hsx
\langle
\theta(\gamma^{-1}\gamma_0) \phi_j, \theta(\gamma^{-1}) \phi_j
\ranglesubEE
\\[15pt]
&=\ 
\frac{1}{\abs{G}} \hsx
\sum\limits_{\gamma \in \Gamma} \hsx
\frac{\abs{G}}{\abs{\Gamma}} \hsx
\langle
\theta(\gamma^{-1}) \theta(\gamma_0) \phi_j, \theta(\gamma^{-1}) \phi_j
\ranglesubEE
\\[15pt]
&=\ 
\frac{1}{\abs{G}} \hsx
\sum\limits_{\gamma \in \Gamma} \hsx
\frac{\abs{G}}{\abs{\Gamma}} \hsx
\langle
\theta(\gamma_0) \phi_j, \phi_j
\ranglesubEE
\\[15pt]
&=\ 
\langle
\theta(\gamma_0) \phi_j, \phi_j
\ranglesubEE.
\end{align*}
Therefore
\allowdisplaybreaks
\begin{align*}
\tr(\rho(\chisubC)) \ 
&=\ 
\frac{\abs{G}}{\abs{\Gamma}} \hsx
\sum\limits_{j=1}^m \hsx
\sum\limits_\ell \hsx
\sum\limits_{\gamma \in C_\ell} \hsx
\langle
\theta(\gamma) \phi_j, \phi_j
\ranglesubEE
\\[15pt]
&=\ 
\frac{\abs{G}}{\abs{\Gamma}} \hsx
\sum\limits_\ell \hsx
\sum\limits_{\gamma \in C_\ell} \hsx
\sum\limits_{j=1}^m \hsx
\langle
\theta(\gamma) \phi_j, \phi_j
\ranglesubEE
\\[15pt]
&=\ 
\frac{\abs{G}}{\abs{\Gamma}} \hsx
\sum\limits_\ell \hsx
\sum\limits_{\gamma \in C_\ell} \hsx
\tr(\theta(\gamma))
\\[15pt]
&=\ 
\frac{\abs{G}}{\abs{\Gamma}} \hsx
\sum\limits_\ell \hsx
\sum\limits_{\gamma \in C_\ell} \hsx
\chisubtheta(\gamma)
\\[15pt]
&=\ 
\frac{\abs{G}}{\abs{\Gamma}} \hsx
\sum\limits_\ell \hsx
\abs{C_\ell} \hsx 
\chisubtheta(C_\ell).
\end{align*}
I.e.: 
\[
\chisubpi(C) 
\ = \ 
\frac{\abs{G}}{\abs{\Gamma}} \hsx 
\sum\limits_\ell \hsx 
\frac{\abs{C_\ell}}{\abs{C}} \hsx 
\chisubtheta(C_\ell).
\]
\\[-.5cm]

[Note: \ 
If $\theta$ is the trivial representation of $\Gamma$ on $\E = \Cx$, then $\chisubtheta = 1_\Gamma$ 
(the function $\equiv 1$) and matters reduce to 
\[
\chisubpi(C) 
\ = \ 
\frac{\abs{G}}{\abs{\Gamma}} \hsx 
\frac{\abs{C \cap \Gamma}}{\abs{C}}.] 
\]
\end{x}
\vspace{0.3cm}

\begin{x}{\small\bf \un{N.B.}} \ 
Take $C = \{e\}$:
\[
\chisubpi(e) 
\ = \ 
\frac{\abs{G}}{\abs{\Gamma}} \hsx 
\chisubtheta(e)
\]
\qquad\qquad $\implies$
\[
\dim \E_{\Gamma, \theta}^G 
\ = \ 
\frac{\abs{G}}{\abs{\Gamma}} \hsx 
\dim \E \qquad \text{(cf. $\#5$).}
\]

\end{x}
\vspace{0.3cm}

\begin{x}{\small\bf LEMMA} \ 
Set $\pi = \Ind_{\Gamma,\theta}^G$ $-$then for any class function $f \in CL(G)$, 
\[
\langle \chisubpi, f \ranglesubG
\ = \ 
\langle \chisubtheta, \restr{f}{\Gamma} \ranglesubGamma . 
\]
\vspace{0.2cm}

PROOF \ 
\allowdisplaybreaks
\begin{align*}
\langle \chisubpi, f \ranglesubG\ 
&=\ 
\frac{1}{\abs{G}} \hsx
\sum\limits_{x \in G} \hsx 
\chisubpi (x) \ov{f(x)}
\\[15pt]
&=\ 
\frac{1}{\abs{G}} \hsx
\sum\limits_{C \in \CON(G)} \hsx
\abs{C} \chisubpi (C) \ov{f(C)}
\\[15pt]
&=\ 
\frac{1}{\abs{G}} \hsx
\sum\limits_{C \in \CON(G)} \hsx
\abs{C} \hsx
\frac{\abs{G}}{\abs{\Gamma}} \hsx
\sum\limits_{\ell} \hsx
\frac{\abs{C_\ell}}{\abs{C}} \hsx
\chisubtheta (C_\ell) \ov{f(C)} 
\\[15pt]
&=\ 
\frac{1}{\abs{\Gamma}} \hsx
\sum\limits_{C \in \CON(G)} \hsx
\sum\limits_{\ell} \hsx
\abs{C_\ell} \chisubtheta (C_\ell) \ov{f(C)}
\\[15pt]
&=\ 
\frac{1}{\abs{\Gamma}} \hsx
\sum\limits_{C_\ell \in \CON(\Gamma)} \hsx
\abs{C_\ell} \chisubtheta (C_\ell) \ov{f(C_\ell)}
\\[15pt]
&=\ 
\frac{1}{\abs{\Gamma}} \hsx
\sum\limits_{\gamma \in \Gamma} \hsx
\chisubtheta (\gamma) \ov{f(\gamma)} 
\\[15pt]
&=\ 
\langle \chisubtheta, \restr{f}{\Gamma} \ranglesubGamma.
\end{align*}
\vspace{0.2cm}

[Note: \ 
One cannot simply quote \S7, \#6 \ldots \ .]
\end{x}
\vspace{0.3cm}

\begin{x}{\small\bf APPLICATION} \ 
Take $f = \chisubPi$ $(\Pi \in \widehat{G})$ and suppose that $\theta$ is irreducible $-$then the multiplicity of $\Pi$ in 
$\Ind_{\Gamma, \theta}^G$ equals the multiplicity of $\theta$ in the restriction of $\Pi$ to $\Gamma$ (cf. \S5, \#5).
\end{x}
\vspace{0.3cm}

\begin{x}{\small\bf THEOREM} \ 
Set $\pi = \Ind_{\Gamma, \theta}^G$ $-$then 
\[
i_{_{\Gamma \ra G}}  \chisubtheta \ = \ \chisubpi.
\]
\vspace{0.2cm}

PROOF \ 
The function 
\[
i_{_{\Gamma \ra G}} \chisubtheta
\]
is a class function on \mG, as is $\chisubpi$, thus it suffices to show that $\forall \ \Pi \in \widehat{G}$, 
\[
\langle i_{_{\Gamma \ra G}} \chisubtheta , \chisubPi \ranglesubG
\ = \ 
\langle \chisubpi, \chisubPi \ranglesubG.
\]
But
\allowdisplaybreaks
\begin{align*}
\langle i_{_{\Gamma \ra G}} \chisubtheta, \chisubPi \ranglesubG\ 
&=\ 
\frac{1}{\abs{G}} \hsx
\sum\limits_{x \in G} \hsx
( i_{_{\Gamma \ra G}} \chisubtheta) (x) 
\ov{\chisubPi(x)}
\\[15pt]
&=\ 
\frac{1}{\abs{G}} \hsx
\frac{1}{\abs{\Gamma}} \hsx
\sum\limits_{x \in G} \hsx
\sum\limits_{y \in G} \hsx
\mathring{\chi}_{_\theta} (y x y^{-1}) \hsx 
\ov{\chisubPi(x)}
\\[15pt]
&=\ 
\frac{1}{\abs{G}}  \hsx
\frac{1}{\abs{\Gamma}} \hsx
\sum\limits_{x \in G} \hsx
\sum\limits_{y \in G} \hsx
\mathring{\chi}_\theta (x) \hsx 
\ov{\chisubPi(y^{-1}xy)}
\\[15pt]
&=\ 
\frac{1}{\abs{G}}  \hsx
\frac{1}{\abs{\Gamma}} \hsx
\sum\limits_{x \in G} \hsx
\sum\limits_{y \in G} \hsx
\mathring{\chi}_\theta (x) \hsx 
\ov{\chisubPi(x)}
\\[15pt]
&=\ 
\frac{1}{\abs{\Gamma}} \hsx
\sum\limits_{x \in G} \hsx
\mathring{\chi}_\theta (x) \hsx
\ov{\chisubPi(x)}
\\[15pt]
&=\ 
\frac{1}{\abs{\Gamma}} \hsx
\sum\limits_{\gamma \in \Gamma} \hsx
\chisubtheta (\gamma) \hsx
\ov{\chisubPi(\gamma)}
\\[15pt]
&=\ 
\langle \chisubtheta, \restr{\chisubPi}{\Gamma} \ranglesubGamma
\\[15pt]
&=\ 
\langle \chisubpi, \chisubPi \ranglesubG\qquad \text{(cf. \#8)}.
\end{align*}
\end{x}
\vspace{0.3cm}

\begin{x}{\small\bf \un{N.B.}} \ 
It is this result that provides the link with the machinery developed in \S7 and \S8.
\end{x}
\vspace{0.3cm}

Suppose that $\Gamma_1  \subset  \Gamma_2 \subset G$ are subgroups.  
Let $(\theta_1,E_1)$ be a unitary representation of $\Gamma_1$ $-$then one can form $\Ind_{\Gamma_1, \theta}^G$.  
On the other hand, one can first form $\theta_2 = \Ind_{\Gamma_1,\theta_1}^{\Gamma_2}$ and then form $\Ind_{\Gamma_2, \theta_2}^G$.
\\[-.25cm]

\begin{x}{\small\bf INDUCTION IN STAGES} \ 
\[
\pi_1 
\ \equiv \ 
\Ind_{\Gamma_1, \theta_1}^G 
\ \approx \ 
\Ind_{\Gamma_2, \theta_2}^G
\ \equiv \ 
\pi_2.
\]

\vspace{0.2cm}

[Apply \S7, \#12: 
\[
\chisubpiOne
\ = \ 
\chisubpiTwo.]
\]
\vspace{0.2cm}

[Note: \ 
Characters determine representations up to equivalence (cf. \S5, \#10).]
\\[-.25cm]
\end{x}

\begin{x}{\small\bf LEMMA} \ 
If $(\theta_1,E_1)$, $(\theta_2,E_2)$ are unitary representations of $\Gamma$, then 
\[
\Ind_{\Gamma, \theta_1 \oplus \theta_2}^G 
\ \approx \ 
\Ind_{\Gamma, \theta_1}^G \ \oplus \ \Ind_{\Gamma, \theta_2}^G .
\]
\end{x}
\vspace{0.3cm}

\begin{x}{\small\bf \un{N.B.}} \ 
Consequently, $\Ind_{\Gamma, \theta}^G$ cannot be irreducible unless $\theta$ itself is irreducible (cf. \S10, \#3).
\\[-.5cm]
\end{x}

Let $G_1$, $G_2$ be finite groups, let $\Gamma_1 \subset G_1$, $\Gamma_2 \subset G_2$ be subgroups.

Put
\[
G 
\ = \ 
G_1 \times G_2, 
\quad 
\Gamma 
\ = \ 
\Gamma_1 \times \Gamma_2.
\]
\\[-1.00cm]

\begin{x}{\small\bf LEMMA} \ 
If
\[
\begin{cases}
\ \chisubOne \ \text{is a character of} \ \Gamma_1  \\
\ \chisubTwo \ \text{is a character of} \ \Gamma_2
\end{cases}
,
\]
then $\chisubOne \chisubTwo$ is a character of $\Gamma$ and 
\[
i_{_{\Gamma \ra G}} \chisubOne \chisubTwo 
\ = \ (i_{_{\Gamma_1 \ra G_1}} \chisubOne) \hsx (i_{_{\Gamma_2 \ra G_2}} \chisubTwo).
\]
\end{x}
\vspace{0.3cm}


\chapter{
$\boldsymbol{\S}$\textbf{10}.\quad  IRREDUCIBILITY OF $\Ind_{\Gamma,\theta}^G$}
\setlength\parindent{2em}
\setcounter{theoremn}{0}
\renewcommand{\thepage}{A II \S10-\arabic{page}}

\qquad Let \mG be a finite group.
\\[-.2cm]

\begin{x}{\small\bf DEFINITION} \ 
Let $(\pi_1,V_1)$, $(\pi_2,V_2)$ be unitary representations of \mG $-$then $\pi_1$ and $\pi_2$ are 
\un{disjoint} 
if they have no common nonzero unitarily equivalent subrepresentations.
\end{x}
\vspace{0.3cm}

\begin{x}{\small\bf LEMMA} \ 
$\pi_1$ and $\pi_2$ are disjoint iff $\chisubpiOne$ and $\chisubpiTwo$ are orthogonal:
\[
\langle
\chisubpiOne, \chisubpiTwo
\ranglesubG
\ = \ 0.
\]
\end{x}
\vspace{0.3cm}

\begin{x}{\small\bf THEOREM} \ 
Let $\Gamma$ be a subgroup of \mG, $(\theta,\E)$ an irreducible unitary representation of $\Gamma$ $-$then 
$\Ind_{\Gamma,\theta}^G$ is irreducible iff for every $x \in G - \Gamma$, the unitary representations
\[
\gamma \ra \theta(\gamma), 
\quad 
\gamma \ra \theta(x^{-1} \gamma x)
\]
of the subgroup
\[
\Gamma(x) 
\ = \ 
\Gamma^x \hsx \cap \hsx \Gamma 
\qquad (\Gamma^x \ = \ x \Gamma x^{-1})
\]
are disjoint.
\\[-.2cm]

PROOF \ 
Set $\pi = \Ind_{\Gamma,\theta}^G$ $-$then on general grounds, $\pi$ is irreducible iff 
$\langle \chisubpi,\chisubpi \ranglesubG = 1$ (cf. \S5, \#11).
\\[-.5cm]

\noindent I.e.: Iff
\[
\langle
i_{_{\Gamma \ra G}} \chisubtheta, i_{_{\Gamma \ra G}} \chisubtheta
\ranglesubG 
\ = \ 
1 \qquad \text{(cf. \S9, \#10)}
\]
or still, iff 
\[
\langle
\chisubtheta, r_{_{G \ra \Gamma}} (i_{_{\Gamma \ra G}} \chisubtheta)
\ranglesubGamma
\ = \ 
1 \qquad \text{(cf. \S7, \#6)}
\]
or still, iff 
\allowdisplaybreaks
\begin{align*}
\langle
\chisubtheta, \sum\limits_{s \in S} \hsx i_{_{\Gamma(s) \ra \Gamma}} (\chisubtheta)_s
\ranglesubGamma \ 
&=\ 
\sum\limits_{s \in S} \hsx
\langle \chisubtheta, i_{_{\Gamma(s) \ra \Gamma}} (\chisubtheta)_s \ranglesubGamma 
\\[15pt]
&=\ 
1 \qquad \text{(cf. \S8, \#3)}.
\end{align*}
Here $S = \Gamma \backslash G / \Gamma$ and it can be assumed that one element of the sum is $s = e$ in which case 
$(\chisubtheta)_s = \chisubtheta$, $\Gamma(s) = \Gamma$, hence
\allowdisplaybreaks
\begin{align*}
\langle \chisubpi,\chisubpi \ranglesubG  \ 
&=\ 
\langle \chisubtheta,\chisubtheta \ranglesubGamma
\ + \ 
\sum\limits_{\substack{s \in \Gamma \backslash G / \Gamma \\ s \notin \Gamma}} \hsx
\langle \chisubtheta, i_{_{\Gamma(s) \ra \Gamma}} (\chisubtheta)_s \ranglesubGamma
\\[11pt]
&=\ 
\langle \chisubtheta,\chisubtheta \ranglesubGamma
\ + \ 
\sum\limits_{\substack{s \in \Gamma \backslash G / \Gamma \\ s \notin \Gamma}} \hsx
\langle r_{_{\Gamma \ra \Gamma(s)}} \chisubtheta, (\chisubtheta)_s \rangle_{\Gamma(s)} 
\qquad \text{(cf. \S7, \#6)}
\\[11pt]
&=\ 
1 
\ + \ 
\sum\limits_{\substack{s \in \Gamma \backslash G / \Gamma \\ s \notin \Gamma}} \hsx
\langle r_{_{\Gamma \ra \Gamma(s)}} \chisubtheta, (\chisubtheta)_s \rangle_{\Gamma(s)} 
\qquad \text{(cf. \S5, \#11)}.
\end{align*}
Each term 
\[
\langle r_{_{\Gamma \ra \Gamma(s)}} \chisubtheta, (\chisubtheta)_s \rangle_{\Gamma(s)}
\]
is nonnegative and per $\Gamma(s)$, 
\[
\begin{cases}
\ r_{_{\Gamma \ra \Gamma(s)}} \chisubtheta \ \text{is the character of} \ \gamma \ra \theta(\gamma)\\[4pt]
\  (\chisubtheta)_s  \hspace{.4cm} \text{is the character of} \ \gamma \ra \theta(s^{-1}\gamma s)
\end{cases}
.
\]
If now $\pi = \Ind_{\Gamma,\theta}^G$ is irreducible, then $\langle \chisubpi,\chisubpi \ranglesubG = 1$, 
thus $\forall \ s \in \Gamma \backslash G / \Gamma$ $(s \notin \Gamma)$, 
\[
r_{_{\Gamma \ra \Gamma(s)}} \chisubtheta 
\quad \text{and} \quad
(\chisubtheta)_s 
\]
are orthogonal.  
Since \mS can be chosen so that it contains any given element of $G - \Gamma$, the disjointness claim is manifest.  
Conversely, the orthogonality of
\[
r_{_{\Gamma \ra \Gamma(s)}} \chisubtheta 
\quad \text{and} \quad
(\chisubtheta)_s 
\]
$\forall \ s \in \Gamma \backslash G / \Gamma$ $(s \notin \Gamma)$ forces 
$\langle \chisubpi, \chisubpi \ranglesubG = 1$.

\end{x}
\vspace{0.3cm}


\chapter{
$\boldsymbol{\S}$\textbf{11}.\quad  BURNSIDE RINGS}
\setlength\parindent{2em}
\setcounter{theoremn}{0}
\renewcommand{\thepage}{A II \S11-\arabic{page}}

\qquad Let \mG be a finite group.
\\[-.2cm]

\begin{x}{\small\bf DEFINITION} \ 
Let $\chi_{_1}, \ldots, \chi_{_t}$ be the characters of the irreducible unitary representations of \mG $-$then the 
\un{character ring}
\index{character ring}\index{$X(G)$} 
$X(G)$ is the free abelian group on generators $\chi_{_1}, \ldots, \chi_{_t}$ under pointwise addition and multiplication with unit 
$1_G$ (cf. \S8, \#5).
\vspace{0.2cm}

[Note: \ 
Recall that 
\[
t 
\ = \ 
\widehat{\abs{G}}
\ = \ 
\abs{\CON(G)}
\ = \ 
\dim CL(G).]
\]
\end{x}
\vspace{0.1cm}

\begin{x}{\small\bf \un{N.B.}} \ 
The pointwise sum or product of two characters is a character and the canonical arrow 
\[
X(G) \  \otimes_\Z \ \Cx \ra CL(G) 
\]
is an isomorphism.
\end{x}
\vspace{0.1cm}

\begin{x}{\small\bf DEFINITION} \ 
An element of $X(G)$ is called a 
\index{virtual character} 
\un{virtual character}.
\end{x}
\vspace{0.1cm}

\begin{x}{\small\bf LEMMA} \ 
A class function $f \in CL(G)$ is a virtual character iff 
$\langle f, \chisubPi \ranglesubG \in \Z$ for all $\Pi \in \widehat{G}$.
\end{x}
\vspace{0.1cm}

\begin{x}{\small\bf REMARK} \ 
The values of a virtual character are algebraic integers (cf. \S6, \#5), hence $X(G)$ is a proper subring of $CL(G)$.
\vspace{0.2cm}

[Note: \ On the other hand, a class function whose values are algebraic integers need not be a virtual character.]
\end{x}
\vspace{0.1cm}

\begin{x}{\small\bf NOTATION} \ 
Let $\sH$ be a collection of subgroups of \mG with the property that
\[
H \in \sH \ \& \ H^\prime \subset H \implies H^\prime \in \sH,
\]
in which case $\sH$ is termed a 
\un{hereditary class}.
\index{hereditary class}
\end{x}
\vspace{0.1cm}

Given $\sH$, let $X(G;\sH)$ be the additive subgroup of $X(G)$ spanned by the 
\[
i_{H \ra G}1_H \qquad (H \in \sH).
\]
\\[-1.25cm]

\begin{x}{\small\bf LEMMA} \ 
$X(G;\sH)$ is a subring of $X(G)$.
\\[-.2cm]

PROOF \ 
Let $H_1$, $H_2 \in \sH$ $-$then the claim is that
\[
(i_{_{H_1 \ra G}} 1_{H_1}) (i_{_{H_2 \ra G}} 1_{H_2}) \in \sH.
\]
Put
\[
\chi \ = \ (i_{_{H_1 \ra G}} 1_{H_1}) 
\]
and write
\allowdisplaybreaks
\begin{align*}
\chi(i_{_{H_2 \ra G}} 1_{H_2}) \ 
&=\ 
i_{_{H_2 \ra G}} ((r_{_{G \ra H_2}} \chi) 1_{H_2}) \qquad \text{(cf. \S7, \#8)}
\\[15pt]
&=\ 
i_{_{H_2 \ra G}} ((r_{_{G \ra H_2}} \chi)).
\end{align*}
Then, thanks to \S8, \#3, there are subgroups $K_1, \ldots, K_r$ of $H_2$ such that 
\[
r_{_{G \ra H_2}} \chi
\ = \ 
\sum\limits_{\ell = 1}^r \hsx 
i_{_{K_\ell \ra H_2}} 1_{K_\ell}.
\]
Therefore
\allowdisplaybreaks
\begin{align*}
i_{_{H_2 \ra G}} ((r_{_{G \ra H_2}} \chi)) \ 
&=\ 
i_{_{H_2 \ra G}} \bigg( 
\sum\limits_{\ell = 1}^r \hsx 
i_{_{K_\ell \ra H_2}} 1_{K_\ell}
\bigg)
\\[11pt]
&=\ 
\sum\limits_{\ell = 1}^r \hsx 
i_{_{H_2 \ra G}}(i_{_{K_\ell \ra H_2}} 1_{K_\ell})
\\[11pt]
&=\ 
\sum\limits_{\ell = 1}^r \hsx 
i_{_{K_\ell \ra G}} 1_{K_\ell} 
\qquad \text{(cf. \S7, \#12)}
\\[11pt]
&
\in X(G;\sH).
\end{align*}
\end{x}
\vspace{0.1cm}

\begin{x}{\small\bf DEFINITION} \ 
$X(G; \sH)$ is the 
\un{Burnside ring}
\index{Burnside ring} 
of \mG associated with the hereditary class $\sH$.
\vspace{0.1cm}

[Note: \ 
It is not a priori evident that $1_G \in X(G;\sH)$.]
\end{x}
\vspace{0.1cm}

\begin{x}{\small\bf CRITERION} \ 
Let \mR be a ring of $\Z$-valued functions on a finite set \mX under pointwise operations.  
Suppose that for each $x \in X$ and for each prime $p$ there exists $f \in R$ such that 
$f(x) \not\equiv 0 \ \modx \hsx p$ $-$then 
$1_X \in R$.
\vspace{0.2cm}

[Attach to each $x \in X$ the ideal
\[
I_x 
\ = \ \{f(x): f \in R\} \subset \Z.
\]
Then, in view of the assumption, $I_x = \Z$ so there exsits $f_x \in R$ such that $f_x(x) = 1$, 
hence
\[
\prod\limits_{x \in X} \hsx (1 - f_x) 
\ = \ 
0.
\]
Now expand the product to get 1 as a sum of elements of \mR.]
\end{x}
\vspace{0.1cm}

Let \mG be a finite group.
\\[-.2cm]

\begin{x}{\small\bf DEFINITION} \ 
Let $p$ be a prime $-$then \mG is a 
\un{$p$-group}
\index{$p$-group} 
if every element $x \in G$ has order a power of $p$.
\vspace{0.2cm}

[Note: \ Every $p$-group is nilpotent.]
\end{x}
\vspace{0.1cm}


\begin{x}{\small\bf LEMMA} \ 
\mG is a $p$-group iff $\abs{G}$ is a power of $p$.
\end{x}
\vspace{0.1cm}

\begin{x}{\small\bf DEFINITION} \ 
Let $p$ be a prime $-$then a subgroup \mP of \mG is a 
\un{Sylow $p$-subgroup}
\index{Sylow $p$-subgroup} 
of \mG if it is a maximal $p$-subgroup of \mG.
\end{x}
\vspace{0.1cm}

\begin{x}{\small\bf THEOREM} \ 
\\[-.2cm]

\qquad \textbullet \quad
Sylow $p$-subgroups exist.
\\[-.2cm]

\qquad \textbullet \quad
All Sylow $p$-subgroups are conjugate.
\\[-.2cm]

\qquad \textbullet \quad 
Every $p$-subgroup is contained in a Sylow $p$-subgroup.
\end{x}
\vspace{0.1cm}

\begin{x}{\small\bf \un{N.B.}} \ 
The number of Sylow $p$-subgroups of \mG is a divisor of $\abs{G}$.
\end{x}
\vspace{0.1cm}

\begin{x}{\small\bf DEFINITION} \ 
Given a prime $p$, a finite group \mH is 
\un{$p$-elementary}
\index{$p$-elementary} 
if it is the direct product of a cyclic group \mC of order prime to $p$ and a $p$-group \mP. 
\vspace{0.2cm}

[Note: \ 
Accordingly, \mC and \mP are normal  subgroups, $C \hsx \cap \hsx P = \{e\}$, and $H = C P$.]
\end{x}
\vspace{0.1cm}

\begin{x}{\small\bf LEMMA} \ 
Subgroups of $p$-elementary groups are again $p$-elementary, hence the $p$-elementary subgroups of \mG constitute a hereditary class $\sE_p(G)$.
\end{x}
\vspace{0.1cm}

\begin{x}{\small\bf DEFINITION} \ 
A finite group \mH is 
\un{elementary}
\index{elementary} 
if it is $p$-elementary for some prime $p$.
\end{x}
\vspace{0.1cm}

\begin{x}{\small\bf NOTATION} \ 
Put 
\[
\sE(G) 
\ = \ 
\bigcup\limits_p \hsx \sE_p(G).
\]
\end{x}
\vspace{0.1cm}

\begin{x}{\small\bf \un{N.B.}} \ 
Since \ $\sE(G)$ \ is a hereditary class, one can form its Burnside ring \ $X(G;\sE(G))$.
\end{x}
\vspace{0.1cm}

\begin{x}{\small\bf DEFINITION} \ 
Given a prime $p$, a group \mH is 
\un{$p$-semielementary}
\index{$p$-semielementary} 
if it is the semidirect product of a cyclic subgroup \mC of order prime to $p$ and a $p$-group \mP.
\vspace{0.2cm}

[Note: \ 
Accordingly, \mC is a normal subgroup, $C \hsx \cap \hsx P = \{e\}$, and $H = C P$.]
\end{x}
\vspace{0.1cm}

\begin{x}{\small\bf LEMMA} \ 
Subgroups of $p$-semielementary groups are again $p$-semielementary, hence the $p$-semielementary subgroups of \mG constitute a hereditary class $\sS \sE_p (G)$.
\end{x}
\vspace{0.1cm}

\begin{x}{\small\bf DEFINITION} \ 
A finite group \mH is 
\un{semielementary}
\index{semielementary} 
if it is $p$-semielementary for some prime $p$.
\end{x}
\vspace{0.1cm}

\begin{x}{\small\bf NOTATION} \ 
Put 
\[
\sS \sE (G) 
\ = \ 
\bigcup\limits_p \ \sS \sE_p (G).
\]
\end{x}
\vspace{0.1cm}

\begin{x}{\small\bf \un{N.B.}} \ 
Since $\sS \sE (G) $ is a hereditary class, one can form its Burnside ring 
$X(G; \sS \sE (G))$.
\end{x}
\vspace{0.1cm}

\begin{x}{\small\bf LEMMA} \ 
\[
1_G \in X(G; \sS \sE (G)),
\]
i.e., there exist integers $a_H(H \in \sS \sE (G))$ such that 
\[
1_G 
\ = \ 
\sum\limits_{H \in \sS \sE (G)} \ 
a_H(i_{_{H \ra G}} 1_H).
\]
\vspace{0.2cm}

PROOF \ 
It suffices to show that the ring $X(G; \sS \sE (G))$ satisfies the assumptions of \#9:  
For every $x \in G$ and for every prime $p$, there exists a group 
$H_{x,p} \equiv H \in \sS \sE (G)$ such that
\[
i_{_{H \ra G}} 1_H (x) \ \not\equiv 0 \ \modx p.
\]
This said, factor the order of $x$ as $p^an$ 
$(p {\not\hspace{.095cm}\mid} \hsx n)$ 
and let $C = \langle x^{p^a} \rangle$ (hence $\abs{C} = n$,
hence is prime to $p$).  
Let \mN be the normalizer of \mC in \mG, let \mP be a Sylow $p$-subgroup of \mN containing $x$, and let 
$H_{x, p} \equiv H = C P$ $-$then \mH is $p$-semielementary and the claim is that 
\[
(i_{_{H \ra G}} 1_H) (x) \ \not\equiv 0 \ \modx p.
\]
By definition,
\[
(i_{_{H \ra G}} 1_H) (x) 
\ = \ 
\frac{1}{\abs{H}} \hsx 
\sum\limits_{\substack{y \in G,\\  yxy^{-1} \in H}} \ 
1_H (y x y^{-1}).
\]
But
\allowdisplaybreaks
\begin{align*}
y x y^{-1} \in H \ 
&\implies 
y C y^{-1} \subset H
\\[15pt]
&\implies 
y C y^{-1} \hsx = \hsx C
\\[15pt]
&\implies 
y \in N.
\end{align*}
Therefore
\[
(i_{_{H \ra G}} 1_H) (x) 
\ = \ 
(i_{_{H \ra N}} 1_H) (x),
\]
the term on the right being the number of left cosets of \mH in \mN fixed by $x$ (cf. \S7, \#11).  
Since \mC is a normal subgroup of \mN and since $C \subset H$, it follows that \mC must fix the left cosets of \mH in \mN.  
Thus the $x$-orbits have cardinality dividing $p^a$, thus each nontrivial $x$-orbit has cardinality divisible by $p$.  
On the other hand, the number of left cosets of \mH in \mN is prime to $p$ ($H$ contains a Sylow $p$-subgroup of \mN).  
Combining these facts then leads to the conclusion that the number of left cosets of \mH in \mN fixed by $x$ is prime to $p$, 
i.e., 
\[
(i_{_{H \ra G}} 1_H) (x)  \ \not\equiv \ 0 \ \modx \hsx p. 
\]
\end{x}
\vspace{0.1cm}

\begin{x}{\small\bf DEFINITION} \ 
\ A 
\un{monomial character}
\index{monomial character} 
of a finite\  group is a \ character of degree 1.
\end{x}
\vspace{0.1cm}


\begin{x}{\small\bf DEFINITION} \ 
A finite group \mH is said to be an 
\un{\mM-group}
\index{\mM-group} 
if each irreducible character of \mH is induced by a monomial character of a subgroup of \mH.
\end{x}
\vspace{0.1cm}

\begin{x}{\small\bf THEOREM} \ 
Suppose that \mH is a finite group which is a semidirect product of an abelian normal subgroup of a nilpotent group 
(in particular, a $p$-group) $-$then \mH is an \mM-group.
\end{x}
\vspace{0.1cm}

\begin{x}{\small\bf APPLICATION} \ 
$p$-elementary groups and $p$-semielementary groups are \mM-groups.
\end{x}
\vspace{0.1cm}


\chapter{
$\boldsymbol{\S}$\textbf{12}.\quad  BRAUER THEORY}
\setlength\parindent{2em}
\setcounter{theoremn}{0}
\renewcommand{\thepage}{A II \S12-\arabic{page}}

\qquad Let \mG be a finite group.
\\[-.2cm]

\begin{x}{\small\bf CHARACTERIZATION OF CHARACTERS} \ 
A class function $f \in CL(G)$ is a virtual character (i.e., belongs to $X(G))$ iff for every $H \in \sE(G)$, 
\[
r_{_{G \ra H}} f \in X(H).
\]
\end{x}
\vspace{0.1cm}

\begin{x}{\small\bf INDUCTION PRINCIPLE} \ 
A class function $f \in CL(G)$ is a virtual character (i.e., belongs to $X(G))$ iff there exist elementary subgroups $H_i$, 
monomial characters $\lambda_i$ of $H_i$, and integers $a_i$ $(1 \leq i \leq n)$ such that 
\[
f 
\ = \ 
\sum\limits_{i=1}^n \hsx 
a_i (i_{_{H_i \ra G}} \lambda_i).
\]
\end{x}
\vspace{0.1cm}

These are the main results.  
Turning to their proofs, let $\sR$ be the ring with unit $1_G$ whose elements are the class functions $f$ on \mG such that
\[
r_{_{G \ra H}} f \in X(H)
\]
for all $H \in \sE(G)$ and let $\sL$ be the subgroup of $X(G)$ spanned over $\Z$ by characters of the form 
$i_{_{H \ra G}} \lambda$, where $\lambda$ is a monomial character for some $H \in \sE(G)$.
\vspace{0.1cm}

\begin{x}{\small\bf LEMMA} \ 
Statements 1 and 2 are equivalent to $\sL = \sR$.
\vspace{0.2cm}

[Note: \ 
Obviously, 
\[
\sL \subset X(G) \subset \sR.]
\]
\end{x}
\vspace{0.1cm}

\begin{x}{\small\bf LEMMA} \ 
$\sL$ is an ideal in $\sR$.
\vspace{0.1cm}

PROOF \ 
Let $\Lambda \in \sL$, say
\[
\Lambda 
\ = \ 
\sum\limits_i \hsx 
a_i (i_{_{H_i \ra G}} \lambda_i),
\]
and let $\psi \in \sR$ $-$then 
\allowdisplaybreaks
\begin{align*}
\psi \Lambda \ 
&=\ 
\sum\limits_i \hsx 
a_i  \hsx \psi (i_{_{H_i \ra G}} \lambda_i)
\\[11pt]
&=\ 
\sum\limits_i \hsx 
a_i (i_{_{H_i \ra G}} ((r_{_{G \ra H_i}} \psi)\lambda_i)) \qquad \text{(cf. \S7, \#8)}.
\end{align*}
Since
\[
r_{_{G \ra H_i}} \psi \in X(H_i), 
\]
there exist integers $b_{i j}$ such that 
\[
r_{_{G \ra H_i}} \psi
\ = \ 
\sum\limits_j \hsx 
b_{i j} \xi_{i j},
\]
$\xi_{i j}$ running through the irreducible characters of $H_i$, hence
\[
\psi \Lambda 
\ = \ 
\sum\limits_{i, j} \hsx 
a_i b_{i j} \hsx (i_{_{H_i \ra G}} \xi_{i j}).
\]
But elementary groups are $M$-groups (cf. \S11, \#29), so $\xi_{i j}$ is induced by a monomial character of some subgroup of $H_i$.  
Taking into account that $\sE(G)$ is a hereditary class, apply \S7, \#12 to conclude that $\psi \Lambda \in \sL$.  
Therefore $\sL$ is an ideal in $\sR$.  
\vspace{0.2cm}

[Note: \ 
Operations in $\sR$ are pointwise and, of course, $\sR$ is commutative.]
\end{x}
\vspace{0.1cm}

Matters thus reduce to showing that $1_G \in \sL$.  To this end, suppose that it were possible to write
\[
1_G 
\ = \ 
\sum\limits_k \hsx 
c_k (i_{_{H_k \ra G}}\chisubk), 
\]
where $c_k \in \Z$ and $\chisubk$ is a character of some proper subgroup $H_k$ of \mG.  
Inductively, it can be assumed that \#2 holds for $H_k$, hence that 
$\chisubk$ can be written as a $\Z$-linear combination of induced monomial characters from elements of $\sE(H_k)$.  
But then $1_G \in \sL$, as desired.
\vspace{0.1cm}

[Note: \ 
Nothing need be done if \mG is elementary to begin with (it being automatic that $1_G \in \sL$).]

\begin{x}{\small\bf LEMMA} \ 
If \mG is not elementary, then $1_G$ can be written as a $\Z$-linear combination of induced characters from proper subgroups of \mG.
\vspace{0.1cm}

\un{Case 1}: \ 
\mG is not semielementary, thus $G \notin \sS\sE(G)$ and the $H \in \sS\sE(G)$ are proper subgroups.  
The contention then follows from \S11, \#25.
\vspace{0.1cm}

\un{Case 2}: \ 
\mG is semielementary: \ $G \in \sS\sE(G)$, say $G = C P$ for some prime $p$.  
Let \mN be the normalizer of \mP in \mG, hence $N = (C \hsx \cap \hsx N) \times P$ is $p$-elementary and it can be assumed that $N \neq G$ (otherwise \mG is elementary and there is nothing to prove).  
Write
\[
i_{_{N \ra G}} 1_N
\ = \ a_0 1_G + 
\sum\limits_{i > 0} \hsx 
a_i \chisubi,
\]
where the $\chisubi \neq 1_G$ are irreducible characters and the $a_i$ are positive integers.
\end{x}
\vspace{0.1cm}

\begin{x}{\small\bf \un{N.B.}} \ 
\allowdisplaybreaks
\begin{align*}
a_0 \ 
&=\ 
\langle i_{_{N \ra G}} 1_N, 1_G \ranglesubG
\\[15pt]
&=\ 
\langle 1_N, r_{_{G \ra N}} 1_G \ranglesubN \qquad \text{(cf. $\S7$, $\#6$)}
\\[15pt]
&=\ 
\langle 1_N, 1_N \ranglesubN
\\[15pt]
&=\ 
1.
\end{align*}
\end{x}
\vspace{0.1cm}

\begin{x}{\small\bf \un{N.B.}} \ 
$\chisubi(e) > 1$ for all $i > 0$.
\\[-.2cm]

[Suppose that $\chisubi(e) = 1$ $(\exists \ i)$.  
Write
\[
\res_{_{G \ra N}} \chisubi 
\ = \ 
c 1_N + \chi
\]
for some character $\chi$ orthogonal to $1_N$ $-$then 
\allowdisplaybreaks
\begin{align*}
c \ 
&=\ 
\langle 1_N, \res_{_{G \ra N}} \chisubi \ranglesubN
\\[6pt]
&=\ 
\langle 1_{_{N \ra G}}1_N, \chisubi \ranglesubG 
\qquad \text{(cf. \S7, \#6)}
\\[6pt]
&=\ 
a_i.
\end{align*}
And
\allowdisplaybreaks
\begin{align*}
1 \ 
&=\ 
\chisubi(e) 
\\[6pt]
&=\ 
a_i + \chisubi(e)
\end{align*}
\qquad \qquad $\implies$
\\[-.9cm]
\[
a_i  
\ = \ 
1
\]
\qquad \qquad $\implies$
\\[-.9cm]
\[
\res_{_{G \ra N}} \chisubi
\ = \ 
1_N.
\]
Recall now that the kernel $K_i$ of $\chisubi$ is the proper normal subgroup of \mG consisting of those $x \in G$ such that 
$\chisubi(x) = \chisubi(e)$ or still, consisting of those $x \in G$ such that $\chisubi(x) = 1$, thus $N \subset K_i$ 
(since $\res_{_{G \ra N}} \chisubi = 1_N$).  
But this is impossible:  $P$ is a Sylow $p$-subgroup of $K_i$, so $G = K_i$ (cf. infra).
\vspace{0.1cm}

[Note: \ 
Let $x \in G$ $-$then both \mP and $x P x^{-1}$ are Sylow $p$-subgroups of $K_i$, hence
\[
k x P x^{-1} k^{-1} 
\ = \ 
P
\]
for some $k \in K_i$ which implies that $kx \in N \subset K_i$, thereby forcing $x \in K_i$, so $G = K_i$.]
\end{x}

Return to the formula
\[
i_{_{N \ra G}} 1_N 
\ = \ 
1_G + 
\sum\limits_{i > 0} \hsx 
a_i \chisubi.
\]
Since $\chisubi(e) > 1$ for all $i > 0$, the $\chisubi$ are not monomial.  
On the other hand, $G = C P$ is semielementary, 
thus is an \mM-group, thus each $\chisubi$ is induced by a monomial character $\lambda_i$ of some proper subgroup $H_i$ of \mG.  
Therefore
\[
1_G 
\ = \ 
i_{_{N \ra G}} 1_N - 
\sum\limits_{i > 0} \hsx 
a_i (i_{_{{H_i} \ra G}} \lambda_i),
\]
which completes the proof of \#5.


\chapter{
$\boldsymbol{\S}$\textbf{13}.\quad  GROUPS OF LIE TYPE}
\setlength\parindent{2em}
\setcounter{theoremn}{0}
\renewcommand{\thepage}{A II \S13-\arabic{page}}

\qquad Let $k$ be a finite field.
\\[-.2cm]

\begin{x}{\small\bf DEFINITION} \ 
A \un{$k$-group} is a linear algebraic group defined over $k$.
\\[-.2cm]

[Note: \ 
A \un{$k$-subgroup} of a $k$-group is a subgroup which is a $k$-group.]
\end{x}
\vspace{0.1cm}

\begin{x}{\small\bf NOTATION} \ 
Given $k$-groups $\uun{A}$, $\uun{B}$, $\uun{C}$, \ldots, denote their group of $k$-rational points 
$\uun{A}(k)$, $\uun{B}(k)$, $\uun{C}(k)$, \ldots, by $A$, $B$, $C$, \ldots \hsx .
\end{x}
\vspace{0.1cm}

Let $\uun{G}$ be a connected reductive $k$-group.
\\[-.2cm]

\begin{x}{\small\bf DEFINITION} \ 
\mG is said to be a \un{group of Lie type}.
\end{x}
\vspace{0.1cm}

\begin{x}{\small\bf \un{N.B.}} \ 
\mG is, of course, finite and it is possible to compute $\abs{G}$ explicitly.
\end{x}
\vspace{0.1cm}

\begin{x}{\small\bf DEFINITION} \ 
A maximal closed connected solvable subgroup of $\uun{G}$ is called a \un{Borel subgroup}.
\\[-.25cm]

[Note: \ 
The conditions ``closed'' and ``connnected'' can be omitted from the definition.]
\end{x}
\vspace{0.1cm}

\begin{x}{\small\bf LEMMA} \ 
\\[-.2cm]

\qquad \textbullet \quad 
Any two Borel subgroups of $\uun{G}$ are conjugate.
\\[-.2cm]

\qquad \textbullet \quad 
Every element of $\uun{G}$ belongs to some Borel subgroup of $\uun{G}$.
\\[-.2cm]

\qquad \textbullet \quad 
Every closed subgroup of $\uun{G}$ containing a Borel subgroup is equal to its own normalizer and is connected.
\\[-.2cm]

\qquad \textbullet \quad 
Any two closed subgroups of $\uun{G}$ containing the same Borel subgroup and conjugate in $\uun{G}$ are equal.
\end{x}
\vspace{0.1cm}

\begin{x}{\small\bf \un{N.B.}}  \ 
Since $k$ is finite, $\uun{G}$ is quasi-split, hence contains a Borel subgroup defined over $k$.
\\[-.25cm]

[Note: \ 
Any two such are $\uun{G}$-conjugate.]
\end{x}
\vspace{0.1cm}

Let $\uun{B}$ be a Borel $k$-subgroup of $\uun{G}$, 
let $\uun{T} \subset \uun{B}$ be a maximal torus of $\uun{G}$ defined over $k$, 
and put
\[
\uun{N} 
\ = \ 
N_{\uun{G}}(\uun{T}).
\]
\\[-.75cm]

\begin{x}{\small\bf LEMMA} \ 
$\uun{N}$ is a $k$-subgroup of $\uun{G}$.
\end{x}
\vspace{0.1cm}

\begin{x}{\small\bf NOTATION} \ 
Set
\[
\uun{W}
\ = \ 
\uun{N} / \uun{T}.
\]
\end{x}
\vspace{0.1cm}

\begin{x}{\small\bf LEMMA} \ 
\[
W 
\ \approx \ 
N / T.
\]

[Note: \ 
\[
\uun{B} \cap \uun{N} = \uun{T}
\ \implies \ 
B \cap N = T.]
\]
\end{x}
\vspace{0.1cm}

\begin{x}{\small\bf LEMMA} \ 
\mW is a finite Coxeter group.
\\[-.2cm]

[Note: \ 
Spelled out, \mW admits a finite system of generators 
$w_1, \ldots, w_\ell$ $(w_i \neq 1$ and $w_i \neq w_j$ for $i \neq j$) 
subject to the relations
\[
w_i^2 = 1, 
\quad 
(w_i w_j)^{m_{i j}} = 1 
\qquad (i \neq j),
\]
where $m_{i j}$ is the order of $w_i w_j$ $(i \neq j)$.]
\end{x}
\vspace{0.1cm}

\begin{x}{\small\bf BRUHAT LEMMA} \ 
\[
G 
\ = \ 
\coprod\limits_{w \in W} \hsx 
B w B.
\]
\end{x}
\vspace{0.1cm}


\begin{x}{\small\bf DEFINITION} \ 
A closed subgroup $\uun{P}$ of $\uun{G}$ is \un{parabolic} if it contains a Borel subgroup of $\uun{G}$.
\end{x}
\vspace{0.1cm}

\begin{x}{\small\bf LEMMA} \ 
Let $\uun{P}_{\hsx 1}$, $\uun{P}_{\hsx 2}$ be parabolic $k$-subgroups of $\uun{G}$ $-$then 
$\uun{P}_{\hsx 1} = \uun{P}_{\hsx 2}$ iff $P_1 = P_2$.
\end{x}
\vspace{0.1cm}

\begin{x}{\small\bf NOTATION} \ 
Given a parabolic $k$-subgroup of $\uun{G}$, denote its unipotent radical by $\uun{U}$.
\\[-.25cm]

[Note: \ 
Recall that $\uun{P}$ is the normalizer of $\uun{U}$.]
\end{x}
\vspace{0.1cm}


\begin{x}{\small\bf DEFINITION} \ 
Let $\uun{P}$ be a parabolic $k$-subgroup of $\uun{G}$ $-$then a closed connected reductive $k$-subgroup 
$\uun{L}$ of $\uun{P}$ is a 
\un{Levi subgroup} of $\uun{P}$ if $\uun{P}$ is the semidirect product 
$\uun{L \hsx U}$ (hence $P = L U$).
\end{x}
\vspace{0.1cm}

\begin{x}{\small\bf LEMMA} \ 
Levi subgroups of $\uun{P}$ exist and any two such are conjugate by a unique element of $U$.  
\end{x}
\vspace{0.1cm}

\begin{x}{\small\bf \un{N.B.}} \ 
\mL is a group of Lie type.
\end{x}
\vspace{0.1cm}

\begin{x}{\small\bf LEMMA} \ 
Let $\uun{P}_{\hsx 1}$, $\uun{P}_{\hsx 2}$ be parabolic $k$-subgroups of $\uun{G}$ $-$then the following conditions are equivalent.
\\[-.2cm]

\qquad \textbullet \quad
$P_1 \cap U_2 \subset U_1$, \quad $P_2 \cap U_1 \subset U_2$
\\[-.2cm]

\qquad \textbullet \quad
$\uun{P}_{\hsx 1}$ and $\uun{P}_{\hsx 2}$ have a common Levi subgroup.
\end{x}
\vspace{0.1cm}

\begin{x}{\small\bf APPLICATION} \ 
\[
U_1 = U_2 
\ \implies \ 
\uun{P}_{\hsx 1} = \uun{P}_{\hsx 2}.
\]

[For under these circumstances, $\uun{P}_{\hsx 1}$ and $\uun{P}_{\hsx 2}$ have a common Levi subgroup $\uun{L}$, thus
\[
P_1 
\ = \ 
L U_1 
\ = \ 
L U_2 
\ = \ 
P_2,
\]
so one can quote $\#14$.]
\end{x}
\vspace{0.1cm}

\begin{x}{\small\bf DEFINITION} \ 
Let $\uun{P}_{\hsx 1}$, $\uun{P}_{\hsx 2}$ be parabolic $k$-subgroups of $\uun{G}$ $-$then
$\uun{P}_{\hsx 1}$ and  $\uun{P}_{\hsx 2}$
are said to be \un{associate} if there exists an $x \in G$ such that 
$\uun{P}_{\hsx 1}$ and $x \uun{P}_{\hsx 2} x^{-1}$ have a common Levi subgroup.
\end{x}
\vspace{0.1cm}

\begin{x}{\small\bf \un{N.B.}} \ 
The relation determined by ``to be associate'' is an equivalence relation on the set of parabolic $k$-subgroups of $\uun{G}$.
\end{x}
\vspace{0.1cm}

\begin{x}{\small\bf LEMMA} \ 
If $\uun{P}_{\hsx 1}$, $\uun{P}_{\hsx 2}$ are not associate, then $\forall \ x \in G$, 
$\uun{P}_{\hsx 1}$, $x \uun{P}_{\hsx 2} x^{-1}$ are not associate.  

[If there exists $x \in G$ such that $\uun{P}_{\hsx 1}$ and $x \uun{P}_{\hsx 2} x^{-1}$ are associate, then there exists 
$y \in G$ such that 
$\uun{P}_{\hsx 1}$ and $y x \uun{P}_{\hsx 2} x^{-1} y^{-1}$ have a common Levi subgroup, thus 
$\uun{P}_{\hsx 1}$ and  $\uun{P}_{\hsx 2}$ are associate, contradiction.]
\end{x}
\vspace{0.1cm}

\begin{x}{\small\bf LEMMA} \ 
Let $\uun{P}_{\hsx 1}$, $\uun{P}_{\hsx 2}$ be parabolic $k$-subgroups of $\uun{G}$.  
Assume: 
$\uun{P}_{\hsx 1}$ and $\uun{P}_{\hsx 2}$ are associate $-$then 
$\abs{P_1} = \abs{P_2}$.

[There is no loss of generality in supposing that 
$\uun{P}_{\hsx 1}$ and  $\uun{P}_{\hsx 2}$ have a common Levi subgroup $\uun{L}$, 
thereby reducing matters to the claim that 
$\abs{U_1} = \abs{U_2}$.]
\end{x}
\vspace{0.1cm}

\begin{x}{\small\bf DESCENT} \ 
Fix a parabolic $k$-subgroup $\uun{P} \subset \uun{G}$ and let $\uun{L} \subset \uun{P}$ be a Levi subgroup 
$-$then 
there is a 1-to-1 correspondence between the set of parabolic $k$-subgroups of $\uun{G}$ contained in $\uun{P}$ 
and the set of parabolic $k$-subgroups of $\uun{L}$.
\\[-.2cm]

\qquad \textbullet \quad
Given a parabolic $k$-subgroup $\uun{P}^\prime \subset \uun{P}$, write
$\uun{P}^\prime = \uun{L}^\prime \hsx \uun{U}^\prime$ and put
$^*\uun{P} = \uun{P}^\prime \cap \uun{L}$ 
$-$then
$^*\uun{P}$ is a parabolic $k$-subgroup of $\uun{L}$ with unipotent radical 
$^*\uun{U} = \uun{U}^\prime \cap \uun{L}$.
\\[-.2cm]

\qquad \textbullet \quad
Given a parabolic $k$-subgroup $^*\uun{P}$ of $\uun{L}$, 
write
$^*\uun{P} = \hsy ^*\uun{L} \hsx ^*\uun{U}$ and put
$\uun{L}^\prime = \hsx ^*\uun{L}$,
$\uun{U}^\prime = \hsx ^*\uun{U} \hsx \uun{U}$ 
$-$then 
$\uun{P}^\prime = \uun{L}^\prime \hsx \uun{U}^\prime$ 
is a parabolic $k$-subgroup of $\uun{G}$ such that 
$\uun{P}^\prime \subset \uun{P}$.  
\\[-.2cm]

The bijection in question is the assignment 
$\uun{P}^\prime \ra \hsx ^*\uun{P}$.
\end{x}
\vspace{0.1cm}

\begin{x}{\small\bf \un{N.B.}} \ 
$\uun{P}^\prime$ and $\uun{P}^{\prime\prime}$ are conjugate by an element of \mG iff $\uun{P}^\prime \cap L$ and 
$\uun{P}^{\prime\prime} \cap L$ are conjugate by an element of \mL.  
\end{x}
\vspace{0.1cm}

\[
\text{APPENDIX}
\]

\qquad {\small\bf LEMMA} \ 
Suppose that 
$\uun{P}_{\hsx 1} = \uun{L}_1 \hsx \uun{U}_1$ 
and 
$\uun{P}_{\hsx 2} = \uun{L}_2 \hsx \uun{U}_2$ are associate 
$-$then $\uun{L}_1$ and $\uun{L}_2$ are conjugate by an element of \mG.
\\[-.2cm]

[Choose $x \in G$ such that $\uun{P}_{\hsx 1}$ and $x \uun{P}_{\hsx 2} x^{-1}$ have a common Levi subgroup $\uun{L}$.  
Choose $u_1 \in U_1$: 
\[
u_1 \uun{L} u_1^{-1} 
\ = \ 
\uun{L}_1.
\]
Choose $xu_2x^{-1} \in \ xU_2 x^{-1}$:
\[
x u_2 x^{-1} \uun{L} x u_2^{-1} x^{-1} 
\ = \ 
x \uun{L}_2 x^{-1}.
\]
Then 
\[
u_2 x^{-1} \uun{L} x u_2^{-1}
\ = \ 
\uun{L}_2
\]
\qquad $\implies$
\[
\uun{L} 
\ = \ 
xu_2^{-1} \uun{L}_2 u_2 x^{-1}
\]
\qquad $\implies$
\allowdisplaybreaks
\begin{align*}
\uun{L}_1 \ 
&=\ 
u_1 \uun{L} u_1^{-1} 
\\[11pt]
&=\ 
u_1 x u_2^{-1} \uun{L}_2 u_2 x^{-1} u_1^{-1}.]
\end{align*}


\chapter{
$\boldsymbol{\S}$\textbf{14}.\quad  HARISH-CHANDRA THEORY}
\setlength\parindent{2em}
\setcounter{theoremn}{0}
\renewcommand{\thepage}{A II \S14-\arabic{page}}

\qquad Let $k$ be a finite field, \un{\un{\mG}} a connected reductive $k$-group.
\\[-.2cm]

\begin{x}{\small\bf DEFINITION} \ 
Let \uun{\mP} be a parabolic $k$-subgroup of \uun{\mG} $-$then \mP is termed a 
\un{cuspidal subgroup}
\index{cuspidal subgroup} 
of \mG.
\end{x}
\vspace{0.1cm}

\begin{x}{\small\bf NOTATION} \ 
Given a cuspidal subgroup $P = L U$ of \mG and an $f \in C(G)$, let
\[
f_P(x) \ = \ \sum\limits_{u \in U} \hsx f(xu) \qquad (x \in G).
\]
\vspace{0.1cm}

[Note: \ If $P = G$, then 
\[
f_G(x) \ = \  f(x) \qquad (x \in G).]
\]
\end{x}
\vspace{0.1cm}

\begin{x}{\small\bf DEFINITION} \ 
Let $f \in C(G)$ $-$then $f$ is said to be a 
\un{cusp form}
\index{cusp form} 
if $f_P = 0$ for all $P \neq G$.
\end{x}
\vspace{0.1cm}

\begin{x}{\small\bf NOTATION} \ 
Write  $^0 C(G)$ for the set of cusp forms and put
\[
^0 CL(G) \ = \ CL(G) \hsx \cap \hsx ^0 C(G).
\]
\end{x}
\vspace{0.1cm}

\begin{x}{\small\bf LEMMA} \ 
 $^0 C(G)$ is a linear subspace of $C(G)$.
\end{x}
\vspace{0.1cm}

\begin{x}{\small\bf LEMMA} \ 
 $^0 C(G)$ is stable under left translations, hence is a left ideal in $C(G)$.
\end{x}
\vspace{0.1cm}

\begin{x}{\small\bf REMARK} \ 
If \uun{\mG} is a torus, then  $^0 C(G) = C(G)$.
\end{x}
\vspace{0.1cm}

\begin{x}{\small\bf NOTATION} \ 
Given $f \in C(G)$, write $f_P \sim 0$ if 
\[
\sum\limits_{\ell \in L} \hsx f_P(x \ell) \hsx \ov{\phi(\ell)} \hsx = \hsx0
\]
for all $\phi \in \hsx ^0 C(L)$ and all $x \in G$.
\vspace{0.1cm}

[Note: \  Bear in mind that \mL is a group of Lie type (cf. \S13, \#19).]
\end{x}
\vspace{0.1cm}


\begin{x}{\small\bf \un{N.B.}} \ 
Matters are independent of the choice of \mL in \mP.
\end{x}
\vspace{0.1cm}

\begin{x}{\small\bf LANGLANDS PRINCIPLE} \ 
If $f_P \sim 0$ for all cuspidal subgroups \mP of \mG (including $P = G$), then $f = 0$.
\vspace{0.2cm}

PROOF \  Proceed by induction on the semisimple $k$-rank $s$ of \uun{\mG}, 
the case $s = 0$ being trivial (because then \uun{\mG} is anisotropic, there is only one \mP, viz. $P = G$, 
and $L = G$, $^0C(L) = C(G) \ldots$).  
So assume that $s$ is positive and let $P = L U$ be for the moment a proper cuspidal subgroup, thus $U \neq \{e\}$ 
and the semisimple $k$-rank of $\uun{L}$ is strictly smaller than that of $\uun{G}$.  
Using now \S13, \#26, let $^*P = \hsx ^*L ^*U$ be a cuspidal subgroup of \mL $-$then 
$P^\prime = L^\prime U^\prime = \hsx ^*L U^\prime$ is a cuspidal subgroup of \mG contained in \mP.  
Freeze $x \in G$ and put $g(\ell) = f_P(x \ell)$ $(\ell \in L)$:
\allowdisplaybreaks
\begin{align*}
g_{^*P}(\ell) \ 
&=\ \sum\limits_{^* u \hsx \in \hsx ^*U} \hsx g(\ell \hsy {^* u})
\\[11pt]
&=\ \sum\limits_{^* u \hsx \in \hsx ^*U} \hsx f_P (x \ell \hsy {^* u})
\\[11pt]
&=\ \sum\limits_{^* u \hsx \in \hsx ^*U} \hsx \sum\limits_{u \hsx \in U} \hsx f(x \ell \hsy {^* u} u)
\\[11pt]
&=\ \sum\limits_{u^\prime \hsx \in U^\prime} \hsx f(x \ell u^\prime)
\\[11pt]
&=\ f_{P^\prime} (x \ell).
\end{align*}
But by assumption,
\[
\sum\limits_{^* \ell \hsy \in \hsx ^*L} \hsx f_{P^\prime}(x \ell \hsy ^* \ell) \ov{\phi(^* \ell)} \hsx = \hsx 0
\]
for all $\phi \in \hsy ^0C(\hsx ^*L)$ or still, 
\[
\sum\limits_{^* \ell \hsy \in \hsx ^*L} \hsx g_{^*P} (\ell ^* \ell) \ov{\phi(^* \ell)} \hsx = \hsx 0
\]
for all $\phi \in \hsy ^0C(\hsx ^*L)$.  
The induction hypothesis then implies that $g = 0$, hence 
\[
f_P(x) \hsx = \hsx g(e) \hsx = \hsx 0.
\]
Therefore $f$ is a cusp form ($x \in G$ being arbitrary), i.e., $f \in \hsx ^0C(G)$.  
Finally, 
\[
f_G \sim 0 \implies \sum\limits_{y \in G} \hsx f(x y) \hsx \ov{\phi(y)} \hsx = \hsx 0
\]
for all $\phi \in \hsx ^0C(G)$ and all $x \in G$.  Take $x = e$ to conclude that 
\[
\sum\limits_{y \hsy \in G} \hsx f(y) \hsx \ov{\phi(y)} \hsx = \hsx 0
\]
for all $\phi \in \hsx ^0C(G)$ and then take $\phi = f$ to conclude that 
\[
\langle f, f \ranglesubG \hsx = \hsx 0 \implies f \hsx = \hsx 0.
\]
\end{x}
\vspace{0.1cm}

\begin{x}{\small\bf NOTATION} \ 
Given a cuspidal subgroup $P = L U$ of \mG, let $C(G;P)$ be the subspace of $C(G)$ consisting of those $f$ such that
\\[-.3cm]

$\qquad (i) \quad f(x u ) \hsx = \hsx f(x) \hspace{1.95cm} (x \in G, u \in U)$

\noindent and 
\\[-.3cm]

$\qquad (ii) \quad \ell \ra f(x \ell) \in \hsx ^0C(L) \qquad (x \in G, \ell \in L).$
\\[-.2cm]

[Note: \ $C(G;P)$ is stable under left translations, hence is a left ideal in $C(G)$.]
\end{x}
\vspace{0.1cm}

\begin{x}{\small\bf EXAMPLE} \ 
\[
C(G;G) \hsx = \hsx ^0C(G).
\]
\end{x}
\vspace{0.1cm}

\begin{x}{\small\bf SUBLEMMA} \ 
Fix \mP $-$then $\forall \ f \in C(G;P)$ and $\forall \ g \in C(G)$, 
\[
\langle f, g_P \ranglesubG \hsx = \hsx \abs{U} \langle f, g \ranglesubP.
\]
\vspace{0.2cm}

PROOF
\allowdisplaybreaks
\begin{align*}
\langle f, g_P \ranglesubG \ 
&=\ \frac{1}{\abs{G}} \hsx \sum\limits_{x \in G} \hsx f(x) \ov{g_P(x)}
\\[11pt]
&=\ \frac{1}{\abs{G}} \hsx \sum\limits_{x \in G} \hsx f(x) \hsx \sum\limits_{u \in U} \ov{g(x u)}
\\[11pt]
&=\ \frac{1}{\abs{G}} \hsx \sum\limits_{u \in U} \hsx \sum\limits_{x \in G} \hsx f(x) \ov{g(x u)}
\\[11pt]
&=\ \frac{1}{\abs{G}} \hsx \sum\limits_{u \in U} \hsx \sum\limits_{x \in G} \hsx f(x u) \ov{g(x u)}
\\[11pt]
&=\  \hsx \sum\limits_{u \in U} \hsx \frac{1}{\abs{G}} \hsx \sum\limits_{x \in G} \hsx f(x) \ov{g(x)}
\\[11pt]
&=\ \sum\limits_{u \in U} \hsx \langle f, g \ranglesubG
\\[11pt]
&=\ \abs{U} \langle f, g \ranglesubG.
\end{align*} 
\end{x}
\vspace{0.1cm}

\begin{x}{\small\bf RAPPEL} \ 
Let $\sH$ be a finite dimensional complex Hilbert space $-$then a subset $M \subset \sH$ is 
\un{total}
\index{total} 
if $M_{\text{$\ell$in}} = \sH$, this being the case iff $M^\perp = \{0\}$.
\vspace{0.1cm}

[Note: \ Subspaces of $\sH$ are necessarily closed \ldots \hsx .]
\\[-.25cm]
\end{x}

Put 
\[
M \hsx = \hsx \bigcup\limits_{P} C(G;P).
\]
\vspace{0.1cm}

\begin{x}{\small\bf LEMMA} \ 
$C(G)$ is spanned by the $f \in M$.
\vspace{0.2cm}

PROOF
It suffices to show that if for some $g \in C(G)$, we have 
\[
\langle f, g \ranglesubG \hsx = \hsx 0
\]
for all $f \in C(G;P)$ and for all cuspidal \mP, then $g = 0$.  
And to this end, it need only be established that $g_P \sim 0$ for all cuspidal \mP (cf. \#10).  
So fix $x \in G$ and let $\phi \in \hsx ^0C(L)$.  
Define $f \in C(G)$ as follows:
\[
\begin{cases}
\ f(y) \hsx = \hsx 0 \hspace{1.5cm} \text{if} \ y \notin x P\\
\ f(x \ell u) \hsx = \hsx \phi(\ell)  \hspace{0.5cm}  (\ell \in L, u \in U)
\end{cases}
.
\]
Then $f \in C(G;P)$, so
\[
0 \hsx = \hsx \langle f, g \ranglesubG \hsx = \hsx \frac{1}{\abs{U}} \langle f, g_P \ranglesubG
\]
$\implies$
\allowdisplaybreaks
\begin{align*}
0 \ 
&=\ \frac{1}{\abs{U}} \hsx \langle g_P, f \ranglesubG
\\[11pt]
&=\ \frac{1}{\abs{U}} \hsx \sum\limits_{y \in G} \hsx g_P(y) \ov{f(y)} 
\\[11pt]
&=\ \frac{1}{\abs{U}} \hsx \sum\limits_{y \in x P} \hsx g_P(y) \ov{f(y)} 
\\[11pt]
&=\ \frac{1}{\abs{U}} \hsx \sum\limits_{\ell, u} \hsx g_P(x \ell u) \ov{f(x \ell u)} 
\\[11pt]
&=\ \frac{1}{\abs{U}} \hsx \sum\limits_{\ell, u } \hsx g_P(x \ell) \ov{\phi(\ell)} 
\\[11pt]
&=\ \frac{1}{\abs{U}} \hsx \sum\limits_{u} \hsx \sum\limits_{\ell \in L} \hsx g_P(x \ell) \ov{\phi(\ell)} 
\\[11pt]
&=\ \sum\limits_{\ell \in L} \hsx g_P(x \ell) \ov{\phi(\ell)}.
\end{align*}
Therefore $g_P \sim 0$ .
\end{x}
\vspace{0.1cm}


\begin{x}{\small\bf CONVENTION} \ 
Cuspidal subgroups $P_1$, $P_2$ are said to be 
\un{associate}
\index{Cuspidal subgroups //associate} 
if this is the case of $\uun{P_1}$, $\uun{P_2}$.
\end{x}
\vspace{0.1cm}

\begin{x}{\small\bf LEMMA} \ 
If $P_1$, $P_2$ are associate, then 
\[
C(G; P_1) \hsx = \hsx C(G; P_2).
\]
\end{x}
\vspace{0.1cm}

\begin{x}{\small\bf LEMMA} \ 
If $P_1$, $P_2$ are not associate, then $C(G; P_1)$, $C(G; P_2)$ are orthogonal.
\end{x}
\vspace{0.1cm}

Let $P_1, \ldots, P_r$ be a set of representatives for the association classes of cuspidal subgroups of \mG.
\vspace{0.1cm}

\begin{x}{\small\bf THEOREM} \ 
There is an orthogonal decomposition
\[
C(G)  \hsx = \hsx \bigoplus\limits_{i = 1}^r \hsx C(G; P_i).
\]
\end{x}
\vspace{0.1cm}

\begin{x}{\small\bf \un{N.B.}} \ 
\#17, \#18 can be established without the use of representation theory but its introduction leads to another approach.
\end{x}
\vspace{0.1cm}

\begin{x}{\small\bf LEMMA} \ 
Let $\Pi \in \widehat{G}$ $-$then $\chi_\Pi$ is a cusp form iff $\forall$ cuspidal $P \neq G$, 
\[
m(\Pi, \Ind_{U, \theta}^G) \hsx = \hsx 0,
\]
where $\theta$ is the trivial representation of \mU on $\E = \Cx$.  
I.e.: Iff
\allowdisplaybreaks
\begin{align*}
\langle \chi_\Pi,  i_{_U \ra _G}\hsx 1_U\ranglesubG  \ 
&=\ 
\langle r_{_G \ra _U} \chi_\Pi, 1_U \rangle_U
\\[11pt]
&=\ \langle 1_U,  r_{_G \ra _U} \chi_\Pi\rangle_U 
\\[11pt]
&=\ m(\theta, \restr{\Pi}{U})
\\[11pt]
&=\ 0.
\end{align*}
\end{x}
\vspace{0.1cm}


\begin{x}{\small\bf LEMMA} \ 
Let $\Pi \in \widehat{G}$ $-$then $\chi_\Pi$ is a cusp form iff $\forall$ cuspidal $P \neq G$, 
\[
\sum\limits_{u \in U} \hsx \Pi(u) \hsx = \hsx 0.
\]
\end{x}
\vspace{0.1cm}

\begin{x}{\small\bf \un{N.B.}} \ 
Let 
\[
V(\Pi)_U \hsx = \hsx \{v \in V(\Pi): \sum\limits_{u \in U} \hsx \Pi(u) v \hsx = \hsx 0\}.
\]
Then $\chi_{_\Pi}$ is a cusp form iff $\forall$ cuspidal $P \neq G$, 
\[
V(\Pi) \hsx = \hsx V(\Pi)_U.
\]
\end{x}
\vspace{0.1cm}

\begin{x}{\small\bf DEFINITION} \ 
Let $\Pi \in \widehat{G}$ $-$then $\Pi$ is said to be in the 
\un{discrete series}
\index{discrete series} 
if its character $\chi_{_\Pi}$ is a cusp form.
\\[-.25cm]
\end{x}

\begin{x}{\small\bf NOTATION} \ 
$^0\widehat{G}$ is the subset of $\widehat{G}$ consisting of those $\Pi$ in the discrete series.
\\[-.25cm]
\end{x}

Given $P = L U$ and $\theta \in \hsx ^0\widehat{L}$, one can lift $\theta$ to \mP and form $\Ind_{P, \theta}^G$ with character
\[
i_{_P \ra _G}\chisubtheta \qquad (\text{cf. $\S9$, $\#10$).}
\]
\\[-.95cm]

\begin{x}{\small\bf THEOREM} \ 
Let $\Pi \in \widehat{G} \hsx - \hsx ^0\widehat{G}$ $-$then there exists a proper cuspidal $P = L U$ and a 
$\theta \in \hsy ^0\widehat{L}$ such that $\Pi$ occurs as a subrepresentation of $\Ind_{P, \theta}^G$:
\[
\langle \chi_{_\Pi}, \chi_{_\pi} \ranglesubG \hsx \neq \hsx 0 \qquad \bigl(\pi = \Ind_{P, \theta}^G\bigr) 
 \qquad (\text{cf. $\S5$, $\#5$).}
\]

PROOF \ 
Proceed by induction on the semisimple $k$-rank $s$ of $\uun{G}$, there being nothing to prove if $s = 0$, 
so assume that $s > 0$ $-$then there exists a proper cuspidal $P = L U$ such that $V(\Pi) \neq V(\Pi)_U$.  
Claim: \ $V(\Pi)_U$ is \mP-invariant: 
$\forall \ \ell_0 \in L$, $\forall \ u_0 \in U$, 
$\forall \ v \in V(\Pi): \sum\limits_{u \in U} \hsx \Pi(u) v = 0$, 
\allowdisplaybreaks
\begin{align*}
\sum\limits_{u \hsy \in U} \hsx \Pi(u) \Pi(\ell_0 u_0) v  \ 
&=\ \sum\limits_{u \in U} \hsx  \Pi(u \ell_0 u_0) v  
\\[11pt]
&=\ \sum\limits_{u \in U} \hsx  \Pi(\ell_0(\ell_0^{-1} u \ell_0) u_0) v  
\\[11pt]
&=\ \Pi(\ell_0) \hsx \bigl( \sum\limits_{u \in U} \hsx  \Pi(\ell_0^{-1} u \ell_0) \bigr) \Pi(u_0) v  
\\[11pt]
&=\ \Pi(\ell_0) \hsx \bigl( \sum\limits_{u \in U} \hsx \Pi(u)\bigr) \Pi(u_0) v
\\[11pt]
&=\ \Pi(\ell_0) \hsx  \sum\limits_{u \in U} \hsx \Pi(u u_0) v
\\[11pt]
&=\ \Pi(\ell_0) \hsx  \sum\limits_{u \in U} \hsx \Pi(u) v
\\[11pt]
&=\ \Pi(\ell_0) 0
\\[11pt]
&=\ 0.
\end{align*}
\end{x}
\vspace{0.1cm}

Consequently, \mP operates on the quotient $V(\Pi) / V(\Pi)_U$.  
Moreover, its restriction to \mU is trivial: $\forall \ u_0 \in U$, $\forall \ v \in V(\Pi) $, 
\allowdisplaybreaks
\begin{align*}
\sum\limits_{u \in U} \hsx \Pi(u) (\Pi(u_0) v - v) \ 
&=\ 
\sum\limits_{u \in U} \hsx \Pi(u u_0) v - \sum\limits_{u \in U} \hsx \Pi(u)  v
\\[11pt]
&=\ \sum\limits_{u \in U} \hsx \Pi(u)  v - \sum\limits_{u \in U} \hsx \Pi(u)  v
\\[11pt]
&=\ 0
\end{align*}
$\implies$
\[
\Pi(u_0)v \ \equiv \ v\hsx  \modx V(\Pi)_U.
\]
On the other hand, while its restriction to \mL need not be irreducible, there is in any event an \mL-invariant subspace $\sV$ of 
$V(\Pi)$ containing $V(\Pi)_U$ such that the quotient representation $\theta$ of \mL on 
\[
V(\Pi) /V(\Pi)_U / \sV / V(\Pi)_U \ \approx \ V(\Pi) / \sV
\]
is irreducible.  
Pass now to $\Ind_{P, \theta}^G$ and note that $\Pi$ occurs as a subrepresentation of $\Ind_{P, \theta}^G$ (see below).  
Accordingly, if $\theta \in \hsx ^0\widehat{L}$, then we are done.  
If, however, $\theta \notin \hsx ^0\widehat{L}$, then, thanks to the induction hypothesis, there exists a proper cuspidal subgroup 
$^*P = \hsx ^*L ^*U$ of \mL and a discrete series representation $^*\theta$ of $^*L$ such that $\theta$ occurs as a subrepresentation of $\Ind_{\hsy ^*P, ^*\theta}^L$.  
Form $P^\prime = L^\prime U^\prime = \hsx ^*L ^* U U$, view $^*\theta$ as a representation of $P^\prime$ trivial on 
$U^\prime = \hsx ^*U U$, and utilize the induction in stages rule (cf. \S9, \#12)
\[
\Ind_{P^\prime, ^*\theta}^G 
\ \approx \
\Ind_{P, \Ind_{\hsy ^*P, ^*\theta}^L}^G
\]
to conclude that $\Pi$, which occurs as a subreprentation of $\Ind_{P, \theta}^G$, must actually occur as a subrepresntation of 
$\Ind_{P^\prime, ^*\theta}^G$:
\allowdisplaybreaks
\begin{align*}
\theta \subset\ \Ind_{\hsy ^*P, ^*\theta}^L \implies \Pi 
&\subset \hsx \Ind_{P, \theta}^G
\\[11pt]
&\subset \hsx \Ind_{P, \Ind_{\hsy ^*P, ^*\theta}^L}^G.
\end{align*}
\vspace{0.2cm}

[Note: \ To confirm that 
\[
I_G(\Pi,\Ind_{P, \theta}^G) \hsx \neq \hsx 0,
\]
define an intertwining operator
\[
T:V(\Pi) \ra E_{P,\theta}^G
\]
by assigning to each $v \in V$ the function 
\[
f_v: G \ra V(\Pi) / \sV
\] 
given by the prescription 
\[
f_v(x) \hsx = \hsx \Pi(x^{-1}) v + \sV.]
\]
\vspace{0.2cm}

This result reduces the problem of describing the elements of $\widehat{G}$ into two parts.
\vspace{0.1cm}

\qquad \textbullet \quad  Isolate the discrete series (Deligne-Lusztig theory).
\vspace{0.1cm}

\qquad \textbullet \quad  Explicate the decomposition of $\Ind_{P, \theta}^G$ and determine its irreducibility (Howlett-Lehrer theory).
\vspace{0.2cm}

We shall pass in silence on the first of these points (for a recent survey, consult 
\url{https://arxiv.org/abs/1404.0861}
) and settle for a summary on the second (cf. \S15).
\vspace{0.1cm}

\begin{x}{\small\bf LEMMA} \ 
The canonical representation of \mG on $C(G;P)$ is equivalent to 
\[
\bigoplus\limits_\theta \hsx \Ind_{P,\theta}^G,
\]
where $\theta$ runs through the elements of $^0\widehat{L}$.
\end{x}
\vspace{0.1cm}

\begin{x}{\small\bf NOTATION} \ 
Given a parabolic $k$-subgroup $\uun{P}$ of $\uun{G}$, let $^0C(P)$ be the subspace of $C(P)$ 
consisting of those $f$ which are invariant to the right under \mU and have the property that the function on 
$P/U$ thereby defined belongs to $^0C(P/U)$.
\end{x}
\vspace{0.1cm}

\begin{x}{\small\bf LEMMA} \ 
Let $\uun{P_1}$, $\uun{P_2}$ be parabolic $k$-subgroups of $\uun{G}$ and let
\[
f_1 \in \hsx ^0C(P_1), \quad f_2 \in \hsx ^0C(P_2).
\]
Then 
\[
\langle r_{_{P_1 \ra P_1 \cap P_2}}f_1, r_{_{P_2 \ra P_2 \cap P_1}} f_2 \rangle_{P_1 \cap P_2}
\hsx = \hsx 
0
\]
unless $\uun{P_1}$ and $\uun{P_2}$ have a common Levi subgroup $\uun{L}$.
\\[-.2cm]

PROOF \ 
Ignoring constant factors (signified by $\overset{\bullet}{=}$), we have 
\begin{align*}
\langle r_{_{P_1 \ra P_1 \cap P_2}}f_1, r_{_{P_2 \ra P_2 \cap P_1}} f_2 \rangle_{P_1 \cap P_2} \hsx
&\overset{\bullet}{=} \hsx \sum_{P_1 \cap P_2} f_1(x) \hsx \ov{f_2(x)}
\\[11pt]
&\overset{\bullet}{=} \hsx \sum_{P_1 \cap P_2 / U_1 \cap U_2} f_1(x) \hsx \ov{f_2(x)}
\\[11pt]
&\overset{\bullet}{=} 
\hsx \sum_{P_1 \cap P_2 / P_1 \cap U_2} \hsx \ov{f_2(x)}
\hsx \sum_{P_1 \cap U_2 / U_1 \cap U_2} \hsx f_1(x u).
\end{align*}
Let $\pi_1: P_1 \ra P_1 / U_1 \approx L_1$ be the canonical projection $-$then 
$^*P = \pi_1(P_1 \cap P_2)$ is a cuspidal subgroup of $L_1$ with unipotent radical 
$^*U = P_1 \cap U_2 / U_1 \cap U_2$.  
Given $x \in P_1 \cap P_2$, write $x = \ell_1 \hsy u_1$ ($\ell_1 \in L_1$, $u_1 \in U_1$), thus
\allowdisplaybreaks
\begin{align*}
f_1(x u) \hsx
&= \hsx f_1(\ell_1 u_1 u) 
\\[11pt]
&= \hsx f_1(\ell_1 u_1 u u_1^{-1}),
\end{align*}
so
\[
\sum_{P_1 \cap U_2 / U_1 \cap U_2} \hsx f_1(x u) 
\hsx \overset{\bullet}{=}  \hsx 
\sum_{^*U} f_1(\ell_1 u) 
\hsx = \hsx 
0
\]
unless $^*U = \{e\}$, i.e., unless
\[
P_1 \hsx \cap \hsx U_2 
\hsx = \hsx
U_1 \hsx \cap \hsx U_2 \subset U_1.
\]
Switching roles leads to 
\[
P_2 \hsx \cap \hsx U_1 
\hsx = \hsx
U_2 \hsx \cap \hsx U_1 \subset U_2.
\]
Therefore the relevant integrals vanish unless $\uun{P_1}$ and $\uun{P_2}$ have a common Levi subgroup 
(cf. \S13, \#20).]

\end{x}
\vspace{0.1cm}

\begin{x}{\small\bf APPLICATION} \ 
Assume: \ 
$\uun{P_1}$ and $\uun{P_2}$ are not associate $-$then 
\[
\langle r_{_{P_1 \ra P_1 \cap P_2}} f_1, r_{_{P_2 \ra P_2 \cap P_1}} f_2 \rangle_{P_1 \cap P_2} \hsx = \hsx 0.
\]
\end{x}
\vspace{0.1cm}

\begin{x}{\small\bf THEOREM} \ 
Let $P_1 = L_1 U_1$, $P_2 = L_2 U_2$ be cuspidal subgroups of \mG.  
Suppose that $\uun{P_1}$ and $\uun{P_2}$ are not associate $-$then 
$\forall \ \theta_1 \in \hsx ^0\widehat{L}_1$, $\forall \ \theta_2 \in \hsx ^0\widehat{L}_2$, 
\[
\pi_1 \hsx = \hsx \Ind_{P_1,\theta_1}^G 
\qquad \text{and} \quad 
\pi_2 \hsx = \hsx \Ind_{P_2,\theta_2}^G 
\]
are disjoint: 
\[
\langle i_{_{P_1 \ra G}} \chisubthetaOne, i_{_{P_2 \ra G}}\chisubthetaTwo \ranglesubG 
\hsx = \hsx 
0 \qquad (\text{cf.} \ \S10, \ \#2).
\]
\vspace{0.2cm}

PROOF \ 
In the notation of \S8, \#4, 
\allowdisplaybreaks
\begin{align*}
\langle i_{_{P_1 \ra G}} \chisubthetaOne, i_{_{P_2 \ra G}}\chisubthetaTwo \ranglesubG \ 
&=\  
\sum\limits_{s \in S} \hsx 
\langle r_{_{P_1 \ra P_2}(s)}  \chisubthetaOne, (\chisubthetaTwo)_s \rangle_{P_2}(s)
\\[11pt]
&=\  
\sum\limits_{s \in S} \hsx 
\langle r_{_{P_1 \ra P_2}(s)} \chisubthetaOne,
r_{_{P_2^s \ra P_2}(s)} (\chisubthetaTwo)^s \rangle_{P_2}(s),
\end{align*}
where
\[
P_2(s) 
\hsx = \hsx
P_2^s \hsx \cap \hsx P_1 \qquad (= \ s P_2 s^{-1} \hsx \cap \hsx P_1).
\]
But $\uun{P_1}$ and $\uun{P_2}$ are not associate, hence $\uun{P_1}$ and $s \uun{P_2} s^{-1}$ are not associate 
(cf. \S13, \#24).  
Therefore each of the terms in the sum $\ds\sum\limits_{s \in S}$ must vanish (cf. \#30).
\\[-.75cm]
\end{x}

\begin{x}{\small\bf NOTATION} \ 
Given a parabolic $k$-subgroup $\uun{P}$ of $\uun{G}$ and a Levi subgroup $\uun{L} \subset \uun{P}$, put
\[
\uun{W}_{\uun{L}}
\  = \ 
\uun{N}_{\uun{G}} (\uun{L}) / \uun{L}.
\]
\end{x}
\vspace{0.1cm}

\begin{x}{\small\bf \un{N.B.}} \ 
If $\uun{L}^\prime$ is another Levi subgroup of $\uun{P}$, then there is a unique $u \in U$ such that 
$\uun{L}^\prime = u \uun{L} u^{-1}$, hence there is a canonical isomorphism 
\[
\uun{W}_{\uun{L}} \ra \uun{W}_{\uun{L}^\prime}.
\]

Set
\[
W_L
\hsx = \hsx
\uun{W}_{\uun{L}}(k) \qquad (= \hsx N_G (\uun{L}) / L).
\]
Then each $w \in W_L$ can be represented by an element $n_w \in N_G (\uun{L})$.
\end{x}
\vspace{0.1cm}

\begin{x}{\small\bf LEMMA} \ 
The arrow
\[
W_L \ra P \backslash G /P
\]
given by 
\[
w \ra P n_w P
\]
is injective.
\end{x}
\vspace{0.1cm}

\begin{x}{\small\bf LEMMA} \ 
$W_L$ operates on $^0C(P)$.
\end{x}
\vspace{0.1cm}

\begin{x}{\small\bf REDUCTION PRINCIPLE} \ 
Let $\uun{P_1}$, $\uun{P_2}$ be parabolic $k$-subgroups of $\uun{G}$ and let 
\[
f_1 \in \hsx ^0CL(P_1), \ f_2 \in \hsx ^0CL(P_2).
\]
Assume: $\uun{P_1}$ and $\uun{P_2}$ have a common Levi subgroup $\uun{L}$ $-$then 
\[
\langle i_{_{P_1 \ra G}}f_1, i_{_{P_2 \ra G}}f_2 \ranglesubG 
\hsx = \hsx 
\sum\limits_{w \in W_L}
\langle r_{_{P_1 \ra L}}f_1, r_{_{P_2 \ra L}}(w \cdot f_2) \ranglesubL. 
\]
\vspace{0.2cm}

PROOF \ 
In the notation of \S8, \#4,
\allowdisplaybreaks 
\begin{align*}
\langle i_{_{P_1 \ra G}}f_1, i_{_{P_2 \ra G}}f_2 \ranglesubG \hsx
&= \hsx \sum\limits_{s \in S} \hsx
\langle r_{_{P_1 \ra P_2}(s)}f_1, (f_2)_s \rangle_{P_2} (s)
\\[11pt]
&= \hsx \sum\limits_{s \in S} \hsx
\langle r_{_{P_1 \ra P_2}(s)}f_1, f_{P_2^s \ra P_2(s)} (f_2)^s \rangle_{P_2} (s), 
\end{align*}
where
\[
P_2(s) 
\hsx = \hsx
P_2^s \hsx \cap \hsx  P_1 
\qquad (= s \hsx P_2 \hsx s^{-1} \hsx \cap \hsx P_1).
\]
The only nonzero terms in the sum are those for which $\uun{P_1}$ and $s \hsx \uun{P_2} \hsx s^{-1} $ have a common Levi subgroup 
$\uun{L}^\prime$ (cf. \#31).  
Choose $u_1 \in U_1$ such that $u_1 \hsx \uun{L}^\prime \hsx u_1^{-1} = \uun{L}$.  
Next
\[
\uun{L}^\prime 
\subset 
s \hsx \uun{P_2} \hsx s^{-1} 
\implies 
s^{-1} \uun{L}^\prime s 
\subset 
\uun{P_2}.
\]
Choose $u_2 \in U_2$ such that $u_2 \hsx \uun{L} \hsx u_2^{-1} = s^{-1} \hsx \uun{L}^\prime \hsx s$, thus
\[
\uun{L}^\prime 
\hsx = \hsx
s \hsy u_2 \hsy \uun{L} \hsy u_2^{-1} \hsy s^{-1}
\]

\qquad\qquad $\implies$
\[
\uun{L} 
\hsx = \hsx
u_1\hsy  \uun{L} ^\prime \hsy  u_1^{-1}
\hsx = \hsx
u_1 \hsy s \hsy u_2 \hsy \uun{L}  \hsy u_2^{-1} \hsy s^{-1} \hsy u_1^{-1}
\]
\qquad\qquad $\implies$
\[
u_1 \hsy s \hsy u_2 \in N_G(\uun{L} ).
\]
On the other hand, 
\[
u_1 \hsy s \hsy u_2 \in P_1 \backslash G / P_2.
\]
Therefore the double cosets $P_1\backslash G /P_2$ that intervene are those containing an element $N_G(\uun{L})$, so
\[
\langle i_{_{P_1 \ra G}} f_1, i_{_{P_2 \ra G}} f_2 \ranglesubG
\hsx = \hsx
\sum\limits_{w \in W_L} \hsx 
\langle r_{_{P_1 \ra P_2}(w)} f_1, r_{_{P_2^w \ra P_2}(w)} (w \cdot f_2) \rangle_{P_2(w)}.
\]
Noting that $\uun{L} = w \uun{L} w^{-1} \subset w \uun{P_2} w^{-1}$ is a Levi subgroup of $w \uun{P_2} w^{-1}$, write
\[
P_2(w)
\hsx = \hsx
P_2^w \hsx \cap \hsx P_1
\hsx = \hsx
P_1 \hsx \cap \hsx w P_2 w^{-1} 
\hsx = \hsx
L \cdot (L \hsx \cap \hsx w U_2 w^{-1}) \cdot (U_1 \hsx \cap \hsx L) \cdot (U_1 \hsx \cap \hsx w U_2 w^{-1})
\]
with uniqueness of expression $-$then 
\[
L \hsx \cap \hsx w U_2 w^{-1} \hsx = \hsx \{e\}, \quad U_1 \hsx \cap \hsx L \hsx = \hsx \{e\}
\]
and 
\[
\langle r_{_{P_1 \ra P_2}(w)} f_1, r_{_{P_2^w \ra P_2}(w)} (w \cdot f_2) \rangle_{P_2(w)} 
\
\hsx = \hsx
\frac{1}{\abs{P_2(w)}} \hsx 
\sum\limits_{x,u} \hsx f_1(x u) \ov{(w \cdot f_2) (x u)},
\]
where the sum runs over all $x \in L$ and all $u \in U_1 \hsx \cap \hsx w U_2 w^{-1}$.  
Since $f_1$ and $w \cdot f_2$ are invariant to the right under $U_1 \hsx \cap \hsx w U_2 w^{-1}$, the above expression equals 
\[
\frac{\abs{U_1 \hsx \cap \hsx w U_2 w^{-1}}}{\abs{P_2(w)}} 
\sum\limits_x \hsx f_1(x) \hsx \ov{(w \cdot f_2) (x)}
\hsx = \hsx
\frac{\abs{U_1 \hsx \cap \hsx w U_2 w^{-1}}}{\abs{P_2(w)}}  
\hsx 
\abs{L} \langle r_{_{P_1 \ra L}} f_1, r_{_{P_2 \ra L}}(w \cdot f_2)\ranglesubL.
\]
And
\[
\frac{\abs{U_1 \hsx \cap \hsx w U_2 w^{-1}}}{\abs{P_2(w)}} \hsx \abs{L}
\hsx = \hsx
\frac{\abs{U_1 \hsx \cap \hsx w U_2 w^{-1}} \cdot \abs{L}}{\abs{L} \cdot \abs{U_1 \hsx \cap \hsx w U_2 w^{-1}}}  
\hsx = \hsx
1.
\]
\end{x}
\vspace{0.1cm}

\begin{x}{\small\bf SUBLEMMA} \ 
Let $\sH$ be a Hilbert space and let $x$, $y \in \sH$.  
Assume: 
\[
\langle x,x \rangle
\hsx = \hsx
\langle x,y \rangle
\hsx = \hsx
\langle y,y \rangle.
\]
Then $x = y$.
\vspace{0.2cm}

PROOF \ 
In fact, 
\allowdisplaybreaks
\begin{align*}
\langle x - y,x - y \rangle \hsx
&= \hsx \langle x,x \rangle + \langle y,y \rangle - \langle x,y \rangle - \langle y,x \rangle
\\[11pt]
&= \hsx \langle x,y \rangle + \langle x,y \rangle - \langle x,y \rangle - \ov{\langle x,y \rangle}
\\[11pt]
&= \hsx \langle x,y \rangle - \ov{\langle x,y \rangle}
\\[11pt]
&= \hsx \langle x,y \rangle - \langle x,y \rangle
\\[11pt]
&= \hsx 0.
\end{align*}
\end{x}
\vspace{0.1cm}

\begin{x}{\small\bf APPLICATION} \ 
If 
\[
r_{_{P_1 \ra L}} f_1 
\hsx = \hsx
r_{_{P_2 \ra L}} f_2,
\]
then
\[
i_{_{P_1 \ra G}}f_1 
\hsx = \hsx
i_{_{P_2 \ra G}} f_2.
\]
\vspace{0.2cm}

[It follows from \#36 that 
\allowdisplaybreaks
\begin{align*}
\langle i_{_{P_1 \ra G}} f_1, i_{_{P_1 \ra G}}f_1 \ranglesubG \hsx
&= \  
\langle i_{_{P_1 \ra G}} f_1, i_{_{P_2 \ra G}} f_2\ranglesubG 
\\[11pt]
&= \  
\langle i_{_{P_2 \ra G}} f_2, i_{_{P_2 \ra G}}f_2 \ranglesubG.]
\end{align*}
\end{x}
\vspace{0.1cm}

\begin{x}{\small\bf NOTATION} \ 
Given a cuspidal subgroup $P = L U$ of \mG and a $\theta \in \hsx ^0\widehat{L}$, let
\[
W_L(\theta) \hsx = \hsx \{w \in W_L: w \cdot \chisubtheta = \chisubtheta\}.
\]
\end{x}
\vspace{0.1cm}

\begin{x}{\small\bf THEOREM} \ 
\[
\langle i_{_{P \ra G}} \chisubtheta, i_{_{P \ra G}} \chisubtheta \ranglesubG \hsx = \hsx \abs{W_L(\theta)}.
\]
\vspace{0.2cm}

[In \#36, take $\uun{P_1} = \uun{P_2} = \uun{P}$ and note that 
\[
\langle r_{_{P \ra L}}\chisubtheta, r_{_{P \ra L}}(w \cdot \chisubtheta) \ranglesubL
\]
equals 1 if $w \cdot \chisubtheta = \chisubtheta$ and equals 0 if $w \cdot \chisubtheta \neq \chisubtheta$.]
\end{x}
\vspace{0.1cm}

Let $\sP$ be the set of parabolic $k$-subgroups of $\uun{G}$.  
Decompose $\sP$ into association classes: $\sP = \coprod \sC$.  
Given $\sC$, take a $\uun{P} \in \sC$ and denote by $\widehat{G}(\sC)$ 
the subset of $\widehat{G}$ comprised of those $\Pi$ which occur as a subrepresentation of
\[
\Ind_{P,\theta}^G
\]
for some $\theta \in \hsx ^0\widehat{L}$.

\begin{x}{\small\bf LEMMA} \ 
$\widehat{G}(\sC)$ is independent of the choice of $\uun{P} \in \sC$.
\vspace{0.2cm}

PROOF \ The theory does not change if $\uun{P}$ is replaced by $x \uun{P} x^{-1}$ $(x \in G)$, so if 
$\uun{P_1}$, $\uun{P_2}$ are associate, then there is no loss of generality in assuming that 
$\uun{P_1}$ and $\uun{P_2}$ have a common Levi subgroup $\uun{L}$, thus 
$
\begin{cases}
\ L \subset P_1 \\
\ L \subset P_2
\end{cases}
.  
$
Given $\theta \in \hsx ^0\widehat{L}$, lift
\[
\begin{cases}
\theta \ \text{to} \ P_1, \  \text{call it} \ \theta_1\\
\theta \ \text{to} \ P_2, \  \text{call it} \ \theta_2\\
\end{cases}
.
\]
Then
\[
\begin{cases}
\ i_{_{P_1 \ra G}} \chisubthetaOne \quad \text{is the character of} \ \Ind_{P_1,\theta_1}^G\\[3pt]
\ i_{_{P_2 \ra G}} \chisubthetaTwo \quad \text{is the character of} \ \Ind_{P_2,\theta_2}^G
\end{cases}
.
\]
But
\[
\begin{cases}
\ r_{_{P_1 \ra L}}\chisubthetaOne \hsx = \hsx \chisubtheta\\[3pt]
\ r_{_{P_2 \ra L}}\chisubthetaTwo \hsx = \hsx \chisubtheta
\end{cases}
\]
\qquad\qquad $\implies$
\[
i_{_{P_1 \ra G}} \chisubthetaOne \hsx = \hsx i_{_{P_2 \ra G}} \chisubthetaTwo \qquad (\text{cf.} \ \#38). 
\]
Therefore
\[
\Ind_{P_1,\theta_1}^G \hsx \approx \hsx \Ind_{P_2,\theta_2}^G.
\]
\end{x}
\vspace{0.1cm}

\begin{x}{\small\bf LEMMA} \ 
If $\sC_1 \neq \sC_2$, then 
\[
\widehat{G}(\sC_1) \hsx \cap \hsx \widehat{G}(\sC_2) \neq \emptyset \qquad (\text{cf.} \ \#30).
\]
\end{x}
\vspace{0.1cm}
Accordingly: 
\vspace{0.1cm}

\begin{x}{\small\bf THEOREM} \ 
There is a disjoint decomposition
\[
\widehat{G} \hsx = \hsx \coprod\limits_\sC \ \widehat{G}(\sC).
\]
\end{x}
\vspace{0.1cm}

\begin{x}{\small\bf NOTATION} \ 
Given $\uun{P} \in \sP$, let $[\uun{P}]$ be the association class to which $\uun{P}$ belongs.
\end{x}
\vspace{0.1cm}

\begin{x}{\small\bf EXAMPLE} \ 
Take $\uun{P} = \uun{G}$ $-$then the elements of $\widehat{G}([G])$ comprise the 
\un{discrete} \un{series}
\index{discrete series} 
for \mG, i.e., $\widehat{G}([G]) = \hsx ^0\widehat{G}$.
\end{x}
\vspace{0.1cm}

\begin{x}{\small\bf EXAMPLE} \ 
Take $\uun{P} = \uun{B}$ $-$then the elements of $\widehat{G}([B])$ comprise the 
\un{principal series}
\index{principal series} 
for \mG.
\end{x}
\vspace{0.1cm}

\begin{x}{\small\bf REMARK} \ 
$W_L$ operates on $^0\widehat{L}$, hence $^0\widehat{L}$ breaks up into $W_L$-orbits.
Let $\theta_1$, $\theta_2 \in \hsx ^0\widehat{L}$ $-$then there are two possibilities.
\vspace{0.2cm}

\qquad \textbullet \quad If $\theta_1$, $\theta_2$ are on the same $W_L$-orbit, then 
\vspace{0.2cm}
\[
\Ind_{P,\theta_1}^G \hsx \approx \hsx \Ind_{P,\theta_2}^G.
\]

\qquad \textbullet \quad If $\theta_1$, $\theta_2$ are not on the same $W_L$-orbit, then 
\[
\Ind_{P,\theta_1}^G \qquad \text{and} \qquad \Ind_{P,\theta_2}^G
\]
are disjoint.

\end{x}
\vspace{0.1cm}


\chapter{
$\boldsymbol{\S}$\textbf{15}.\quad  HOWLETT-LEHRER THEORY}
\setlength\parindent{2em}
\setcounter{theoremn}{0}
\renewcommand{\thepage}{A II \S15-\arabic{page}}

\qquad In view of \S14, \#40, 
\[
\langle i_{_{P \ra G}} \chisubtheta, i_{_{P \ra G}} \chisubtheta \ranglesubG \ = \ \abs{W_L(\theta)}.
\]
And on general grounds 
(cf. \S5, \#11), $\Ind_{P,\theta}^G$ is irreducible iff 
\[
\langle i_{_{P \ra G}} \chisubtheta, i_{_{P \ra G}} \chisubtheta \ranglesubG \ = \ 1.
\]

\begin{x}{\small\bf DEFINITION} \ 
$\theta$ is 
\un{unramified}
\index{unramified} 
if $\abs{W_L(\theta)} = 1$.
\end{x}
\vspace{0.3cm}

\begin{x}{\small\bf THEOREM} \ 
$\Ind_{P,\theta}^G$  is irreducible iff $\theta$ is unramified.
\end{x}
\vspace{0.3cm}

To discuss the decomposability of $\Ind_{P,\theta}^G$ , note that $\Pi \in \widehat{G}$ occurs as a subrepresentation of 
$\Ind_{P,\theta}^G$  iff 
\[
\langle \chisubPi, i_{_{P\ra G}} \chisubtheta \ranglesubG \hsx \neq \hsx 0.
\]
\vspace{0.3cm}

\begin{x}{\small\bf LEMMA} \ 
There is a one-to-one correspondence between the $\Pi \in \widehat{G}$ such that
\[
\langle \chisubPi, i_{_{P \ra G}} \chisubtheta \ranglesubG \neq 0
\] 
and the irreducible representations $\rho$ of
\[
I_G(\Ind_{P,\theta}^G, \Ind_{P,\theta}^G).
\]
And if $\Pi \longleftrightarrow \rho$, then 
\[
\chi_{_\rho}(1) \hsx = \hsx \langle \chisubPi, i_{_{P_\ra G}} \chisubtheta \ranglesubG,
\]
the positive integer on the right being the multiplicity
\[
m (\Pi,\Ind_{P,\theta}^G)
\]
of $\Pi$ in $\Ind_{P,\theta}^G$.
\end{x}
\vspace{0.3cm}


\begin{x}{\small\bf THEOREM} \ 
The semisimple algebra
\[
I_G(\Ind_{P,\theta}^G, \Ind_{P,\theta}^G)
\]
is isomorphic to the  semisimple algebra
\[
C(W_L(\theta)).
\]
\end{x}
\vspace{0.3cm}

The irreducible components of $\Ind_{P,\theta}^G$ are therefore parameterized by the elements of $W_L(\theta)$:  
If $\omega \in W_L(\theta)$ and if $\Pi(\omega) \in \widehat{G}$ is the irreducible component of 
$\Ind_{P,\theta}^G$ corresponding to $\omega$, then 
\[
\langle \chi_{_{\Pi(\omega)}}, i_{_{P\ra G}} \chisubtheta \ranglesubG 
\ = \ 
\chi_{_\omega}(1),
\]
the dimension of the representation space of $\omega$.


\chapter{
$\boldsymbol{\S}$\textbf{16}.\quad  MODULE LANGUAGE}
\setlength\parindent{2em}
\setcounter{theoremn}{0}
\renewcommand{\thepage}{A II \S16-\arabic{page}}

\qquad Let \mG be a finite group, $\Gamma \subset G$ a subgroup.  
View $C(G)$ as a left $C(G)$-module and as a right $C(\Gamma)$-module.
\\[-.2cm]

\begin{x}{\small\bf CONSTRUCTION} \ 
Let $\theta : \Gamma \ra \GL(\E)$ be a representation of $\Gamma$ $-$then the tensor product
\[
C(G) \otimes_{C(\Gamma)} \hsx \E
\]
is a left $C(G)$-module or, equivalently, a representation, the 
\un{representation $\Ind_{\Gamma, \theta}^G$ of \mG} \un{induced by $\theta$}.
\index{representation $\Ind_{\Gamma, \theta}^G$ of \mG induced by $\theta$}
\end{x}
\vspace{0.3cm}

\begin{x}{\small\bf \un{N.B.}} \ 
The left action is given by
\[
\bigg( \hsx \sum\limits_{x \in G} \hsx f(x) \delta_x \bigg) \bigl(\delta_y \hsx \otimes \hsx X \bigr) 
\hsx = \hsx
\sum\limits_{x \in G} \delta_{x y} \hsx \otimes \hsx f(x) X \qquad (X \in \E)
\]
and from the definitions, $\forall \ \gamma \in \Gamma$, 
\[
\delta_x \hsx \delta_\gamma \hsx \otimes \hsx X
\hsx = \hsx
\delta_x \hsx \otimes \hsx \theta(\gamma) X \qquad (X \in \E).
\]
\end{x}
\vspace{0.3cm}

\begin{x}{\small\bf LEMMA} \ 
Write
\[
G 
\hsx = \hsx
\coprod\limits_{k = 1}^n \hsx x_k \Gamma.
\]
Then as a vector space
\[
\Ind_{\Gamma, \theta}^G 
\hsx = \hsx
\bigoplus\limits_{k = 1}^n \ \bigl(\delta_{x_k} \hsx \otimes \hsx \E \bigr).
\]
\vspace{0.3cm}

[Note: \ The summand
\allowdisplaybreaks
\begin{align*}
\delta_{x_k} \hsx \otimes \hsx \E \ 
&=\  \{\delta_{x_k} \hsx \otimes \hsx X : X \in \E\}\\
&\approx\ \E\bigl(\delta_{x_k} \hsx \otimes \hsx X \longleftrightarrow X \bigr)
\end{align*}
is the transform of $\delta_e \hsx \otimes \hsx \E  \approx \E$ under the action of $\delta_{x_k}$: 
\[
\delta_{x_k} (\delta_e \otimes \hsx \ X)  \hsx 
\hsx = \hsx
\delta_{x_k} \hsx \otimes \hsx X.]
\]
\end{x}
\vspace{0.3cm}

The following result justifies the notation and the terminology.
\vspace{0.3cm}

\begin{x}{\small\bf THEOREM} \ 
Set $\pi = \Ind_{\Gamma, \theta}^G$ $-$then 
\[
i_{{_\Gamma \ra G}} \chi_{_\theta} \hsx = \hsx \chi_\pi \qquad \text{(cf.} \ \S9, \ \#10).
\]

PROOF \ 
Let $X_1, \ldots, X_d$ be a basis for $\E$ and define $\theta_{i j}(\gamma)$ by 
\[
\theta(\gamma) X_j
\hsx = \hsx
\sum\limits_i \hsx \theta_{i j}(\gamma) X_i.
\]
Equip $C(G) \hsx \otimes_{C(\Gamma)} \hsx \E$ with the basis 
\[
\{\delta_{x_1} \otimes X_1, \ldots ,\delta_{x_1} \otimes X_d, 
\delta_{x_2} \otimes X_1, \ldots ,\delta_{x_2} \otimes X_d, 
\ldots, 
\delta_{x_n} \otimes X_1, \ldots ,\delta_{x_n} \otimes X_d\}
\]
and write $xx_k = x_\ell \gamma$ $-$then 
\allowdisplaybreaks
\begin{align*}
\pi(x) (\delta_{x_k} \hsx \otimes \hsx X_j) \
&=\ \delta_{x x_k} \hsx \otimes \hsx X_j
\\[11pt]
&=\ \delta_{x_\ell \gamma} \hsx \otimes \hsx X_j
\\[11pt]
&=\ \delta_{x_\ell} \delta_\gamma \hsx \otimes \hsx X_j
\\[11pt]
&=\ \delta_{x_\ell} \hsx \otimes \hsx \theta (\gamma) X_j
\\[11pt]
&=\ \delta_{x_\ell} \hsx \otimes \hsx \sum\limits_i \theta_{i j}(\gamma) X_i
\\[11pt]
&=\ \sum\limits_i \theta_{i j} (\gamma) \delta_{x_\ell} \hsx \otimes \hsx X_i
\\[11pt]
&=\ \sum\limits_i \theta_{i j} (x_\ell^{-1} x x_k) \delta_{x_\ell} \hsx \otimes \hsx X_i.
\end{align*}
Define $\mathring{\theta}$ on \mG by 
$\mathring{\theta}(\gamma) = [\theta_{i j}(\gamma)]$ \ $(\gamma \in \Gamma)$ 
and 
$\mathring{\theta}(x) =0_d$ if $x \notin \Gamma$ $(0_d$ 
the zero 
$d$-by-$d$ matrix), thus the block matrix representing $\pi(x)$ is 
\[
\begin{pmatrix}
\mathring{\theta}(x_1^{-1} x x_1) & \mathring{\theta}(x_1^{-1} x x_2)  &\ldots & \mathring{\theta}(x_1^{-1} x x_n) 
\\[8pt]
\mathring{\theta}(x_2^{-1} x x_1) & \mathring{\theta}(x_2^{-1} x x_2)  &\ldots & \mathring{\theta}(x_2^{-1} x x_n) 
\\[8pt]
\vdots &\vdots &\vdots &\vdots\\
\mathring{\theta}(x_n^{-1} x x_1) & \mathring{\theta}(x_n^{-1} x x_2)  &\ldots & \mathring{\theta}(x_n^{-1} x x_n) 
\\[8pt]
\end{pmatrix}
\]
Taking the trace
\allowdisplaybreaks
\begin{align*}
\chi_\pi(x) \ 
&=\ 
\tr(\pi(x))
\\[11pt]
&=\ \sum\limits_{k = 1}^n \hsx \tr (\mathring{\theta}(x_k^{-1} x x_k))
\\[11pt]
&=\ \sum\limits_{k = 1}^n \hsx \mathring{\chi_\theta} (x_k^{-1} x x_k)
\\[11pt]
&=\ (i_{_{\Gamma \ra G}} \chi_{_\theta})(x)
\end{align*}
finishes the proof.
\end{x}
\vspace{0.3cm}

\begin{x}{\small\bf NOTATION} \ 
$\MOD(\Gamma)$ is the category of left $C(\Gamma)$-modules and $\MOD(G)$ is the category of left $C(G)$-modules. 
\vspace{0.2cm}

[Note: \ All data is over $\Cx$ and finite dimensional.]
\end{x}
\vspace{0.3cm}

\begin{x}{\small\bf \un{N.B.}} \ 
Morphisms are intertwining operators.
\end{x}
\vspace{0.3cm}

\begin{x}{\small\bf SCHOLIUM} \ 
The assignment
\[
(\theta,\E) \ra \Ind_{\Gamma, \theta}^G
\]
defines a functor 
\[
\MOD(\Gamma) \ra \MOD(G).
\]
\end{x}
\vspace{0.3cm}

\begin{x}{\small\bf NOTATION} \ 
Given a representation $(\pi,V)$ of \mG, denote its restriction to $\Gamma$ by $\Res_{\Gamma,\pi}^G$.
\end{x}
\vspace{0.3cm}

\begin{x}{\small\bf SCHOLIUM} \ 
The assignment
\[
(\pi,V) \ra \Res_{\Gamma,\pi}^G
\]
defines a functor 
\[
\MOD(G) \ra \MOD(\Gamma).
\]
\end{x}
\vspace{0.3cm}

Here now are the fundamental formalities.
\vspace{0.3cm}

\begin{x}{\small\bf LEMMA} \ 
\[
I_G \bigl(\Ind_{\Gamma, \theta}^G , (\pi,V)\bigr) 
\hsx \approx \hsx 
I_\Gamma \bigl((\theta,E) , \Res_{\Gamma,\theta}^G\bigr). 
\]
\\[-1.25cm]
\end{x}

\begin{x}{\small\bf SLOGAN} \ 
The restriction functor is a right adjoint for the induction functor.
\\[-.5cm]
\end{x}

\begin{x}{\small\bf LEMMA} \ 
\[
I_G\bigl((\pi,V), \Ind_{\Gamma, \theta}^G \bigr) 
\hsx \approx \hsx 
I_\Gamma\bigl(\Res_{\Gamma,\pi}^G, (\theta,\E) \bigr). 
\]
\\[-1.25cm]
\end{x}


\begin{x}{\small\bf SLOGAN} \ 
The restriction functor is a left adjoint for the induction functor.
\end{x}
\vspace{0.3cm}

Moving on:
\vspace{0.3cm}

\begin{x}{\small\bf DEFINITION} \ 
Let $\theta:\Gamma \ra \GL(\E)$ be a representation of $\Gamma$ $-$then 
\[
\Inv_\Gamma(\E) 
\hsx = \hsx 
\{X \in \E : \theta(\gamma) X = X \ \forall \ \gamma \in \Gamma\}
\]
is the set of 
\un{$\Gamma$-invariants}
\index{$\Gamma$-invariants} 
per $\E$.
\end{x}
\vspace{0.3cm}

\begin{x}{\small\bf DEFINITION} \ 
Let $\theta:\Gamma \ra \GL(\E)$ be a representation of $\Gamma$ $-$then 
\[
\CoInv_\Gamma(\E) 
\hsx = \hsx 
\E / I_\Gamma \E
\]
is the set of 
\un{$\Gamma$-coinvariants}
\index{$\Gamma$-coinvariants} 
per $\E$.
\vspace{0.2cm}

[Note: \ $I_\Gamma \subset C(\Gamma)$ is the augmentation ideal, thus $I_\Gamma \E$ stands for the set of all finite sums 
$\sum\limits_i \theta(\gamma_i) X_i$ $(\delta_{\gamma_i} \in I_\Gamma, X_i \in \E)$.]
\end{x}
\vspace{0.3cm}

Specialize and assume that \mG is a group of Lie type (cf. \S13, \#3).
\vspace{0.3cm}

\begin{x}{\small\bf NOTATION} \ 
Given a cuspidal subgroup $P = L U$ of \mG, 
\[
\Inf_{L,P} : \MOD(L) \ra \MOD(P)
\]
is the inflation functor.
\vspace{0.2cm}

[In other words, given a representation $(\theta,\E)$ of \mL, $\Inf_{L,P} (\theta)$ is the lift of $\theta$ to \mP, 
i.e., $\E$ viewed as a left $C(P)$-module with trivial \mU-action.]
\end{x}
\vspace{0.3cm}

\begin{x}{\small\bf DEFINITION} \ 
The composite
\[
\Ind_{P,-}^G \hsx \circ \hsx \Inf_{L,P}
\]
defines a functor 
\[
R_{L,P}^G :  \MOD(L) \ra \MOD(G)
\]
termed 
\un{Harish-Chandra induction}.
\index{Harish-Chandra induction}
\end{x}
\vspace{0.3cm}

\begin{x}{\small\bf THEOREM} \ 
If $P_1 = L U_1$, $P_2 = L U_2$ are cuspidal subgroups of \mG, then the functors 
\[
\begin{cases}
\ R_{L,P_1}^G\\
\ R_{L,P_2}^G
\end{cases}
\]
are naturally isomorphic.
\vspace{0.2cm}

[Note: \ 
Accordingly, the left $C(G)$-module isomorphism class of $R_{L,P}^G (\theta,\E)$ depends only on $\theta$ 
(it being independent of the particular cuspidal subgroup $P = L U$).]
\end{x}
\vspace{0.3cm}

\begin{x}{\small\bf LEMMA} \ 
\allowdisplaybreaks
\begin{align*}
I_G(R_{L,P}^G (\theta,\E), (\pi,V)) \ 
&\approx\ 
I_P(\Inf_{L,P} \theta, \Res_{P,\pi}^G)
\\[11pt]
&\approx\ 
I_L((\theta,\E), \Inv_U (\Res_{P,\pi}^G)).
\end{align*}

[Note: \ 
For any left $C(P)$-module \mM, the set $\Inv_U(M)$ is canonically a left $C(L)$-module.]
\end{x}
\vspace{0.1cm}

\begin{x}{\small\bf SLOGAN} \ 
The composite of restriction followed by the taking of invariants is a right adjoint for Harish-Chandra induction.
\end{x}
\vspace{0.3cm}

\begin{x}{\small\bf LEMMA} \ 
\allowdisplaybreaks
\begin{align*}
I_G((\pi,V),R_{L,P}^G (\theta,\E)) \ 
&\approx\ 
I_P(\Res_{P,\pi}^G,\Inf_{L,P} (\theta))
\\[11pt]
&\approx\ 
I_L( \CoInv_U (\Res_{P,\pi}^G), (\theta,\E)).
\end{align*}
\vspace{0.1cm}

[Note: \ 
For any left $C(P)$-module \mM, the set $\CoInv_U(M)$ is canonically a left $C(L)$-module.]
\end{x}
\vspace{0.1cm}

\begin{x}{\small\bf SLOGAN} \ 
The composite of restriction followed by the taking of coinvariants is a left adjoint for Harish-Chandra induction.
\end{x}
\vspace{0.1cm}

\begin{x}{\small\bf SUBLEMMA} \ 
For any left $C(P)$-module \mM,
\[
\Inv_U(M) 
\hsx \approx \hsx
\CoInv_U(M).
\]
\end{x}
\vspace{0.1cm}

\begin{x}{\small\bf SCHOLIUM} \ 
The left and right adjoint of Harish-Chandra induction are naturally isomorphic.
\end{x}
\vspace{0.1cm}

\begin{x}{\small\bf DEFINITION} \ 
Harish-Chandra restriction $^*R_{L,P}^G$ is the left and right adjoint of Harish-Chandra induction.
\end{x}
\vspace{0.1cm}

\begin{x}{\small\bf LEMMA} \ 
\[
^*R_{L,P}^G ((\pi,V))
\hsx = \hsx 
e_U V,
\]
where
\[
e_U 
\hsx = \hsx 
\frac{1}{\abs{U}} \hsx \sum\limits_{u \in U} \hsx \pi(u).
\]
\end{x}
\vspace{0.1cm}

\begin{x}{\small\bf THEOREM} \ 
If $P_1 = L U_1$, $P_2 = L U_2$ are cuspidal subgroups of \mG, then the functors
\[
\begin{cases}
\ ^*RG_{L,P_1}\\
\ ^*RG_{L,P_2}
\end{cases}
\]
are naturally isomorphic.
\vspace{0.2cm}

[Note: \ 
Accordingly, the left $C(L)$-module class of $^*R^G_{L,P} (\pi,V)$ depends only on $\pi$ 
(it being independent of the particular cuspidal parabolic subgroup $P = L U$).]

\end{x}
\vspace{0.5cm}

\[
\text{APPENDIX}
\]
\\[-1.25cm]

Let $P = L U$ be a cuspidal subgroup of \mG.
\\[-.2cm]

\qquad {\small\bf DEFINITION} \ 
Given $\phi \in CL(L)$, define $\widetilde{\phi} \in C(P)$ by the rule
\[
\widetilde{\phi}(\ell u) \hsx = \hsx \phi(\ell).
\]
\\[-.2cm]

\qquad {\small\bf LEMMA} \ 
$\widetilde{\phi}$ is a class function, i.e., 
\[
\widetilde{\phi} \in CL(P).
\]

PROOF \ 
The claim is that $\forall \ p \in P$, $\forall \ p_1 \in P$, 
\[
\widetilde{\phi} (p p_1 p^{-1}) \hsx = \hsx \widetilde{\phi} (p_1).
\]
Write $p = \ell u$, $p_1 = \ell_1 u_1$ $-$then
\allowdisplaybreaks
\begin{align*}
\widetilde{\phi} (p p_1 p^{-1}) \ 
&=\ \widetilde{\phi} (\ell u \ell_1 u_1 u^{-1} \ell^{-1})
\\[11pt]
&=\ \widetilde{\phi} (\ell u \ell_1 u_1 \ell^{-1} \ell u^{-1} \ell^{-1})
\\[11pt]
&=\ \widetilde{\phi} (\ell u \ell_1 u_1 \ell^{-1} v) \qquad\qquad (v = \ell u^{-1} \ell^{-1} \in U)
\\[11pt]
&=\ \widetilde{\phi} (\ell u \ell_1 \ell^{-1} \ell u_1 \ell^{-1} v)
\\[11pt]
&=\ \widetilde{\phi} (\ell u \ell_1 \ell^{-1} v_1 v) \qquad\qquad (v_1 = \ell u_1 \ell^{-1} \in U)
\\[11pt]
&=\ \widetilde{\phi} (\ell (\ell_1 \ell^{-1}) (\ell_1 \ell^{-1})^{-1} u(\ell_1 \ell^{-1}) v_1 v)
\\[11pt]
&=\ \widetilde{\phi} (\ell \ell_1 \ell^{-1} v_2 v_1 v) \qquad\qquad (v_2 = (\ell_1 \ell^{-1})^{-1} u (\ell_1 \ell^{-1}) \in U)
\\[11pt]
&=\ {\phi} (\ell \ell_1 \ell^{-1})
\\[11pt]
&=\ {\phi} (\ell_1)
\\[11pt]
&=\ \widetilde{\phi} (p_1).  
\end{align*}

Thus there is an arrow 
\[
CL(L) \ra CL(P) \ra CL(G),
\]
namely
\[
\phi \ra \widetilde{\phi} \ra i_{_{P \ra G}} \widetilde{\phi}.
\]
On the other hand, there is an arrow
\[
CL(G) \ra CL(L),
\]
namely
\[
\psi \ra \restr{\psi_P}{L} \hsx \equiv \hsx r_{_{G \ra L}} \psi_P.
\]
\vspace{0.2cm}

[Note: \ 
$\forall \ \ell \in L$, $\forall \ \ell_1 \in L$, 
\allowdisplaybreaks
\begin{align*}
\sum\limits_{u \in U} \hsx \psi(\ell \ell_1 \ell^{-1} u) \ 
&=\ \sum\limits_{u \in U} \hsx \psi(\ell \ell_1 \ell^{-1} u \ell \ell^{-1})
\\[11pt]
&=\ \sum\limits_{u \in U} \hsx \psi(\ell_1 \ell^{-1} u \ell)
\\[11pt]
&=\ \sum\limits_{u \in U} \hsx \psi(\ell_1 u).]
\end{align*}
\vspace{0.3cm}

\qquad {\small\bf LEMMA} \ 
Let $\phi \in CL(L)$, $\psi \in CL(G)$ $-$then
\[
\langle i_{_{P \ra G}} \widetilde{\phi}, \psi \rangle_G
\hsx = \hsx
\langle \phi, r_{_{G \ra L}} \psi_P \rangle_L.
\]
\\[-.2cm]

PROOF 
\\[-1.5cm]

\allowdisplaybreaks
\begin{align*}
\langle i_{_{P \ra G}} \widetilde{\phi}, \psi \rangle_G \ 
&=\ 
\frac{1}{\abs{G}} \hsx \sum\limits_{x \in G} \hsx (i_{_{P \ra G}} \widetilde{\phi})(x)  \ov{\psi(x)}
\\[11pt]
&=\ 
\frac{1}{\abs{G}} \hsx \sum\limits_{x \in G} \hsx\hsx
\frac{1}{\abs{P}} \hsx\hsx \sum\limits_{y \in G} \hsx 
\mathring{\widetilde{\phi}}(y x y^{-1}) \ov{\psi(x)}
\\[11pt]
&=\ 
\frac{1}{\abs{G}} \hsx \sum\limits_{y \in G} \hsx\hsx
\frac{1}{\abs{P}} \hsx\hsx \sum\limits_{x \in G} \hsx 
\mathring{\widetilde{\phi}}(y x y^{-1}) \ov{\psi(x)}
\\[11pt]
&=\ 
\frac{1}{\abs{G}} \hsx \sum\limits_{y \in G} \hsx\hsx
\frac{1}{\abs{P}} \hsx\hsx \sum\limits_{x \in G} \hsx \mathring{\widetilde{\phi}}(x) \ov{\psi(y^{-1} x y)}
\\[11pt]
&=\ 
\frac{1}{\abs{G}} \hsx \sum\limits_{y \in G} \hsx\hsx
\frac{1}{\abs{P}} \hsx\hsx \sum\limits_{p \in P} \hsx \widetilde{\phi}(p) \ov{\psi(y p y^{-1})}
\\[11pt]
&=\ 
\frac{1}{\abs{G}} \hsx \sum\limits_{y \in G} \hsx\hsx
\frac{1}{\abs{P}} \hsx\hsx \sum\limits_{\ell, u} \hsx \widetilde{\phi}(\ell u) \ov{\psi(y \ell u y^{-1})}
\\[11pt]
&=\ 
\frac{1}{\abs{G}} \hsx \sum\limits_{y \in G} \hsx\hsx
\frac{1}{\abs{P}} \hsx\hsx \sum\limits_{\ell, u} \hsx \phi(\ell) \ov{\psi(y \ell u y^{-1})}
\\[11pt]
&=\ 
\frac{1}{\abs{G}} \hsx \sum\limits_{y \in G} \hsx\hsx
\frac{1}{\abs{P}} \hsx\hsx \sum\limits_{\ell, u} \hsx \phi(\ell) \ov{\psi(\ell u)}
\\[11pt]
&=\ 
\frac{1}{\abs{P}} \hsx \sum\limits_{\ell, u} \hsx \phi(\ell) \ov{\psi(\ell u)}
\\[11pt]
&=\ 
\frac{1}{\abs{L} \hsx \abs{U}} \hsx \sum\limits_{\ell, u} \hsx \phi(\ell) \ov{\psi(\ell u)}
\\[11pt]
&=\ 
\frac{1}{\abs{L}} \hsx \sum\limits_{\ell} \hsx 
\phi(\ell) \hsx 
\frac{1}{\abs{U}} \hsx \sum\limits_{u \in U} \hsx  \ov{\psi(\ell u)}
\\[11pt]
&=\ 
\frac{1}{\abs{L}} \hsx \sum\limits_{\ell} \hsx 
\phi(\ell) \hsx
\ov{\psi_P(\ell)}
\\[11pt]
&=\ 
\langle \phi, r_{_{G \ra L}} \psi_P \rangle_L.
\end{align*}


\chapter{
$\boldsymbol{\S}$\textbf{1}.\quad  ORBITAL SUMS}
\setlength\parindent{2em}
\setcounter{theoremn}{0}
\renewcommand{\thepage}{A III \S1-\arabic{page}}

\qquad Let \mG be a finite group.
\\[-.2cm]

\begin{x}{\small\bf DEFINITION} \ 
Given $f \in C(G)$ and $\gamma \in G$, put
\[
\sO(f, \gamma) \hsx = \hsx \sum\limits_{x \in G} \hsx f(x \gamma x^{-1}),
\]
the 
\un{orbital sum}
\index{orbital sum} 
of $f$ at $\gamma$.
\end{x}
\vspace{0.1cm}

\begin{x}{\small\bf LEMMA} \ 
The function $\sO(f)$ defined by the assignment
\[
\gamma \ra \sO(f,\gamma)
\]
is a class function on \mG, i.e., is an element of $CL(G)$.
\end{x}
\vspace{0.1cm}

\begin{x}{\small\bf LEMMA} \ 
There is an expansion
\[
\sO(f,\gamma) 
\hsx = \hsx 
\sum\limits_{\Pi \in \widehat{G}} \hsx \tr(\Pi^*(f)) \hsy \chisubPi(\gamma),
\]
where 
\[
\Pi^*(f)
\hsx = \hsx 
\sum\limits_{x \in G} f(x) \hsy \Pi^*(x).
\]
\vspace{0.2cm}

PROOF \ 
Since $\sO(f)$ is a class function, $\forall \ \gamma \in G$, 
\[
\sO(f,\gamma) 
\hsx = \hsx 
\sum\limits_{\Pi \in \widehat{G}} \ \langle \sO(f), \chisubPi \ranglesubG \hsx \chisubPi(\gamma) 
\qquad (\text{cf. II,} \ \S4, \#17).
\]
But
\allowdisplaybreaks
\begin{align*}
\langle \sO(f), \chisubPi \ranglesubG \ 
&=\ \frac{1}{\abs{G}}\hsx  \sum\limits_{\gamma \in G} \hsx \sO(f,\gamma) \hsx \ov{\chisubPi(\gamma)}
\\[11pt]
&=\ \frac{1}{\abs{G}}\hsx  \sum\limits_{\gamma \in G} \hsx 
\sum\limits_{x \in G} f(x \gamma x^{-1}) \hsx \ov{\chisubPi(\gamma)}
\\[11pt]
&=\ \frac{1}{\abs{G}}\hsx  \sum\limits_{x \in G} \hsx 
\sum\limits_{\gamma \in G} f(x \gamma x^{-1}) \hsx \ov{\chisubPi(\gamma)}
\\[11pt]
&=\ \frac{1}{\abs{G}}\hsx  \sum\limits_{x \in G} \hsx \sum\limits_{\gamma \in G} \hsx 
f(\gamma) \hsx  \ov{\chisubPi(x^{-1} \gamma x)}
\\[11pt]
&=\ \frac{1}{\abs{G}}\hsx  \sum\limits_{x \in G} \hsx \sum\limits_{\gamma \in G} \hsx 
f(\gamma)   \hsx \ov{\chisubPi(\gamma)}
\\[11pt]
&=\ \frac{1}{\abs{G}}\abs{G} \hsx  \sum\limits_{\gamma \in G} \hsx   f(\gamma) \hsx \ov{\chisubPi(\gamma)}
\\[11pt]
&=\ \sum\limits_{\gamma \in G} \hsx f(\gamma) \hsx \ov{\chisubPi(\gamma)}
\\[11pt]
&=\ \sum\limits_{x \in G} \hsx f(x) \hsx \ov{\chisubPi(x)}
\\[11pt]
&=\ \tr(\Pi^*(f)).
\end{align*}
\\[-.75cm]

[Note: \ 
Recall that 
\[
\chisubPiStar \hsx = \hsx \ov{\chisubPi} \qquad \text{(cf. II, \S4, \#4)}.]
\]
\end{x}
\vspace{0.1cm}

\begin{x}{\small\bf \un{N.B.}}  \ 
In terms of the Fourier transform, 
\[
\Pi^*(f) 
\hsx = \hsx
\widehat{f}(\Pi^*)
\implies 
\tr(\Pi^*(f))
\hsx = \hsx
\tr(\widehat{f}(\Pi^*)).]
\]
\end{x}
\vspace{0.1cm}


\chapter{
$\boldsymbol{\S}$\textbf{2}.\quad  THE LOCAL TRACE FORMULA}
\setlength\parindent{2em}
\setcounter{theoremn}{0}
\renewcommand{\thepage}{A III \S2-\arabic{page}}

\qquad Let \mG be a finite group.
\\[-.2cm]

\begin{x}{\small\bf NOTATION} \ 
Denote by $\pi_{L,R}$ the representation of $G \times G$ on $C(G)$ given by
\[
(\pi_{L,R} (x_1, x_2) f) (x)
\hsx = \hsx
f(x_1^{-1} x x_2) \qquad (\text{cf. II, \S1, \#14}).
\]
\vspace{0.2cm}

Define a linear bijection 
\[
T : C(G) \ra C(G \times G / G)
\]
via the prescription 
\[
Tf(x_1, x_2) 
\hsx = \hsx
f(x_1 x_2^{-1}).
\]
\end{x}
\vspace{0.1cm}

\begin{x}{\small\bf \un{N.B.}} \ 
Embed \mG diagonally into $G \times G$ $-$then $\forall \ x \in G$, 
\allowdisplaybreaks
\begin{align*}
T f ((x_1,x_2) (x,x)) \ 
&=\ 
Tf (x_1 \hsy x, x_2 \hsy x) \
\\[11pt]
&=\ 
f(x_1 \hsy x \hsy x^{-1} \hsy x_2^{-1})
\\[11pt]
&=\ 
f(x_1 \hsy x_2^{-1})
\\[11pt]
&=\ 
T f (x_1,x_2).
\end{align*}
\end{x}
\vspace{0.1cm}

\begin{x}{\small\bf NOTATION}\ 
Set 
\[
L_{G \times G / G} 
\hsx = \hsx 
\Ind_{G,\theta}^{G \times G},
\]
where $\theta$ is the trivial representation of \mG on $\E = \Cx$.
\end{x}
\vspace{0.1cm}

\begin{x}{\small\bf LEMMA} \ 
\[
T \in I_{G \times G} (\pi_{L,R}, L_{G \times G / G}).
\]
\vspace{0.2cm}

PROOF \ 
$\forall \ x_1, x_2 \in G$, $\forall \ f \in C(G)$,
\allowdisplaybreaks
\begin{align*}
(T \pi_{L,R} (x_1,x_2) f) (y_1,y_2)\
&=\ (\pi_{L,R} (x_1,x_2) f) (y_1 y_2^{-1})
\\[11pt]
&=\ f(x_1^{-1} y_1 y_2^{-1} x_2).
\end{align*}
And
\allowdisplaybreaks
\begin{align*}
(L_{G \times G / G} (x_1, x_2) Tf) (y_1,y_2) \ 
&=\ T f ( (x_1,x_2)^{-1} (y_1,y_2)) 
\\[11pt]
&=\ T f ( (x_1^{-1},x_2^{-1}) (y_1,y_2)) 
\\[11pt]
&=\ T f(x_1^{-1} y_1, x_2^{-1} y_2)
\\[11pt]
&=\ f(x_1^{-1} y_1 y_2^{-1} x_2).
\end{align*}
\end{x}
\vspace{0.1cm}

\begin{x}{\small\bf \un{N.B.}}  \ 
\mT is unitary: \ $\forall \ f, \ g \in C(G)$, 
\[
\langle Tf, Tg \ranglesubGxG
\hsx = \hsx 
\langle f, g \ranglesubG.
\]

[By definition,
\allowdisplaybreaks
\begin{align*}
\langle Tf, Tg \ranglesubGxG \ 
&=\ \frac{1}{\abs{G \times G}} \hsx \sum\limits_{(x_1, x_2) \in G \times G} \hsx 
Tf(x_1,x_2) \ov{Tg(x_1,x_2)}
\\[11pt]
&=\ 
\frac{1}{\abs{G \times G}} \hsx \sum\limits_{(x_1, x_2) \in G \times G} \hsx 
f(x_1 x_2^{-1}) \hsx \ov{g(x_1 x_2^{-1})} 
\\[11pt]
&=\ 
\frac{1}{\abs{G \times G}} \hsx \sum\limits_{x_1 \in G} \hsx 
\sum\limits_{x_2 \in G}\hsx  f(x_1 x_2^{-1}) \hsx \ov{g(x_1 x_2^{-1})} 
\\[11pt]
&=\ 
\frac{1}{\abs{G \times G}} \hsx 
\sum\limits_{x_1 \in G} \hsx 
\sum\limits_{x_2 \in G}\hsx  f(x_1 x_2) \hsx \ov{g(x_1 x_2)} 
\\[11pt]
&=\ \frac{1}{\abs{G \times G}} \hsx 
\sum\limits_{x_1 \in G} \hsx 
\sum\limits_{x_2 \in G}\hsx  f(x_2) \hsx \ov{g(x_2)} 
\\[11pt]
&=\ \frac{\abs{G}}{\abs{G \times G}} \hsx 
\sum\limits_{x \in G} \hsx 
f(x) \hsx \ov{g(x)} 
\\[11pt]
&=\ \frac{1}{\abs{G}} \hsx 
\sum\limits_{x \in G} 
\hsx f(x) \hsx \ov{g(x)} 
\\[11pt]
&=\ \langle f, g \ranglesubG \hsy.]
\end{align*}

\end{x}
\vspace{0.1cm}

\begin{x}{\small\bf NOTATION} \ 
Given an $x \in G$, write $C(x)$ for its conjugacy class and $G_x$ for its centralizer (cf. II, \S4, \#10).
\end{x}
\vspace{0.1cm}

\begin{x}{\small\bf EXAMPLE} \ 
$\forall \ f \in C(G)$ and $\forall \ \gamma \in G$, 
\[
\sO(f,\gamma) 
\ = \ 
\sum\limits_{x \in G} \hsx f(x \gamma x^{-1}) 
\ = \ 
\abs{G_\gamma} \hsx 
\sum\limits_{x \in C(\gamma)} \hsx f(x).
\]
\end{x}
\vspace{0.1cm}

\begin{x}{\small\bf LEMMA} \ 
Abbreviate $\chisubpiLR$ to $\chisubLR$ $-$then 
\[
\chisubLR (x_1,x_2) 
\hsx = \hsx 
\abs{\{x \in G : \delta_x (x_1^{-1} x x_2) \neq 0\}}
\]
\[
\hspace{3.25cm}
= \ 
\begin{cases}
\ \abs{G_x} \qquad (x = x_1 = x_2)\\
\ 0 \hspace{1.35cm} (C(x_1) \neq C(x_2))
\end{cases}
.
\]
\vspace{0.2cm}

[Work instead with the character of $L_{G \times G / G}$ and apply II, \S7, \#11.]
\end{x}
\vspace{0.1cm}

Given $f_1$, $f_2 \in C(G)$, define $f \in C(G \times G)$ by 
\[
f(x_1, x_2) \hsx = \hsx f_1(x_1) f_2(x_2), 
\]
and let 
\[
\pi_{L,R}(f) 
\ = \ 
 \sum\limits_{x_1 \in G} \hsx 
 \sum\limits_{x_2 \in G} \hsx
 f_1(x_1) f_2(x_2) 
 \pi_{L,R}(x_1, x_2).
\]
Then $\forall \ \phi \in C(G)$, 
\allowdisplaybreaks
\begin{align*}
(\pi_{L,R} (f) \phi) (x) \ 
&=\ \sum\limits_{x_1 \in G} \hsx \sum\limits_{x_2 \in G} \hsx 
f_1(x_1) f_2(x_2) \phi(x_1^{-1} x x_2)\\
&=\ \sum\limits_{y \in G} \hsx K_f(x,y) \phi(y), 
\end{align*}
where 
\[
K_f(x,y) 
\hsx = \hsx
\sum\limits_{z \in G} \hsx f_1(x z) f_2(z y).
\]
Therefore $\pi_{L,R} (f)$ is an integral operator on $C(G)$ (a.k.a. $L^2(G) \ldots$) with kernel $K_f(x,y)$.
\vspace{0.1cm}

\begin{x}{\small\bf LEMMA} \ 
The $\sqrt{\abs{G}} \ \delta_x$ $(x \in G)$ constitute an orthonormal basis for $C(G)$.
\end{x}
\vspace{0.1cm}

\begin{x}{\small\bf LEMMA} \ 
$\forall \ f = f_1 f_2$, 
\[
\tr(\pi_{L,R}(f))
\ = \ 
\sum\limits_{x \in G} \hsx
K_f(x,x).
\]
\vspace{0.2cm}

PROOF \
In fact, 
\allowdisplaybreaks
\begin{align*}
\tr(\pi_{L,R}(f)) \ 
&=\ 
\sum\limits_{x \in G} \
\langle \pi_{L,R}(f) \hsx \sqrt{\abs{G}} \hsx \delta_x, \sqrt{\abs{G}} \hsx \delta_x  \ranglesubG
\\[15pt]
&=\ 
\abs{G} \
\sum\limits_{x \in G} \
\frac{1}{\abs{G}} \
\sum\limits_{y \in G} \
(\pi_{L,R} (f) \delta_x) (y) \delta_x(y)
\\[15pt]
&=\ 
\sum\limits_{x \in G} \
(\pi_{L,R} (f) \delta_x) (x)
\\[15pt]
&=\ 
\sum\limits_{x \in G} \
\sum\limits_{x_1 \in G} \
\sum\limits_{x_2 \in G} \
f_1(x_1) f_2(x_2) \delta_x (x_1^{-1} x x_2)
\\[15pt]
&=\ 
\sum\limits_{x \in G} \hsx
\sum\limits_{z \in G} \hsx
f_1(x z) f_2(z y)
\\[15pt]
&=\ 
\sum\limits_{x \in G} \hsx
K_f(x,x).
\end{align*}
\end{x}
\vspace{0.1cm}

Enumerate the elements of $\CON(G)$, say
\[
\CON(G) 
\ = \ 
\{C_1, \ldots, C_n\}.
\]
For each $i$, fix a $\gamma_i \in C_i$ $(1 \leq i \leq n)$.
\\[-.2cm]

\begin{x}{\small\bf LEMMA} \ 
$\forall \ f = f_1 f_2$, 
\[
\sum\limits_{x \in G} \hsx K_f(x,x)
\ = \ 
\sum\limits_{i=1}^n \
\frac{1}{\abs{G_{\gamma_i}}} \hsx 
\sO(f_1, \gamma_i) \sO(f_2, \gamma_i).
\]
\vspace{0.2cm}

PROOF \ 
Start with the LHS:
\allowdisplaybreaks
\begin{align*}
\sum\limits_{x \in G} \hsx K_f(x,x) \ 
&=\
\sum\limits_{x \in G} \
\sum\limits_{z \in G} \
f_1(xz) f_2(zx)
\\[15pt]
&=\
\sum\limits_{x \in G} \
\sum\limits_{y \in G} \
f_1(y) f_2(x^{-1} y x)
\\[15pt]
&=\
\sum\limits_{y \in G} \
F(y),
\end{align*}
where
\[
F(y) 
\ = \ 
\sum\limits_{x \in G} \
f_1(y) f_2(x^{-1} y x ).
\]
Using now \S4, \#2 below, we have 
\[
\sum\limits_{y \in G} \
F(y)
\ = \ 
\sum\limits_{i=1}^n \
\frac{1}{\abs{G_{\gamma_i}}} \hsx \sO(F,\gamma_i).
\]
And
\allowdisplaybreaks
\begin{align*}
\sO(F,\gamma_i) \ 
&=\
\sum\limits_{x \in G} \
F(x \gamma_i x^{-1})
\\[15pt]
&=\ 
\sum\limits_{x \in G} \
\sum\limits_{y \in G} \
f_1(x \gamma_i x^{-1}) f_2(y^{-1} x \gamma_i x^{-1} y)
\\[15pt]
&=\ 
\sum\limits_{x \in G} \
f_1(x \gamma_i x^{-1}) \
\sum\limits_{y \in G} \
f_2(y^{-1} x \gamma_i x^{-1} y)
\\[15pt]
&=\ 
\sum\limits_{x \in G} \
f_1(x \gamma_i x^{-1}) \
\sum\limits_{y \in G} \
f_2(y^{-1} \gamma_i  y)
\\[15pt]
&=\ 
\sum\limits_{x \in G} \
f_1(x \gamma_i x^{-1}) \
\sum\limits_{y \in G} \
f_2(y\gamma_i  y^{-1} )
\\[15pt]
&=\ 
\sO(f_1, \gamma_i) \hsx \sO(f_2, \gamma_i).
\end{align*}

\end{x}
\vspace{0.1cm}

\begin{x}{\small\bf LEMMA} \ 
$\forall \ f = f_1 f_2$, 
\[
\sum\limits_{x \in G} \hsx K_f(x,x)
\ = \ 
\sum\limits_{\Pi \in \widehat{G}} \hsx
\tr( \widehat{f}_1 (\Pi)) \hsx \tr( \widehat{f}_2 (\Pi^*)).
\]
\vspace{0.2cm}

PROOF \ 
Write
\allowdisplaybreaks
\begin{align*}
\sum\limits_{x \in G} \hsx K_f(x,x) \ 
&=\ 
\sum\limits_{x \in G} \
\sum\limits_{y \in G} \
f_1(y) f_2(x^{-1} y x)
\\[15pt]
&=\ 
\sum\limits_{y \in G} \
f_1(y) 
\sum\limits_{x \in G} \
f_2(x y x^{-1})
\\[15pt]
&=\ 
\sum\limits_{y \in G} \
f_1(y) \sO(f_2,y)
\\[15pt]
&=\ 
\sum\limits_{y \in G} \
f_1(y)
\sum\limits_{\Pi \in \widehat{G}} \
\tr( \Pi^* (f_2)) \chisubPi(y) 
\qquad \text{(cf. \S1, \#3)}
\\[15pt]
&=\ 
\sum\limits_{\Pi \in \widehat{G}} \
\bigg(
\sum\limits_{y \in G} \hsx
f_1(y) \chisubPi(y) \bigg)\tr( \Pi^* (f_2)\bigr)
\\[15pt]
&=\ 
\sum\limits_{\Pi \in \widehat{G}} \
\tr( \Pi(f_1)) \tr( \Pi^* (f_2))
\\[15pt]
&=\ 
\sum\limits_{\Pi \in \widehat{G}} \hsx
\tr( \widehat{f}_1 (\Pi)) \tr( \widehat{f}_2 (\Pi^*)).
\end{align*}
\end{x}
\vspace{0.1cm}

\begin{x}{\small\bf DEFINITION} \ 
Given $f = f_1 f_2$, the 
\un{local trace formula}
\index{local trace formula} 
is the relation
\[
\sum\limits_{\Pi \in \widehat{G}} \hsx
\tr( \widehat{f}_1 (\Pi)) \tr( \widehat{f}_2 (\Pi^*))
\hsx = \hsx 
\sum\limits_{i=1}^n \hsx 
\frac{1}{\abs{G_{\gamma_i}}} \hsx 
\sO(f_1, \gamma_i) \hsx \sO(f_2,\gamma_i).
\]
\end{x}
\vspace{0.1cm}

\begin{x}{\small\bf EXAMPLE} \ 
Suppose that $f_1 = f_2$ is real valued, call it $\phi$ $-$then
\allowdisplaybreaks
\begin{align*}
\tr( \widehat{\phi} (\Pi^*)) \ 
&=\ 
\sum\limits_{x \in G} \hsx
\phi(x) \chisubPiStar(x)
\\[15pt]
&=\ 
\sum\limits_{x \in G} \hsx
{\phi(x)} \ov{\chisubPi(x)}
\\[15pt]
&=\ 
\sum\limits_{x \in G} 
\ov{\phi(x)} \ \ov{\chisubPi(x)}
\\[15pt]
&=\ 
\ov{
\sum\limits_{x \in G} \hsx
\phi(x) \chisubPi(x)
}\\[15pt]
&=\ 
\ov{\tr(\widehat{\phi}(\Pi))}.
\end{align*}

Therefore
\allowdisplaybreaks
\begin{align*}
\sum\limits_{\Pi \in \widehat{G}} \ \tr(\widehat{\phi}(\Pi)) \ov{\tr(\widehat{\phi}(\Pi))} \ 
&=\ 
\sum\limits_{\Pi \in \widehat{G}} \
\abs{\tr(\widehat{\phi}(\Pi))}^2\\[15pt]
&=\ 
\sum\limits_{i=1}^n \ 
\frac{1}{\abs{G_{\gamma_i}}} 
\sO(\phi,\gamma_i)^2.
\end{align*}
\vspace{0.2cm}

[Note: \ 
Specialize and take $f_1 = f_2 = \delta_e$ $-$then 
\allowdisplaybreaks
\begin{align*}
\tr(\widehat{\delta}_e(\Pi)) \ 
&=\ 
\chisubpi(e)\\
&=\ d_\Pi
\end{align*}
and
\[
\sO(\delta_e,\gamma_i) 
\hsx = \hsx 
0 \qquad (\gamma_i \neq e)
\]
while 
\allowdisplaybreaks
\begin{align*}
\sO(\delta_e,e) \ 
&=\ 
\sum\limits_{x \in G} \hsx
\delta_e(x e x^{-1})\\[15pt]
&=\ 
\sum\limits_{x \in G} \hsx
\delta_e(e)\\[15pt]
&=\ 
\abs{G}.
\end{align*}
Consequently, 
\[
\sum\limits_{\Pi \in \widehat{G}} \ 
d_\Pi^2
\ = \ 
\frac{1}{\abs{G}} \hsx \abs{G}^2
\ = \  
\abs{G} \qquad \text{(cf. II, \S3, \#5 and II, \S5, \#9)}.]
\]
\end{x}
\vspace{0.1cm}

From the definitions, 
\[
\begin{cases}
\ \tr(\widehat{f}_1 (\Pi)) 
\hsx = \hsx 
\sum\limits_{x \in G} \hsx f_1(x) \chisubPi(x) 
\hsx = \hsx 
\abs{G} \langle f_1, \chisubPiStar \ranglesubG\\[15pt]
\ \tr(\widehat{f}_2 (\Pi^*)) 
\hsx = \hsx 
\sum\limits_{x \in G} \hsx f_2(x) \chisubPiStar(x) 
\hsx = \hsx 
\abs{G} \langle f_2, \chisubPi \ranglesubG
\end{cases}
.
\]
Therefore
\allowdisplaybreaks
\begin{align*}
\sum\limits_{\Pi \in \widehat{G}} \tr(\widehat{f}_1 (\Pi)) \tr(\widehat{f}_2 (\Pi^*)) \ 
&=\ 
\abs{G}^2 \
\sum\limits_{\Pi \in \widehat{G}} \
\langle f_1, \chisubPiStar \ranglesubG
\langle f_2, \chisubPi \ranglesubG
\\[15pt]
&=\ 
\abs{G}^2 \
\sum\limits_{\Pi \in \widehat{G}} \
\langle f_1, \chisubPi \ranglesubG
\langle f_2, \chisubPiStar \ranglesubG.
\end{align*}
\vspace{0.1cm}

\begin{x}{\small\bf \un{N.B.}} \ 
Assume in addition that $f_1$ and $f_2$ are class functions.  Write
\[
\begin{cases}
\ f_1(x) \hsx = \hsx \sum\limits_{\Pi \in \widehat{G}} \hsx \langle f_1, \chisubPi \ranglesubG \chisubPi (x) \\[15pt]
\ \ov{f_2(x)} \hsx = \hsx \sum\limits_{\Pi \in \widehat{G}} \hsx \langle \ov{f_2}, \chisubPi \ranglesubG \chisubPi (x) 
\end{cases}
\qquad\qquad \text{(cf. II, \S4, \#17)}.
\]
Then
\allowdisplaybreaks
\begin{align*}
\langle f_1, \ov{f_2} \ranglesubG \ 
&=\  \sum\limits_{\Pi \in \widehat{G}} \ 
\langle f_1, \chisubPi \ranglesubG
\ov{\langle \ov{f_2}, \chisubPi \ranglesubG}
\qquad \text{(first orthogonality relations)}
\\[11pt]
&=\  \sum\limits_{\Pi \in \widehat{G}} \ 
\langle f_1, \chisubPi \ranglesubG
\langle f_2, \chisubPiStar \ranglesubG.
\end{align*}
On the other hand, 
\allowdisplaybreaks
\begin{align*}
\langle f_1, \ov{f_2} \ranglesubG \ 
&=\ \frac{1}{\abs{G}} \ \sum\limits_{x \in G} \hsx f_1(x) \ov{\ov{f_2(x)}}
\\[11pt]
&=\ 
\frac{1}{\abs{G}} \hsx
\sum\limits_{i = 1}^n \
\sum\limits_{x \in C_i} \ f_1(x) f_2(x)
\\[11pt]
&=\ 
\frac{1}{\abs{G}} \ \sum\limits_{i = 1}^n \ \abs{C_i} \hsx f_1(\gamma_i) f_2(\gamma_i)
\\[11pt]
&=\ \sum\limits_{i = 1}^n \ \frac{\abs{C_i}}{\abs{G}} \hsx f_1(\gamma_i) f_2(\gamma_i)
\\[11pt]
&=\ \sum\limits_{i = 1}^n \ \frac{1}{\abs{G_{\gamma_i}}} \hsx f_1(\gamma_i) f_2(\gamma_i)
\end{align*}

\qquad\qquad$\implies$
\allowdisplaybreaks
\begin{align*}
\abs{G}^2 \langle f_1, \ov{f_2} \ranglesubG \ 
&=\ \sum\limits_{i = 1}^n \hsx \frac{1}{\abs{G_{\gamma_i}}} \ 
\abs{G} \ f_1(\gamma_i) \abs{G} \hsx f_2(\gamma_i)\\
&=\ \sum\limits_{i = 1}^n \hsx \frac{1}{\abs{G_{\gamma_i}}} \ 
\sO(f_1, \gamma_i) \hsx \sO(f_2, \gamma_i).
\end{align*}
\\[-.2cm]

The irreducible representations of $G \times G$ are the outer tensor products
\[
\Pi_1 \hsx \un{\otimes} \hsx \Pi_2 \quad (\Pi_1, \Pi_2 \in \widehat{G}) \qquad \text{(cf. II, \S5, \#13)}.
\]
Moreover, 
\[
\chisubPiOneuXTwo 
\ = \ 
\chisubPiOne \hsx \chisubPiTwo.
\]

Consider now the direct sum decomposition 
\[
L_{G \times G / G}
\hsx = \hsx 
\bigoplus\limits_{\Pi_1, \Pi_2 \in \widehat{G}} \hsx 
m(\Pi_1 \hsx \un{\otimes} \hsx \Pi_2, L_{G \times G / G}) \Pi_1 \hsx \un{\otimes} \hsx \Pi_2.
\]
Then 
\[
\tr(L_{G \times G / G}(f))
\ = \ 
\sum\limits_{\Pi_1, \Pi_2 \in \widehat{G}} 
\
m(\Pi_1 \hsx \un{\otimes} \hsx \Pi_2, L_{G \times G / G}) 
\
\tr(\Pi_1(f_1)) \hsx \tr(\Pi_2(f_2)).
\]
I.e. (cf. \#4):
\[
\tr(\Pi_{L,R} (f))
\ = \ 
\sum\limits_{\Pi_1, \Pi_2 \in \widehat{G}} 
\
m(\Pi_1 \hsx \un{\otimes} \hsx \Pi_2, L_{G \times G / G}) 
\hsx
\tr(\Pi_1(f_1)) \hsx \tr(\Pi_2(f_2)).
\]
I.e. (cf. \#12):
\[
\sum\limits_{\Pi \in \widehat{G}} \tr (\Pi(f_1)) \hsx \tr(\Pi^*(f_2)) 
\hsx = \hsx 
\sum\limits_{\Pi_1, \Pi_2 \in \widehat{G}}  \hsx
m(\Pi_1 \hsx \un{\otimes} \hsx \Pi_2, L_{G \times G / G}) 
\hsx
\tr(\Pi_1(f_1)) \hsx \tr(\Pi_2(f_2)).
\]

Therefore, thanks to I, \S3, \#9, 
\[
m(\Pi_1 \hsx \un{\otimes} \hsx \Pi_2, L_{G \times G / G})
\]
must vanish unless $\Pi_1 = \Pi$, $\Pi_2 = \Pi^*$, in which case the coefficient is equal to 1.
\end{x}
\vspace{0.1cm}

\begin{x}{\small\bf SCHOLIUM} \ 

\[
\Pi_{L,R} 
\ \approx \
\bigoplus\limits_{\Pi \in \widehat{G}} \hsx \Pi \hsx \un{\otimes} \hsx \Pi^*.
\]
\end{x}


\chapter{
$\boldsymbol{\S}$\textbf{3}.\quad  THE GLOBAL PRE-TRACE FORMULA}
\setlength\parindent{2em}
\setcounter{theoremn}{0}
\renewcommand{\thepage}{A III \S3-\arabic{page}}

\qquad Let \mG be a finite group, $\Gamma \subset G$ a subgroup.
\\[-.2cm]

\begin{x}{\small\bf NOTATION} \ 
Set 
\[
L_{G / \Gamma} 
\hsx = \hsx
\Ind_{\Gamma, \theta}^G,
\]
where $\theta$ is the trivial representation of $\Gamma$ on $\E = \Cx$.
\vspace{0.2cm}

[Note: \ 
Accordingly, $\chisubtheta = 1_\Gamma$ and $\E_{\Gamma, \theta}^G = C(G / \Gamma)$.]
\end{x}
\vspace{0.1cm}

\begin{x}{\small\bf EXAMPLE} \ 
In the special case when $\Gamma = \{e\}$, $L_{G / \Gamma} = L$, the left translation representation of \mG on $C(G)$ 
(cf. II, \S1, \#12).
\end{x}
\vspace{0.1cm}

\begin{x}{\small\bf \un{N.B.}} \ 
The pair $(G \hsx\times\hsx G, G)$ figuring in \S2 is an instance of the overall setup.
\end{x}
\vspace{0.1cm}

Given $f \in C(G)$, $\phi \in C(G/\Gamma)$, we have
\allowdisplaybreaks
\begin{align*}
\bigl(L_{G / \Gamma}  (f) \phi \bigr) (x) 
&=\ \sum\limits_{y \in G} \ 
f(y) \hsx \bigl(L_{G / \Gamma}  (y) \phi \bigr) (x) 
\\[15pt]
&=\ \sum\limits_{y \in G} \ 
f(y)  \phi (y^{-1} x) 
\\[15pt]
&=\ 
\sum\limits_{y \in G} \ 
f(x y^{-1})  \hsx \phi (y) 
\\[15pt]
&=\ \sum\limits_{y \in G} \ 
f(x y^{-1})  \hsx 
\frac{1}{\abs{\Gamma}} \ \sum\limits_{\gamma \in \Gamma} \ \phi(y \gamma) 
\\[15pt]
&=\ \sum\limits_{y \in G} \ 
\frac{1}{\abs{\Gamma}} \ \sum\limits_{\gamma \in \Gamma} \ f(x y^{-1})  \hsx \phi(y \gamma) 
\\[15pt]
&=\  
\sum\limits_{\gamma \in \Gamma} \ \frac{1}{\abs{\Gamma}} \
\sum\limits_{y \in G} \hsx f(x y^{-1})  \ \phi(y \gamma) 
\\[15pt]
&=\ 
\sum\limits_{\gamma \in \Gamma} \hsx \frac{1}{\abs{\Gamma}} \
\sum\limits_{y \in G} \hsx f(x \gamma y^{-1})  \hsx \phi(y) 
\\[15pt]
&=\ 
\sum\limits_{y \in G} \ \bigg(
\frac{1}{\abs{\Gamma}} \ 
\sum\limits_{\gamma \in \Gamma} \ 
f(x \gamma y^{-1})  \bigg) \hsx \phi(y) 
\\[15pt]
&=\ \sum\limits_{y \in G} \ K_f(x,y) \hsx \phi(y),
\end{align*}
where
\[
 K_f(x,y)
 \ = \ 
 \frac{1}{\abs{\Gamma}} \ \sum\limits_{\gamma \in \Gamma} \ f(x \gamma y^{-1}).
\]
\vspace{0.1cm}

To summarize:
\\[-.2cm]

\begin{x}{\small\bf LEMMA} \ 
$\forall \ f \in C(G)$, $\forall \ \phi \in C(G/\Gamma)$, $\forall \ x \in G$, 
\[
\bigl(
L_{G / \Gamma} \hsx (f) \phi \bigr) (x) 
\ = \ 
\sum\limits_{y \in G} \
K_f(x,y) \hsx \phi(y), 
\]
where 
\[
K_f(x,y)
\hsx = \hsx 
\frac{1}{\abs{\Gamma}} \
\sum\limits_{\gamma \in \Gamma} \
f(x \gamma y^{-1}).
\]
\end{x}
\vspace{0.1cm}

Write
\[
G
\ = \ 
\coprod\limits_{k=1}^n \
x_k \Gamma.
\]
Then for any $f \in C(G)$, 
\[
\sum\limits_{x \in G} \
f(x) 
\ = \ 
\sum\limits_{k=1}^n \
\sum\limits_{\gamma \in \Gamma} \
f(x_k \gamma),
\]
thus for any $\phi \in C(G/\Gamma)$, 
\[
\sum\limits_{x \in G} \
\phi(x) 
\ = \ 
\abs{\Gamma} \
\sum\limits_{k=1}^n \
\phi(x_k).
\]
\\[-.2cm]

\begin{x}{\small\bf RAPPEL} \ 
(cf. II, \S9, \#2) \ The Hilbert space structure on $C(G/\Gamma)$ is defined by the inner product
\allowdisplaybreaks
\begin{align*}
\langle \phi, \psi \ranglesubtheta \ 
&=\ 
\frac{1}{\abs{G}} \
\sum\limits_{x \in G} \hsx
\phi(x) \ov{\psi(x)}\\[15pt]
&=\ 
\frac{\abs{\Gamma}}{\abs{G}} \
\sum\limits_{k=1}^n \hsx
\phi(x_k) \ov{\psi(x_k)}.
\end{align*}
\end{x}
\vspace{0.1cm}

\begin{x}{\small\bf NOTATION} \ 
Define functions $\delta_k \in C(G/\Gamma)$ by the rule 
\[
\delta_k(x_\ell \hsy \gamma) 
\hsx = \hsx 
\delta_{k \ell} 
\qquad (1 \leq k, \ \ell \leq n).
\]
\end{x}
\vspace{0.1cm}

\begin{x}{\small\bf LEMMA} \ 
The
\[
\Delta_k 
\ = \ 
\bigg( \frac{\abs{G}}{\abs{\Gamma}} \bigg)^{1/2} \delta_k
\]
constitute an orthonormal basis for $C(G/\Gamma)$.
\vspace{0.2cm}

PROOF \ 
A given $\phi \in C(G / \Gamma)$ admits the decomposition
\[
\phi
\ = \ 
\sum\limits_{k = 1}^n \ \phi(x_k) \hsx \delta_k.
\]
In addition, 
\allowdisplaybreaks
\begin{align*}
\langle \Delta_k, \Delta_\ell \ranglesubtheta \ 
&=\ \frac{\abs{\Gamma}}{\abs{G}} \ \sum\limits_{j = 1}^n \hsx \Delta_k(x_j) \hsx \Delta_\ell(x_j) \\[15pt]
&=\ \frac{\abs{\Gamma}}{\abs{G}} \ \frac{\abs{G}}{\abs{\Gamma}} \\[15pt]
&=\ 1
\end{align*}
if $k = \ell$ and is 0 otherwise.
\end{x}
\vspace{0.1cm}

\begin{x}{\small\bf LEMMA} \ 
$\forall \ f \in C(G)$,
\[
\tr \bigl(L_{G / \Gamma} (f) \bigr)
\ = \
\sum\limits_{x \in G} \hsx 
K_f(x,x).
\]
\vspace{0.2cm}

PROOF \ 
In fact, 
\allowdisplaybreaks
\begin{align*}
\tr \bigl(L_{G / \Gamma} (f) \bigr) \ 
&=\ 
\sum\limits_{k=1}^n \hsx 
\langle L_{G / \Gamma}(f) \Delta_k, \Delta_k \ranglesubtheta
\\[15pt]
&=\ 
\sum\limits_{k=1}^n \ 
\frac{\abs{\Gamma}}{\abs{G}} \
\frac{\abs{G}}{\abs{\Gamma}} \
\sum\limits_{\ell=1}^n \ 
\bigl(L_{G / \Gamma} (f) \delta_k\bigr) (x_\ell) \delta_k (x_\ell)
\\[15pt]
&=\ 
\sum\limits_{k=1}^n \ 
\bigl(L_{G / \Gamma} (f) \delta_k\bigr) (x_k)
\\[15pt]
&=\ 
\sum\limits_{k=1}^n \ 
\sum\limits_{y \in G} \ 
f(x_k y^{-1}) \delta_k(y)
\\[15pt]
&=\ 
\sum\limits_{k=1}^n \ 
\sum\limits_{\ell=1}^n \ 
\sum\limits_{\gamma \in \Gamma} \
f(x_k \hsx \gamma^{-1} \hsx x_\ell^{-1}) \hsx \delta_k(x_\ell \gamma)
\\[15pt]
&=\ 
\sum\limits_{k=1}^n \ 
\sum\limits_{\ell=1}^n \ 
\sum\limits_{\gamma \in \Gamma} \
f(x_k \hsx \gamma^{-1} \hsx x_\ell^{-1}) \hsx \delta_{k \ell}
\\[15pt]
&=\ 
\sum\limits_{k=1}^n \ 
\sum\limits_{\gamma \in \Gamma} \
f(x_k \hsx \gamma^{-1} \hsx x_k^{-1})
\\[15pt]
&=\ 
\sum\limits_{k=1}^n \ 
\sum\limits_{\gamma \in \Gamma} \
f(x_k \hsx \gamma \hsx x_k^{-1})
\\[15pt]
&=\ 
\frac{1}{\abs{\Gamma}} \
\sum\limits_{k=1}^n \ 
\abs{\Gamma}
f(x_k \hsx \gamma \hsx x_k^{-1})
\\[15pt]
&=\ 
\frac{1}{\abs{\Gamma}} \
\sum\limits_{k=1}^n \ 
\sum\limits_{\eta \in \Gamma} \
\sum\limits_{\gamma \in \Gamma} \
f(x_k \hsx \eta \hsx \gamma  \hsx \eta^{-1} \hsx x_k^{-1})
\\[15pt]
&=\ 
\frac{1}{\abs{\Gamma}}  \
\sum\limits_{x \in G} \
\sum\limits_{\gamma \in \Gamma} \
f(x \hsx \gamma \hsx x^{-1})
\\[15pt]
&=\ 
\sum\limits_{x \in G} \
\frac{1}{\abs{\Gamma}}
\sum\limits_{\gamma \in \Gamma} \
f(x \hsx \gamma \hsx x^{-1})
\\[15pt]
&=\ 
\sum\limits_{x \in G} \
K_f(x,x).
\end{align*}

\end{x}
\vspace{0.1cm}

\begin{x}{\small\bf EXAMPLE} \ 
Take $\Gamma = G$ $-$then $\forall \ f \in C(G)$, 
\allowdisplaybreaks
\begin{align*}
\tr \bigl(L_{G / G} (f) \bigr) \ 
&=\ 
\sum\limits_{x \in G} \
\frac{1}{\abs{G}} \
\sum\limits_{y \in G} \
f(x \hsx y \hsx x^{-1})\\[15pt]
&=\ 
\frac{1}{\abs{G}} \
\sum\limits_{x \in G} \
\sum\limits_{y \in G} \
f(x \hsx y \hsx x^{-1})\\[15pt]
&=\ 
\frac{1}{\abs{G}} \
\sum\limits_{x \in G} \
\sum\limits_{y \in G} \
f(y)\\[15pt]
&=\ 
\frac{\abs{G}}{\abs{G}} \
\sum\limits_{y \in G} \
f(y)\\[15pt]
&=\ 
\sum\limits_{x \in G} \ f(x).
\end{align*}
\end{x}
\vspace{0.1cm}

\begin{x}{\small\bf EXAMPLE} \ 
Fix $C \in \CON(G)$ and $x \in C$ $-$then
\[
\abs{C}\chisubLsubGmodGamma (x) 
\hsx = \hsx 
\frac{\abs{G}}{\abs{\Gamma}} \hsx \abs{C \cap \Gamma} \qquad \text{(cf. II, \S9, \#6)}.
\]
\vspace{0.2cm}

[Work with $f = \chisubC$, thus
\allowdisplaybreaks
\begin{align*}
\tr(L_{G / \Gamma} (\chisubC)) \ 
&=\ 
\sum\limits_{y \in G} \ 
\chisubC (y) \chisubLsubGmodGamma (y)\\[15pt]
&=\ 
\sum\limits_{y \in C} \ 
\chisubC (y) \chisubLsubGmodGamma (y)\\[15pt]
&=\ 
\abs{C} \chisubLsubGmodGamma (x) .
\end{align*}

Meanwhile
\allowdisplaybreaks
\begin{align*}
\tr(L_{G / \Gamma} (\chisubC)) \ 
&=\ 
\sum\limits_{y \in G} \hsx 
K_{\chisubC} (y,y)\\[15pt]
&=\ 
\sum\limits_{y \in G} \ 
\frac{1}{\abs{\Gamma}} \
\sum\limits_{\gamma \in \Gamma} \ 
\chisubC(y \gamma y^{-1})\\[15pt]
&=\ 
\sum\limits_{y \in G} \ 
\frac{1}{\abs{\Gamma}} \
\sum\limits_{\gamma \in \Gamma} \ 
\chisubC(\gamma)\\[15pt]
&=\ 
\frac{\abs{G}}{\abs{\Gamma}}
\abs{C \cap \Gamma}.
\end{align*}

On general grounds, there is a direct sum decomposition
\[
L_{G/\Gamma}
\ = \
\bigoplus\limits_{\Pi \in \widehat{G}} \hsx m(\Pi, L_{G/\Gamma})  \hsx \Pi.
\]

[Note: \ 
\[
m(\Pi, L_{G/\Gamma}) \hsx \neq \hsx 0
\]
iff the restriction of $\Pi$ to $\Gamma$ contains the trivial representation $\theta$ of $\Gamma$ on $\E = \Cx$ 
(cf. II, \S9, \#9) (but see below (cf. \#14)).]
\end{x}
\vspace{0.1cm}

\begin{x}{\small\bf SCHOLIUM} \ 
$\forall \ f \in C(G)$,
\[
\tr(L_{G / \Gamma}(f)) 
\ = \
\sum\limits_{\Pi \in \widehat{G}} \ m(\Pi, L_{G/\Gamma})  \hsx \tr(\widehat{f}(\Pi)).
\]
\vspace{0.2cm}

[Note: \ 
Explicated, 
\[
\tr(\widehat{f}(\Pi)
\ = \
\sum\limits_{x \in G} \hsx 
f(x) \hsx \chisubPi(x) 
\hsx = \hsx
\tr(\Pi(f)).]
\]
\end{x}
\vspace{0.1cm}

\begin{x}{\small\bf DEFINITION} \ 
Given $f \in C(G)$, the 
\un{global pre-trace formula}
\index{global pre-trace formula} 
is the relation
\allowdisplaybreaks
\begin{align*}
\sum\limits_{\Pi \in \widehat{G}} \hsx m(\Pi, L_{G/\Gamma}) \hsx \tr(\widehat{f}(\Pi))
&=\ \sum\limits_{x \in G} \hsx K_f(x,x) 
\\[15pt]
&=\ 
\sum\limits_{x \in G} \ \frac{1}{\abs{\Gamma}} \ \sum\limits_{\gamma \in \Gamma} \hsx f(x \gamma x^{-1}).
\end{align*}
\end{x}
\vspace{0.1cm}
\begin{x}{\small\bf APPLICATION} \ 
Take $\Gamma = \{e\}$ $-$then
\[
\sum\limits_{\Pi \in \widehat{G}} \ m(\Pi, L_{G/\Gamma}) \hsx \tr(\widehat{f}(\Pi))
\]
becomes
\[
\sum\limits_{\Pi \in \widehat{G}} \ d_\Pi \hsx \tr(\widehat{f}(\Pi))
\]
while
\[
\sum\limits_{x \in G} \ \frac{1}{\abs{\Gamma}} \ \sum\limits_{\gamma \in \Gamma} \ f(x \gamma x^{-1})
\]
becomes
\[
\abs{G} f(e).
\]
I.e.:
\[
f(e) 
\hsx = \hsx 
\frac{1}{\abs{G}} \ \sum\limits_{\Pi \in \widehat{G}} \ d_\Pi  \hsx \tr(\widehat{f}(\Pi)),
\]
the so-called ``Plancherel theorem'' for \mG.
\end{x}
\vspace{0.1cm}

\begin{x}{\small\bf APPLICATION} \ 
Fix $\Pi_0 \in \widehat{G}$ and take $f = \ov{\chisubPiZero}$.
\\[-.2cm]

\qquad \textbullet \ $\Pi \neq \Pi_0$
\allowdisplaybreaks
\begin{align*}
\implies \ \tr(\Pi(\ov{\chisubPiZero})) \ 
&=\ 
\sum\limits_{x \in G} \hsx
\ov{\chisubPiZero(x)} \chisubPi(x) 
\\[15pt]
&=\ 
0.
\end{align*}
\vspace{0.1cm}

\qquad \textbullet \ $\Pi = \Pi_0$
\allowdisplaybreaks
\begin{align*}
\implies \ \tr(\Pi_0(\ov{\chisubPiZero})) \ 
&=\ 
\sum\limits_{x \in G} \
\ov{\chisubPiZero(x)} \hsx \chisubPiZero(x)\\[15pt]
&=\ 
\abs{G}.
\end{align*}
Therefore
\[
\sum\limits_{\Pi \in \widehat{G}} \ 
m(\Pi, L_{G / \Gamma}) \hsx \tr(\Pi(\ov{\chisubPiZero}))
\]
reduces to 
\[
\abs{G} \hsx m(\Pi_0,L_{G / \Gamma}).
\]
On the other hand,
\allowdisplaybreaks
\begin{align*}
\sum\limits_{x \in G} \
K_f(x,x)
&=\ 
\sum\limits_{x \in G} \
\frac{1}{\abs{\Gamma}} \
\sum\limits_{\gamma \in \Gamma} \
\ov{\chisubPiZero (x \gamma x^{-1})}
\\[15pt]
&=\ 
\sum\limits_{x \in G} \
\frac{1}{\abs{\Gamma}} \
\sum\limits_{\gamma \in \Gamma} \
\ov{\chisubPiZero (\gamma)}
\\[15pt]
&=\ 
\frac{\abs{G}}{\abs{\Gamma}} \
\sum\limits_{\gamma \in \Gamma} \
\ov{\chisubPiZero (\gamma)}
\\[15pt]
&=\ 
\abs{G} \ 
\langle 1_\Gamma, \restr{\chisubPiZero}{\Gamma} \ranglesubGamma
\\[15pt]
&=\ 
\abs{G} \ 
m(\theta, \restr{\Pi_0}{\Gamma}).
\end{align*}
So
\[
\abs{G} \hsx m(\Pi_0,L_{G / \Gamma})
\ = \
\abs{G} \hsx m(\theta, \restr{\Pi_0}{\Gamma}) 
\]
\qquad\qquad $\implies$
\[
m(\Pi_0,L_{G / \Gamma})
\hsx = \hsx
m(\theta, \restr{\Pi_0}{\Gamma}) .
\]
\vspace{0.2cm}

[Note: \ 
As above, $\theta$ is the trivial representation of $\Gamma$ on $\E = \Cx$.]
\end{x}
\vspace{0.1cm}

\begin{x}{\small\bf \un{N.B.}} \ 
Take $\Gamma = \{e\}$ $-$then
\[
m(\theta, \restr{\Pi_0}{\Gamma}) 
\hsx = \hsx
d_{\Pi_0},
\]
hence
\[
m(\Pi_0,L_{G / \Gamma})
\hsx = \hsx
d_{\Pi_0} \qquad \text{(cf. II, \S5, \#8)}.
\]
\end{x}
\vspace{0.1cm}


\chapter{
$\boldsymbol{\S}$\textbf{4}.\quad  THE GLOBAL TRACE FORMULA}
\setlength\parindent{2em}
\setcounter{theoremn}{0}
\renewcommand{\thepage}{A III \S4-\arabic{page}}

\qquad Let \mG be a finite group, $\Gamma \subset G$ a subgroup.
\\[-.2cm]

\begin{x}{\small\bf NOTATION} \ 
For any $\gamma \in \Gamma$, 
\[
\begin{cases}
\ G_\gamma \ = \ \text{centralizer of $\gamma$ in \mG}\\
\ \Gamma_\gamma \ = \ \text{centralizer of $\gamma$ in $\Gamma$}\\
\end{cases}
.
\]

Given an $f \in C(G)$, we have 
\[
\tr(L_{G/\Gamma} (f))
\ = \ 
\sum\limits_{x \hsy \in \hsy G} \hsx
K_f(x,x),
\]
where 
\[
K_f(x,x)
\ = \ 
\frac{1}{\abs{\Gamma}} \ 
\sum\limits_{\gamma \hsy \in \hsy \Gamma} \
f(x \hsy \gamma \hsy x^{-1})
\qquad \text{(cf. \S3, \#8)}.
\]
\\[-.2cm]

Enumerate the elements of $\CON(\Gamma)$, say
\[
\CON(\Gamma) 
\ = \ 
\{C_1, \ldots, C_n\}.
\]
For each $i$, fix a $\gamma_i \in C_i$ $(1 \leq i \leq n)$.
\end{x}
\vspace{0.1cm}

\begin{x}{\small\bf LEMMA} \ 
$\forall \ f \in C(G)$, 
\[
\sum\limits_{x \hsy \in \hsy G} \hsx K_f(x,x) 
\ = \ 
\sum\limits_{i=1}^n \ \frac{1}{\abs{\Gamma_{\gamma_i}}} \hsx \sO(f,\gamma_i).
\]
\vspace{0.2cm}

PROOF \ 
Write
\[
\Gamma
\ = \ 
\coprod\limits_{k} \ \gamma_{i,k} \hsx \Gamma_{\gamma_i}.
\]
Then $\forall \ x \in G$, 
\allowdisplaybreaks
\begin{align*}
K_f(x,x) \ 
&=\ \frac{1}{\Gamma} \hsx 
\sum\limits_{\gamma \hsy \in \hsy \Gamma} \hsx f(x \hsy \gamma \hsy x^{-1})
\\[15pt]
&=\ \frac{1}{\Gamma} \ 
\sum\limits_{i=1}^n \ 
\sum\limits_{\gamma \in C_i} \hsx f(x \hsy \gamma \hsy x^{-1})
\\[15pt]
&=\ \frac{1}{\Gamma} \ 
\sum\limits_{i=1}^n \ 
\sum\limits_{k} \ 
f(x \hsy \gamma_{i,k} \hsy \gamma_i  \hsy\gamma_{i,k}^{-1} \hsy x^{-1}).
\end{align*}
Therefore
\allowdisplaybreaks
\begin{align*}
\sum\limits_{x \hsy \in \hsy G} \hsx K_f(x,x) \ 
&=\ 
\frac{1}{\abs{\Gamma}} \ 
\sum\limits_{i=1}^n \ 
\sum\limits_{k} \
\sum\limits_{x \hsy \in \hsy G} \
f(x \hsy \gamma_{i,k} \hsy \gamma_i  \hsy\gamma_{i,k}^{-1} \hsy x^{-1})
\\[15pt]
&=\ 
\frac{1}{\abs{\Gamma}} \ 
\sum\limits_{i=1}^n \ 
[\Gamma:\Gamma_{\gamma_i}]
\sum\limits_{x \hsy \in \hsy G} \
f(x \hsy \gamma_i  \hsy x^{-1}).
\end{align*}
Write
\[
G 
\ = \ 
\coprod\limits_k x_{i,k} \ \Gamma_{\gamma_i}.
\]
Then
\allowdisplaybreaks
\begin{align*}
\sum\limits_{x \hsy \in \hsy G} \ K_f(x,x) \ 
&=\ 
\frac{1}{\abs{\Gamma}} \
\sum\limits_{i=1}^n \ 
[\Gamma:\Gamma_{\gamma_i}] \
\sum\limits_k \
\sum\limits_{\eta_i \in \Gamma_{\gamma_i}} \
f(x_{i,k} \hsy \eta_i \hsy \gamma_i \hsy \eta_i^{-1}  \hsy x_{i,k}^{-1})
\\[15pt]
&=\ 
\frac{1}{\abs{\Gamma}} \
\sum\limits_{i=1}^n \ 
[\Gamma:\Gamma_{\gamma_i}][\Gamma_{\gamma_i}] \
\sum\limits_k \
f(x_{i,k} \hsy \gamma_i \hsy x_{i,k}^{-1})
\\[15pt]
&=\ 
\sum\limits_{i=1}^n \ 
\frac{[\Gamma:\Gamma_{\gamma_i}][\Gamma_{\gamma_i}]}{\abs{\Gamma}} \
\sum\limits_k \
f(x_{i,k} \hsy \gamma_i \hsy x_{i,k}^{-1})
\\[15pt]
&=\ 
\sum\limits_{i=1}^n \
\sum\limits_k \ f(x_{i,k} \gamma_i x_{i,k}^{-1}).
\end{align*}
Write
\[
\begin{cases}
\ds\ G \ = \ \coprod\limits_\ell \ y_{i,\ell} \hsy G_{\gamma_i}
\\[15pt]
\ds\ G_{\gamma_i} \ = \ \coprod\limits_m \ z_{i,m} \hsy \Gamma_{\gamma_i}
\end{cases}
\]
\qquad\qquad $\implies$ 
\[
G \ = \ \coprod\limits_\ell \ \coprod\limits_m \ y_{i,\ell}  \hsy z_{i,m} \hsy \Gamma_{\gamma_i}.
\]
Then
\allowdisplaybreaks
\begin{align*}
\sum\limits_{x \hsy \in \hsy G} \ K_f(x,x) \ 
&=\ 
\sum\limits_{i=1}^n \ 
\sum\limits_\ell \
\sum\limits_m \
f(y_{i,\ell} \hsy z_{i,m} \hsy \gamma_i \hsy z_{i,m}^{-1} \hsy y_{i,\ell}^{-1})
\\[15pt]
&=\ 
\sum\limits_{i=1}^n \ 
[G_{\gamma_i}:\Gamma_{\gamma_i}] 
\sum\limits_\ell \
f(y_{i,\ell} \hsy \gamma_i \hsy y_{i,\ell}^{-1})
\\[15pt]
&=\ 
\sum\limits_{i=1}^n \ 
\abs{\frac{G_{\gamma_i}}{\Gamma_{\gamma_i}}} \
\sum\limits_\ell \ 
f(y_{i,\ell} \hsy \gamma_i \hsy y_{i,\ell}^{-1})
\\[15pt]
&=\ 
\sum\limits_{i=1}^n \ 
\frac{1}{\abs{\Gamma_{\gamma_i}}} \hsx \sO(f,\gamma_i) .
\end{align*}
\end{x}
\vspace{0.1cm}

\begin{x}{\small\bf \un{N.B.}} \ 
$\forall \ \gamma$, 
\[
\sO(f,\gamma) 
\ = \
\abs{G_\gamma} \ \sum\limits_{x \in G / G_\gamma} \ f(x \hsy \gamma \hsy x^{-1}), 
\]
the sum on the right being taken over a set of representatives for the left cosets of $G_\gamma$ in \mG.
\end{x}
\vspace{0.1cm}

\begin{x}{\small\bf EXAMPLE} \ 
Take $\Gamma = G$ $-$then $\forall \ f \in C(G)$, 
\allowdisplaybreaks
\begin{align*}
\sum\limits_{i=1}^n \ 
\frac{1}{\abs{G_{\gamma_i}}} \
\sO(f,\gamma_i)
&=\ 
\sum\limits_{i=1}^n \ 
\frac{1}{\abs{G_{\gamma_i}}} \ 
\sum\limits_{x \hsy \in \hsy G} \ 
f(x \hsy \gamma_i \hsy x^{-1})
\\[15pt]
&=\ 
\sum\limits_{i=1}^n \ 
\frac{1}{\abs{G_{\gamma_i}}}\hsx \abs{G_{\gamma_i}}
\sum\limits_{y_i \in G/G_{\gamma_i}} \ 
f(y_i \hsy \gamma_i \hsy y_i^{-1})
\\[15pt]
&=\ 
\sum\limits_{i=1}^n \ 
\sum\limits_{y_i \hsy \in \hsy G/G_{\gamma_i}} \ 
f(y_i \hsy \gamma_i \hsy y_i^{-1})
\\[15pt]
&=\ 
\sum\limits_{i=1}^n \ 
\sum\limits_{y \hsy \in \hsy C_i} \ 
f(y)
\\[15pt]
&=\ 
\sum\limits_{x \hsy \in \hsy G} \ 
f(x) 
\qquad \text{(cf. \S3, \#9)}.
\end{align*}
\end{x}
\vspace{0.1cm}

\begin{x}{\small\bf EXAMPLE} \ 
Suppose that $f \in CL(G)$ $-$then 
\allowdisplaybreaks
\begin{align*}
\sum\limits_{x \hsy \in \hsy G} \ K_f(x,x) 
&=\ 
\sum\limits_{x \hsy \in \hsy G} \ 
\frac{1}{\abs{\Gamma}} \ 
\sum\limits_{\gamma \hsy \in \hsy \Gamma} \ 
f(x \hsy \gamma \hsy x^{-1})
\\[15pt]
&=\ 
\sum\limits_{x \hsy \in \hsy G} \ 
\frac{1}{\abs{\Gamma}} \ 
\sum\limits_{\gamma \hsy \in \hsy \Gamma} \ 
f(\gamma)\\[15pt]
&=\ 
\frac{\abs{G}}{\abs{\Gamma}} \ 
\sum\limits_{\gamma \hsy \in \hsy \Gamma} \ 
f(\gamma).
\end{align*}
In the other direction,
\allowdisplaybreaks
\begin{align*}
\sum\limits_{i=1}^n \ 
\frac{1}{\abs{\Gamma_{\gamma_i}}} \
\sO(f,\gamma_i) 
&=\ 
\sum\limits_{i=1}^n \hsx 
\frac{1}{\abs{\Gamma_{\gamma_i}}} \
\sum\limits_{x \hsy \in \hsy G} \hsx 
f(x \hsy \gamma_i \hsy x^{-1}) \\[15pt]
&=\ 
\abs{G} \ 
\sum\limits_{i=1}^n \ 
\frac{f(\gamma_i)}{\abs{\Gamma_{\gamma_i}}}.
\end{align*}
Therefore
\[
\frac{1}{\abs{\Gamma}} \ \sum\limits_{\gamma \hsy \in \hsy \Gamma} \ f(\gamma)
\ = \ 
\sum\limits_{i=1}^n \ 
\frac{f(\gamma_i)}{\abs{\Gamma_{\gamma_i}}}.
\]
\end{x}
\vspace{0.1cm}

\begin{x}{\small\bf DEFINITION} \ 
Given $f \in C(G)$, the 
\un{global trace formula}
\index{global trace formula}  
is the relation
\[
\sum\limits_{\Pi \in \widehat{G}} \hsx m(\Pi, L_{G /\Gamma}) \hsy \tr(\widehat{f}(\Pi)) 
\ = \ 
\sum\limits_{i=1}^n \hsx \frac{1}{\abs{\Gamma_{\gamma_i}}} \hsx \sO(f,\gamma_i) 
\qquad \text{(cf. \S3, \#12)}.
\]

\end{x}
\vspace{0.1cm}

\begin{x}{\small\bf EXAMPLE} (POISSON SUMMATION) \ 
Take \mG abelian and identify $\widehat{G}$ with the character group of 
$G: \Pi \longleftrightarrow \chi$, hence
\[
\widehat{f}(\chi) 
\ = \ 
\sum\limits_{x \hsy \in \hsy G} \hsx 
f(x) \chi(x).
\]
Consider now the sum
\[
\sum\limits_{\chi \in \widehat{G}} \hsx 
m(\chi,L_{G / \Gamma}) \hsx \widehat{f}(\chi).
\]
Let $\Gamma^\perp = \{\chi \in \widehat{G}: \chi(\gamma) = 1 \ \forall \ \gamma \in \Gamma\}$ $-$then
\[
\begin{cases}
\ \chi \in \Gamma^\perp \implies m(\chi,L_{G / \Gamma}) = 1\\[11pt]
\ \chi \notin \Gamma^\perp \implies m(\chi,L_{G / \Gamma}) = 0
\end{cases}
\qquad \text{(cf. \S3, \#14)}.
\]
Therefore matters on the ``spectral side'' reduce to 
\[
\sum\limits_{\chi \hsy \in \hsy \Gamma^\perp} \hsx \widehat{f}(\chi) .
\]
And on the ``geometric side'', 
\allowdisplaybreaks
\begin{align*}
\sum\limits_{i=1}^n \ 
\frac{1}{\abs{\Gamma_{\gamma_i}}} \
\sO(f,\gamma_i) 
&=\ 
\frac{1}{\abs{\Gamma}} \
\sum\limits_{\gamma \hsy \in \hsy \Gamma} \ 
\sO(f,\gamma)
\\[15pt]
&=\ 
\frac{1}{\abs{\Gamma}} \
\sum\limits_{\gamma \hsy \in \hsy \Gamma} \ 
\sum\limits_{x \hsy \in \hsy G} \hsx 
f(x \gamma x^{-1})
\\[15pt]
&=\ 
\frac{1}{\abs{\Gamma}} \
\sum\limits_{\gamma \hsy \in \hsy \Gamma} \ 
\abs{G} f(\gamma)
\\[15pt]
&=\ 
\frac{\abs{G}}{\abs{\Gamma}} \
\sum\limits_{\gamma \hsy \in \hsy \Gamma} \ 
f(\gamma).
\end{align*}
Therefore
\[
\frac{1}{\abs{G}} \
\sum\limits_{\chi \hsy \in \hsy \Gamma^\perp} \ 
\widehat{f}(\chi)
\ = \ 
\frac{1}{\abs{\Gamma}} \
\sum\limits_{\gamma \hsy \in \hsy \Gamma} \ 
f(\gamma).
\]
\vspace{0.1cm}

Each element $\zeta$ in the center $Z(\Gamma)$ of $\Gamma$ determines a one-element conjugacy class $\{\zeta\}$.
\end{x}
\vspace{0.1cm}

\begin{x}{\small\bf DEFINITION} \ 
The 
\un{central contribution}
\index{central contribution} 
to the global trace formula is the subsum
\[
\sum\limits_{\zeta \in Z(\Gamma)} \ \frac{1}{\abs{\Gamma_\zeta}} \ \sO(f,\zeta)
\]
of 
\[
\sum\limits_{i=1}^n \ \frac{1}{\abs{\Gamma_{\gamma_i}}} \ \sO(f,\gamma_i).
\]
Accordingly, 
\allowdisplaybreaks
\begin{align*}
\sum\limits_{\zeta \in Z(\Gamma)} \hsx \frac{1}{\abs{\Gamma_\zeta}} \ \sO(f,\zeta) \ 
&=\ 
\frac{1}{\abs{\Gamma}} \
\sum\limits_{\zeta \in Z(\Gamma)} \hsx \sO(f,\zeta)
\\[15pt]
&=\ 
\frac{1}{\abs{\Gamma}} \
\sum\limits_{\zeta \in Z(\Gamma)} \
\sum\limits_{x \hsy \in \hsy G} \hsx f(x \zeta x^{-1})\\[15pt]
&=\ 
\frac{1}{\abs{\Gamma}} \
\sum\limits_{\zeta \in Z(\Gamma)} \ 
\abs{G} f(\zeta)\\[15pt]
&=\ 
\frac{\abs{G}}{\abs{\Gamma}} \
\sum\limits_{\zeta \in Z(\Gamma)} \ f(\zeta).
\end{align*}
\end{x}
\vspace{0.1cm}


%% file: _B.tex
\chapter{
$\boldsymbol{\S}$\textbf{1}.\quad  UNITARY REPRESENTATI0NS}
\setlength\parindent{2em}
\setcounter{theoremn}{0}
\renewcommand{\thepage}{B I \S1-\arabic{page}}

\qquad Let \mG be a compact group.
\\

\begin{x}{\small\bf NOTATION} \ 
$\td_G$ is normalized Haar measure on \mG: 
\[
\int\limits_G \hsx 
1 \hsx  \td_G(x) 
\ = \ 
1.
\]
\end{x}
\vspace{0.3cm}

\begin{x}{\small\bf LEMMA} \ 
\[
L^1(G) 
\ \supset \
L^2(G) 
\ \supset \
C(G) 
\]
and 

\qquad \textbullet \ $\forall \ f \in L^2(G)$, $\norm{f}_2 \geq \norm{f}_1$.
\\[-.2cm]

\qquad \textbullet \ $\forall \ f \in C(G)$, $\norm{f}_2 \leq \norm{f}_\infty$.
\\[-.2cm]
\end{x}

\begin{x}{\small\bf \un{N.B.}} \ 
The convolution operator 
\[
*: L^2(G) \times  L^2(G) \ra C(G)
\]
is given by 
\begin{align*}
(f * g) (x) \ 
&=\ 
\int\limits_G \hsx 
f(x y^{-1}) g(y) \td_G(y)
\\[15pt]
&=\ 
\int\limits_G \hsx 
f(y) g(y^{-1} x) \td_G(y).
\end{align*}
\\[-.5cm]
\end{x}

\begin{x}{\small\bf DEFINITION} \ 
A 
\un{unitary representation} 
\index{unitary representation} 
of \mG on a Hilbert space $\sH$ is a homomorphism $\pi:G \ra \UN (\sH)$ from \mG to the unitary group 
$\UN(\sH)$ of $\sH$ such that $\forall \ a \in \sH$, the map 
\[
x \ra \pi(x) a
\]
of \mG into $\sH$ is continuous.
\end{x}
\vspace{0.3cm}


\begin{x}{\small\bf DEFINITION} \ 
\\

\qquad \textbullet \  \ The left translation representation of \mG on $L^2(G)$ is the prescription
\[
L(x) f(y) 
\ = \ 
f(x^{-1} y).
\]

\qquad \textbullet \ \ The right translation representation of \mG on $L^2(G)$ is the prescription
\[
R(x) f(y) 
\ = \ 
f(y x).
\]

[Note: \ 
Both \mL and \mR are unitary.]
\\[-.2cm]
\end{x}

\begin{x}{\small\bf \un{N.B.}} \ 
There is also a unitary representation $\pisubLR$ of $G \times G$ on $L^2(G)$, namely
\[
(\pisubLR (x_1, x_2) f)(x) 
\ = \ 
f(x_1^{-1} x x_2).
\]
\end{x}
\vspace{0.3cm}

\begin{x}{\small\bf DEFINITION} \ 
A unitary representation $\pi$ of \mG on a Hilbert space $\sH \neq \{0\}$ is 
\un{irreducible}
\index{unitary representation // irreducible} 
if the only closed subspaces of $\sH$ which are invariant under $\pi$ are $\{0\}$ and $\sH$.
\end{x}
\vspace{0.3cm}

\begin{x}{\small\bf THEOREM} \ 
Let $\pi$ be a unitary representation of \mG $-$then $\pi$ is the Hilbert space direct sum of finite dimensional irreducible unitary representations.
\end{x}
\vspace{0.3cm}

\begin{x}{\small\bf APPLICATION} \ 
Every irreducible unitary representation of \mG is finite dimensional.
\end{x}
\vspace{0.3cm}

\begin{x}{\small\bf NOTATION} \ 
$\widehat{G}$ is the set of unitary equivalence classes of irreducible unitary representations of \mG.
\vspace{0.2cm}

[Note: \ 
Generically, $\Pi \in \widehat{G}$ with representation space $V(\Pi)$ and $\td_\Pi = \dim V(\Pi)$ is its dimension.]
\end{x}
\vspace{0.3cm}


\begin{x}{\small\bf \un{N.B.}} \ \ 
Let $\pi$ be a unitary representation of \mG $-$then there exist cardinal numbers $n_\Pi \ (\Pi \in \widehat{G})$ such that 
\[
\pi 
\ = \ 
\widehat{\bigoplus\limits_{\Pi \in \widehat{G}}} \ 
n_\Pi \hsx \Pi.
\]
\end{x}
\vspace{0.3cm}

\begin{x}{\small\bf EXAMPLE} \
Take $\pi = L$ $-$then 
\[
L 
\ = \ 
\widehat{\bigoplus\limits_{\Pi \in \widehat{G}}} \ 
\td_\Pi \hsx \Pi.
\]
\vspace{0.2cm}

[Note: \ 
There is also an analog of \mA, III, \S2, \#16.]
\end{x}
\vspace{0.3cm}

\begin{x}{\small\bf THEOREM} \ 
$\forall \ x \in G$ $(x \neq e)$, $\exists$ an irreducible unitary representation $\Pi$ such that $\Pi(x) \neq \id$ (Gelfand-Raikov).
\end{x}
\vspace{0.3cm}

\begin{x}{\small\bf APPLICATION} \ 
\[
\bigcap\limits_{\Pi \in \widehat{G}}
\Ker \Pi 
\ = \ 
\
\{e\}.
\]
\end{x}
\vspace{0.3cm}

\begin{x}{\small\bf LEMMA} \ 
Given $\Pi \in \widehat{G}$, suppose that $A \in \Hom(V(\Pi), V(\Pi))$ has the property that $\forall \ x \in G$, 
\[
A \hsx \Pi (x) 
\ = \ 
\Pi (x) \hsx A.
\]
Then \mA is a scalar multiple of the identity (Schur), call it $\lambda_A$.
\end{x}
\vspace{0.3cm}


\chapter{
$\boldsymbol{\S}$\textbf{2}.\quad  EXPANSION THEORY}
\setlength\parindent{2em}
\setcounter{theoremn}{0}
\renewcommand{\thepage}{B I \S2-\arabic{page}}

\qquad Let \mG be a compact group.
\\[-.2cm]

\begin{x}{\small\bf DEFINITION} \ 
Let $\pi$ be a finite dimensional unitary representation of \mG $-$then its \un{character} is the function 
\[
\chisubPi:G \ra \Cx
\]
given by the prescription
\[
\chisubPi(x) 
\ = \ 
\tr(\pi(x)) \qquad (x \in G).
\]
\end{x}

\vspace{0.3cm}

\begin{x}{\small\bf DEFINITION} \ 
The character of an irreducible unitary representation is called an \un{irreducible character}.
\end{x}
\vspace{0.3cm}

\begin{x}{\small\bf LEMMA} \ 
Let $\Pi_1$, $\Pi_2 \in \widehat{G}$ and suppose that $\Pi_1 \neq \Pi_2$ $-$then
\[
\langle \chisubPiOne, \chisubPiTwo\rangle 
\ = \ 0.
\]
\end{x}
\vspace{0.3cm}

\begin{x}{\small\bf LEMMA} \ 
Let $\Pi \ \in \widehat{G}$ $-$then
\[
\langle \chisubPi, \chisubPi\rangle 
\ = \ 1.
\]
\end{x}
\vspace{0.3cm}

\begin{x}{\small\bf DEFINITION} \ 
A continuous complex valued function $\phi$ on \mG is of \un{positive type} if for all 
$x_1, \ldots, x_n \in G$ and $\lambda_1, \ldots, \lambda_n \in \Cx$, 
\[
\sum\limits_{i, j = 1}^n \ \lambda_i \ov{\lambda_j} \phi(x_i^{-1} x_j) 
\ \geq \ 
0.
\]
\end{x}
\vspace{0.3cm}

\begin{x}{\small\bf \un{N.B.}} \ 
The sum of two functions of positive type is of positive type and a positive scalar multiple of a functon of 
positive type is of positive type.
\end{x}
\vspace{0.3cm}

\begin{x}{\small\bf LEMMA} \ 
If $\pi:G \ra \UN(\sH)$ is a unitary representation and if $a \in \sH$, then 
\[
\phi(x) 
\ = \ 
\langle \pi(x) a, a \rangle \qquad (x \in G)
\]
is of positive type.
\\[-.25cm]

[Note: \ 
\[
\norm{\phi}_\infty 
\ = \ 
\langle a, a \rangle.]
\]
\end{x}
\vspace{0.3cm}

\begin{x}{\small\bf EXAMPLE} \ 
$\forall \ \Pi \in \widehat{G}$, $\chisubPi$ is of positive type.
\\[-.25cm]

[Fix an orthonormal basis $v_1, \ldots, v_n$ in $V(\Pi)$ $-$then 
\[
\chisubPi(x) 
\ = \ 
\langle \Pi(x) v_1, v_1 \rangle + \cdots + \langle \Pi(x) v_n, v_n \rangle,
\]
from which the assertion.]
\end{x}
\vspace{0.3cm}

\begin{x}{\small\bf NOTATION} \ 
Given $\Pi \in \widehat{G}$ and $f \in L^2(G)$, put
\[
\Pi(f) 
\ = \ 
\int_G \hsx f(x) \Pi(x) \hsx \td_G(x).
\]
\end{x}
\vspace{0.3cm}

\begin{x}{\small\bf LEMMA} \ 
$\forall \ f_1, \ f_2 \in L^2(G)$, 
\[
\Pi(f_1 * f_2) 
\ = \ 
\Pi(f_1) \circ \Pi(f_2).
\]

\end{x}
\vspace{0.3cm}

\begin{x}{\small\bf NOTATION} \ 
Given $f \in L^2(G)$, define $f^* \in L^2(G)$ by
\[
f^*(x) 
\ = \ 
\ov{f(x^{-1})} \qquad (=\ov{\widecheck{f}(x)}).
\]
\end{x}
\vspace{0.3cm}

\begin{x}{\small\bf LEMMA} \ 
$\forall \ f \in L^2(G)$, $\forall \ v_1, \ v_2 \in V(\Pi)$, 
\[
\langle \Pi(f)v_1, v_2 \rangle 
\ = \ 
\langle v_1, \Pi(f^*) v_2 \rangle,
\]
i.e., 
\[
\Pi(f)^* 
\ = \ 
\Pi(f^*).
\]

PROOF \ 

\allowdisplaybreaks
\begin{align*}
\langle \Pi(f)v_1, v_2\rangle \ 
&=\ 
\int_G \hsx f(x) \langle \Pi(x)v_1, v_2\rangle \hsx \td_G(x)
\\[11pt]
&=\
\int_G \hsx \langle v_1,\ov{f(x)} \hsx \Pi(x^{-1}) v_2\rangle \hsx \td_G(x) 
\\[11pt]
&=\ 
\int_G \hsx \langle v_1,\ov{f(x^{-1})} \hsx \Pi(x) v_2\rangle \hsx \td_G(x) 
\\[11pt]
&=\ 
\int_G \hsx \langle v_1,f^*(x)  \hsx\Pi(x) v_2\rangle \hsx \td_G(x) 
\\[11pt]
&=\ 
\langle v_1, \int_G \hsx f^*(x) \hsx \Pi(x)v_2 \hsx \td_G(x)\rangle
\\[11pt]
&=\
\langle v_1, \Pi(f^*)v_2\rangle.
\end{align*}

\end{x}
\vspace{0.3cm}

\begin{x}{\small\bf THEOREM} \ 
Let $f \in L^2(G)$ $-$then
\[
\int_G \hsx \abs{f(x)}^2 \hsx \td_G(x) 
\ = \ 
\sum\limits_{\Pi \in \widehat{G}} \hsx \td_\Pi \hsx \tr(\Pi(f)\Pi(f)^*).
\]
\end{x}
\vspace{0.3cm}

\begin{x}{\small\bf THEOREM} \ 
Let $f \in L^2(G)$ $-$then
\[
f \ = \ 
\sum\limits_{\Pi \in \widehat{G}} \ \td_\Pi(f * \chisubPi),
\]
the series converging in $L^2(G)$.
\end{x}
\vspace{0.3cm}

\begin{x}{\small\bf THEOREM} \ 
Let
\[
f \in \spanx_\Cx (L^2(G) * L^2(G)) \subset C(G).
\]
Then
\[
f(e) 
\ = \ 
\sum\limits_{\Pi \in \widehat{G}} \hsx \td_\Pi \hsx \tr(\Pi(f)).
\]


PROOF \ 
Put $f = f_1 * f_2$:

\allowdisplaybreaks
\begin{align*}
f(e) \ 
&=\ 
\int_G \hsx f_1(x^{-1}) f_n(x) \hsx \td_G(x)
\\[11pt]
&=\ 
\int_G \hsx \ov{\ov{f_1(x^{-1})}} \hsx f_2(x) \ \hsx \td_G(x)
\\[11pt]
&=\ 
\int_G \hsx f_2(x) \ov{\ov{f_1(x^{-1})}} \hsx \td_G(x) 
\\[11pt]
&=\ 
\int_G \hsx f_2(x) \ov{f_1^*(x)} \hsx \td_G(x)
\\[11pt]
&=\ 
\langle f_2, f_1^* \rangle
\\[11pt]
&=\ 
\sum\limits_{\Pi \in \widehat{G}} \ \td_\Pi \hsx \tr(\Pi(f_1 * f_2))
\\[11pt]
&=\ 
\sum\limits_{\Pi \in \widehat{G}} \ \td_\Pi \hsx \tr(\Pi(f)).
\end{align*}

[Note: \ 
This is the so-called ``Plancherel theorem'' for \mG (cf. A, III, $\S3$, $\#13$).
\end{x}
\vspace{0.3cm}

\begin{x}{\small\bf \un{N.B.}} \ 
The foregoing may fail if $f$ is only assumed to be continuous 
(e.g., take $G = \bS^1 \ldots$).
\end{x}
\vspace{0.3cm}

\begin{x}{\small\bf DEFINITION} \ 
A function $f \in L^2(G)$ is said to be an \un{$L^2$ class function} if
\[
f(x) 
\ = \ 
f(y x y^{-1})
\]
for almost all $x$ and all $y$.
\end{x}
\vspace{0.3cm}

\begin{x}{\small\bf \un{N.B.}} \ 
$\forall \ \Pi \in \widehat{G}$, $\chisubPi$ is an $L^2$ class function.
\end{x}
\vspace{0.3cm}

\begin{x}{\small\bf THEOREM} \ 
Suppose that $f \in L^2(G)$ is an $L^2$ class function $-$then
\[
f 
\ = \ 
\sum\limits_{\Pi \in \widehat{G}} \ \langle f,\chisubPi\rangle \chisubPi,
\]
the series converging in $L^2(G)$, and 
\[
\norm{f}^2 
\ = \ 
\sum\limits_{\Pi \in \widehat{G}} \hsx \abs{ \langle f,\chisubPi\rangle }^2.
\]
\end{x}
\vspace{0.3cm}

\begin{x}{\small\bf SCHOLIUM} \ 
The $\{\chisubPi: \Pi \in \widehat{G}\}$ constitute an orthonormal basis for the set of $L^2$ class functions.
\end{x}
\vspace{0.3cm}

\begin{x}{\small\bf NOTATION} \ 
Write $C(G)_\fin(L)$ for the set of \mG-finite functions in $C(G)$ per \mL:
\[
f \in C(G)_\fin(L) 
\ \Leftrightarrow \ 
\dim\{L(x) f: x \in G\}_{\ell\text{in}} 
\ < \ 
\infty.
\]
\end{x}
\vspace{0.3cm}

\begin{x}{\small\bf NOTATION} \ 
Write $C(G)_\fin(R)$ for the set of \mG-finite functions in $C(G)$ per \mR:
\[
f \in C(G)_\fin(R) 
\ \Leftrightarrow \ 
\dim\{R(x) f: x \in G\}_{\ell\text{in}}
\ < \ 
\infty.
\]
\end{x}
\vspace{0.3cm}

\begin{x}{\small\bf LEMMA} \ 
\[
C(G)_\fin(L) 
\ = \ C(G)_\fin(R).
\]
\end{x}
\vspace{0.3cm}

\begin{x}{\small\bf NOTATION} \ 
Write $C(G)_\fin$ unambiguously for the \mG-finite functions per either action.
\end{x}
\vspace{0.3cm}

Recalling $\S1$, $\#6$, $\pisubLR$ operates on $C(G)_\fin$ and it turns out that 
\[
C(G)_\fin
\ \approx \ 
\bigoplus\limits_{\Pi \in \widehat{G}} \hsx V(\Pi^*) \otimes V(\Pi).
\]
Here the identification sends an element
\[
v^* \otimes v \in V(\Pi^*) \otimes V(\Pi)
\]
to 
\[
f_{v^* \otimes v} \in C(G)_\fin,
\]
where
\[
f_{v^* \otimes v}(x) 
\ = \ 
v^*(\Pi(x^{-1})v).
\]

[Note: \ 
\[
L^2(G) 
\ \approx \ 
\widehat{\bigoplus\limits_{\Pi \in \widehat{G}}} \ V(\Pi^*) \otimes V(\Pi).]
\]
\\[-.2cm]

\begin{x}{\small\bf THEOREM} \ 
 $C(G)_\fin$ is dense in $C(G)$.
\end{x}
\vspace{0.3cm}

\begin{x}{\small\bf THEOREM} \ 
 $C(G)_\fin$ is dense in $L^2(G)$.
\end{x}
\vspace{0.3cm}

\begin{x}{\small\bf DEFINITION} \ 
A function $f \in C(G)$ is said to be a \un{continuous class function} if $f(x) = f(y x y^{-1})$ for all 
$x$, $y \ \in G$ (written $f \in CL(G)$).
\end{x}
\vspace{0.3cm}

\begin{x}{\small\bf EXAMPLE} \ 
$\forall \ \Pi \in \widehat{G}$, $\chisubPi$ is a continuous class function: $\chisubPi \in CL(G)$.
\end{x}
\vspace{0.3cm}

\begin{x}{\small\bf THEOREM} \ 
The span of the $\chisubPi$ $(\Pi \in \widehat{G})$ equals the set of continuous class functions in $C(G)_\fin$.
\end{x}
\vspace{0.3cm}

\begin{x}{\small\bf THEOREM} \ 
The span of the $\chisubPi$ $(\Pi \in \widehat{G})$ is dense in the set of continuous class functions.
\end{x}
\vspace{0.3cm}


\chapter{
$\boldsymbol{\S}$\textbf{3}.\quad  STRUCTURE THEORY}
\setlength\parindent{2em}
\setcounter{theoremn}{0}
\renewcommand{\thepage}{B I \S3-\arabic{page}}

\qquad Let \mG be a compact group.
\\[-.2cm]

\begin{x}{\small\bf NOTATION} \ 
$G^0 \subset G$ is the connected component of the identity of \mG.
\end{x}
\vspace{0.3cm}

\begin{x}{\small\bf LEMMA} \ 
$G^0$ is a closed normal subgroup of \mG.
\end{x}
\vspace{0.3cm}

\begin{x}{\small\bf LEMMA} \ 
The quotient $G  / G^0$ is compact and totally disconnected.
\end{x}
\vspace{0.3cm}

\begin{x}{\small\bf DEFINITION} \ 
A topological group possessing a neighborhood of the identity which does not contain a nontrivial subgroup is said to be a group with 
\un{no small subgroups}.
\index{no small subgroups}
\end{x}

\begin{x}{\small\bf RAPPEL} \ 
A Lie group has no small subgroups.
\end{x}
\vspace{0.3cm}

\begin{x}{\small\bf THEOREM} \ 
The following conditions on a compact group \mG are equivalent.

\qquad \textbullet \ \ 
\mG is a Lie group.
\vspace{0.2cm}

\qquad \textbullet \ \ 
\mG has no small subgroups.
\vspace{0.2cm}

\qquad \textbullet \ \ 
\mG has a faithful finite dimensional representation.
\end{x}
\vspace{0.3cm}

\begin{x}{\small\bf REMARK} \ 
Every compact group is the projective limit of compact Lie groups.
\end{x}
\vspace{0.3cm}

Let \mG be a compact Lie group.
\vspace{0.3cm}

\begin{x}{\small\bf \un{N.B.}} \ 
Every finite group (discrete topology) is a compact Lie group.
\end{x}
\vspace{0.3cm}

\begin{x}{\small\bf EXAMPLE} \ 
The product $\ds\prod\limits_{n=1}^\infty \hsx \SU(n)$ is a compact group but it is not a Lie group.
\end{x}
\vspace{0.3cm}

\begin{x}{\small\bf EXAMPLE} \ 
The $p$-adic integers 
\[
\Z_p 
\ = \ 
\lim\limits_{\substack{\la\\n \geq 1}} \hsx
\bigl(\Z / p^n \Z \bigr)
\]
are a compact group but they are not a Lie group.
\end{x}
\vspace{0.3cm}

\begin{x}{\small\bf DEFINITION} \ 
A 
\un{torus}
\index{torus} 
is a compact Lie group which is isomorphic to $\R^n / \Z^n \hsx \approx \hsx (\R / \Z)^n$ for some $n \geq 0$.
\vspace{0.2cm}

[Note: \
The nonnegative integer $n$ is called the 
\un{rank}
\index{rank (of a torus)} 
of the torus.]
\end{x}
\vspace{0.3cm}

\begin{x}{\small\bf THEOREM} \ 
Every compact abelian Lie group is isomorphic to the product of a torus and a finite abelian group.
\end{x}
\vspace{0.3cm}

\begin{x}{\small\bf DEFINITION} \ 
A compact Lie group is 
\un{topologically cyclic}
\index{topologically cyclic} 
if it contains an element whose powers are dense.
\end{x}
\vspace{0.3cm}

\begin{x}{\small\bf LEMMA} \ 
Every torus \Torus is topologically cyclic.
\vspace{0.2cm}

[Note: \ 
There are infinitely many topologically cyclic elements in \Torus and their totality has full measure in any Haar measure on \Torus.]
\end{x}
\vspace{0.3cm}

\begin{x}{\small\bf THEOREM} \ 
A compact Lie group is topologically cyclic iff it is isomorphic to the product of a torus and a finite cyclic group.
\end{x}
\vspace{0.3cm}

Let \mG be a compact Lie group, $\fg$ its Lie algebra.
\index{$\fg$}
\vspace{0.3cm}

\begin{x}{\small\bf LEMMA} \ 
$G^0$ is an open normal subgroup of \mG.
\end{x}
\vspace{0.3cm}

Therefore the compact quotient $G / G^0$ is discrete, hence is a finite group, the 
\un{group of} \un{components}
\index{group of components} 
of \mG.
\vspace{0.3cm}

\begin{x}{\small\bf NOTATION} \ 
$Z(G)$ is the center of \mG and $Z(G)^0 \subset Z(G)$ is the connected component of the identity element of $Z(G)$. 
\index{$Z(G)^0$} 
\end{x}
\vspace{0.3cm}

\begin{x}{\small\bf \un{N.B.}} \ 
In general, $Z(G)$ is not connected (consider $\SU(3))$.
\end{x}
\vspace{0.3cm}

\begin{x}{\small\bf THEOREM} \ 
Assume that \mG is connected $-$then $Z(G)^0$ is a compact abelian Lie subgroup of \mG and its Lie algebra is the center 
$Z(\fg)$
of $\fg$, i.e., the ideal
\[
\{X \in \fg : [X,Y] 
\ = \ 
0 \ \forall \ Y \in \fg\}.
\]
\end{x}
\vspace{0.3cm}

\begin{x}{\small\bf DEFINITION} \ 
\vspace{0.2cm}

\qquad \textbullet \ \
A Lie algebra is 
\un{simple}
\index{simple} 
if it is noncommutative and has no proper nontrivial ideals.
\vspace{0.2cm}

\qquad \textbullet \ \
A Lie algebra is 
\un{semisimple}
\index{semisimple} 
if it is noncommutative and has no proper nontrivial commutative ideals.
\vspace{0.2cm}

\qquad \textbullet \ \ 
A Lie algebra is 
\un{reductive}
\index{reductive} 
if it is the direct sum of an abelian Lie algebra and a semisimple Lie algebra.
\\[-.2cm]
\end{x}

\begin{x}{\small\bf \un{N.B.}} \ 
A Lie group is simple, semisimple, or reductive if this is the case of its Lie algebra.
\\[-.2cm]
\end{x}

\begin{x}{\small\bf LEMMA} \ 
A semisimple Lie algebra has a trivial center (it being a commutative ideal).
\\[-.2cm]
\end{x}

\begin{x}{\small\bf LEMMA} \ 
A semisimple Lie algebra can be decomposed as a finite direct sum of simple ideals.
\\[-.2cm]
\end{x}

\begin{x}{\small\bf DEFINITION} \ 
If \mG and \mH are Lie groups and if \mH is a subgroup of \mG, then \mH is a 
\un{Lie subgroup}
\index{Lie subgroup} 
of \mG if the arrow $H \ra G$ of the inclusion is continuous.
\\[-.25cm]

[Note: \ 
If \mG is a Lie group and if \mH is a closed subgroup of \mG, then \mH is a Lie group.]
\end{x}
\vspace{0.3cm}

\begin{x}{\small\bf \un{N.B.}} \ 
A Lie subgroup of a compact Lie group needn't be compact nor carry the relative topology.
\end{x}
\vspace{0.3cm}

\begin{x}{\small\bf THEOREM} \ 
Let \mG be a compact Lie group and let \mH be a semisimple connected Lie subgroup of \mG $-$then as a subset of \mG, \mH is closed, and as a Lie subgroup of \mG, \mH carries the relative topology.
\end{x}
\vspace{0.3cm}

\begin{x}{\small\bf NOTATION} \ 
\vspace{0.2cm}

\qquad \textbullet \ \ $\fz (\fg)$ is the center of $\fg$.
\vspace{0.2cm}

\qquad \textbullet \ \ $\fg_\ssx$ is the ideal in $\fg$ spanned by $[\fg,\fg]$.
\end{x}
\vspace{0.3cm}

\begin{x}{\small\bf LEMMA} \ 
$\fg_\ssx$ is a semisimple Lie algebra.
\end{x}
\vspace{0.3cm}

\begin{x}{\small\bf THEOREM} \ 
Let \mG be a compact Lie group $-$then 
\[
\fg 
\ = \ 
\fz(\fg) \hsx \oplus \hsx \fg_\ssx,
\]
thus $\fg$ is reductive or still, \mG is reductive.
\end{x}
\vspace{0.3cm}

\begin{x}{\small\bf NOTATION} \ 
$G_\ssx$ is the analytic subgroup of \mG corresponding to $\fg_\ssx$.
\end{x}
\vspace{0.3cm}

\begin{x}{\small\bf NOTATION} \ 
$G^*$ is the commutator subgroup of \mG, i.e., the subgroup of \mG generated by the 
\[
x y x^{-1} y^{-1} \qquad (x, y \in G).
\]

[Note: \ 
$G^*$ is necessarily normal.]
\end{x}
\vspace{0.3cm}

\begin{x}{\small\bf THEOREM} \ 
Assume that \mG is connected $-$then $G^*$ is a compact connected Lie subgroup of \mG with Lie algebra $\fg_\ssx$, so 
$G^* = G_\ssx$, hence is semisimple.
\end{x}
\vspace{0.3cm}

\begin{x}{\small\bf THEOREM} \ 
Assume that \mG is connected $-$then \mG is the commuting product $Z(G)^0 G_\ssx$.
\end{x}
\vspace{0.3cm}

\begin{x}{\small\bf THEOREM} \ 
Assume that \mG is connected $-$then 
\[
G 
\ \approx \ 
(Z(G)^0  \hsx \times \hsx G_\ssx) / \Delta,
\]
where 
\[
\Delta 
\ \approx \ 
Z(G)^0  \hsx \cap \hsx G_\ssx
\]
is embedded in $Z(G)^0 \times G_\ssx$ via the arrow $z \ra (z^{-1}, z)$.
\vspace{0.2cm}

[Note: \ 
Spelled out, there is an exact sequence
\[
\{1\} 
\ra 
Z(G)^0  \hsx \cap \hsx G_\ssx
\overset{\iota}{\lra} 
Z(G)^0  \hsx \times \hsx G_\ssx 
\overset{\mu}{\lra} 
G 
\ra \{1\},
\]
where
\[
\iota(z) 
\ = \ 
(z^{-1}, z), 
\quad 
\mu(z,x) 
\ = \  zx.]
\]
\end{x}
\vspace{0.3cm}

\begin{x}{\small\bf \un{N.B.}} \ 
Structurally, $Z(G)^0$ is a torus and
\[
Z(G_\ssx) 
\ = \ 
Z(G) \hsx \cap \hsx G_\ssx
\]
is a finite abelian group.
\end{x}
\vspace{0.3cm}

\begin{x}{\small\bf SCHOLIUM} \ 
Assume that \mG is connected $-$then \mG is semisimple iff $Z(G)$ is finite.
\vspace{0.2cm}

[Note: \ 
Here is another way to put it: \mG is semisimple iff $G = G_\ssx$ or still, 
iff $G = G^*$.  
To see that connectedness is essential, consider the 8 element quaternion group 
$\{\pm 1, \pm i, \pm j, \pm k\}$ $-$then its commutator group is $\{\pm 1\}$.]
\end{x}
\vspace{0.3cm}

\begin{x}{\small\bf EXAMPLE} \ 
The center of $G/Z(G)$ is trivial, so $G/Z(G)$ (which is connected) is semisimple.
\end{x}
\vspace{0.3cm}

There are simple ideals $\fh_i \subset \fg_\ssx$ such that 
\[
\fg_\ssx 
\ = \ 
\bigoplus\limits_{i = 1}^r \ \fh_i
\]
with $[\fh_i,\fh_j] = 0$ for $i \neq j$ and such that the span of $[\fh_i,\fh_i] = \fh_i$.
\vspace{0.3cm}

Put $H_i = \exp \fh_i$.
\vspace{0.3cm}

\begin{x}{\small\bf LEMMA} \ 
$H_i$ is a compact connected normal Lie subgroup of $G_\ssx$ and its Lie algebra is $\fh_i$ (hence $H_i$ is simple).
\end{x}
\vspace{0.3cm}

\begin{x}{\small\bf LEMMA} \ 
A proper compact normal Lie subgroup of $H_i$ is necessarily discrete, finite, and central.
\end{x}
\vspace{0.3cm}

\begin{x}{\small\bf LEMMA} \ 
There is a decomposition
\[
G_\ssx 
\ = \ 
H_1 \cdots H_r,
\]
where $H_i$ and $H_j$ commute ($i \neq j$).
\end{x}
\vspace{0.3cm}

\begin{x}{\small\bf \un{N.B.}} \ 
The differential of the arrow 
\[
H_1 \times \cdots \times H_r \ra G_\ssx
\]
defined by the rule
\[
(x_1, \ldots, x_r) \ra x_1 \cdots x_r
\]
is the identity map, thus its kernel $\Delta$ is discrete and normal, thus finite and central as well, so 
\[
G_\ssx 
\ \approx \ 
(H_1 \times \cdots \times H_r) /  \Delta.
\]
\end{x}
\vspace{0.75cm}

\[
\text{APPENDIX}
\]
\vspace{0.05cm}

Let \mG be a compact connected Lie group.
\\

\qquad{\small\bf DEFINITION} \ 
\mG is 
\un{tall}
\index{tall} 
if for each positive integer $n$, there are but finitely many elements of $\widehat{G}$ of degree $n$.
\\[-.2cm]

\qquad{\small\bf THEOREM} \ 
\mG is semisimple iff \mG is tall.
\\[-.2cm]

\qquad{\small\bf REMARK} \ 
If \mG is not semisimple, then \mG possesses infinitely many nonisomorphic irreducible representations of degree 1.

\chapter{
$\boldsymbol{\S}$\textbf{4}.\quad  MAXIMAL TORI}
\setlength\parindent{2em}
\setcounter{theoremn}{0}
\renewcommand{\thepage}{B I \S4-\arabic{page}}

\qquad Let \mG be a connected Lie group, $\mathfrak{g}$ its Lie algebra.

\begin{x}{\small\bf LEMMA} \ 
Every connected abelian subgroup $A \subset G$ is contained in a maximal connected abelian subgroup $T \subset G$.
\end{x}
\vspace{0.3cm}

\begin{x}{\small\bf \un{N.B.}} \ 
\mT is compact
\vspace{0.2cm}

[In fact, $\ov{T}$ is connected and abelian.]
\end{x}
\vspace{0.3cm}

\begin{x}{\small\bf DEFINITION}  \ 
A 
\un{maximal torus}
\index{maximal torus} 
$\Torus \subset G$ is a maximal connected abelian subgroup of \mG.
\vspace{0.2cm}

[Note: \
$\Torus$ is a torus $\ldots$ .]
\end{x}
\vspace{0.3cm}

\begin{x}{\small\bf THEOREM} \ 
Assume that \mG is connected and let $\Torus_1 \subset G$, $\Torus_2 \subset G$ be maximal tori $-$then 
$\exists \ x \in G$ such that $x \Torus_1 x^{-1} = \Torus_2$.
\end{x}
\vspace{0.3cm}

\begin{x}{\small\bf THEOREM} \ 
Assume that \mG is connected and let $\Torus \subset G$ be a maximal torus $-$then 
\[
G 
\ = \ 
\bigcup\limits_{x \in G} \hsx 
x \Torus x^{-1}.
\]
\end{x}
\vspace{0.3cm}

\begin{x}{\small\bf APPLICATION} \ 
The exponential map $\exp : \fg \ra G$ is surjective.
\vspace{0.2cm}

[Every element of \mG belongs to a maximal torus and the exponential map of a torus is surjective.]
\end{x}
\vspace{0.3cm}

\begin{x}{\small\bf LEMMA} \ 
Assume that \mG is connected and let $\Torus \subset G$ be a maximal torus $-$then the centralizer of $\Torus$ in \mG is $\Torus$ itself.
\end{x}
\vspace{0.3cm}


\begin{x}{\small\bf APPLICATION} \ 
The center of \mG is contained in $\Torus$, i.e., $Z(G) \subset \Torus$.
\\[-.5cm]

[Note: \
More is true, viz.
\[
Z(G) 
\ = \ 
\bigcap\limits_\Torus, 
\]
the intersection being taken over all maximal tori in \mG.]
\\[-.2cm]
\end{x}

\begin{x}{\small\bf LEMMA} \ 
Assume that \mG is connected and let $\Torus \subset G$ be a maximal torus $-$then $\Torus$ is a maximal abelian subgroup.
\\[-.2cm]
\end{x}

\begin{x}{\small\bf REMARK} \ 
A maximal abelian subgroup need not be a maximal torus.
\\[-.5cm]

[In $\SO(3)$, there is a maximal abelian subgroup which is isomorphic to $(\Z / 2\Z)^2$, hence is not a maximal torus.]
\\[-.2cm]
\end{x}
\vspace{0.3cm}

\begin{x}{\small\bf NOTATION} \ 
Given a torus $\Torus \subset G$, let $N(\Torus)$ be its normalizer in \mG.
\index{$N(\Torus)$ }
\\[-.2cm]
\end{x}

\begin{x}{\small\bf LEMMA} \ 
The quotient $N(\Torus) / \Torus$ is finite iff $\Torus$ is a maximal torus.
\\[-.2cm]
\end{x}

Let \mG be a compact connnected Lie group, $\Torus \subset G$ a maximal torus.
\\

\begin{x}{\small\bf DEFINITION} \ 
The 
\un{Weyl group}
\index{Weyl group} 
of $\Torus$ in \mG is the quotient 
\[
W 
\ = \ 
N(\Torus) / \Torus.
\]
\end{x}

\begin{x}{\small\bf \un{N.B.}} \ 
Different choices of $\Torus$ give rise to isomorphic Weyl groups.
\vspace{0.3cm}

Fix a maximal torus $\Torus \subset G$ $-$then $N(\Torus)$ operates on $\Torus$ by conjugation:
\[
\begin{cases}
\ N(\Torus) \hsx \times \hsx \Torus \ra \Torus\\
\ (n,t) \ra n \Torus n^{-1}
\end{cases}
.
\]
Since $\Torus$ operates trivially on itself, there is an induced operation of the Weyl group:
\[
W \hsx \times \hsx \Torus \ra \Torus.
\]

[Note: \
The action of \mW is on the left, thus the orbit space is denoted by $W \backslash \Torus$.]
\\[-.2cm]
\end{x}

\begin{x}{\small\bf LEMMA} \ 
The canonical homomorphism $W \ra \Aut \Torus$ is injective.
\\[-.2cm]
\end{x}

\begin{x}{\small\bf LEMMA} \ 
Two elements of $\Torus$ are conjugate in \mG iff they lie on the same orbit under the action of \mW.
\\[-.2cm]
\end{x}

\begin{x}{\small\bf RAPPEL} \ 
Let \mG be a compact group and let \mX be a Hausdorff topological space on which \mG operates to the left $-$then the action arrow 
\[
G \hsx \times \hsx X \ra X
\]
is a closed map.  
Equip the orbit space $G \backslash X$ with the quotient topology and let $\pi:X \ra G \backslash X$ be the projection.  
Then:
\\[-.2cm]

\qquad \textbullet \quad 
$G \backslash X$ is a Hausdorff space.
\\[-.2cm]

\qquad \textbullet \quad 
\mX is compact iff $G \backslash X$ is compact.
\\[-.2cm]

\qquad \textbullet \quad
$\pi:X \ra G \backslash X$ is open, closed, and proper.
\\[-.2cm]
\end{x}

\begin{x}{\small\bf EXAMPLE} \ 
$W \backslash \Torus$ is a compact Hausdorff space.
\\[-.2cm]
\end{x}

\begin{x}{\small\bf NOTATION} \ 
$\CON(G)$ is the set of conjugacy classes of \mG.
\\[-.2cm]
\end{x}

Geometrically, $\CON(G)$ is the orbit space under the action of \mG on itself via inner automorphisms:
\[
\begin{cases}
\ G \hsx \times G \hsx \ra G\\
\ (x, y) \ra x y x^{-1}
\end{cases}
.
\]
It carries the quotient topology per the projection $G \ra \CON(G)$ under which it is a compact Hausdorff space.
\\[-.2cm]

\begin{x}{\small\bf RAPPEL} \ 
A one-to-one continuous map from a compact Hausdorff space \mX onto a Hausdorff space \mY is a homeomorphism.
\\[-.2cm]
\end{x}

\begin{x}{\small\bf THEOREM} \ 
The arrow
\[
W\backslash \Torus \ra \CON(G)
\]
which sends the $W$-orbit $W t$ of $t \in \Torus$ to the conjugacy class of $t \in \Torus$ in \mG is a (well defined) homeomorphism.

[The map is injective (cf. $\#16$), continuous (see below), and surjective (cf. $\#5$), so $\#20$ is applicable.]
\\[-.3cm]

[Note: \
To check the continuity of the arrow
\[
W \backslash \Torus \ra \CON(G),
\]
bear in mind that $W \backslash \Torus$ has the quotient topology, thus it suffices to check the continuity of the composition
\[
\Torus \ra W \backslash \Torus \ra \CON(G).
\]
But this map is just the restriction to $\Torus$ of the arrow
\[
G \ra \CON(G).]
\]
\end{x}

\begin{x}{\small\bf NOTATION} \ 
\\[-.2cm]

\qquad \textbullet \quad Given $f \in C(\Torus)$ and $w \in W$, $w \cdot f$ is the function in $C(\Torus)$ defined by the rule
\[
(w \cdot f) (t) 
\ = \ 
f(n^{-1} t n) 
\qquad (w = n \Torus).
\]
\\[-.5cm]

\qquad \textbullet \quad Given $f \in C(G)$ and $x \in G$, $x \cdot f$ is the function in $C(G)$ defined by the rule
\[
(x \cdot f) (y) 
\ = \ 
f(x^{-1} y x) .
\]
\end{x}
\vspace{0.3cm}

\begin{x}{\small\bf \un{N.B.}} \ 
These rules define operations
\[
\begin{cases}
\ W \hsx \times \hsx C(\Torus) \ra C(\Torus) \\
\ G \hsx \times \hsx C(G) \ra C(G)
\end{cases}
\]
with associated invariants
\[
\begin{cases}
\ C( W \backslash \Torus) \hsx = \hsx C(\Torus)^W\\
\ CL(G) \hsx = \hsx C(G)^G
\end{cases}
.
\]
\vspace{0.2cm}

[Note: \ 
$CL(G)$ is the subspace of $C(G)$ comprised of the continuous class functions (cf. \S2, \#27) or still, the space $C(\CON(G))$.]
\end{x}
\vspace{0.3cm}

\begin{x}{\small\bf LEMMA} \ 
The arrow 
\[
f \ra \restr{f}{\Torus}
\]
of restriction defines an isomorphism
\[
CL(G) \ra C(\Torus)^W.
\]
\end{x}
\vspace{0.3cm}


\chapter{
$\boldsymbol{\S}$\textbf{5}.\quad  REGULARITY}
\setlength\parindent{2em}
\setcounter{theoremn}{0}
\renewcommand{\thepage}{B I \S5-\arabic{page}}

\qquad Let \mG be a connected Lie group with lie algebra $\mathfrak{g}$.  
Consider the polynomial
\[
\det((t + 1) - \Ad(x)) \ = \ 
\sum\limits_{i = 0}^n \ D_i(x) \hsx t^i \qquad (x \in G),
\]
where $t$ is an indeterminate and $n = \dim G$.  
The $D_i$ are real analytic functions on \mG and $D_n = 1$.  
Let $\ell$ be the smallest positive integer such that $D_\ell \neq 0$ $-$then $\ell$ is called the 
\un{rank}
\index{rank} 
of \mG and an element $x \in G$ is said to be 
\un{singular}
\index{singular} 
or 
\un{regular}
\index{regular} 
according to whether $D_\ell(x) = 0$ or not.

\begin{x}{\small\bf NOTATION} \ 
$G^\sreg$ is the set of regular elements in \mG.
\\[-.25cm]
\end{x}

\begin{x}{\small\bf LEMMA} \ 
$G^\sreg$ is an open, dense subset of \mG while its complement, the set of singular elements, is a set of Haar measure zero 
(right or left).
\\[-.25cm]
\end{x}

\begin{x}{\small\bf \un{N.B.}} \ 
$G^\sreg$ is inner automorphism invariant and stable under multiplication by elements from the center of \mG.
\\[-.25cm]
\end{x}

From this point forward, assume that \mG is a compact connected Lie group.
\\[-.25cm]

\begin{x}{\small\bf LEMMA} \ 
The set of singular elements in \mG is a finite union of submanifolds of \mG, each of dimension $\leq \dim G - 3$.
\\[-.25cm]

[Note: \ Therefore $G^\sreg$ is path connected.]
\end{x}
\vspace{0.3cm}

\begin{x}{\small\bf RAPPEL} \ 
The fundamental group of a connected Lie group is abelian.
\\[-.25cm]
\end{x}

Fix a maximal torus $\Torus$.
\\[-.2cm]


\begin{x}{\small\bf LEMMA} \ 
The quotient $G/\Torus$ is simply connected.
\\[-.25cm]
\end{x}

\begin{x}{\small\bf LEMMA} \ 
The induced map $\pi_1(\Torus) \ra \pi_1(G)$ is surjective.
\\[-.2cm]

PROOF \ 
Consider the exact sequence
\[
\pi_1(\Torus) \ra \pi_1(G) \ra \pi_1(G/\Torus)
\]
arising from the fibration $\Torus \ra G \ra G/\Torus$.
\\[-.25cm]
\end{x}

\begin{x}{\small\bf THEOREM} \ 
$\pi_1(G)$ is a finitely generated abelian group.
\\[-.25cm]

[Note: \ 
If \mG is semisimple, then $\pi_1(G)$ is finite, thus its universal covering group $\widetilde{G}$ is compact.]
\end{x}
\vspace{0.3cm}

\begin{x}{\small\bf LEMMA} \ 
An element $x \in G$ is regular iff $x$ lies in a unique maximal torus.  
\\[-.5cm]
\end{x}

Put
\[
\Torus^\sreg 
\ = \ 
\Torus \cap G^\sreg.
\]
\\[-.5cm]

\begin{x}{\small\bf THEOREM} \ 
\[
G^\sreg 
\ = \ 
\bigcup\limits_{x \in G} \hsx x \Torus^\sreg x^{-1}.
\]
\end{x}

\begin{x}{\small\bf THEOREM} \ 
The map
\[
\mu:G/\Torus \times \Torus^\sreg \ra G^\sreg
\]
that sends
\[
(x \Torus, t) 
\quad \text{to} \quad 
x \hsy t \hsy x^{-1}
\]
is a surjective, $\abs{W}$-to-one local diffeomorphism.
\\[-.2cm]

[To verify the ``$\abs{W}$-to-one'' claim, observe that $\forall \ w \in W$ $(w = n\Torus)$, 
\allowdisplaybreaks
\begin{align*}
\\mu(x \hsy n^{-1} \hsy \Torus, n \hsy t \hsy n^{-1}) \ 
&=\ 
x \hsy n^{-1} \cdot n \hsy t \hsy n^{-1} \cdot n \hsy x^{-1} 
\\[11pt]
&=\ 
x \hsy t \hsy x^{-1}
\\[11pt]
&=\ 
\\mu (x \hsy \Torus,t),
\end{align*}
hence
\[
\abs{\\mu^{-1}(x t x^{-1})} 
\ \geq \ 
\abs{W}.
\]
In the opposite direction, suppose that
\[
x \hsy t \hsy x^{-1}
\ = \ 
y \hsy s \hsy y^{-1} 
\qquad (t, s \in \Torus^\sreg).
\]
Then there is a $w \in W$ such that
\[
s 
\ = \ 
n t n^{-1} 
\qquad (w = n \Torus) \qquad \text{(cf. $\S4$, $\#16$)}
\] 
from which 
\[
x \hsy t \hsy x^{-1}
\ = \ 
y \hsy n \hsy t \hsy n^{-1} \hsy y^{-1},
\]
so $x^{-1} \hsy y \hsy n \in G_t$, the centralizer of $t$ in \mG.  
But 
\[
t \in \Torus^\sreg 
\implies 
G_t^0 
\ = \ 
\Torus
\]
which implies that conjugation by $x^{-1} y n$ preserves $\Torus$ 
($G_t^0$ being the identity component of $G_t$), i.e., 
\[
n^\prime 
\ \equiv \ 
x^{-1} y n \in \N
\]

\qquad $\implies$
\allowdisplaybreaks
\begin{align*}
(y \hsy \Torus , s) \ 
&=\ 
(x \hsy (x^{-1} \hsy y \hsy n)n^{-1}\Torus,n \hsy t \hsy n^{-1})
\\[11pt]
&=\ 
(x \hsy n^\prime n^{-1}  \Torus, n \hsy t \hsy n^{-1})
\\[11pt]
&=\ 
(x \hsy n^\prime \hsy n^{-1} \Torus, n(n^{-1} \hsy y^{-1} \hsy x) t (x^{-1} \hsy y \hsy n) n^{-1})
\\[11pt]
&=\ 
(x(n^\prime n^{-1}) \Torus, (n^\prime \hsy  n^{-1})^{-1} t(n^\prime \hsy n^{-1}))
\\[11pt]
&\in \\mu^{-1}(x t x^{-1}).]
\end{align*}

\end{x}


[Note: \ 
$G_t^0$ is a compact connected Lie group and $\Torus \subset G_t^0$ is a maximal torus.  
If $\Torus \neq G_t^0$, $\exists \ z \in G_t^0$ : $z \Torus z^{-1} \neq \Torus$ 
(cf. $\S4$, $\#5$ (applied to $G_t^0$)).  
But then 
\[
t 
\ = \ 
z t z^{-1} \in z\Torus z^{-1},
\]
contradicting the regularity of $t$ (cf. $\#9$).]
\vspace{0.3cm}

Let $\fg$ be the Lie algebra of \mG, $\ft$ the Lie algebra of $\Torus$.  
Since \mG is compact, there is a positive definite symmetric bilinear form on $\fg$ which is invariant under the adjoint representation:
\[
\Ad : G \ra \Aut \fg.
\]
Denote by $\fg/\ft$ the orthogonal complement of $\ft$ in $\fg$ $-$then $\fg/\ft$ is stable under $\Ad \hsx \Torus$, 
which gives rise to an induced action
\[
\Ad_{G/\Torus} : \Torus \ra \Aut \fg/\ft.
\]
Denoting by $I_{G/\Torus}$ the identity map $\fg/\ft \ra \fg/\ft$, one may then attach to each $t \in \Torus$ the endomorphism
\[
\Ad_{G/\Torus}\big(t^{-1}\big) - I_{G/\Torus}
\]
of $\fg/\ft$.
\\[-.2cm]

\begin{x}{\small\bf LEMMA} \ 
The determinant of 
\[
\Ad_{G/\Torus}\big(t^{-1}\big) - I_{G/\Torus}
\]
is positive on the subset of $\Torus$ comprised of the topologically cyclic elements.
\end{x}
\vspace{0.3cm}

\begin{x}{\small\bf INTEGRATION FORMULA} \ 
For any continuous function $f \in C(G)$,
\[
\int_G \hsx f(x) \hsx \td_G(x) 
\ = \ 
\frac{1}{W} \hsx 
\int_\Torus \hsx 
\bigg[
\det\big(\Ad_{G/\Torus} (t^{-1}) - I_{G/\Torus}\big) 
\int_G \hsx
f(x t x^{-1}) \hsx \td_G(x)
\bigg]
\hsx \td_\Torus(t).
\]
\\[-.7cm]

[Note: \ 
$\td_G(x)$ is normalized Haar measure on \mG and $\td_\Torus(t)$ is normalized Haar measure on $\Torus$.]
\end{x}
\vspace{0.3cm}

\begin{x}{\small\bf SCHOLIUM} \ 
For any continuous class function $f \in CL(G)$,

\[
\int_G \hsx f(x) \hsx \td_G(x) 
\ = \ 
\frac{1}{\abs{W}} \hsx 
\int_\Torus \hsx \det\big(\Ad_{G/\Torus} (t^{-1}) - I_{G/\Torus}\big) f(t) \hsx \td_\Torus(t).
\]
\end{x}
\vspace{0.3cm}

\[
\text{APPENDIX}
\]

Consider the polynomial
\[
\det(t - \ad(X)) 
\ = \ 
\sum\limits_{i = 0}^n \hsx 
d_i(X) t^i 
\qquad (X \in \fg),
\]
where $t$ is an indeterminate and $n = \dim \fg$.  
The $d_i$ are polynomial functions on $\fg$ and $d_n = 1$.  
Let $\ell$ be the smallest positive integer such that $d_\ell \neq 0$ $-$then $\ell$ is called the \un{rank} of $\fg$ 
and an element $X \in \fg$ is said to be \un{singular or regular} according to whether $d_\ell(X) = 0$ or not.
\\[-.2cm]

{\small\bf \un{N.B.}} \ 
The rank of $\fg$ equals the rank of \mG, both being equal to the dimension of $\ft$.
\\[-.2cm]

{\small\bf NOTATION} \ 
$\fg^\sreg$ is the set of regular elements in $\fg$.
\\[-.2cm]

{\small\bf LEMMA} \ 
$\fg^\sreg$ is an open, dense subset of $\fg$.
\\[-.2cm]

{\small\bf NOTATION} \ 
$\fG \equiv \Int \fg$ is the adjoint group of $\fg$.
\\[-.2cm]

[Note: \ 
Recall that the arrow
\[
\Ad : G \ra \fG
\]
is surjective with kernel $Z(G)$, so
\[
G/Z(G) 
\ \approx \ 
\fG.]
\]

Put
\[
\ft^\sreg 
\ = \ 
\ft \cap \fg^\sreg.
\]

{\small\bf THEOREM} \ 
\[
\fg^\sreg 
\ = \ 
\bigcup\limits_{x \in \fG} \hsx x\big(\ft^\sreg\big).
\]

\chapter{
$\boldsymbol{\S}$\textbf{6}.\quad  WEIGHTS AND ROOTS}
\setlength\parindent{2em}
\setcounter{theoremn}{0}
\renewcommand{\thepage}{B I \S6-\arabic{page}}

\qquad Let \mG be a compact connected semisimple Lie group, $\Torus \subset G$ a maximal torus.  
Denote their respective Lie algebras by $\mathfrak{g}$, $\mathfrak{t}$ and let 
$\mathfrak{g}_\Cx$, $\mathfrak{t}_\Cx$ stand for their complexifications.  

Suppose that $(\pi,V)$ is a representation of \mG $-$then \mV can be equipped with a \mG-invariant inner product, thus rendering matters unitary.
\vspace{0.3cm}

\begin{x}{\small\bf LEMMA} \ 
 $\td\pi$ is skew-adjoint on $\fg$ (hence self-adjoint on $\sqrt{-1} \fg$).
 \\[-.25cm]
 
[Given $X \in \fg$, apply $\ds{\frac{\td}{\td t}\ds |_{t = 0}}$ to
\[
 \langle \pi(exp \hsy t) v_1, \pi(\exp \hsy t X) v_2) \rangle
 \ = \ 
 \langle v_1, v_2 \rangle
\]
 to get
\[
\langle \td \pi (X) v_1, v_2 \rangle + \langle v_1, \td\pi(X) v_2 \rangle = 0.]
\]
\\[-.5cm]
\end{x}

\begin{x}{\small\bf \un{N.B.}} \ 
$\forall \ X \in \fg$, 
\[
\pi(\exp \hsy X) 
\ = \ 
e^{\td\pi(X)}.
\]
\\[-.5cm]
\end{x}

\begin{x}{\small\bf LEMMA} \ 
\mV is simultaneously diagonalizable under the action of $\ft_\Cx$.
\\[-.25cm]

[This is because 
\[
\{\td\pi(H): H \in \ft_\Cx\}
\]
 is a commuting family of normal operators.]
 \\[-.25cm]
\end{x}

Consequently, there is a finite set $\Phi(V) \subset \ft_\Cx^* - \{0\}$, the elements of which being the 
\un{weights} of \mV, such that
\[
V 
\ = \ 
V^0 \bigoplus\limits_{\lambda \in \Phi(V)} \hsx V^\lambda,
\]
where
\[
V^0 
\ = \ 
\{v \in V: \td\pi(H)v = 0\} \qquad (H \in \ft_\Cx)
\]
 and
\[
V^\lambda 
\ = \ 
\{v \in V: \td\pi(H)v = \lambda(H) v\} \qquad (H \in \ft_\Cx).
\]
 \\[-.2cm]

\begin{x}{\small\bf LEMMA} \ 
Fix a $\lambda \in \Phi(V)$ $-$then $\restr{\lambda}{\ft}$ is purely imaginary and 
$\restr{\lambda}{\sqrt{-1} \hsy \ft}$ is purely real.
\end{x}
\vspace{0.3cm}

\begin{x}{\small\bf \un{N.B.}} \ 
Given $t \in \Torus$, choose $H \in \ft$ such that $t = \exp H$ $-$then $\forall \ v \in V^\lambda$, 
\[
\pi(t) v 
\ = \ 
\pi(\exp H) v 
\ = \ 
e^{\td \pi(H)} v 
\ = \ 
e^{\lambda(H)}v.
\]
\end{x}
\vspace{0.3cm}

\begin{x}{\small\bf RAPPEL} \ 
Denote by $I_x$ the inner automorphism $y \ra x\hsy y \hsy x^{-1}$ attached to $x \in G$ $-$then 
the \un{adjoint representation} of \mG is the homomorphism 
$\Ad:G \ra \Aut \fg$ defined by the rule
\[
\Ad(x) 
\ = \ 
(\td I_x)_e
\]
 and the \un{adjoint representation} of $\fg$ is the homomorphism 
$\ad:\fg \ra \End \fg$ defined by the rule
\[
\ad(X) 
\ = \ 
(\td \Ad)_e(X).
\]
\end{x}
\vspace{0.3cm}

\begin{x}{\small\bf \un{N.B.}} \ 
$\forall \ X, Y \in \fg$, 
\[
\ad(X) Y 
\ = \ 
[X,Y].
\] 
\end{x}

\qquad \textbullet \quad 
For each $x \in G$, extend the domain of $\Ad(x)$ from $\fg$ to $\fg_\Cx$ by complex linearity.
\\[-.2cm]

\qquad \textbullet \quad 
For each $X \in \fg$, extend the domain of $\ad(X)$ from $\fg$ to $\fg_\Cx$ by complex linearity.


\begin{x}{\small\bf LEMMA} \ 
$(\Ad, \fg_\Cx)$ is a representation of \mG with differential $(\ad,\fg_\Cx)$.
\end{x}
\vspace{0.3cm}

Take now $V = \fg_\Cx$, let $\pi = \Ad$, and abbreviate $(\fg_\Cx)^\alpha$ to $\fg^\alpha$ ($\alpha \in \Phi(\fg_\Cx)$) 
$-$then $\fg^0 = \ft_\Cx$ and there is a weight space decomposition
\[
\fg_\Cx 
\ = \ 
\fg^0 \hsx \bigoplus\limits_{\alpha \in \Phi(\fg_\Cx)} \hsx \fg^\alpha.
\]

\begin{x}{\small\bf TERMINOLOGY} \ 
The elements $\alpha \in \Phi(\fg_\Cx)$ are called the \un{roots} of the pair $(\fg_\Cx, \ft_\Cx)$.
\end{x}
\vspace{0.3cm}

\begin{x}{\small\bf \un{N.B.}} \ 
\[
\fg^\alpha 
\ = \ 
\{X \in \fg_\Cx : [H,X] = \alpha(H) \hsy X (H  \in \ft_\Cx)\}.
\]
\end{x}
\vspace{0.3cm}

\begin{x}{\small\bf LEMMA} \ 
$\forall \ \alpha \in \Phi(\fg_\Cx)$, $\forall \ \lambda \in \Phi(V) \cup \{0\}$, 
\[
\td \pi (\fg^\alpha) V^\lambda \subset V^{\alpha + \lambda}.
\]
 
 PROOF \ 
 Let $H \in \ft_\Cx$, $X_\alpha \in \fg^\alpha$, $v_\lambda \in V^\lambda$ $-$then
 \allowdisplaybreaks
 \begin{align*}
 \td \pi(H) \td \pi(X_\alpha) \hsy v_\lambda \ 
  &=
 (\td \pi(X_\alpha) \td \pi(H)  + [\td \pi(H), \td \pi(X_\alpha) ]) \hsy v_\lambda
 \\[11pt]
&=
 (\td \pi(X_\alpha) \td \pi(H)  + \td \pi([H,X_\alpha]) \hsy v_\lambda
 \\[11pt]
 &=
 (\td \pi(X_\alpha) \hsy \td \pi(H)  +\alpha(H) \hsy \td \pi(X_\alpha) \hsy v_\lambda 
 \\[11pt]
 &=
(\lambda(H) + \alpha(H) \td \pi(X_\alpha) \hsy v_\lambda 
\end{align*}
\qquad $\implies$ 
\[
\td \pi(X_\alpha) v_\lambda  \in V^{\alpha + \lambda}.
\]

[Note: \ 
Take $\lambda = 0$ to see that
\[
\td \pi(\fg^\alpha) V^0  \subset V^\alpha.
\]
\\[-.5cm]
\end{x}

\begin{x}{\small\bf APPLICATION} \ 
$\forall \ \alpha$, $\beta \in \Phi(\fg_\Cx) \cup \{0\}$, 
\[
[\fg^\alpha, \fg^\beta] 
\ \subset \ 
\fg^{\alpha + \beta}.
\]
\\[-1.20cm]
\end{x}

\begin{x}{\small\bf LEMMA} \ 
Let $\langle \ , \ \rangle$ be an $\Ad \hsx G$ invariant inner product on $\fg_\Cx$ $-$then for all 
$\alpha$, $\beta \in \Phi(\fg_\Cx) \cup \{0\}$,  
\[
\langle\fg^\alpha, \fg^\beta\rangle
\ = \ 
0 \quad \text{if} \quad \alpha + \beta \neq 0.
\]
\\[-1.2cm]
\end{x}

\begin{x}{\small\bf LEMMA} \ 
$\forall \ \alpha \in \Phi(\fg_\Cx)$, $\dim \fg^\alpha = 1$ and the only multiples of $\alpha$ in $\Phi(\fg_\Cx)$ are 
$\pm\alpha$.
\\[-.25cm]
\end{x}

\begin{x}{\small\bf NOTATION} \ 
$\sigma : \fg_\Cx \ra \fg_\Cx$ is the map that sends $Z = X + \sqrt{-1} \hsx Y$ to 
$\ov{Z} = X - \sqrt{-1} \hsx Y$ $(X, \ Y \in \fg)$.
\\[-.25cm]
\end{x}

\begin{x}{\small\bf LEMMA} \ 
$\sigma$ is an $\R$-linear involution which preserves the bracket, i.e., 
\[
\sigma \big([Z_1, Z_2]\big) 
\ = \ 
[\sigma Z_1, \sigma Z_2] \qquad (Z_1, Z_2 \in \fg_\Cx).
\]
\\[-1.20cm]
\end{x}

\begin{x}{\small\bf \un{N.B.}} \ 
$\forall \ \alpha \in \Phi(\fg_\Cx)$, 
\[
\sigma \fg^\alpha 
\ = \ 
\fg^{-\alpha}.
\]
\\[-1.20cm]
\end{x}

\begin{x}{\small\bf RAPPEL} \ 
The \un{Killing form} of $\fg_\Cx$ is the bilinear form 
$B: \fg_\Cx \times \fg_\Cx \ra \Cx$ gvien by
\[
B(Z_1, Z_2) 
\ = \ 
\tr(\ad(Z_1) \circ \ad(Z_2)).
\]
\end{x}


\begin{x}{\small\bf PROPERTIES} \ 
\\[-.3cm]

\qquad \textbullet \quad 
$\forall \ x \in G$, $\forall \ Z_1, Z_2 \in \fg_\Cx$, 
\[
B(\Ad(x) Z_1, \Ad(x) Z_2) 
\ = \ 
B(Z_1, Z_2).
\]

\qquad \textbullet \quad 
$\forall \ Z, Z_1, Z_2 \in \fg_\Cx$, 
\[
B(\ad(Z) Z_1, Z_2) 
\ = \ 
-B(Z_1, \ad(Z)Z_2).
\]
\end{x}

\begin{x}{\small\bf \un{N.B.}} \ 
The prescription
\[
\langle Z_1, Z_2 \rangle_\sigma 
\ = \ 
-B(Z_1, \sigma Z_2)
\]
is an $\Ad$ $G$ invariant inner product on $\fg_\Cx$.
\end{x}
\vspace{0.3cm}

Every $\alpha \in \Phi(\fg_\Cx)$ is determined by its restriction to either $\ft$ or $\sqrt{-1} \ \ft$, so
$\alpha$ can be viewed as an element of $(\sqrt{-1} \ \ft)^*$ (purely real) or of $\ft^*$ (purely imaginary).
\\[-.2cm]

\begin{x}{\small\bf CONSTRUCTION} \ 
\mB induces an isomorphism between $\sqrt{-1} \ \ft$ and $(\sqrt{-1} \ \ft)^*$ as follows:  
Given $\lambda \in (\sqrt{-1} \ \ft)^*$, define $H_\lambda \in \sqrt{-1} \ \ft$ by the relation
\[
\lambda(H) 
\ = \ 
B(H, H_\lambda) \qquad (H \in \sqrt{-1} \ \ft).
\]

[Note: \ 
\mB is negative definite on $\ft \times \ft$, hence \mB is a real inner product on the real vector space 
$\sqrt{-1} \ \ft$ and for $\lambda_1$, $\lambda_2 \in (\sqrt{-1} \ \ft)^*$, one writes
\[
B(\lambda_1, \lambda_2) 
\ = \ 
B(H_{\lambda_1}, H_{\lambda_2}).]
\]
\end{x}
\vspace{0.3cm}

\begin{x}{\small\bf DEFINITION} \ 
The vector $H_\alpha \in \sqrt{-1} \ \ft$ is called the \un{root vector} associated with $\alpha$.
\end{x}
\vspace{0.3cm}

\begin{x}{\small\bf LEMMA} \ 
The roots span $(\sqrt{-1} \ \ft)^*$ and the root vectors span $\sqrt{-1} \ \ft$.
\end{x}
\vspace{0.3cm}


\begin{x}{\small\bf LEMMA} \ 
Let $X_\alpha \in \fg^\alpha$, $X_{-\alpha} \in \fg^{-\alpha}$ $-$then
\[
[X_\alpha,X_{-\alpha}] 
\ = \ 
B(X_\alpha, X_{-\alpha}) H_\alpha.
\]

PROOF \ 
First of all, 
\[
[\fg^\alpha, \fg^{-\alpha}] 
\subset 
\fg^{\alpha - \alpha} 
\ = \ 
\fg^0 
\ = \ 
\ft_\Cx \qquad \text{(cf. $\#12$)},
\]
thus
\[
[X_\alpha,X_{-\alpha}] \in \ft_\Cx.
\]
Proceeding, $\forall \ H \in \ft_\Cx$, 
\allowdisplaybreaks
\begin{align*}
B([X_\alpha, X_{-\alpha}],H) \ 
&=\ 
-B([X_{-\alpha}, X_\alpha],H)
\\[11pt]
&=\ 
-B(\ad(X_{-\alpha})X_\alpha,H)
\\[11pt]
&=\ 
B(X_{-\alpha},\ad(X_{-\alpha})H)
\\[11pt]
&=\ 
B(X_\alpha,[X_{-\alpha},H])
\\[11pt]
&=\ 
-B(X_\alpha,[H,X_{-\alpha}])
\\[11pt]
&=\ 
-B(X_\alpha, -\alpha(H) X_{-\alpha})
\\[11pt]
&=\ 
\alpha(H) B(X_\alpha,X_{-\alpha})
\\[11pt]
&=\ 
B(H,H_\alpha) B(X_\alpha,X_{-\alpha})
\\[11pt]
&=\ 
B(H_\alpha,H) B(X_\alpha, X_{-\alpha})
\\[11pt]
&=\ 
B(B(X_\alpha, X_{-\alpha})H_\alpha,H)
\end{align*}
\qquad $\implies$ 
\[
[X_\alpha,X_{-\alpha}] 
\ = \ 
B(X_\alpha,X_{-\alpha}) H_\alpha.
\]
\end{x}
\vspace{0.3cm}


\begin{x}{\small\bf NOTATION} \ 
Put
\[
h_\alpha 
\ = \ 
2 \frac{H_\alpha}{B(H_\alpha,H_\alpha)}.
\]
Then $\alpha(h_\alpha) = 2$.
\end{x}
\vspace{0.3cm}

\begin{x}{\small\bf \un{N.B.}} \ 
$\forall \ \lambda \in (\sqrt{-1} \hsx t)^*$, 
\allowdisplaybreaks
\begin{align*}
\lambda(h_\alpha) \ 
&=\ 
\lambda \bigg(2 \frac{H_\alpha}{B(H_\alpha,H_\alpha)}\bigg)
\\[11pt]
&=\ 
\lambda \bigg(2 \frac{H_\alpha}{B(\alpha,\alpha)}\bigg)
\\[11pt]
&=\ 
2 \frac{\lambda (H_\alpha)}{B(\alpha,\alpha)}
\\[11pt]
&=\ 
2 \frac{B(H_\alpha, H_\lambda)}{B(\alpha,\alpha)}
\\[11pt]
&=\ 
2 \frac{B(H_\lambda, H_\alpha)}{B(\alpha,\alpha)}
\\[11pt]
&=\ 
2 \frac{B(\lambda, \alpha)}{B(\alpha,\alpha)}
\end{align*}
and analogously, $\forall \ H \in \sqrt{-1} \ \ft$, 
\[
\alpha(H) 
\ = \ 
2 \frac{B(H, h_\alpha)}{B(h_\alpha, h_\alpha)}.
\]
\end{x}
\vspace{0.3cm}

\begin{x}{\small\bf NORMALIZATION} \ 
Scale the data and choose $e_\alpha \in \fg^\alpha$, $f_\alpha \in \fg^{-\alpha}$ such that
\[
[e_\alpha, f_\alpha] 
\ = \ h_\alpha,
\]

hence 
\[
\begin{cases}
\ [h_\alpha, e_\alpha] \ = \ 2e_\alpha \\[3pt]
\ [h_\alpha, f_\alpha] \ = \ -2f_\alpha
\end{cases}
.
\]
Consequently, 
\[
\spanx_\Cx \{h_\alpha, e_\alpha, \ f_\alpha\} 
\ \approx \ 
\fs\fl(2,\Cx),
\]
where
\[
h_\alpha \longleftrightarrow h 
\ = \ 
\begin{pmatrix}
1 & \hspace{0.25cm} 0\\
0 & -1
\end{pmatrix}
, \quad
e_\alpha \longleftrightarrow e 
\ = \ 
\begin{pmatrix}
0 & \ 1\\
0 & \ 0
\end{pmatrix}
, \quad
f_\alpha \longleftrightarrow f 
\ = \ 
\begin{pmatrix}
0 & \ 0\\
1 & \ 0
\end{pmatrix}
.
\]
\end{x}
\vspace{0.3cm}

\begin{x}{\small\bf \un{N.B.}} \ 
Under this correspondence,
\begin{align*}
\fs\fu(2) \ 
&\approx \ 
\spanx_\R \{\sqrt{-1} \ h_\alpha, e_\alpha - f_\alpha, \sqrt{-1}\ (e_\alpha + f_\alpha) \}
\\[11pt]
&\equiv \ 
\fs_\alpha,
\end{align*}
where
\[
\sqrt{-1} \ h_\alpha \longleftrightarrow \sqrt{-1} \ h
\ = \ 
\begin{pmatrix}
\sqrt{-1} & \hspace{0.25cm} 0\\
0 &-\sqrt{-1}
\end{pmatrix}
\]
and
\[
e_\alpha - f_\alpha \longleftrightarrow e - f 
\ = \ 
\begin{pmatrix}
\hspace{0.25cm} 0 & 1\\
-1 &0
\end{pmatrix}
, \quad
\sqrt{-1} \ (e_\alpha + f_\alpha) \longleftrightarrow \sqrt{-1} \ (e + f)
\ = \
\begin{pmatrix}
\hspace{0.25cm} 0 & \sqrt{-1}\\
\sqrt{-1} &0
\end{pmatrix}
.
\]
\end{x}
\vspace{0.3cm}

\begin{x}{\small\bf LEMMA} \ 
The analytic subgroup $S_\alpha$ of \mG with Lie algebra $\fs_\alpha$ is compact and isomorphic to 
$\SU(2)$ or $\SU(2)/\Z_2$.
\end{x}
\vspace{0.3cm}

\begin{x}{\small\bf LEMMA} \ 
Let $(\pi,V)$ be a unitary representation of \mG $-$then $\forall \ \lambda \in \Phi(V)$, $\lambda(h_\alpha) \in \Z$.
\\[-.2cm]

PROOF \ 
In $\SU(2)$, $e^{2 \pi \sqrt{-1} \hsx h} = I$.  
This said, let $\phi_\alpha:\SU(2) \ra G$ be the arrow realizing the preceding setup and consider $\pi \circ \phi_\alpha$: 

\allowdisplaybreaks
\begin{align*}
I \ 
&=\ 
\pi \big(\phi_\alpha\big(e^{2 \pi \hsx \sqrt{-1} \hsx h}\big)\big)
\\[11pt]
&=\ 
\pi \big(e^{2 \pi \hsx \td \phi_\alpha (\sqrt{-1} \hsx h)} \big)
\\[11pt]
&=\ 
\pi \big(e^{2 \pi \hsx \sqrt{-1} \hsx h_\alpha} \big)
\\[11pt]
&=\ 
e^{2 \pi \hsx \sqrt{-1} \hsx \td \pi (h_\alpha)}.
\end{align*}
On the other hand, $\forall \ v \in V^\lambda$, 
\allowdisplaybreaks
\begin{align*}
v \ 
&=\ 
e^{2 \pi \hsx \sqrt{-1} \hsx \td \pi (h_\alpha)} \hsx v
\\[11pt]
&=\ 
e^{2 \pi \hsx \sqrt{-1} \hsx \lambda (h_\alpha)} \hsx v \qquad \text{(cf. $\#5$).}
\end{align*}
Therefore $\lambda(h_\alpha) \in \Z$.
\end{x}
\vspace{0.3cm}


\chapter{
$\boldsymbol{\S}$\textbf{7}.\quad  LATTICES}
\setlength\parindent{2em}
\setcounter{theoremn}{0}
\renewcommand{\thepage}{B I \S7-\arabic{page}}

\qquad Let \mV be a finite dimensional vector space over $\R$.
\vspace{0.3cm}

\begin{x}{\small\bf DEFINITION} \ 
A 
\un{lattice}
\index{lattice} 
in \mV is an additive subgroup $L \subset V$ such that 
\\[-.5cm]

\qquad \textbullet \quad 
\mL is closed;
\vspace{0.2cm}

\qquad \textbullet \quad
\mL is discrete;
\vspace{0.2cm}

\qquad \textbullet \quad 
\mL spans $V$.
\\[-.25cm]
\end{x}

\begin{x}{\small\bf EXAMPLE} \ 
$\Z^n$ is a lattice in $\R^n$.
\\[-.25cm]
\end{x}

\begin{x}{\small\bf DEFINITION} \ 
A 
\un{basis}
\index{basis (of a lattice)} 
for a lattice $L \subset V$ is a set $\{e_1, \ldots, e_n\} \subset L$ ($n = \dim V$) such that
\[
L 
\ = \ 
\bigg\{ \sum\limits_{i = 1}^n \  k_i \hsx e_i : k_i \in \Z \bigg\}.
\]
\\[-.75cm]
\end{x}

\begin{x}{\small\bf LEMMA} \ 
Every lattice has a basis.
\end{x}
\vspace{0.3cm}

\begin{x}{\small\bf DEFINITION} \ 
If \mL, \mK are lattices in \mV, then \mL is a \un{sublattice} of \mK if \mL is a subset of \mK.
\\[-.25cm]
\end{x}

\begin{x}{\small\bf LEMMA} \ 
If \mL is a sublattice of \mK, then $K/L$ is a finite group \mG.  
Moreover, there is a one-to-one correspondence between subgroups $H \subset G$ and the lattices 
$L \subset M \subset K$, viz.
\[
\pi(M) \ =\ H
\quad \text{and} \quad 
M \ = \ \pi^{-1}(H),
\]
where $\pi:K \ra K/L$ is the projection.
\\[-.25cm]
\end{x}

\begin{x}{\small\bf NOTATION} \ 
Given a lattice $L \subset V$, let
\[
L^* 
\ = \ 
\{v^* \in V^* : v^* (x) \in \Z \ \forall \ x \in L\}.
\]
\\[-1.25cm]
\end{x}


\begin{x}{\small\bf LEMMA} \ 
$L^*$ is a lattice in $V^*$, the 
\un{dual}
\index{dual (lattice)} 
of \mL.
\\[-.25cm]
\end{x}

Let $\{e_1, \ldots, e_n\}$ be a basis for a lattice $L \subset V$.  
Define $\{f_1, \ldots, f_n\}$ by
\[
f_j(e_i) 
\ = \ 
\delta_{i j}.
\]
\\[-1cm]

\begin{x}{\small\bf LEMMA} \ 
$\{f_1, \ldots, f_n\}$ is a basis for $L^*$.
\\[-.25cm]
\end{x}

\begin{x}{\small\bf APPLICATION} \ 
\[
L^{**} 
\ \approx \ 
L.
\]
\\[-.75cm]

[In fact, the condition
\[
f_j (e_i) 
\ = \ 
\delta_{i j}
\]
is symmetric in $f$ and $e$.]
\\[-.25cm]
\end{x}

\begin{x}{\small\bf LEMMA} \ 
Suppose that \mL is a sublattice of \mK $-$then $K^* \subset L^*$ and 
\[
L^* / K^* 
\ \approx \ 
\widehat{K / L}.
\]

PROOF \ 
The first point is obvious.  
As for the second, define a homomorphism $\rho:L^* \ra \widehat{K / L}$ by stipulating that 
\[
\rho(\ell^*) (x + L) 
\ = \ 
\exp (2 \pi \hsx \sqrt{-1} \ \ell^*(x)).
\]
Then the kernel of $\rho$ is $K^*$, so $\rho$ induces an injection $L^* / K^* \ra \widehat{K/L}$, thus 
\[
\abs{L^* / K^*} 
\ \leq \
\abs{\widehat{K/L}}
\ = \ 
\abs{K/L}.
\]
But then by duality, 

\[
\abs{L^* / K^*} 
\ \geq \ 
\abs{K^{**} / L^{**}}
\ = \ 
\abs{K/L}.
\]
\\[-1.25cm]
\end{x}

Let \mG be a compact connected semisimple Lie group, $\Torus \subset G$ a maximal torus.
\\[-.25cm]

\begin{x}{\small\bf CONVENTION} \ 
Identify $(\sqrt{-1} \ \ft)^{**}$ with $\sqrt{-1} \ \ft$ and let \mL be a lattice in $(\sqrt{-1} \ \ft)^*$ $-$then 
its dual is the lattice $L^* \subset \sqrt{-1} \ \ft$ specified by the prescription
\[
\{H \in \sqrt{-1} \ \ft : \lambda (H) \in \Z \ \forall \ \lambda \in L\}.
\]
\\[-1.25cm]
\end{x}

\begin{x}{\small\bf DEFINITION} \ 
The \un{root lattice} is the lattice $L_{\rt}$ in $(\sqrt{-1} \ \ft)^*$ generated by the $\alpha \in \Phi(\fg_\Cx)$.
\\[-.25cm]
\end{x}

\begin{x}{\small\bf DEFINITION} \ 
The \un{weight lattice} is the lattice $L_{\wt}$ in $(\sqrt{-1} \ \ft)^*$ given by 
\[
\{\lambda \in (\sqrt{-1} \ \ft)^* : \lambda(h_\alpha) \in \Z \ \forall \ \alpha \in \Phi(\fg_\Cx)\}. 
\]
\\[-1.25cm]
\end{x}

\begin{x}{\small\bf LEMMA} \ 
$L_\rt$ is a sublattice of $L_\wt$.
\\[-.25cm]
\end{x}

Given a character $\chi : \Torus \ra \bS^1$, there is a commutative diagram
\[
\begin{tikzcd}[sep=large]
{\ft} 
\ar{d}[swap]{\exp}  
\arrow[rr,"\td_\chi"] 
&&{\sqrt{-1} \ \R} 
\ar{d}{\exp}
\\
{\Torus} 
\ar{rr}[swap]{\chi} 
&&{\bS^1}
\end{tikzcd}
\]
and the arrow $\chi \ra \td_\chi$ implements an identification of $\widehat{\Torus}$ with the lattice
\[
\td\widehat{\Torus} 
\ \equiv \ 
\{\lambda \in (\sqrt{-1} \ \ft)^* : \restr{\lambda}{\exp^{-1} (e)} \subset 2 \pi \sqrt{-1} \ \Z\}.
\]
Here
\[
\td_\chi \in \Hom_\R (\ft,\sqrt{-1} \ \R)
\]
which we shall view as an element of 
\[
\Hom_\R (\sqrt{-1} \ \ft, \R)
\]
by writing
\[
\td \lambda (\sqrt{-1} \ H) 
\ = \ 
\sqrt{-1} \ \td \lambda (H) \qquad (H \in \ft).
\]

[Note: \ 
$\sqrt{-1} \ \R$ is the Lie algebra of $\bS^1$, the exponential map $\exp: \sqrt{-1} \ \R \ra \bS^1$ 
being the usual exponential function $\sqrt{-1} \ \theta \ra e^{\sqrt{-1} \ \theta}$.]
\\[-.2cm]

\begin{x}{\small\bf LEMMA} \ 
$L_\rt$ is a sublattice of $\td \widehat{\Torus}$ and $\td \widehat{\Torus}$ is a sublattice of $L_\wt$.
\\[-.25cm]
\end{x}

\begin{x}{\small\bf THEOREM} \ 
\\[-.25cm]

\qquad \textbullet \quad 
$Z(G) \ \approx \ \td\widehat{\Torus} /  L_\rt$.  
\\[-.25cm]

\qquad \textbullet \ 
$\pi_1(G) \quad \approx \ L_\wt / \td\widehat{\Torus}$.
\end{x}


\chapter{
$\boldsymbol{\S}$\textbf{8}.\quad  WEYL CHAMBERS AND WEYL GROUPS}
\setlength\parindent{2em}
\setcounter{theoremn}{0}
\renewcommand{\thepage}{B I \S8-\arabic{page}}

\qquad Let \mG be a compact connected semisimple Lie group, $\Torus \subset G$ a maximal torus, 
$\Phi(\mathfrak{g}_\Cx)$ the roots of the pair $(\mathfrak{g}_\Cx, \ft_\Cx)$.
\\[-.25cm]

\begin{x}{\small\bf DEFINITION} \ 
A subset $\Psi$ of $\Phi(\fg_\Cx)$ is a \un{simple system} of roots if it is a vector space basis for 
$(\sqrt{-1} \ \ft)^*$ and has the property that every root can be written as a linear combination
\[
\sum\limits_{\alpha \in \Psi} \ n_\alpha \hsx \alpha,
\]
where the $n_\alpha$ are integers all of the same sign.
\\[-.25cm]
\end{x}

\begin{x}{\small\bf DEFINITION} \ 
The elements in a simple system of roots are said to be \un{simple}.
\\[-.25cm]
\end{x}

\begin{x}{\small\bf \un{N.B.}} \ 
Simple systems exist (cf. infra).
\\[-.25cm]
\end{x}

\begin{x}{\small\bf CONSTRUCTION} \ 
Let $\Psi$ be a simple system of roots.
\\[-.2cm]

\qquad \textbullet \quad 
The \un{positive roots} per $\Psi$ is the set
\[
\Phi^+ 
\ = \ 
\{\beta \in \Phi(\fg_\Cx) : \beta = \sum\limits_{\alpha \in \Psi} \ n_\alpha \hsx \alpha \qquad (n_\alpha \in \Z_{\geq 0})\}.
\]

\qquad \textbullet \quad 
The \un{negative roots} per $\Psi$ is the set
\[
\Phi^- 
\ = \ 
\{\beta \in \Phi(\fg_\Cx) : \beta = \sum\limits_{\alpha \in \Psi} \ n_\alpha \hsx \alpha \qquad (n_\alpha \in \Z_{\leq 0})\}.
\]

Accordingly, 
\[
\Phi(\fg_\Cx)
\ = \ 
\Phi^+ \hsx \coprod \hsx \Phi^-.
\]
\\
\end{x}


\begin{x}{\small\bf DEFINITION} \ 
\\[-.2cm]

\qquad \textbullet \quad 
The connected components of 
\[
(\sqrt{-1} \ \ft)^* \ - \bigcup\limits_{\alpha \in \Phi(\fg_\Cx)} \hsx \alpha^\perp
\]
are called the \un{Weyl chambers} of $(\sqrt{-1} \ \ft)^*$.
\\[-.2cm]

\qquad \textbullet \quad 
The connected components of 
\[
\sqrt{-1} \ \ft  \ - \bigcup\limits_{\alpha \in \Phi(\fg_\Cx)} \hsx h_\alpha^\perp
\]
are called the \un{Weyl chambers} of $\sqrt{-1} \ \ft$.
\\[-.25cm]
\end{x}

\begin{x}{\small\bf DEFINITION} \ 
\\[-.2cm]

\qquad \textbullet \quad 
If $C \subset (\sqrt{-1} \ \ft)^*$ is a Weyl chamber, then $\alpha \in \Phi(\fg_\Cx)$ is said to be 
\un{\mC-positive} if $B(C,\alpha) > 0$ and \un{\mC-negative} if $B(C,\alpha) < 0$.
\\[-.2cm]

\qquad \textbullet \quad 
If $C \subset \sqrt{-1} \ \ft$ is a Weyl chamber, then $\alpha \in \Phi(\fg_\Cx)$ is said to be 
\un{\mC-positive} if $B(C,h_\alpha) > 0$ and \un{\mC-negative} if $B(C,h_\alpha) < 0$.
\\[-.25cm]
\end{x}

\begin{x}{\small\bf DEFINITION} \ 
\\[-.2cm]

\qquad \textbullet \quad 
If $C \subset (\sqrt{-1} \ \ft)^*$ is a Weyl chamber and if $\alpha$ is \mC-positive, then $\alpha$ is 
\un{decomposable}w.r.t. \mC if there exists $\beta$, $\gamma \in \Phi(\fg_\Cx)$ such that $\alpha = \beta + \gamma$ 
(otherwise, $\alpha$ is \un{indecomposable} w.r.t. $C$).
\\[-.2cm]

\qquad \textbullet \quad 
If $C \subset \sqrt{-1} \ \ft$ is a Weyl chamber and if $\alpha$ is \mC-positive, then $\alpha$ is 
\un{decomposable} w.r.t. \mC if there exists $\beta$, $\gamma \in \Phi(\fg_\Cx)$ such that $\alpha = \beta + \gamma$ 
(otherwise, $\alpha$ is \un{indecomposable} w.r.t. $C$).
\\[-.25cm]
\end{x}

\begin{x}{\small\bf NOTATION} \ 
\\[-.25cm]

\qquad \textbullet \quad 
Given a Weyl chamber $C \subset (\sqrt{-1} \ \ft)^*$, let $\Psi(C)$ be the subset of $\Phi(\fg_\Cx)$
comprised of those $\alpha$ which are \mC-positive and indecomposable.
\\[-.2cm]

\qquad \textbullet \quad 
Given a Weyl chamber $C \subset \sqrt{-1} \ \ft$, let $\Psi(C)$ be the subset of $\Phi(\fg_\Cx)$
comprised of those $\alpha$ which are \mC-positive and indecomposable.
\\[-.25cm]
\end{x}

\begin{x}{\small\bf LEMMA} \ 
In either case, $\Psi(C)$ is a simple system of roots.
\\[-.25cm]
\end{x}

\begin{x}{\small\bf NOTATION} \ 
\\[-.2cm]

\qquad \textbullet \quad
Given a simple system of roots $\Psi$, let 
\[
C(\Psi) 
\ = \ 
\{\lambda \in (\sqrt{-1} \ \ft)^* : B(\lambda, \alpha) > 0 \ \forall \ \alpha \in \Psi\}.
\]

\qquad \textbullet \quad
Given a simple system of roots $\Psi$, let 
\[
C(\Psi) 
\ = \ 
\{H \in \sqrt{-1} \ \ft : B(H, h_\alpha) > 0 \ \forall \ \alpha \in \Psi\}.
\]
\\[-1.25cm]
\end{x}

\begin{x}{\small\bf LEMMA} \ 
In either case, $C(\Psi) $ is a Weyl chamber.
\\[-.25cm]
\end{x}

\begin{x}{\small\bf THEOREM} \ 
\\[-.2cm]

\qquad \textbullet \quad 
There is a one-to-one correspondence between the simple systems of roots and the Weyl chambers of $(\sqrt{-1} \ \ft)^*$:
\[
\begin{cases}
\ \Psi \ra C(\Psi) \\[3pt]
\ C \ra \Psi(C)
\end{cases}
.
\]

\qquad \textbullet \quad 
There is a one-to-one correspondence between the simple systems of roots and the Weyl chambers of $\sqrt{-1} \ \ft$:
\[
\begin{cases}
\ \Psi \ra C(\Psi) \\[3pt]
\ C \ra \Psi(C)
\end{cases}
.
\]
\\[-1.25cm]
\end{x}


The Weyl group $W = N(\Torus)/\Torus$ operates via $\Ad$ on $\sqrt{-1} \ \ft$ and $(\sqrt{-1} \ \ft)^*$.
\\[-.2cm]

\begin{x}{\small\bf LEMMA} \ 
The action of \mW on $\sqrt{-1} \ \ft$ and $(\sqrt{-1} \ \ft)^*$ is faithful, i.e., $w \in W$ acts trivially iff $w$ 
is the identity element.
\\[-.25cm]

PROOF \ 
Suppose that $\Ad(n)$ $(n \in \N)$ is the identity element on $\ft$ and consider the commutative diagram
\[
\begin{tikzcd}[sep=huge]
{\ft} 
\ar{d}[swap]{\text{\small{exp \ }}}  
\ar{rr}{\text{\small{$\Ad(n)$}}}
&&{\ft} 
\ar{d}{\text{\small{\ exp}}}
\\
{\Torus} 
\ar{rr}[swap]{\text{\small{$I_n$}}}
&&{\Torus}
\end{tikzcd}
.
\]
Then
\[
\exp \ft 
\ = \ 
\Torus
\]
and $\forall \ X \in \ft$, 
\[
I_n(\exp X) 
\ = \ 
n(\exp X) n ^{-1}
\ = \ 
\exp(\Ad(n)X) 
\ = \ 
\exp X.
\]
Therefore $n$ centralizes $\Torus$, hence $n \in \Torus$ (cf. $\S4$, $\#7$), i.e., $n$ represents the identity element of \mW.
\end{x}
\vspace{0.3cm}

\begin{x}{\small\bf LEMMA} \ 
\mW preserves $\Phi(\fg_\Cx)$ and $w h_\alpha = h_{w \alpha}$ $(w \in W)$.
\end{x}
\vspace{0.3cm}

\begin{x}{\small\bf NOTATION} \ 
\\[-.2cm]

\qquad \textbullet \quad 
Given $\alpha \in \Phi(\fg_\Cx)$, define
\[
r_\alpha : (\sqrt{-1} \ \ft)^* \ra  (\sqrt{-1} \ \ft)^*
\]
by
\begin{align*}
r_\alpha(\lambda) \ 
&= \ 
\lambda - 2 \ \frac{B(\lambda,\alpha)}{B(\alpha,\alpha)} \ \alpha 
\\[8pt]
&= \ 
\lambda - \lambda(h_\alpha) \alpha.
\end{align*}


\qquad \textbullet \quad 
Given $\alpha \in \Phi(\fg_\Cx)$, define
\[
r_{h_\alpha} : \sqrt{-1} \ \ft \ra  \sqrt{-1} \ \ft
\]
by
\begin{align*}
r_{h_\alpha}(H) \
&= \ 
h - 2 \ \frac{B(H,h_\alpha)}{B(h_\alpha,h_\alpha)} \ h_\alpha 
\\[8pt]
&= \ 
H - \alpha(H) h_\alpha.
\end{align*}

[Note: \ 
Geometrically, $r_\alpha$ is the reflection of $(\sqrt{-1} \ \ft)^*$ across the hyperplane perpendicular to $\alpha$ 
and $r_{h_\alpha}$ is the reflection of $\sqrt{-1} \ \ft$ across the hyperplane perpendicular to $h_\alpha$.]
\end{x}
\vspace{0.3cm}

\begin{x}{\small\bf NOTATION} \ 
Depending on the context, $W(\Phi(\fg_\Cx))$ is the group generated by
\[
\{r_\alpha : \alpha \in \Phi(\fg_\Cx)\} 
\quad \text{or} \quad 
\{r_{h_\alpha} : \alpha \in \Phi(\fg_\Cx)\}.
\]
\\[-1.25cm]
\end{x}

\begin{x}{\small\bf \un{N.B.}} \ 
$W(\Phi(\fg_\Cx))$ operates on $\ft^*$ and $\ft$ (extension by complex linearity).
\\[-.25cm]
\end{x}

\begin{x}{\small\bf LEMMA} \ 
$\forall \ \alpha \in \Phi(\fg_\Cx)$, $\exists \ n_\alpha \in N(\Torus)$ such that the action of $n_\alpha$ on  $(\sqrt{-1} \ \ft)^*$ 
is given by $r_\alpha$ and the action of $n_\alpha$ on $\sqrt{-1} \ \ft$ is given by $r_{h_\alpha}$.
\\[-.25cm]
\end{x}

\begin{x}{\small\bf THEOREM} \ 
\\[-.2cm]

\qquad \textbullet \quad 
Per $(\sqrt{-1} \ \ft)^*$, $W \approx W(\Phi(\fg_\Cx))$.
\\[-.2cm]

\qquad \textbullet \quad 
Per $\sqrt{-1} \ \ft$, $W \approx W(\Phi(\fg_\Cx))$.
\\[-.2cm]

[Note: \ 
It follows from $\#18$ that in either case,
\[
W(\Phi(\fg_\Cx)) \subset W,
\]
so the crux is the reversal of this.]
\\[-.25cm]
\end{x}


\begin{x}{\small\bf LEMMA} \ 
\mW operates simply transitively on the set of Weyl chambers in $(\sqrt{-1} \ \ft)^*$ or $\sqrt{-1} \ \ft$.

[Note: \ 
In other words, there is exactly one element of the Weyl group mapping a given Weyl chamber onto another one.]
\\[-.25cm]
\end{x}

\begin{x}{\small\bf \un{N.B.}} \ 
It is a corollary that $\abs{W}$ is the cardinality of the set of Weyl chambers.
\\[-.25cm]
\end{x}

\begin{x}{\small\bf EXAMPLE} \ 
Given a Weyl chamber \mC (be it in $(\sqrt{-1} \ \ft)^*$ or $\sqrt{-1} \ \ft$), there exists a unique element 
$w^{\hsx \stickfigure} \in W$ 
which maps \mC to its negative $-C$, hence 
$w^{\hsx \stickfigure}  \Psi(C) = \ -\Psi(C)$. 
\\[-.5cm]

[Note: \ 
In general, $-e \notin W$.]
\\[-.25cm]
\end{x}

\begin{x}{\small\bf THEOREM} \ 
Let
\[
C \subset (\sqrt{-1} \ \ft)^* 
\quad \text{or} \quad 
C \subset \sqrt{-1} \ \ft
\]
be a Weyl chamber $-$then its closure $\ov{C}$ is a fundamental domain for the action of \mW, i.e., 
$\ov{C}$ meets each \mW orbit exactly once.
\\[-.25cm]
\end{x}

Fix a Weyl chamber $C \subset (\sqrt{-1} \ \ft)^*$ and thereby determine the simple system $\Psi(C)$, hence $\Phi^+$.
\\[-.2cm]

\begin{x}{\small\bf NOTATION} \ 
$W(C)$ is the subgroup of $W(\Phi(\fg_\Cx))$ generated by the 
$r_\alpha$ $(\alpha \in \Psi(C))$.
\\[-.25cm]
\end{x}

\begin{x}{\small\bf LEMMA} \ 
\[
W(C) 
\ = \ 
W(\Phi(\fg_\Cx)).
\]
\\[-1.25cm]
\end{x}

\begin{x}{\small\bf NOTATION} \ 
Given $w \in W(\Phi(\fg_\Cx))$, let $\ell(w)$ be the smallest $k$ such that 
$w$ can be factored as a product 
$r_{\alpha_1} \cdots r_{\alpha_k}$, where the $\alpha_i \in \Psi(C)$ (set $\ell(w) = 0$ if $w = e$).  

[Note: \ 
$\ell(w)$ is referred to as the \un{length} of $w$.]
\end{x}
\vspace{0.3cm}

\begin{x}{\small\bf LEMMA} \ 
$\ell(w)$  is the number of $\alpha \in \Phi^+$ such that $w \alpha \in \Phi^-$. 
\end{x}

\begin{x}{\small\bf APPLICATION} \ 
If $w \hsx \Phi^+ = \Phi^+$, then $w = e$.
\end{x}

\begin{x}{\small\bf \un{N.B.}} \ 
The assignment
\[
w \ra \det(w) = \big(-1\big)^{\ell(w)} \in \{\pm1\}
\]
is a character of \mW.
\end{x}
\vspace{0.3cm}

\begin{x}{\small\bf LEMMA} \ 
If $\lambda \in L_\wt$, then $\forall \ w \in W$, $\lambda - w\lambda \in L_\rt$.
\\[-.2cm]

PROOF \ 
This is obvious if $w = r_\alpha$ for some $\alpha \in \Psi(C)$.  
In general, $w =r_{\alpha_1} \cdots r_{\alpha_k}$ $(k = \ell(w))$ and one can write
\[
\lambda - w \lambda 
\ = \ 
(\lambda - r_k(\lambda)) + (r_k(\lambda) - r_{k-1}(r_k(\lambda)) + \cdots .
\]
Let $\alpha_1, \ldots, \alpha_\ell$ be an enumeration of the elements of $\Psi(C)$.
\\[-.25cm] 

[Note: \ 
Recall that $\ell$ is the rank of \mG or still, the dimension of $\Torus$ or still, 
the dimension of $\sqrt{-1} \ \ft$ or still, the dimension of $(\sqrt{-1} \ \ft)^*$.]
\end{x}
\vspace{0.3cm}

\begin{x}{\small\bf DEFINITION} \ 
The \un{fundamental weights} are the $\omega_i \in L_\wt$ per the prescription
\[
2 \ \frac{B(\omega_i,\alpha_j)}{B(\alpha_j,\alpha_j)}
\ = \ 
\delta_{i j} \qquad (1 \leq i, j \leq \ell).
\]

\end{x}
\vspace{0.3cm}

\begin{x}{\small\bf LEMMA} \ 
The set $\{\omega_1, \ldots, \omega_\ell\}$ is a basis for $L_\wt$.
\end{x}


\begin{x}{\small\bf DEFINITION} \ 
A weight $\lambda \in L_\wt$ is said to be \un{dominant} if $B(\lambda, \alpha) \geq 0$ for all $\alpha \in \Psi(C)$.
\end{x}

\begin{x}{\small\bf \un{N.B.}} \ 
To say that $\lambda \in L_\wt$ is dominant amounts to saying that $\lambda \in \ov{C}$ (the closure of \mC).
\end{x}
\vspace{0.3cm}

\begin{x}{\small\bf LEMMA} \ 
A weight $\lambda \in L_\wt$ is dominant iff it is a linear combination with nonnegative integral coefficients of the $\omega_i$.
\end{x}

\begin{x}{\small\bf NOTATION} \ 
Put 
\[
\rho 
\ = \ 
\frac{1}{2} \hsx \sum\limits_{\alpha \in \Phi^+} \hsx \alpha.
\]
\end{x}

\begin{x}{\small\bf \un{N.B.}} \ 
Ultimately, $\rho$ depends on the choice of \mC.
\end{x}
\vspace{0.3cm}

\begin{x}{\small\bf LEMMA} \ 
$\forall \ w \in W$, 
\[
w \rho 
\ = \ 
\rho - \sum\limits_{\alpha \in \Phi^+, w^{-1} \alpha \in \Phi^-} \hsx \alpha.
\]
\end{x}

\begin{x}{\small\bf APPLICATION} \ 
$\forall \ \alpha \in \Psi(C)$, 
\[
r_\alpha(\rho) 
\ = \ 
\rho - \alpha.
\]

[Note: \ 
$\forall \ \alpha \in \Psi(C)$, 
\[
r_\alpha(\Phi^+ - \{\alpha\}) 
\ = \ 
\Phi^+ - \{\alpha\}.]
\]
\end{x}

\begin{x}{\small\bf LEMMA} \ 
\[
\rho 
\ = \ 
\omega_1 + \cdots +\omega_\ell.
\]
PROOF \ 
Given $\alpha_i \in \Psi(C)$,
\begin{align*}
\alpha_i  \ 
&=\ 
\rho - r_{\alpha_i }(\rho)
\\[11pt]
&=\ 
\rho - \bigg(\rho - 2 \ \frac{B(\rho,\alpha_i)}{B(\alpha_i,\alpha_i)} \bigg)
\\[11pt]
&=\
2 \ \frac{B(\rho,\alpha_i)}{B(\alpha_i,\alpha_i)} \ \alpha_i
\end{align*}

\qquad $\implies$
\[
2 \ \frac{B(\rho,\alpha_i)}{B(\alpha_i,\alpha_i)} \ \alpha_i
\ = \ 
1 
\implies
\rho \in L_\wt \qquad \text{(see below)}.
\]
Now write
\[
\rho 
\ = \ 
n_1 \hsx \omega_1 + \cdots + n_\ell \hsx \omega_\ell.
\]
Then
\begin{align*}
1 \ 
&=\ 
2 \ \frac{B(\rho,\alpha_j)}{B(\alpha_j,\alpha_j)}
\\[11pt]
&=\ 
2 \ \frac{B\big(\sum\limits_i \hsx n_i \omega_i, \alpha_j\big)}{B(\alpha_j,\alpha_j)}
\\[11pt]
&=\ 
\sum\limits_i \  n_i \  2 \  \frac{B(\omega_i,\alpha_j)}{B(\alpha_j,\alpha_j)}
\\[11pt]
&=\ 
\sum\limits_i \hsx n_i \hsx \delta_{i j} 
\\[11pt]
&=\ n_i \quad \implies 1 = n_i.
\end{align*}
Therefore
\[
\rho 
\ = \ 
\omega_1 + \cdots +\omega_\ell.
\]
\end{x}
\vspace{0.3cm}

\begin{x}{\small\bf \un{N.B.}} \ 
It follows that $\rho$ is a dominant weight.
\end{x}
\vspace{0.3cm}

\[
\text{APPENDIX}
\]

{\small\bf LEMMA} \ 
Suppose that $\lambda \in (\sqrt{-1} \ \ft)^*$ has the property that
\[
2 \ \frac{B(\lambda,\alpha_i)}{B(\alpha_i,\alpha_i)} \ \in \ \Z \qquad (i = 1, \ldots, \ell).
\]
Then $\lambda \in L_\wt$. 
\\[-.25cm]

PROOF \ 
It is a question of showing that $\forall \ \alpha \in \Phi^+$, 
\[
\lambda(h_\alpha) 
\ = \ 
2 \ \frac{B(\lambda,\alpha)}{B(\alpha,\alpha)} \in \Z.
\]
To this end, let 
$\ds\alpha = \sum\limits_{i = 1}^n \hsx n_i \hsx \alpha_i \in \Phi^+$ and proceed by induction on 
$\ds\abs{\alpha} = \sum\limits_{i = 1}^\ell \hsx n_i$, the \un{level} of $\alpha$.  
The case $\abs{\alpha} = 1$ is the hypothesis, so assume that the assertion is true for all levels $< \abs{\alpha}$.  
Choose $\alpha_i$ such that $B(\alpha,\alpha_i) > 0$, hence
\[
\beta 
\ = \ 
r_{\alpha_i} (\alpha) 
\ = \ 
\alpha - 2 \ \frac{B(\alpha,\alpha_i)}{B(\alpha_i,\alpha_i)} \alpha_i
\]
is positive and has level $< \abs{\alpha}$, thus
\begin{align*}
2 \ \frac{B(\lambda,\alpha)}{B(\alpha,\alpha)} \ 
&=\ 
2 \  \frac{B(r_{\alpha_i}(\lambda),\beta)}{B(\beta,\beta)} 
\\[11pt]
&=\ 
2 \ \frac{B(\lambda,\beta)}{B(\beta,\beta)} 
- 2 \ \frac{B(\lambda,\alpha_i)}{B(\alpha_i,\alpha_i)} 
\ 2 \ \frac{B(\alpha_i,\beta)}{B(\beta,\beta)} 
\end{align*}
is an integer.

\chapter{
$\boldsymbol{\S}$\textbf{9}.\quad  DESCENT}
\setlength\parindent{2em}
\setcounter{theoremn}{0}
\renewcommand{\thepage}{B I \S9-\arabic{page}}

\qquad Let \mG be a compact connected semisimple Lie group, $\Torus \subset G$ a maximal torus, 
$\Phi(\mathfrak{g}_\Cx)$ the roots of the pair $(\mathfrak{g}_\Cx, \mathfrak{t}_\Cx)$, 
$C \subset (\sqrt{-1} \mathfrak{t})^* $ a Weyl chamber, 
$\Psi \ (= \Psi(C))$ the simple system of roots thereby determined, and 
$\Phi^+$ $(\Phi^-)$ the positive (negative) roots per $\Psi$.
\\

\begin{x}{\small\bf RAPPEL} \ 
Given a character $\chi:\Torus \ra \Sphere^1$, there is a commutative diagram
\[
\begin{tikzcd}[sep=huge]
{\ft}  
\ar{d}[swap]{\text{\small{exp \ }}}  
\ar{rr}{\text{\small{$\td_\chi$}}}
&&{\sqrt{-1} \ \R} 
\ar{d}{\text{\small{\ exp}}}
\\
{\Torus} 
\ar{rr}[swap]{\text{\small{$\chi$}}}
&&{\Sphere^1}
\end{tikzcd}
\]
and the arrow $\chi \ra \td_\chi$ implements an identification of $\widehat{\Torus}$ with the lattice
\[
\td \widehat{\Torus} 
\ \equiv \ 
\{\lambda \in (\sqrt{-1} \ \ft)^* : \restr{\lambda}{\exp^{-1}(e)} \ \subset \ 2 \pi \sqrt{-1} \ \Z\}.
\]
\\[-1.25cm]
\end{x}

\begin{x}{\small\bf \un{N.B.}}  \ 
$\td \widehat{\Torus}$ is a sublattice of $L_\wt$ and 
\[
\pi_1(G) 
\ \approx \ 
L_\wt /\td \widehat{\Torus} 
\qquad \text{(cf. $\S7$, $\#17$).}
\]

[Note: \ 
Therefore $L_\wt = \td\widehat{\Torus}$ iff \mG is simply connected.]
\\[-.25cm]
\end{x}

\begin{x}{\small\bf NOTATION} \ 
Each $\lambda \in \td \widehat{\Torus}$ determines a character $\xi_\lambda \in \td \widehat{\Torus}$ such that 
\[
\xi_\lambda (\exp H) 
\ = \ 
e^{\lambda (H)} 
\qquad (H \in \ft).
\]
\\[-1.25cm]
\end{x}

\begin{x}{\small\bf DEFINITION} \ 
A function $f:\ft \ra \Cx$ \un{descends to $\Torus$} if it factors through the exponential map, 
i.e., if $f(H + Z) = f(H)$ $\forall \ H \in \ft$ and $\forall \ Z \in \ft$ such that $\exp Z = e$.
\\[-.25cm]
\end{x}


If $f : \ft \ra \Cx$ descends to $\Torus$, then there is a function $F:\Torus \ra \Cx$ such that
\[
F(\exp H) 
\ = \ 
f(H) \
\qquad (H \in \ft).
\]

\begin{x}{\small\bf EXAMPLE} \ 
Given $\lambda \in \td \widehat{\Torus}$, the function $H \ra e^{\lambda(H)}$ descends to $\Torus$ $(F = \xi_\lambda)$.  
\\[-.25cm]
\end{x}

\begin{x}{\small\bf EXAMPLE} \ 
Put 
\[
\rho 
\ = \ 
\frac{1}{2} \hsx \sum\limits_{\alpha \in \Phi^+} \alpha 
\qquad \text{(cf. $\S8$, $\#36$).} 
\]
Then $\forall \ w \in W$,
\[
w \rho - \rho \in L_\rt  \ \subset \ \td \widehat{\Torus} 
\qquad \text{(cf. $\S8$, $\#38$)},
\]
thus the function
\[
H \ra e^{(w \rho - \rho)(H)} 
\]
descends to $\Torus$ $(F = \xi_{w \rho - \rho})$.
\\[-.25cm]
\end{x}

\begin{x}{\small\bf \un{N.B.}} \ 
It is not claimed nor is it true in general that the function $H \ra e^{\rho(H)}$ descends to $\Torus$.
\\[-.25cm]
\end{x}

\begin{x}{\small\bf DEFINITION} \ 
$\Delta : \ft \ra \Cx$ is the function
\[
H \ra \prod\limits_{\alpha \in \Phi^+} \hsx \big(e^{\alpha(H)/2} - e^{-\alpha(H)/2}\big) 
\qquad (H \in \ft).
\]

[Note: \ 
$\alpha/2$ need not belong to $L_\wt$.]
\\[-.25cm]
\end{x}

\begin{x}{\small\bf LEMMA} \ 
\[
\Delta 
\ = \ 
e^\rho \hsx 
\prod\limits_{\alpha \in \Phi^+} \hsx \big(1 - e^{-\alpha} \big).
\]
Therefore $\Delta$ descends to $\Torus$ iff $e^\rho$ descends to $\Torus$.
\\[-.25cm]
\end{x}


\begin{x}{\small\bf LEMMA} \ 
$\abs{\Delta}^2$ descends to $\Torus$.
\\[-.25cm]

PROOF \ 
$\forall \ H \in \ft$, 
\allowdisplaybreaks
\begin{align*}
\abs{\Delta(H)}^2 \ 
&=\ 
\Delta(H) \ov{\Delta(H)} 
\\[11pt]
&=\ 
e^{\rho(H)} \hsx 
\prod\limits_{\alpha \in \Phi^+} \
\big(1 - e^{-\alpha(H)}\big) 
\
\ov{e^\rho(H) \hsx \prod\limits_{\alpha \in \Phi^+} \hsx\big(1 - e^{-\alpha(H)}\big) }
\\[11pt]
&=\ 
e^{\rho(H)} \hsx 
\prod\limits_{\alpha \in \Phi^+} \
\big(1 - e^{-\alpha(H)}\big) 
\hsx 
e^{-\rho(H)}
\ 
\ov{\prod\limits_{\alpha \in \Phi^+} \hsx \big(1 - e^{-\alpha(H)}\big)}
\\[11pt]
&=\ 
\prod\limits_{\alpha \in \Phi^+} \
(1 - e^{-\alpha(H)}) 
\ 
\ov{\big(1 - e^{-\alpha(H)}\big)}
\\[11pt]
&=\ 
\prod\limits_{\alpha \in \Phi^+} \
\abs{1 - e^{-\alpha(H)}}^2,
\end{align*}
which  descends to $\Torus$.
\\[-.25cm]
\end{x}

\begin{x}{\small\bf LEMMA} \ 
$\forall \ t \in \Torus$, 
\[
\det \big(\Ad_{G/\Torus} (t^{-1}) - \tI_{G/\Torus}\big) 
\ = \ 
\abs{\Delta(t)}^2.
\]

PROOF \ 
The complexification of $\fg/\ft$ is the direct sum of the $\fg^\alpha$ on which $t \in \Torus$ 
acts by $\xi_\alpha(t)$ in the adjoint representation, so
\allowdisplaybreaks
\begin{align*}
\det(\Ad_{G/\Torus} (t^{-1}) - \tI_{G/\Torus}) \ 
&=\ 
\prod\limits_{\alpha \in \Phi(\fg_\Cx)} \hsx 
(\xi_{\alpha}(t^{-1}) - 1)
\\[11pt] 
&=\ 
\prod\limits_{\alpha \in \Phi(\fg_\Cx)} \hsx 
(1 - \xi_{-\alpha}(t))
\\[11pt] 
&=\ 
\prod\limits_{\alpha \in \Phi^+} \hsx 
\abs{1 - \xi_{-\alpha}(t)}^2.
\end{align*}

[Note: \ 
The number of roots is even.]

\end{x}
\vspace{0.3cm}

\begin{x}{\small\bf INTEGRATION FORMULA} \ 
For any continuous function $f \in C(G)$,
\[
\int_G \hsx f(x) \hsx \td_G(x)
\ = \ 
\frac{1}{\abs{W}} \hsx 
\int_\Torus \hsx \abs{\Delta (t)}^2  \ 
\int_G \hsx
f(x t x^{-1})
\hsx \td_G(x)
\hsx \td_\Torus(t) 
\qquad \text{(cf. $\S5$, $\#13$)}.
\]
\\[-1.25cm]
\end{x}

\begin{x}{\small\bf SCHOLIUM} \ 
For any continuous class function $f \in CL(G)$,
\[
\int_G \hsx f(x) \hsx \td_G(x)
\ = \ 
\frac{1}{\abs{W}} \hsx 
\int_\Torus \hsx \abs{\Delta(t)}^2 f(t) \hsx \td_\Torus(t) 
\qquad \text{(cf. $\S5$, $\#14$)}.
\]
\\[-1.25cm]
\end{x}

\begin{x}{\small\bf REMARK} \ 
Let $t \in \Torus$ $-$then $t \in \Torus^\reg$ iff
\[
\abs{\Delta(t)}^2 \ \neq \ 0
\]
or still, iff
\[
\prod\limits_{\alpha \in \Phi^+} \hsx \abs{1 - \xi_{-\alpha}(t)}^2 \ \neq \ 0.
\]
\\[-1.25cm]
\end{x}

\begin{x}{\small\bf \un{N.B.}} \ 
Let $H \in \ft$ $-$then
\[
\abs{\Delta(e^H)}^2 
\ = \ 
2^{\abs{\Phi(\fg_\Cx)}} \hsx 
\prod\limits_{\alpha \in \Phi^+} \hsx 
\sin^2 \bigg(\frac{\alpha(H)}{2 \sqrt{-1}}\bigg).
\]

[Note: \ 
Bear in mind that $\alpha(H) \in \sqrt{-1} \ \R$.]
\\[-.25cm]
\end{x}

\begin{x}{\small\bf NOTATION} \ 
Let 
\[
\Xi 
\ = \ 
\{H \in \ft: \forall \ \alpha \in \Phi(\fg_\Cx), \hsx \alpha(H) \notin 2 \pi \sqrt{-1} \ \Z\}.
\]
\\[-1.25cm]
\end{x}


\begin{x}{\small\bf LEMMA} \ 
$\Xi$ is open and dense in $\ft$.  Moreover, 
\[
\exp \Xi 
\ = \ 
\Torus^\reg.
\]
\\[-1.25cm]
\end{x}

\begin{x}{\small\bf RAPPEL} \ 
The inclusion $\Torus \ra G$ induces a bijection between the orbits of \mW in $\Torus$ and the conjugacy classes of \mG 
(cf. $\S4$, $\#16$).  
Consequently, the class functions on \mG are the ``same thing'' as \mW-invariant functions on $\Torus$.
\\[-.25cm]
\end{x}

\begin{x}{\small\bf NOTATION} \ 
Given $\lambda \in \td \widehat{\Torus}$, define $\gche_\lambda : \Xi \ra \Cx$ by setting
\[
\gche_\lambda(H) 
\ = \ 
\frac{\sum\limits_{w \in W} \hsx \det(w) e^{w(\lambda + \rho)(H)}}{\Delta(H)}
\qquad (H \in \Xi).
\]
\\[-1.25cm]
\end{x}

\begin{x}{\small\bf LEMMA} \ 
$\forall \ w \in W$, 
\[
w(\Delta) 
\ = \ 
(-1)^{\ell(w)} \Delta.
\]
\\[-1.25cm]
\end{x}

Recalling that $\det(w) = (-1)^{\ell(w)} $, (cf. $\S8$, $\#29$), it therefore follows that $\gche_\lambda$ is a \mW-invariant function on $\Xi$.

Next, $\forall \ H \in \Xi$, 
\[
\gche_\lambda(H) 
\ = \ 
\frac{\sum\limits_{w \in W} \hsx \det(w) e^{(w(\lambda + \rho) - \rho) (H)}}
{\prod\limits_{\alpha \in \Phi^+} \hsx \big(1 - e^{-\alpha(H)}\big)}.
\]
Since
\[
e^{(w(\lambda + \rho) - \rho)(H)} 
\ = \ 
e^{w \lambda(H)} e^{(w \rho - \rho)(H)},
\]
the numerator of this fraction descends to $\Torus$ (cf. $\#5$, $\#6$).  
The same also goes for the denominator which is nonzero on $\Xi$.  
Accordingly, $\gche_\lambda$ descends to a \mW-invariant function on $\Torus^\reg$, hence extends to a class function on $G^\reg$ 
(cf. $\S5$, $\#10$), denoted still by $\gche_\lambda$.
\vspace{0.3cm}


\chapter{
$\boldsymbol{\S}$\textbf{10}.\quad  CHARACTER THEORY}
\setlength\parindent{2em}
\setcounter{theoremn}{0}
\renewcommand{\thepage}{B I \S10-\arabic{page}}

\qquad Let \mG be a compact connected semisimple Lie group, $\Torus \subset G$ 
a maximal torus, and maintain the assumptions/notation of \S9.
\\[-.25cm]

\begin{x}{\small\bf THEOREM} \ 
Suppose given a $\Pi \in \widehat{G}$ $-$then there is a $\lambda_\Pi \in \td\widehat{\Torus}$ subject to 
$\lambda_\Pi + \rho \in C$ such that $\forall \ x \in G^\reg$, 
\[
\chisubPi(x) 
\ = \ 
\gche_{\lambda_\Pi} (x).
\]
\\[-1.25cm]
\end{x}

The proof proceeds by a series of lemmas.
\\[-.2cm]

\begin{x}{\small\bf NOTATION} \ 
Given $\gamma \in C$, define $A_\gamma:\ft \ra \Cx$ by
\[
A_\gamma(H) 
\ = \ 
\sum\limits_{w \in W} \hsx \det(w) e^{w \gamma} (H).
\]
\\[-1.25cm]
\end{x}

Rephrased, the claim becomes the assertion that
\[
\chisubPi(\exp H) \hsx \Delta(H) 
\ = \ 
A_{\lambda_\Pi + \rho} (H) \qquad (H \in \Xi)
\]
for some $\lambda_\Pi \in \td\widehat{\Torus}$ subject to $\lambda_\Pi + \rho \in C$.
\\

\begin{x}{\small\bf NOTATION} \ 
$\td\widehat{\Torus}(C)$ is the subset of $\td\widehat{\Torus}$ consisting of those $\lambda$ such that 
$\lambda + \rho \in C$, say $\td\widehat{\Torus}(C) = \{\lambda_k\}$.
\\[-.25cm]

[Note: \ 
 It turns out that $\td\widehat{\Torus}(C) = \td\widehat{\Torus} \cap \hsx \ov{C}$ (cf. $\#9$).]
 \\[-.25cm]
\end{x}

\begin{x}{\small\bf LEMMA} \ 
There exist integers $m_k$ such that $\forall \ H \in \Xi$, 
\[
\chisubPi(\exp H) \hsx \Delta(H) 
\ = \ 
\sum\limits_{k} \hsx m_k A_{\lambda_k + \rho} (H).
\]


[Note: \ 
The point of departure is the fact that $\restr{\chisubPi}{\Torus}$ decomposes as a finite sum
\[
\sum\limits_{\lambda \in \td\widehat{\Torus}} \ n_\lambda \hsx \xi_\lambda \qquad (n_\lambda \in \Z_{\geq 0}).]
\]
Proceeding, 
\allowdisplaybreaks
\begin{align*}
1 \ 
&=\ 
\int_G \ \abs{\chisubPi(x)}^2 \hsx \td_G(x)
\\[11pt]
&=\ 
\frac{1}{\abs{W}} \hsx 
\int_\Torus \ 
\abs{\Delta(t)}^2 \abs{\chisubPi(t)}^2 \hsx \td_\Torus(t) \qquad \text{(cf. $\S9$, $\#13$).}
\end{align*}
\\[-0.75cm]
\end{x}

\begin{x}{\small\bf \un{N.B.}} \ 
The function
\[
\big|\hsx \sum\limits_{k} \hsx m_k A_{\lambda_k + \rho }\hsx\big|^2
\]
descends to $\Torus$ (because $\abs{\Delta}^2$ descends to $\Torus$ (cf. $\S9$, $\#10$)).
\\[-.2cm]

Therefore
\[
1 
\ = \ 
\frac{1}{\abs{W}} \hsx 
\int_\Torus \hsx 
\big|\hsx \sum\limits_{k} \hsx m_k A_{\lambda_k + \rho}\hsx\big|^2
\hsx \td_\Torus(t).
\]
\\[-1.25cm]
\end{x}

\begin{x}{\small\bf LEMMA} \ 
The function 
\[
A_{\lambda_k + \rho} \ \ov{A_{\lambda_{k^\prime} + \rho}}
\ = \ 
\big(e^{-\rho} A_{\lambda_k + \rho}\big) \hsx \ov{\big(e^{-\rho} A_{\lambda_{k^\prime} + \rho}\big)}
\]
descends to $\Torus$ (cf. $\S9$, $\#6$).
\\[-.2cm]

Therefore
\allowdisplaybreaks
\begin{align*}
\frac{1}{\abs{W}} \hsx 
\int_\Torus \hsx 
A_{\lambda_k + \rho} \ \ov{A_{\lambda_{k^\prime} + \rho}}
\hsx \td_\Torus(t) \ 
&=\ 
\frac{1}{\abs{W}} \hsx 
\int_\Torus \ 
\big(e^{-\rho} A_{\lambda_k + \rho}\big) \  \ov{\big(e^{-\rho} A_{\lambda_{k^\prime} + \rho}\big)}
\hsx \td_\Torus(t) 
\\[11pt]
&=\ 
\frac{1}{\abs{W}} \ 
\sum\limits_{w, w^\prime \in W} \ 
\det(w w^\prime) \hsx
\int_\Torus \hsx 
\xi_{w(\lambda_k + \rho) - \rho} 
\ 
\xi_{-w^\prime(\lambda_{k^\prime} + \rho) + \rho}
\ \td_\Torus(t).
\end{align*}
And
\allowdisplaybreaks
\begin{align*}
\int_\Torus \ 
\xi_{w(\lambda_k + \rho) - \rho} &
\ 
\xi_{-(w^\prime(\lambda_{k^\prime} + \rho) - \rho)} 
\ \td_\Torus(t) 
\ = \ 1
\\[11pt]
&\iff 
w(\lambda_k + \rho) - \rho
\ = \ 
w^\prime(\lambda_{k^\prime} + \rho) - \rho
\\[11pt]
&\iff 
w(\lambda_k + \rho)
\ = \ 
w^\prime(\lambda_{k^\prime} + \rho)
\\[11pt]
&\iff 
w = w^\prime \quad \text{and} \quad k = k^\prime
\end{align*}
but is zero otherwise.

Therefore
\[
\frac{1}{\abs{W}} \hsx 
\int_\Torus \hsx 
A_{\lambda_k + \rho} \ \ov{A_{\lambda_k + \rho}}
\hsx \td_\Torus(t) 
\ = \ 
\begin{cases}
\ 1 \ \text{if} \ k = k^\prime \\[3pt]
\ 0 \ \text{if} \ k \neq k^\prime
\end{cases}
.
\]

Matters then reduce to the equation
\[
1 
\ = \ 
\sum\limits_k \hsx m_k^2.
\]
However, the $m_k \in \Z$, hence all but one are zero.  
Consequently, there is a $\lambda_\Pi \in \td\widehat{\Torus}$ subject to 
$\lambda_\Pi + \rho \in C$ such that $\forall \ H \in \Xi$, 
\[
\chisubPi(\exp H) \hsx \Delta(H) 
\ = \ 
\pm A_{\lambda_\Pi + \rho}(H).
\]
\\[-1.25cm]
\end{x}

\begin{x}{\small\bf LEMMA} \ 
The $A_\gamma$ $(\gamma \in C)$ are linearly independent over $\Z$.


[Given $\gamma$, $\gamma^\prime \in C$, 

\[
\langle A_\gamma, A_{\gamma^\prime} \rangle
\ = \ 
\begin{cases}
\ 1 \quad \text{if} \ \gamma = \gamma^\prime \\[3pt]
\ 0 \quad \text{if} \ \gamma \neq \gamma^\prime
\end{cases}
,
\]
the inner product $\langle \ , \ \rangle$ being by definition the multiplicity of the ``zero weight'' in 
\[
\frac{1}{\abs{W}} \hsx 
\bigg[ \sum\limits_{w \in W} \hsx \det(w) e^{w\gamma} \bigg]
\bigg[ \sum\limits_{w^\prime \in W} \hsx \det(w^\prime) e^{-w^\prime\gamma^\prime} \bigg]
\ = \ 
\frac{1}{\abs{W}} \hsx 
\bigg[ \sum\limits_{w,w^\prime \in W} \hsx \det(ww^\prime) e^{w\gamma -w^\prime\gamma^\prime}.
\]
But
\allowdisplaybreaks
\begin{align*}
w \gamma  -w^\prime\gamma^\prime = 0 
&\implies 
\gamma = w^{-1} w^\prime \gamma^\prime 
\\[11pt]
&\implies
w = w^\prime
\\[11pt]
&\implies
\gamma = \gamma^\prime, 
\end{align*}
so the number of solutions is $\abs{W}$ if $\gamma = \gamma^\prime$ and is zero otherwise.]
\\[-.25cm]
\end{x}

\begin{x}{\small\bf APPLICATION} \ 
The linear function $\lambda_\Pi + \rho \in C$ is unique.
\\[-.25cm]
\end{x}

\begin{x}{\small\bf LEMMA} \ 
Let $\lambda \in \td\widehat{\Torus}$ $-$then
\[
\lambda + \rho \in C 
\iff 
\lambda \in \ov{C}.
\]

PROOF \ 
$\forall \ \alpha_i \in \Psi(C)$, 
\[
2 \ \frac{B(\rho, \alpha_i)}{B(\alpha_i, \alpha_i)} 
\ = \ 
1 
\qquad \text{(cf. $\S8$, $\#40$)}
\]
and
\[
2 \ \frac{B(\lambda, \alpha_i)}{B(\alpha_i, \alpha_i)} \ \in \ \Z 
\qquad (\lambda \in \td\widehat{\Torus} \subset L_\wt).
\]

The stated equivalence then follows upon writing 
\allowdisplaybreaks
\begin{align*}
2 \  \frac{B(\lambda + \rho, \alpha_i)}{B(\alpha_i, \alpha_i)} \
&=\ 
2 \ \frac{B(\rho, \alpha_i)}{B(\alpha_i, \alpha_i)} 
 + 
2 \ \frac{B(\lambda, \alpha_i)}{B(\alpha_i, \alpha_i)}
\\[11pt]
&=\ 
1 + 
2 \ \frac{B(\lambda, \alpha_i)}{B(\alpha_i, \alpha_i)}.
\end{align*}
\\[-.25cm]
\end{x}

\begin{x}{\small\bf APPLICATION} \ 
\[
\lambda_\Pi + \rho \in C 
\implies 
\lambda_\Pi \in \ov{C}.
\]
\\[-1.25cm]
\end{x}

Return now to the expression 
\[
\chisubPi(\exp H) \Delta (H) 
\ = \ 
\pm A_{\lambda_\Pi + \rho} (H) 
\]
valid for $H \in \Xi$, the objective then being to establish that it is the plus sign which obtains.
\\[-.2cm]

\begin{x}{\small\bf LEMMA} \ 
$\forall \ H \in \ft$, 
\[
\Delta(H) 
\ = \ 
\sum\limits_{w \in W} \hsx 
\det(w) e^{w \rho(H)}.
\]

[Note: \ 
There is no vicious circle here in that the formula can be derived by direct (albeit somewhat tedious) manipulation, 
the derivation being independent of the preceding considerations (but consistent with the final outcome).]
\\[-.25cm]
\end{x}

From this it follows that $\forall \ H \in \Xi$, 
\allowdisplaybreaks
\begin{align*}
\pm \chisubPi(\exp H) \ 
&=\ 
\frac
{\sum\limits_{w \in W} \hsx \det(w) e^{w(\lambda_\Pi + \rho)(H)}}
{\sum\limits_{w \in W} \hsx \det(w) e^{w \rho(H)}}
\\[11pt]
&=\ 
\gche_{\lambda_\Pi} (H).
\end{align*}

\begin{x}{\small\bf NOTATION} \ 
Define $H_\rho \in \sqrt{-1} \ \ft$ by the relation
\[
\rho(H) 
\ = \ 
B(H, H_\rho) \quad (H \in \sqrt{-1} \ \ft) \qquad \text{(cf. $\S6$, $\#21$).}
\]
\\[-1.25cm]
\end{x}

\begin{x}{\small\bf LEMMA} \ 
$\sqrt{-1} \ t H_\rho \in \Xi$ for small positive $t$.
\\[-.25cm]
\end{x}

\begin{x}{\small\bf LEMMA} \ 
\[
\lim\limits_{t \downarrow 0} \hsx \chisubPi(\exp \sqrt{-1} \ t H_\rho) 
\ = \ 
\td_\Pi.
\]

[For $\restr{\chisubPi}{\Torus}$ is continuous and $\td_\Pi = \chisubPi(e)$.]
\\[-.25cm]
\end{x}

\begin{x}{\small\bf APPLICATION} \ 
\[
\pm \hsx  \td_\Pi 
\ = \ 
\lim\limits_{t \downarrow 0} \hsx \gche_{\lambda_\Pi} (\sqrt{-1} \ t H_\rho).
\]
\\[-1.25cm]
\end{x}

\begin{x}{\small\bf SUBLEMMA} \ 
$\forall \ w \in W$, 
\allowdisplaybreaks
\begin{align*}
w(\lambda_\Pi + \rho) (\sqrt{-1} \ t \hsx H_\rho)\ 
&=\ 
\sqrt{-1} \ t \hsx (\lambda_\Pi + \rho)(w^{-1} H_\rho)
\\[11pt]
&=\ 
\sqrt{-1} \ t \hsx B(H_{\lambda_\Pi + \rho}, w^{-1}H_\rho)
\\[11pt]
&=\ 
\sqrt{-1} \ t \hsx B(w H_{\lambda_\Pi + \rho}, H_\rho)
\\[11pt]
&=\ 
\sqrt{-1} \ t \rho(w H_{\lambda_\Pi + \rho})
\\[11pt]
&=\ 
(w^{-1} \rho)(\sqrt{-1} \ t \hsx H_{\lambda_\Pi + \rho}).
\end{align*}
\\[-.25cm]
\end{x}


\begin{x}{\small\bf LEMMA} \ 
\[
\lim\limits_{t \downarrow 0} \gche_{\lambda_\Pi} \hsx \big(\sqrt{-1} \ t \hsx H_\rho\big) 
\ = \ 
\frac{\prod\limits_{\alpha \in \Phi^+} B(\lambda_\Pi + \rho, \alpha)}{\prod\limits_{\alpha \in \Phi^+} B(\rho, \alpha)}.
\]

PROOF \ 
Write
\allowdisplaybreaks
\begin{align*}
\sum\limits_{w \in W} \ 
\det(w) e^{w(\lambda_\Pi + \rho)\big(\sqrt{-1} \ t \hsx H_\rho\big)} \ 
&=\ 
\sum\limits_{w \in W} \ 
\det(w) e^{(w^{-1}\rho)(\sqrt{-1} \ t \hsx H_{\lambda_\Pi +\rho})}
\\[11pt]
&=\ 
\sum\limits_{w \in W} \ 
\det(w^{-1}) e^{(w^{-1}\rho)(\sqrt{-1} \ t \hsx H_{\lambda_\Pi +\rho})}
\\[11pt]
&=\ 
\sum\limits_{w \in W} \ 
\det(w) 
e^{(w\rho)(\sqrt{-1} \ t \hsx H_{\lambda_\Pi +\rho})}
\\[11pt]
&=\ 
\Delta(\sqrt{-1} \ t \hsx H_{\lambda_\Pi +\rho})
\\[11pt]
&=\ 
\prod\limits_{\alpha \in \Phi^+} \
\bigg(
e^{\alpha(\sqrt{-1} \ t \hsx H_{\lambda_\Pi +\rho})/2}
-
e^{-\alpha(\sqrt{-1} \ t \hsx H_{\lambda_\Pi +\rho})/2}
\bigg)
\\[11pt]
&=\ 
\prod\limits_{\alpha \in \Phi^+} \hsx
\sqrt{-1} \ t \hsx \alpha \hsx (H_{\lambda_\Pi +\rho}) + o(1)
\\[11pt]
&=\ 
\big(\sqrt{-1} \ t\big)^{\abs{\Phi^+}} \hsx
\prod\limits_{\alpha \in \Phi^+} \hsx
B(\lambda_\Pi + \rho, \alpha) + o(1).
\end{align*}
\\[-.25cm]
\end{x}
Analogously, 
\[
\sum\limits_{w \in W} \hsx 
\det(w) 
e^{(w\rho)(\sqrt{-1} \ t H_\rho)}
\ = \ 
\big(\sqrt{-1} \ t\big)^{\abs{\Phi^+}} \hsx
\prod\limits_{\alpha \in \Phi^+} \hsx
B(\rho,\alpha) + o(1).
\]
Taking the limit as $t \downarrow 0$ then finishes the proof.
\\[-.2cm]

\begin{x}{\small\bf \un{N.B.}} \ 
Both $\rho$ and $\lambda_\Pi + \rho$ belong to \mC, thus $\alpha \in \Phi^+$, 
\[
B(\rho,\alpha) > 0 
\quad \text{and} \quad 
B(\lambda_\Pi + \rho,\alpha) > 0,
\]
so
\[
\lim\limits_{t \downarrow 0} \gche_{\lambda_\Pi} \hsx \big(\sqrt{-1} \ t H_\rho\big) 
\ > \ 
0.
\]
\\[-1.25cm]
\end{x}

\begin{x}{\small\bf APPLICATION} \ 
\[
\td_\Pi
\ = \ 
\lim\limits_{t \downarrow 0} \gche_{\lambda_\Pi} \hsx \big(\sqrt{-1} \ t H_\rho\big).
\]
I.e.: The plus sign prevails.
\\[-.25cm]
\end{x}

\begin{x}{\small\bf SCHOLIUM} \ 
\[
\td_\Pi 
\ = \ 
\frac
{\prod\limits_{\alpha \in \Phi^+} \hsx B(\lambda_\Pi + \rho, \alpha)}
{\prod\limits_{\alpha \in \Phi^+} \hsx B(\rho, \alpha)}.
\]
\\[-1.25cm]
\end{x}

\begin{x}{\small\bf LEMMA} \ 
The arrow from $\widehat{G}$ to $\td\widehat{T} \cap \ov{C}$ that sends $\Pi$ to $\lambda_\Pi$ is well-defined 
(cf. $\#8$) and injective.
\\[-.2cm]

PROOF \ 
Given $\Pi_1$, $\Pi_2 \in \widehat{G}$, suppose that $\lambda_{\Pi_1} = \lambda_{\Pi_2}$ $-$then 
$\lambda_{\Pi_1} + \rho = \lambda_{\Pi_2} + \rho$, hence
\[
\gche_{\lambda_{\Pi_1}} 
\ = \ 
\gche_{\lambda_{\Pi_2}},
\]
which implies that 
$\chisubPiOne = \chisubPiTwo$ on $G^\reg$ or still, by continuity, 
$\chisubPiOne = \chisubPiTwo$ on \mG, so 
$\Pi_1 = \Pi_2$.
\\[-.25cm]
\end{x}

\begin{x}{\small\bf LEMMA} \ 
The arrow from $\widehat{G}$ to $\td\widehat{T} \cap \ov{C}$ that sends $\Pi$ to $\lambda_\Pi$ is surjective.
\\[-.2cm]

PROOF \ 
Fix a $\lambda \in \td\widehat{T} \cap \ov{C}$ $-$then
\allowdisplaybreaks
\begin{align*}
\int_G \hsx 
\abs{\gche_\lambda (x)}^2 \td_G(x) \ 
&=\ 
\frac{1}{\abs{W}} \hsx \int_{\Torus^\reg} \
\abs{\Delta(t)}^2 \ \abs{\gche_\lambda (t)}^2
\ 
\td_\Torus(t)
\\[11pt]
&=\ 
\frac{1}{\abs{W}} \hsx \int_\Torus \
A_{\lambda + \rho} \ \ov{A_{\lambda + \rho}} 
\ 
\td_\Torus(t)
\\[11pt]
&=\ 
1.
\end{align*}
Therefore $\gche_\lambda$ is an $L^2$ class function (cf. $\S2$, $\#17$).  
Now fix a  $\Pi_0 \in \widehat{G}$: 
\\[-1.25cm]

\allowdisplaybreaks
\begin{align*}
\langle \gche_\lambda, \chisubPiZero \rangle \ 
&=\ 
\int_G \ 
\gche_\lambda (x) \ \ov{\chisubPiZero(x)} 
\
\td_G(x)
\\[11pt]
&=\ 
\frac{1}{\abs{W}} \hsx 
\int_{\Torus^\reg} \
\abs{\Delta(t)}^2 
\gche_\lambda (t) \ \ov{\chisubPiZero(t)} 
\
\td_\Torus(t)
\\[11pt]
&=\ 
\frac{1}{\abs{W}} \hsx 
\int_{\Torus^\reg} \
\abs{\Delta(t)}^2 
\gche_\lambda (t) \ \ov{\gche_{\lambda_{\Pi_0}}(t)} 
\
\td_\Torus(t)
\\[11pt]
&=\ 
\begin{cases}
1 \quad \text{if} \quad \lambda_{\Pi_0} = \lambda \\
0 \quad \text{if} \quad \lambda_{\Pi_0} \neq \lambda
\end{cases}
\end{align*}

\qquad $\implies$
\allowdisplaybreaks
\begin{align*}
\gche_\lambda \ 
&=\ 
\sum\limits_{\Pi \in \widehat{G}} \hsx 
\langle \gche_\lambda, \chisubPi \rangle_{\chisubPi} \qquad \text{(cf. $\S2$, $\#19$)}
\\[11pt]
&=\ 
\chisubPi
\end{align*}
for a unique $\Pi \in \widehat{G}$ with $\lambda_\Pi = \lambda$.
\\[-.25cm]
\end{x}


\begin{x}{\small\bf SCHOLIUM} \ 
\[
\begin{cases}
\ \widehat{G} \longleftrightarrow \td\widehat{T} \cap \ov{C} \\[3pt]
\ \Pi \longleftrightarrow \lambda_\Pi
\end{cases}
.
\]

[Note: \ 
\mW operates on $\td\widehat{T}$ and $\td\widehat{T} \cap \ov{C}$ is a fundamental domain for this action 
(cf. $\S8$, $\#23$), hence $\widehat{G}$ is parameterized by the orbits of \mW in $\td\widehat{\Torus}$.]
\\[-.25cm]
\end{x}

\begin{x}{\small\bf \un{N.B.}} \ 
\[
\lambda_{\Pi^*} 
\ = \ 
-w ^{\stickfigure} \lambda_\Pi 
\qquad \text{(cf. $\S8$, $\#22$).}
\]
\\[-1.25cm]
\end{x}

\begin{x}{\small\bf REMARK} \ 
It is clear that if $\lambda_\Pi = 1_G$, then $\lambda_\Pi = 0$. 
\\[-.25cm]
\end{x}

\begin{x}{\small\bf LEMMA} \ 
In the restriction of $\chisubPi$ to $\Torus$, $\xi_{\lambda_\Pi}$ occurs with multiplicity 1.
\\[-0.75cm]
\end{x}

\[
\text{APPENDIX}
\]

There are two directions in which the theory can be extended.
\\[-.2cm]

\qquad \textbullet \quad 
Drop the assumption that \mG is semisimple and work with an arbitrary compact connected Lie group.
\\[-.2cm]

\qquad \textbullet \quad 
Drop the assumption that \mG is connected and work with an arbitrary compact Lie group.
\\[-.2cm]

As regards the first point, no essential difficulties are encountered.  
As regards the second point, however, there are definitely some subtleties (see Chapter 1 of D. Vogan's book 
``Unitary Representations of Reductive Lie Groups''.)
\\[-.2cm]

{\small\bf NOTATION} \ 
Let \mG be a compact semisimple Lie group, $\Torus \subset G$ a maximal torus, 
$C \subset \sqrt{-1} \ \ft$ a Weyl chamber and let
\[
N_G(C) 
\ = \ 
\{x \in G: \Ad(x) C \subset C\}.
\]
\\[-1.25cm]

{\small\bf LEMMA} \ 
\[
N_G(C) \cap G^0 
\ = \ 
\Torus, 
\quad 
N_G(C) G^0 
\ = \ 
G.
\]
Therefore
\[
G/G^0 
\ \approx \ 
N_G(C) /\Torus.
\]
\\[-1.25cm]

{\small\bf \un{N.B.}} \ 
Each element of \mG is conjugate to an element of $N_G(C)$.

\chapter{
$\boldsymbol{\S}$\textbf{11}.\quad  THE INVARIANT INTEGRAL}
\setlength\parindent{2em}
\setcounter{theoremn}{0}
\renewcommand{\thepage}{B I \S11-\arabic{page}}

\qquad Let \mG be a compact connected semisimple Lie group, $T \subset G$ a maximal torus etc.
\\[-.25cm]

\begin{x}{\small\bf NOTATION} \ 
Set
\[
\pi 
\ = \ 
\prod\limits_{\alpha \in \Phi^+} \hsx \alpha
\ \equiv \
\prod\limits_{\alpha > 0} \hsx \alpha.
\]
\\[-1.25cm]
\end{x}

\begin{x}{\small\bf LEMMA} \ 
$\pi$ is a homogeneous polynomial of degree $r$ $(= \abs{\Phi^+})$ and $\forall \ w \in W$, 
\[
w \hsx \pi 
\ = \ 
\det(w) \hsx \pi.
\]
\\[-1.25cm]
\end{x}

\begin{x}{\small\bf LEMMA} \ 
If $p$ is a homogeneous polynomial such that $\forall \ w \in W$, 
\[
w \hsx p 
\ = \ 
\det(w) \hsx p,
\]
then $p$ can be written as $\pi P$, where \mP is a homogeneous \mW-invariant polynomial.
\end{x}
\vspace{0.3cm}

\begin{x}{\small\bf \un{N.B.}} \ 
$P \hsx = \hsx 0$ if $\deg p \hsx < \hsx r$ and $P  \hsx= \hsx C$ (a constant) if $\deg p \hsx = \hsx r$.
\\[-.25cm]
\end{x}

\begin{x}{\small\bf DEFINITION} \ 
Given $f \in C^\infty(\fg)$ and $H \in \ft$, put
\[
\phi_f(H) 
\ = \ 
\pi(H) \  
\int_G \  f(\Ad(x) \hsx H) \ \td_G(x),
\]
the \un{invariant integral} of $f$ at \mH.
\\[-.25cm]
\end{x}

\begin{x}{\small\bf FUNCTIONAL EQUATION} \ 
$\forall \ w \in W$ $(w = n \hsx \Torus)$, 
\allowdisplaybreaks
\begin{align*}
\phi_f(w H) \ 
&=\ 
\pi(w H) \
\int_G \ f(\Ad(x) \hsx w \hsx H) 
\ 
\td_G(x)
\\[11pt]
&=\ 
\det(w)\hsx \pi(H) \ \int_G \
f(\Ad(x) \Ad(n) \hsx H)
\ 
\td_G(x)
\\[11pt]
&=\ 
\det(w) \pi(H) \ \int_G \
f(\Ad(x n) \hsx H)
\ 
\td_G(x)
\\[11pt]
&=\ 
\det(w) \pi(H) \ \int_G \
f(\Ad(x) H)
\ 
\td_G(x)
\\[11pt]
&=\ 
\det(w) \hsx \phi_f(H).
\end{align*}
\\[-.75cm]
\end{x}

\begin{x}{\small\bf LEMMA} \ 
\[
f \in C^\infty(\fg) 
\implies 
\phi_f \in C^\infty(\ft).
\]
\\[-1.25cm]
\end{x}

\begin{x}{\small\bf LEMMA} \ 
\[
f \in C_c^\infty(\fg) 
\implies 
\phi_f \in C_c^\infty(\ft).
\]
\\[-1.25cm]
\end{x}

\begin{x}{\small\bf LEMMA} \ 
\[
f \in C(\fg) 
\implies 
\phi_f \in C(\ft).
\]

[If $D \in \fP(\fg)$ is a polynomial differential operator, then there exists a finite number of elements 
$D_1, \ldots, D_p \in \fP(\fg)$ and analytic functions $a_1, \ldots, a_p$ on \mG such that $\forall \ x \in G$, 
\[
\Ad(x) \hsx D 
\ = \ 
\sum\limits_{i = 1}^p \hsx a_i(x) \hsx D_i.
\]

[Note: \
An automorphism of $\fg$ extends to an automorphism of $\fP(\fg)$.]
\end{x}
\vspace{0.3cm}

\begin{x}{\small\bf NOTATION} \ 
Set 
\[
\widetilde{\pi} 
\ = \ 
\prod\limits_{\alpha \in \Phi^+} \ H_\alpha
\ \equiv \ 
\prod\limits_{\alpha > 0} \ H_\alpha.
\]
\\[-1.25cm]
\end{x}

\begin{x}{\small\bf \un{N.B.}} \ 
$\partial(\widetilde{\pi}) \hsx (\pi) $ is a constant (explicated infra).

[The point is that $\pi$ is a homogeneous polynomial of degree $r$ and $\partial(\widetilde{\pi})$ 
is a polynomial differential operator of degree $r$.]
\\[-.25cm]
\end{x}

\begin{x}{\small\bf RAPPEL} \ 
For the record,
\[
\partial H_\alpha(f)|_H
\ = \ 
\frac{\td}{\td t} \ f(H + t H_\alpha)|_{t = 0}.
\]
[In particular, if $f$ is linear, then 
\[
\partial H_\alpha(f)|_H
\ = \ 
f(H_\alpha),
\]
a constant.]
\\[-.2cm]

Put
\[
F(H) 
\ = \ 
\int_G \ f(\Ad (x) \hsx H) \  \td_G(x) \qquad (H \in \ft).
\]
Then
\allowdisplaybreaks
\begin{align*}
(\partial(\widetilde{\pi}) \circ \pi) F|_{H = 0} \ 
&=\ 
F(H; \partial(\widetilde{\pi}) \circ \pi)|_{H = 0}
\\[11pt]
&=\ 
\partial(\widetilde{\pi}) \hsx (\pi) \hsx F(0)
\\[11pt]
&=\ 
\partial(\widetilde{\pi})\hsx (\pi) \hsx f(0).
\end{align*}
\\[-.75cm]
\end{x}

\begin{x}{\small\bf THEOREM} \ 
\[
\partial(\widetilde{\pi})(\pi)
\ = \ 
\abs{W} \ \prod\limits_{\alpha > 0} \ B(\rho, \alpha).
\]
\\[-.2cm]

PROOF \ 
The sum 
\[
\sum\limits_{w \in W} \ \det(w) \hsx (w \rho)^k
\]
is a homogeneous polynomial of degree $k$ which transforms according to the determinant per the action of \mW, 
hence vanishes if $0 \leq k < r$ but if $k = r$, 
\[
\frac{1}{r!} \hsx 
\sum\limits_{w \in W} \ \det(w) \hsx (w \rho)^k
\ = \ 
C(\rho) \pi
\]
for some constant $C(\rho)$ (cf. $\#4$).  
To calculate $C(\rho)$, note that $\rho^r$ is a homogeneous polynomial of degree $r$, thus 
$\partial(\widetilde{\pi})(\rho)^r$ is a constant, so
\allowdisplaybreaks
\begin{align*}
\partial(\widetilde{\pi})(\det(w) \hsx (w \rho)^r) \ 
&=\ 
w(\partial(\widetilde{\pi})) \hsx (\rho^r) 
\\[11pt]
&=\ 
\partial(\widetilde{\pi}) \hsx (\rho)^r
\\[11pt]
&=\ 
\prod\limits_{\alpha > 0} \ \partial(H_\alpha) (\rho)^r
\\[11pt]
&=\ 
r! \ \prod\limits_{\alpha > 0} \ B(\rho,\alpha).
\end{align*}
Therefore, on the one hand, 
\allowdisplaybreaks
\begin{align*}
\partial(\widetilde{\pi}) \ \bigg(\frac{1}{r!} \ 
\sum\limits_{w \in W} \  \det(w) \hsx (w \rho)^r\bigg) \ 
&=\ 
\frac{1}{r!} \
\sum\limits_{w \in W} \ 
\partial(\widetilde{\pi}) (\det(w)\hsx  (w \rho)^r)
\\[11pt]
&=\ 
\frac{1}{r!} \
\abs{W} r! \hsx \hsx \prod\limits_{\alpha > 0} \ B(\rho,\alpha)
\\[11pt]
&=\ 
\abs{W} \ \prod\limits_{\alpha > 0} \ B(\rho,\alpha),
\end{align*}
while on the other
\[
\partial(\widetilde{\pi}) \bigg(\frac{1}{r!} \ 
\sum\limits_{w \in W} \ \det(w) \hsx (w \rho)^r\bigg) 
\ = \ 
C(\rho) \hsx \partial(\widetilde{\pi}) \hsx (\pi).
\]
Consequently, 
\[
\abs{W} \hsx \prod\limits_{\alpha > 0} \ B(\rho, \alpha) 
\ = \ 
C(\rho) \hsx \partial(\widetilde{\pi}) \hsx (\pi)
\]
\qquad $\implies$
\[
\frct{1}{r!} \hsx \sum\limits_{w \in W} \ \det(w) \hsx (w \rho)^r 
\ = \ 
\frct{\abs{W} \ \prod\limits_{\alpha > 0} \ 
B(\rho, \alpha) }{\partial(\widetilde{\pi})(\pi)} \ \pi.
\]
Let $H = \sqrt{-1} \  t H_\rho$ (cf. $\S10$, $\#13$) and write 
$\lim\limits_{H \ra 0}$ in place of $\lim\limits_{t \downarrow 0}$:
\allowdisplaybreaks
\begin{align*}
1 \ 
&=\ 
\frct{\Delta(H)}{\Delta(H)} 
\\[11pt]
&=\ 
\frct
{\sum\limits_{w \in W} \ \det(w) \hsx e^{w \rho(H)}}
{\prod\limits_{\alpha > 0} \ \big(e^{\alpha(H)/2}  - e^{-\alpha(H)/2}\big)}
\\[11pt]
&=\ 
\lim\limits_{H \ra \hsy 0} 
\frct
{\sum\limits_{w \in W} \ \det(w) \hsx e^{w \rho(H)}}
{\prod\limits_{\alpha > 0} \ \big(e^{\alpha(H)/2}  - e^{-\alpha(H)/2}\big)}
\\[11pt]
&=\ 
\lim\limits_{H \ra \hsy 0} 
\frct
{\sum\limits_{w \in W} \ \det(w) \hsx e^{w \rho(H)}}
{e^{-\rho(H)} \hsx \prod\limits_{\alpha > 0} \ \big(e^{\alpha(H)} - 1\big)}
\\[11pt]
&=\ 
\lim\limits_{H \ra \hsy 0} \ 
\bigg[
\frct{e^{\rho(H)}}{\prod\limits_{\alpha > 0} \hsx \frac{e^{\alpha(H) - 1}}{\alpha(H)}}
\hsx \times \hsx
\frct{\sum\limits_{w \in W} \hsx  \det(w) \hsx e^{w \rho(H)}}{\pi(H)}
\bigg]
\\[11pt]
&=\ 
\lim\limits_{H \ra \hsy 0} \ 
\frct
{\sum\limits_{w \in W} \ \det(w) \hsx e^{w \rho(H)}}
{\pi(H)}
\end{align*}
which upon expansion of the exponentials equals
\[
\lim\limits_{H \ra 0}\  (C(\rho) + o(1) 
\ = \ 
C(\rho)
\]
\qquad $\implies$
\[
1 
\ = \ 
C(\rho)
\ = \ 
\frct{\abs{W} \ \prod\limits_{\alpha > 0} \ B(\rho, \alpha)}
{\raisebox{-.1cm}{$\partial(\widetilde{\pi})(\pi)$}}
\]
\qquad $\implies$
\[
\partial(\widetilde{\pi}) \hsx (\pi)
\ = \ 
\abs{W} \hsx \prod\limits_{\alpha > 0} \ B(\rho, \alpha).
\]
\\[-1.25cm]
\end{x}

\begin{x}{\small\bf APPLICATION} \ 
Given $f \in C^\infty(\fg)$, 
\allowdisplaybreaks
\begin{align*}
\phi_f(0; \partial(\widetilde{\pi})) \ 
&=\ 
(\partial(\widetilde{\pi}) \hsx \phi_f)(0)
\\[11pt]
&=\ 
\big(\abs{W} \ \prod\limits_{\alpha > 0} \ B(\rho, \alpha) \big) \hsx f(0).
\end{align*}
\end{x}


\chapter{
$\boldsymbol{\S}$\textbf{12}.\quad  PLANCHEREL}
\setlength\parindent{2em}
\setcounter{theoremn}{0}
\renewcommand{\thepage}{B I \S12-\arabic{page}}

\qquad Keeping to the overall setup of $\S11$, assume in addition that \mG is simply connected, so
\[
L_\wt 
\ = \ 
\td \widehat{\Torus} 
\qquad \text{(cf. $\S7$, $\#17$)}
\]
and $e^\rho$ descends to $\Torus$, so does $\Delta$, thus
\[
\Delta(t) 
\ = \ 
\xi_\rho(t) \ 
\prod\limits_{\alpha > 0} \  (1 - \xi_\alpha(t^{-1})) 
\qquad (t \in \Torus).
\]
\\[-1.25cm]

\begin{x}{\small\bf NOTATION} \ 
Put
\[
\sW 
\ = \ 
L_\wt, 
\quad 
\sW^+ 
\ = \ 
L_\wt \cap \ov{C}.
\]
\\[-1.25cm]
\end{x}

\begin{x}{\small\bf \un{N.B.}} \ 
The elements of $\sW^+$ are the dominant weights (cf. $\S8$, $\#34$).
\\[-.25cm]
\end{x}

\begin{x}{\small\bf NOTATION} \ 
Given $\Lambda \in \sW^+$, $\Pi_\Lambda$ is the irreducible unitary representation of \mG associated with $\Lambda$, $\chisubLambda$ 
its character, 
\[
\td_\Lambda 
\ = \ 
\frct{\prod\limits_{\alpha > 0} \hsx B(\Lambda + \rho, \alpha)}{\prod\limits_{\alpha > 0} \  B(\rho, \alpha)}
\]
its dimension (cf. $\S10$, $\#20$).
\\[-.25cm]
\end{x}

\begin{x}{\small\bf \un{N.B.}}\ 
On $\Torus^\reg$, 
\[
\chisubLambda(t) \Delta(t) 
\ = \ 
\sum\limits_{w \in W} \  \det(w) \hsx \xi_{w(\Lambda + \rho)} (t).
\]
It is wellknown that
\[
C^\infty(G) * C^\infty(G) \  \subset \  C^\infty(G),
\]
so on the basis of $\S2$, $\#15$, the Plancherel theorem is in force:
\allowdisplaybreaks
\begin{align*}
f(e) \ 
&=\ 
\sum\limits_{\Pi \in \widehat{G}} \  \td_\Pi \tr(\Pi(f))
\\[11pt]
&=\ 
\sum\limits_{\Pi \in \widehat{G}} \  \td_\Pi 
\ 
\int_G \hsx f(x) \hsy \chisubPi(x) \hsx \td_G(x)
\end{align*}
or still, 
\[
f(e)
\ = \ 
\sum\limits_{\Lambda \in \sW^+} \  \td_\Lambda 
\ 
\int_G \  f(x) \hsy \chisubLambda(x) \hsx \td_G(x).
\]
Our objective now will be to give another proof of this relation which is independent of the factorization theory for 
$C^\infty(G)$ but hinges instead on the result formulated in $\S11$, $\#13$.
\\[-.25cm]
\end{x}

\begin{x}{\small\bf NOTATION} \ 
Given $f \in C^\infty(G)$ and $t \in \Torus$, put
\[
F_f (t)
\ = \ 
\Delta(t) \hsx
\int_G \  f(x \hsy t \hsy x^{-1}) \hsx \td_G(x),
\]
the \un{invariant integral} of $f$ at $t$.
\\[-.25cm]
\end{x}

\begin{x}{\small\bf LEMMA} \ 
\[
F_f \hsx \in \hsx C^\infty(\Torus).
\]
\\\\[-1.25cm]

Owing to $\S9$, $\#12$, 
\[
\int_G \  f(x) \hsx \td_G(x)
\ = \ 
\frct{1}{\abs{W}} \  
\int_\Torus \  \abs{\Delta(t)}^2 \hsx
\int_G \  f(x \hsy  t \hsy x^{-1}) \hsx \td_G(x)
\hsx \td_\Torus(t)
\]
which equals
\[
\frct{1}{\abs{W}} \  
\int_\Torus \  \ov{\Delta(t)} \hsx \Delta(t) 
\int_G \  f(x \hsy t \hsy x^{-1}) \hsx \td_G(x)
\hsx \td_\Torus(t)
\] 
or still, 
\[
\frct{1}{\abs{W}} \  
\int_\Torus \  \ov{\Delta(t)} \hsx F_f(t) 
\hsx \td_\Torus(t)
\]
or still, 
\[
\frct{(-1)^r}{\abs{W}} \hsx
\int_\Torus \  \Delta(t) \hsy F_f(t) \hsx \td_\Torus(t) \qquad (r = \abs{\Phi^+}).
\]

Therefore
\[
\sum\limits_{\Lambda \in \sW^+} \td_\Lambda \hsx
\int_G \hsx f(x) \chisubLambda(x) \td_G(x) 
\ = \ 
\frac{(-1)^r}{\abs{W}} \ 
\sum\limits_{\Lambda \in \sW^+} \ \td_\Lambda 
\int_\Torus \  \Delta(t) \hsx  \chisubLambda(t)  \hsx  F_f(t) \hsx \td_\Torus(t).
\]
\\[-1.25cm]
\end{x}

\begin{x}{\small\bf LEMMA} \ 
$\forall w \in W$, 
\[
\det(w) \  
\prod\limits_{\alpha > 0} \  
B(\Lambda + \rho,\alpha) 
\ = \ 
\prod\limits_{\alpha > 0} \  
B(w(\Lambda + \rho),\alpha).
\]

Proceeding, 
\allowdisplaybreaks
\begin{align*}
\frct{(-1)^r}{\abs{W}} \ 
&\sum\limits_{\Lambda \in \sW^+}\hsx \td_\Lambda  \  
\int_\Torus \  \Delta(t) \chisubLambda(t) F_f(t) \hsx \td_\Torus(t) 
\\[11pt]
&=\ 
\frct{(-1)^r}{\abs{W}} \ 
\sum\limits_{\Lambda \in \sW^+} \ \td_\Lambda  \  
\int_\Torus \  
\sum\limits_{w \in W} \  \det(w) \xi_{w(\Lambda + \rho)}(t) \hsx 
F_f(t) \hsx \td_\Torus(t) \ 
\\[11pt]
&=\ 
\frct{(-1)^r}{\abs{W} \ 
\prod\limits_{\alpha > 0}\hsx  B(\rho,\alpha)} \ 
\sum\limits_{\Lambda \in \sW^+} \
\sum\limits_{w \in W} \ 
\int_\Torus \ 
\det(w) 
\prod\limits_{\alpha > 0}\hsx  B(\Lambda + \rho,\alpha) \ 
\xi_{w(\Lambda + \rho)}(t)
F_f(t) \hsx \td_\Torus(t)
\\[11pt]
&=\ 
\frct{(-1)^r}{\abs{W} \ 
\prod\limits_{\alpha > 0} \  
B(\rho,\alpha)}  \ 
\sum\limits_{\Lambda \in \sW^+} \
\sum\limits_{w \in W} \ 
\int_\Torus \
\prod\limits_{\alpha > 0}\  
B(w(\Lambda + \rho),\alpha) \hsx
\xi_{w(\Lambda + \rho)}(t)
F_f(t) \hsx \td_\Torus(t)
\\[11pt]
&=\ 
\frct{(-1)^r}{\abs{W} \
\prod\limits_{\alpha > 0}\  B(\rho,\alpha)} \ 
\sum\limits_{\lambda \in \sW}\
\int_\Torus \
\prod\limits_{\alpha > 0} \  
B(\lambda,\alpha)  \hsx
\xi_\lambda(t)
F_f(t) \hsx \td_\Torus(t),
\end{align*}
the $\lambda \in \sW$ for which $\prod\limits_{\alpha > 0}\hsx  B(\lambda,\alpha)  = 0$ making no contribution.
\\[-.25cm]
\end{x}

\begin{x}{\small\bf REMARK} \ 
The elements $\lambda \in \sW$ such that $w \lambda \neq \lambda$ when $w \neq e$ $(w \in W)$ 
are in a one-to-one correspondence with the pairs $(\Lambda, w) \in \sW^+ \times W$ via the arrow 
$(\Lambda,w) \ra w(\Lambda + \rho)$.
\\[-.25cm]
\end{x}

To isolate $f(e)$, put $\widecheck{f}(X) = f(\exp X)$ $(X \in \fg)$ $-$then 
$\widecheck{f} \in C^\infty(\fg)$ and $\forall \ H \in \ft$,
\allowdisplaybreaks
\begin{align*}
F_f(\exp H) \ 
&=\ 
\Delta(\exp H) \
\int_G \ f(x(\exp H) x^{-1}) 
\hsx \td_G(x)
\\[11pt]
&=\ 
\Delta(\exp H) \
\int_G \ f(\exp \Ad(x)H)) 
\hsx \td_G(x)
\\[11pt]
&=\ 
\Delta(\exp H) \
\int_G \
\widecheck{f}(\Ad(x)H)
\hsx \td_G(x).
\end{align*}
\\[-.75cm]

\begin{x}{\small\bf LEMMA} \ 
Let $\lambda$ be a linear function on $\ft_\tc$ $-$then there exists a unique $\Ad \hsx G$ invariant analytic function $\Gamma_\lambda$ on $\fg$ 
such that $\forall \ H \in \ft$, 
\[
\Gamma_\lambda(H) \hsx \pi(H) 
\ = \ 
\sum\limits_{w \in W} \ \det(w) \hsx e^{w \lambda(H)}.
\]
\\[-1.25cm]
\end{x}

\begin{x}{\small\bf APPLICATION} \ 
Take $\lambda = \rho$ $-$then there exists a unique $\Ad \hsx G$ invariant analytic function $\Gamma_\rho$ on $\fg$ such that 
$\forall \ H \in \ft$, 
\allowdisplaybreaks
\begin{align*}
\Gamma_\rho(H) \hsx \pi(H) \ 
&=\ 
\sum\limits_{w \in W} \ 
\det(w) \hsx e^{w \rho (H)} 
\\[11pt]
&=\ 
\Delta(H) \qquad \text{(cf. $\S10$, $\#11$).}
\end{align*}
Therefore
\allowdisplaybreaks
\begin{align*}
\Delta(\exp H) \hsx \int_G \hsx \widecheck{f} (\Ad(x)H) \ \td_G(x)
&=\ 
\Gamma_\rho(H) \pi(H) \ 
\int_G \ \widecheck{f} (\Ad(x)H) \ \td_G(x)
\\[11pt]
&=\ 
\pi(H) \Gamma_\rho(H)  \ 
\int_G \hsx \widecheck{f} (\Ad(x)H) \ \td_G(x)
\\[11pt]
&=\ 
\pi(H)\ 
\int_G \ \Gamma_\rho(H) \hsx  \widecheck{f} (\Ad(x)H) \ \td_G(x)  
\\[11pt]
&=\ 
\pi(H) \ 
\int_G \ \Gamma_\rho(\Ad(x)H) \hsx \widecheck{f} (\Ad(x)H) \ \td_G(x)  
\\[11pt]
&=\ 
\phi_{\Gamma_\rho \widecheck{f}} \hsx (H).
\end{align*}
Summary: $\forall \hsx H \in \ft$, 
\[
F_f(\exp H) 
\ = \ 
\phi_{\Gamma_\rho \widecheck{f}} \hsx (H).
\]
\\[-1.25cm]
\end{x}

\begin{x}{\small\bf SUBLEMMA} \ 
In $\{H: \pi(H) \pi \neq 0\}$, 
\allowdisplaybreaks
\begin{align*}
\Gamma_\rho(0) \ 
&=\ 
\lim\limits_{H \ra \hsy 0} \ 
\frct{\Delta(H)}{\pi(H)} 
\\[11pt]
&=\ 
\lim\limits_{H \ra \hsy 0} \ 
\frct{\prod\limits_{\alpha > 0} \  \bigg(e^{\alpha(H)/2} - e^{-\alpha(H)/2}\bigg)}{\prod\limits_{\alpha > 0} \ \alpha(H)}
\\[11pt]
&=\ 
\lim\limits_{H \ra \hsy 0} \ 
\prod\limits_{\alpha > 0} \ 
\bigg[
\frct{e^{\alpha(H)/2} - e^{-\alpha(H)/2}}{\alpha(H)}
\bigg]
\\[11pt]
&=\ 
1.
\end{align*}

Next
\[
F_f(\exp H) 
\ = \ 
\phi_{\Gamma_\rho \widecheck{f}} \hsx (H) 
\]
\qquad $\implies$
\[
\partial(\widetilde{\pi}) F_f (\exp H) 
\ = \ 
\partial(\widetilde{\pi}) \hsx \phi_{\Gamma_\rho \widecheck{f}} \hsx (H) 
\]
\qquad $\implies$
\allowdisplaybreaks
\begin{align*}
(\partial(\widetilde{\pi}) F_f \circ \exp) (0) \ 
&=\ 
\partial(\widetilde{\pi}) \phi_{\Gamma_\rho \widecheck{f}}(0)
\\[11pt]
&=\ 
\big(\abs{W} \  \prod\limits_{\alpha > 0} \ 
B(\rho,\alpha) \big)( \Gamma_\rho \widecheck{f}) (0) \qquad \text{(cf, $\S11$, $\#14$).}
\end{align*}
And
\[
(\Gamma_\rho \hsx \widecheck{f})  (0)
\ = \ 
\Gamma_\rho (0) \hsx \widecheck{f} (0)
\ = \ 
\widecheck{f} (0)
\ = \ 
f(e).
\]

Therefore
\[
f(e) 
\ = \ 
\frct{1}{\abs{W} \hsx \prod\limits_{\alpha > 0} \ 
B(\rho,\alpha)} 
\ 
\lim\limits_{H \ra \hsy 0} \hsx F_f(\exp H; \partial(\widetilde{\pi})).
\]
\end{x}
\vspace{0.3cm}

\begin{x}{\small\bf NOTATION} \ 
Given $\lambda \in \sW$, put
\[
\widehat{F}_f(\lambda) 
\ = \ 
\int_\Torus \ F_f(t) \hsx \xi_\lambda(t) \ \td_\Torus(t),
\]
the Fourier transform of $F_f$.
\\[-.25cm]
\end{x}

\begin{x}{\small\bf \un{N.B.}} \ 
Assume that the Haar measure on $\widehat{\Torus}$ is normalized so that Fourier inversion is valid 
(thus each $\lambda \in \sW$ is assigned mass 1).
\\[-.2cm]

Write
\allowdisplaybreaks
\begin{align*}
\lim\limits_{H \ra \hsy 0} \hsx F_f(\exp H; \partial(\widetilde{\pi})) \ 
&=\ 
\int_{\widehat{\Torus}} \ 
\widehat{F}_f(\lambda) \ 
\lim\limits_{H \ra \hsy 0} \
\xi_{-\lambda} (\exp H;\partial(\widetilde{\pi})
\td_{\widehat{\Torus}}(\lambda)
\\[11pt]
&=\ 
\int_{\widehat{\Torus}} \ 
\widehat{F}_f(\lambda)
\lim\limits_{H \ra \hsy 0} \
\partial(\widetilde{\pi}) e^{-\lambda(H)} \ 
\td_{\widehat{\Torus}}(\lambda)
\\[11pt]
&=\ 
(-1)^r \ 
\int_{\widehat{\Torus}} \ 
\widehat{F}_f(\lambda)
\prod\limits_{\alpha > 0} \ 
B(\lambda,\alpha)  \
\td_{\widehat{\Torus}}(\lambda)
\\[11pt]
&=\ 
(-1)^r \ 
\int_{\widehat{\Torus}} \ 
\prod\limits_{\alpha > 0} \ 
B(\lambda,\alpha) \widehat{F}_f(\lambda)
\ 
\td_{\widehat{\Torus}}(\lambda)
\\[11pt]
&=\ 
(-1)^r \ 
\int_{\widehat{\Torus}} \hsx 
\prod\limits_{\alpha > 0} \ B(\lambda,\alpha) \ 
\bigg(
\int_\Torus \  
F_f(t) \hsx \xi_\lambda(t) \hsx \td_\Torus(t) \bigg) 
\ 
\td_{\widehat{\Torus}}(\lambda)
\\[11pt]
&=\ 
(-1)^r \ \sum\limits_{\lambda \in \sW} \ 
\int_\Torus \ B(\lambda, \alpha) \hsx \xi_\lambda(t) F_f(t) \  \td_\Torus(t).
\end{align*}

Therefore
\allowdisplaybreaks
\begin{align*}
f(e) \ 
&=\ 
\frct{(-1)^r}{\abs{W} \prod\limits_{\alpha > 0}  \ B(\rho,\alpha) }
\sum\limits_{\lambda \in \sW} \ 
\int_\Torus \ 
B(\lambda, \alpha) \hsx \xi_\lambda(t) F_f(t) \hsx \td_\Torus(t)
\\[11pt]
&=\ 
\sum\limits_{\Lambda \in \sW^+} \ \td_\Lambda \hsx \int_G \hsx f(x) \chisubLambda (x) \  \td_G(x), 
\end{align*}
the relation at issue.
\end{x}


\chapter{
$\boldsymbol{\S}$\textbf{13}.\quad  DETECTION}
\setlength\parindent{2em}
\setcounter{theoremn}{0}
\renewcommand{\thepage}{B I \S13-\arabic{page}}

\qquad Let \mG be a compact group.
\\[-.25cm]

\begin{x}{\small\bf DEFINITION} \ 
The \un{character ring} $X(G)$ is the free abelian group generated by the irreducible characters of \mG 
(i.e., by the $\chisubPi$ $(\Pi \in \widehat{G}))$ under pointwise addition and multiplication with unit $1_G$.
\\[-.25cm]
\end{x}

\begin{x}{\small\bf DEFINITION} \ 
An element of $X(G)$ is called a \un{virtual character}.
\\[-.25cm]
\end{x}

\begin{x}{\small\bf NOTATION} \ 
$CL(G)$ is the subspace of $C(G)$ comprised of continuous class functions. (cf. $\S2$, $\#27$).
\\[-.25cm]
\end{x}

\begin{x}{\small\bf LEMMA} \ 
A class function $f \in CL(G)$ is a virtual character of \mG iff 
\[
\langle f, \chisubPi \rangle 
\ = \ 
\int_G \ f(x) \hsx \ov{\chisubPi(x)} \ \td_G(x) \in \Z
\]
for all $\Pi \in \widehat{G}$.

PROOF \ 
The condition is obviously necessary.  
As for the sufficiency, we have
\[
\norm{f}^2 
\ = \ 
\sum\limits_{\Pi \in \widehat{G}} \ 
\langle f, \chisubPi \rangle^2 \qquad \text{(cf. $\S2$, $\#19$),} 
\]
hence
\[
\langle f, \chisubPi \rangle 
\ = \ 
0
\]
for all but finitely many $\chisubPi$, say $\chisubPiOne, \ldots \chisubPin$ and then
\[
f 
\ = \ 
\sum\limits_{i = 1}^n \ \langle f, \chisubPii \rangle \chisubPii \qquad \text{(ibid.)}.
\]

[Note: \ 
A priori, this is an equality in the $\tL^2$-sense, hence is valid almost everywhere.  
But both sides are continuous, thus equality is valid everywhere.]
\\[-.25cm]
\end{x}


Let \mG be a compact connected Lie group.

\begin{x}{\small\bf NOTATION} \ 
$CL^\infty(G)$ is the set of $C^\infty$ class functions.
\\[-.25cm]
\end{x}

\begin{x}{\small\bf RAPPEL} \ 
The characters of \mG belong to $CL^\infty(G)$.
\\[-.25cm]
\end{x}

\begin{x}{\small\bf \un{N.B.}} \ 
Therefore $X(G)$ is a subring of the ring of $C^\infty$ functions on \mG.
\\[-.25cm]
\end{x}

\begin{x}{\small\bf REMARK} \ 
Per $\#4$, suppose that $f \in CL(G)$ has the property that
\[
\langle f, \chisubPi \rangle \in \Z
\]
for all $\Pi \in \widehat{G}$ $-$then it follows after the fact that $f \in CL^\infty(G)$.
\\[-.25cm]
\end{x}

Let $\Torus \subset G$ be a maximal torus and assign to the symbol $X(\Torus)$ the obvious interpretation.
\\[-.25cm]

\begin{x}{\small\bf RAPPEL} \ 
The arrow
\[
f \ra \restr{f}{\Torus}
\]
of restriction defines an isomorphism
\[
CL(G) \ra C(\Torus)^W \qquad \text{(cf. $\S4$, $\#24$).}
\]
\\[-1.25cm]
\end{x}

\begin{x}{\small\bf APPLICATION} \ 
Restriction to $\Torus$ induces an injective homomorphism
\[
X(G) \ra X(\Torus)^W.
\]
\\[-1.25cm]
\end{x}

Take a $\phi \in X(\Torus)^W$ and let $f \in CL(G)$ be the class function that restricts to $\phi$.
\\[-.25cm]

\begin{x}{\small\bf LEMMA}\ 
$f$ is a virtual character of \mG, i.e., $f \in X(G)$.
\\[-.25cm]

PROOF \ 
With $\#4$ in view, write
\allowdisplaybreaks
\begin{align*}
\langle f, \chisubPi \rangle \ 
&=\ 
\int_G \ f(x) \hsx \ov{\chisubPi(x)} \ \td_G(x)
\\[11pt]
&=\ 
\frac{1}{\abs{W}} \ \int_\Torus \ 
\abs{\Delta(t)}^2  \hsx \phi(t) \hsx \ov{\chisubPi(t)} \ \td_\Torus(t)
\\[11pt]
&=\ 
\frac{1}{\abs{W}} \ \int_\Torus \ 
\Delta(t) \hsx \phi(t) \hsx \ov{\Delta(t) \hsx \chisubPi(t)} \ \td_\Torus(t)
\\[11pt]
&=\ 
\frac{1}{\abs{W}} \ (\abs{W} \Z) \ = \  \Z.
\end{align*}
\\[-.75cm]
\end{x}

\begin{x}{\small\bf SCHOLIUM} \ 
\[
X(G) 
\ \approx \ 
X(\Torus)^W.
\]
\\[-1.25cm]
\end{x}

\begin{x}{\small\bf \un{N.B.}} \ 
Rephrased, a continuous class function $f:G \ra \Cx$ is a virtual character of \mG iff $\restr{f}{\Torus}$ is a virtual character of $\Torus$.
\\[-.25cm]
\end{x}

\begin{x}{\small\bf THEOREM} \ 
Let $f \in CL(G)$ $-$then $f \in X(G)$ iff its restriction to every finite elementary subgroup of \mG is a virtual character.
\\[-.25cm]

PROOF \ 
To establish the nontrivial assertion, let $H \subset G$ be a finite subgroup $-$then
the assumption on $f$ coupled with A, II, $\S12$, $\#1$ implies that $\restr{f}{H} \in X(H)$.  
Matters can thus be reinforced, the assumption on $f$ becoming that its restriction to every finite subgroup of \mG is a virtual character and, 
thanks to what has been said above, one might just as well work with $\Torus$ rather than \mG.  
Choose a sequence $H_1 \subset H_2 \subset \cdots$ of finite subgroups of $\Torus$ whose union is dense in $\Torus$ $-$then 
$\forall \ \chi \in X(\Torus)$, 

\allowdisplaybreaks
\begin{align*}
\langle f, \hsx \chi \ranglesubTorus \ 
&=\ 
\int_\Torus \ f \ov{\chi} 
\\[11pt]
&=\ 
\lim\limits_{n \ra \infty} \ \bigg(\frac{1}{\abs{H_n}} \
\sum\limits_{h \in H_n} \ f(h) \hsx \ov{\chi(h)}\bigg)
\\[11pt]
&=\ 
\lim\limits_{n \ra \infty} \hsx 
\langle f, \chi \ranglesubHn.
\end{align*}
But 
\[
\restr{f}{H_n}  \in X(H_n)
\]
\qquad\qquad $\implies$
\[
\langle f, \chi \ranglesubHn \in \Z
\]
\qquad\qquad $\implies$
\[
\langle f, \chi \ranglesubTorus \in \Z.
\]
\end{x}


\chapter{
$\boldsymbol{\S}$\textbf{14}.\quad  INDUCTION}
\setlength\parindent{2em}
\setcounter{theoremn}{0}
\renewcommand{\thepage}{B I \S14-\arabic{page}}

\qquad Let \mG be a finite group, $\Gamma \subset G$ a subgroup.
\\[-.25cm]

\begin{x}{\small\bf RAPPEL} \ 
There is an arrow
\[
i_{\Gamma \ra G} : CL(\Gamma) \ra CL(G)
\]
which sends characters of $\Gamma$ to characters of \mG (cf. A, II, $\S9$, $\#10$), thus induces an arrow 
\[
X(\Gamma) \ra X(G).
\]
\\[-1.25cm]
\end{x}

\begin{x}{\small\bf \un{N.B.}} \ 
If 
\[
G \ = \ \coprod\limits_{k = 1}^n \ x_k \Gamma
\]
and if $\phi \in CL(\Gamma)$ is a class function, then
\[
\big(i_{\Gamma \ra G} \hsx \phi \big)(x) 
\ = \ 
\sum\limits_{k = 1}^n \ \overset{\circ}{\phi} (x_k^{-1} \hsx x \hsx x_k) \qquad \text{(cf. A, II, $\S7$, $\#10$),}
\]
i.e., 
\[
\big(i_{\Gamma \ra G} \hsx \phi \big)(x) 
\ = \ 
\sum\limits_{k, x_k^{-1} \hsy x \hsy x_k \in \Gamma} \ \phi(x_k^{-1} \hsx x \hsx x_k).
\]
\\[-.25cm]
\end{x}

Let \mG be a compact Lie group, $\Gamma \subset G$ a closed Lie subgroup.
\\[-.25cm]

\begin{x}{\small\bf NOTATION} \ 
Given an $x \in G$, write $(G/\Gamma)^x$ for the fixed point set of the action of $x$ on $(G/\Gamma)$.
\\[-.25cm]
\end{x}

\begin{x}{\small\bf LEMMA} \ 
A coset $y \Gamma$ in $G/\Gamma$  lies in $(G/\Gamma)^x$  iff $y^{-1} \hsx x \hsx y \in \Gamma$.
\\[-.25cm]
\end{x}


\begin{x}{\small\bf LEMMA} \ 
If cosets $y_1 \Gamma$, $y_2 \Gamma$ lie in the same connected component of $(G/\Gamma)^x$, 
then $y_2 \Gamma = y y_1 \Gamma$ for some $y$ in the centralizer of $x$.
\\[-.25cm]
\end{x}

\begin{x}{\small\bf \un{N.B.}} \ 
If $\phi \in CL(\Gamma)$ is a class function and if $y_2 = y \hsy y_1 \hsy \gamma$, then 
\begin{align*}
\phi(y_2^{-1} \hsy x \hsy y_2) \ 
&=\ 
\phi(\gamma^{-1} \hsy y_1^{-1} \hsy y^{-1} \hsy x \hsy y \hsy y_1 \gamma)
\\[11pt]
&=\ 
\phi (y_1^{-1} \hsy y^{-1} \hsy x \hsy y \hsy y_1)
\\[11pt]
&=\ 
\phi (y_1^{-1} \hsx x \hsx y_1).
\end{align*}
\\[-1.25cm]
\end{x}

Let $C_1, \ldots, C_m$ be the connected components of $(G/\Gamma)^x$, thus
\[
(G/\Gamma)^x
\ = \ 
\coprod\limits_{j = 1}^m \ C_j,
\]
let $\chi(C_j)$ be the Euler characteristic of $C_j$, and fix elements 
\[
y_1 \Gamma \in C_1, \ldots, y_m \Gamma \in C_m.
\]
\\[-1.25cm]

\begin{x}{\small\bf NOTATION} \ 
Given a class function $\phi \in CL(\Gamma)$, put
\[
\big(i_{\Gamma \ra G} \hsx \phi\big)(x) 
\ = \ 
\sum\limits_{j = 1}^m \ \chi(C_j) \phi(y_j^{-1} \hsy x \hsx y_j).
\]
\\[-1.25cm]
\end{x}

\begin{x}{\small\bf LEMMA} \ 
\[
i_{\Gamma \ra G} \hsx \phi \in CL(G),
\]
the \un{induced} class function.
\\[-.25cm]
\end{x}

\begin{x}{\small\bf \un{N.B.}} \ 
Therefore
\[
i_{\Gamma \ra G} : CL(\Gamma) \ra CL(G).
\]
\\[-1.25cm]
\end{x}

\begin{x}{\small\bf REMARK} \ 
The definition of $i_{\Gamma \ra G} \hsx \phi$ is independent of the choice of representatives $y_j\Gamma$ for the components of 
$(G/\Gamma)^x$ but it is not quite obvious that $i_{\Gamma \ra G} \hsx \phi$ is continuous.
\\[-.25cm]
\end{x}

\begin{x}{\small\bf RECONCILIATION} \ 
Take the case when \mG and $\Gamma \subset G$ are finite.  
Write
\[
G \ = \ \coprod\limits_{k = 1}^n \ x_k \Gamma.
\]
Then, as recalled in $\#2$, 
\[
\big(i_{\Gamma \ra G} \hsx \phi\big)(x) 
\ = \ 
\sum\limits_{k, x_k^{-1} \hsy x \hsy x_k \in \Gamma} \  \phi(x_k^{-1} \hsy x \hsy x_k)
\]
which, in view of $\#4$, is equal to 
\[
\sum\limits_{\substack{y \Gamma \in G/\Gamma\\ x y \Gamma = y\Gamma}}
\phi(y^{-1} x y)
\]
or still, is equal to 
\[
\sum\limits_{y\Gamma \in (G/\Gamma)^x} \ \phi(y^{-1} \hsy x \hsy y).
\]
But here the $C_j$ are points, say
\[
C_j \ = \ \{y_j \Gamma\} 
\quad (\implies (G/\Gamma)^x = \{\{y_1 \Gamma\}, \dots, \{y_m \Gamma\}\},
\]
so $\chi(C_j) = 1$, thus
\begin{align*}
\sum\limits_{j = 1}^m \ \chi(C_j) \hsx \phi(y_j^{-1} \hsy x \hsy y_j) \ 
&=\ 
\sum\limits_{j = 1}^m \ \phi(y_j^{-1} x y_j)
\\[11pt]
&=\ 
\sum\limits_{y \Gamma \in (G/\Gamma)^x} \ \phi(y^{-1} \hsy x \hsy y).
\end{align*}
\\[-.75cm]
\end{x}


\begin{x}{\small\bf RAPPEL} \ 
A compact connected Lie group of positive dimension has zero Euler characteristic, so the connected components of a 
compact Lie group of positive dimension have zero Euler characteristic.
\\[-.25cm]
\end{x}

\begin{x}{\small\bf EXAMPLE} \ 
Take $\Gamma = \{e\}$, let $\phi = 1_\Gamma$, and assume that $\dim G > 0$ $-$then
$(G/\Gamma)^x$ is empty if 
$x \neq e$, hence for such $x$,
\[
\big(i_{\Gamma \ra G} \hsx \phi\big) (x) 
\ = \ 
0, 
\]
but if $x = e$, then $(G/\Gamma)^e = G$ and 
\begin{align*}
\big(i_{\Gamma \ra G} \hsx \phi\big) (e) \ 
&=\ 
\sum\limits_{j = 1}^m \ \chi(C_j)  \hsx \phi(y_j^{-1}e  \hsy y_j)
\\[11pt]
&=\ 
\bigg(\sum\limits_{j = 1}^m \ \chi(C_j)\bigg)  \hsx \phi(e)
\\[11pt]
&=\ 
0.
\end{align*}
Therefore
\[
i_{\Gamma \ra G}  \hsx \phi
\ = \ 0.
\]
\\[-1.25cm]
\end{x}

\begin{x}{\small\bf DEFINITION} \ 
A closed subgroup \mH of \mG is \un{generic} if it is topologically cyclic and of finite index in its normalizer.

[Note: \ 
Let \mG be a compact connected Lie group, $\Torus \subset G$ a maximal torus $-$then $\Torus$ is generic.]
\\[-.25cm]
\end{x}

\begin{x}{\small\bf DEFINITION} \ 
An element $x \in G$ is \un{generic} if it generates a generic subgroup of \mG.
\\[-.25cm]

[Note: \ 
Let \mG be a compact connected Lie group $-$then a generic element is necessarily regular.]
\\[-.25cm]
\end{x}

\begin{x}{\small\bf LEMMA} \ 
The generic elements are dense in \mG.
\\[-.25cm]
\end{x}

\begin{x}{\small\bf THEOREM} \ 
Suppose that $x \in G$ is generic $-$then 
\[
\abs{(G/\Gamma)^x} 
\ > \ 
\infty
\]
and 
\[
\big(i_{\Gamma \ra G} \hsx \phi\big) (x) 
\ = \ 
\sum\limits_{y \Gamma \in (G/\Gamma)^x} \hsx \phi(y^{-1} \hsy x \hsy y).
\]
\\[-1.25cm]
\end{x}

\begin{x}{\small\bf EXAMPLE} \ 
Take $\Gamma = \{e\}$, let $\phi = 1_\Gamma$, and assume that $\dim G > 0$ $-$then 
at every generic element of \mG,
\[
\big(i_{\Gamma \ra G} \hsx \phi\big)(x) 
\ = \ 
0,
\]
hence by continuity (in conjuction with $\#16$), 
\[
i_{\Gamma \ra G} \hsx \phi
\ = \ 
0 \qquad \text{(cf. $\#13$).}
\]
\\[-1.25cm]
\end{x}

Let \mG be a compact Lie group, let $\Gamma_1$, $\Gamma_2 \subset G$ be closed Lie subgroups, and let
\[
G 
\ = \ 
\bigcup\limits_{s \in S} \ \Gamma_1 s \Gamma_2
\]
be a double coset decomposition of \mG.
\\[-.2cm]

\begin{x}{\small\bf \un{N.B.}} \ 
\[
\Gamma_1 \backslash G/\Gamma_2 
\]
is the orbit space per the action of $\Gamma_1$ by left translation on $G/\Gamma_2$.
\\[-.25cm]
\end{x}

Write
\[
\Gamma_1 \backslash G/\Gamma_2 
\ = \ 
\coprod\limits_{s \in S} \ U_s,
\]
where each $U_s$ is a connected component of one orbit type for the action of $\Gamma_1$ on
$G/\Gamma_2$.  
Fix elements $x_s \in G$ such that $\Gamma_1 x_s \Gamma_2 \in U_s$ and for each $s$ let
\[
\phi_s:CL(\Gamma_2) \ra CL(\Gamma_1) 
\]
denote the following composite: 
Take a $\phi \in CL(\Gamma_2)$ and form $\phi^s \equiv \phi \circ I_{x_s^{-1}}$ (a class function on $x_s \hsy \Gamma_2 \hsy x_s^{-1}$), 
then restrict $\phi^s$ to $\Gamma_2(s) \equiv x_s \hsy \Gamma_2 \hsy x_s^{-1} \cap \Gamma_1$, 
call it $\phi_s$, and finally apply $i_{\Gamma_2(s) \ra \Gamma_1}$.  
I.e.:
\[
\phi_s(\phi) 
\ = \ 
i_{\Gamma_2(s) \ra \Gamma_1} \phi_s.
\]
\\[-.2cm]

\begin{x}{\small\bf THEOREM} \ 
As maps from $CL(\Gamma_2)$ to $CL(\Gamma_1)$, 
\[
r_{G \ra \Gamma_1} \circ i_{\Gamma_2 \ra G} 
\ = \ 
\sum\limits_{s \in S} \ \chi^\#(U_s) \phi_s,
\]
where for each $s \in S$, 
\[
\chi^\#(U_s) 
\ = \ 
\chi(\ov{U_s}) - \chi(\ov{U_s} - U_s).
\]
\\[-1.25cm]
\end{x}

\begin{x}{\small\bf \un{N.B.}} \ 
When \mG and $\Gamma_1, \Gamma_2 \subset G$ are finite, matters reduce to A, II, $\S8$, $\#3$.
\\[-.2cm]

Here is a sketch of the proof.
\\[-.2cm]

\qquad 1. \quad 
Fix a class function $\phi \in CL(\Gamma_2)$ and a $\gamma_1 \in \Gamma_1$.
\\[-.2cm]

\qquad 2. \quad  
Let $\widetilde{U}_s \subset G/\Gamma_2$ denote the inverse image of $U_s$ under the projection to 
$\Gamma_1\backslash G/\Gamma_2$, thus
\[
G/\Gamma_2 
\ = \ 
\coprod\limits_{s \in S} \ \widetilde{U}_s.
\]
\\[-.75cm]

\qquad 3. \ \ 
Let $C_1, \ldots, C_m$ be the connected components of $\big(G/\Gamma_2\big)^{\gamma_1}$, thus
\[
\big(G/\Gamma_2\big)^{\gamma_1}
\ = \ 
\coprod\limits_{j = 1}^m \ C_j.
\]
\\[-.75cm]

\qquad 4. \ \ 
For each pair $(s,j)$, put
\[
V_{s,j} 
\ = \ 
\big(\Gamma_1 \cdot x_s \Gamma_2) \cap C_j \ \subset \ \widetilde{U_s} \cap C_j  \ \subset \ G/\Gamma_2.
\]
\\[-.75cm]

\qquad 5. \ \ 
The arrows
\[
V_{s,j}  \ra \widetilde{U_s} \cap C_j  \ra U_s
\]
are a fibration sequence, hence by the multiplicativity of the Euler characteristic, 
\[
\chi(C_j) 
\ = \ 
\sum\limits_{s \in S} \  \chi(V_{s,j} ) \chi^\#(U_s).
\]
\\[-.75cm]

\qquad 6. \ \ 
Fix elements $\gamma_{s, j} \in \Gamma_1$ such that
\[
\gamma_{s,j} \hsx x_s \Gamma_2 \in V_{s,j}.
\]
Then in particular, 
\[
\gamma_{s,j} \hsx x_s \hsx \Gamma_2 \in \big(G/\Gamma_2\big)^{\gamma_1}
\]
\qquad $\implies$
\[
x_s^{-1} \hsx \gamma_{s,j}^{-1} \hsx \gamma_1 \hsx \gamma_{s,j} \hsx  x_s \in \Gamma_2
\]
\qquad $\implies$
\[
\gamma_{s,j}^{-1} \hsx \gamma_1 \hsx \gamma_{s,j} \in x_s \hsx \Gamma_2 \hsx x_s^{-1},
\]
the domain of $\phi^s$.
\\[-.2cm]

\qquad 7. \ \ 
From the definitions,
\begin{align*}
(r_{G \ra \Gamma_1} (i_{\Gamma_2 \ra G} \hsx \phi) ) (\gamma_1) \ 
&=\ 
\sum\limits_{j = 1}^m  \  \chi(C_j) \hsx \phi(y_j^{-1} \hsy \gamma_1\hsy  y_j)
\\[11pt]
&=\ 
\sum\limits_{s \in S} \ \sum\limits_{j = 1}^m \  
\chi(V_{s,j}) \chi^\#(U_s) \hsx \phi(x_s^{-1} \hsy \gamma_{s,j}^{-1} \hsy \gamma_1 \hsy \gamma_{s,j} \hsy x_s)
\\[11pt]
&=\ 
\sum\limits_{s \in S} \  \chi^\#(U_s) \  \sum\limits_{j = 1}^m  \  
\chi(V_{s,j}) \hsx  \phi^s(\gamma_{s,j}^{-1} \hsx \gamma_1 \hsx \gamma_{s,j}).
\end{align*}
\\[-.75cm]

\qquad 8. \ \ 
The isotropy subgroup of the action of $\Gamma_1$ on $x_s \hsx \Gamma_2 \in G/\Gamma_2$ is
\[
\Gamma_2(s) 
\ = \ 
x_s \hsx \Gamma_2  \hsx x_s^{-1} \cap \Gamma_1.
\]
And
\begin{align*}
\big(\Gamma_1/\Gamma_2(s)\big)^{\gamma_1} \ 
&\approx \ 
(\Gamma_1 \cdot x_s \hsx \Gamma_2)^{\gamma_1}
\\[11pt]
&=\ 
\coprod\limits_{j = 1}^m \  V_{s,j} 
\\[11pt]
&\subset \ 
G/\Gamma_2.
\end{align*}
\\[-.2cm]

\qquad 9. \ \ 
Given $s \in S$, 
\[
\sum\limits_{j = 1}^m \  \chi(V_{s,j}) \hsx \phi^s \hsx \big(\gamma_{s,j}^{-1} \hsx \gamma_1 \hsx \gamma_{s,j}\big)
\ = \ 
\big(i_{\Gamma_2(s) \ra \Gamma_1} \phi^s\big) (\gamma_1).
\]
\\[-.2cm]

\qquad 10. \ \ 
Therefore
\begin{align*}
\big(r_{G \ra \Gamma_1} \hsx \big(i_{\Gamma_2 \ra G} \hsx \phi\big)\big) (\gamma_1)\ 
&=\ 
\sum\limits_{s \in S} \  \chi^\#(U_s) (i_{\Gamma_2(s) \ra \Gamma_1} \hsx \phi^s\big) (\gamma_1)
\\[11pt]
&=\ 
\sum\limits_{s \in S} \  \chi^\#(U_s) \hsx \phi_s(\phi) (\gamma_1),
\end{align*}
the contention.
\\[-.25cm]
\end{x}

\begin{x}{\small\bf THEOREM} \ 
The arrow 
\[
i_{\Gamma \ra G}:CL(\Gamma) \ra CL(G)
\]
sends virtual characters to virtual characters, thus induces an arrow
\[
X(\Gamma) \ra X(G).
\]

PROOF \ 
Recall first that this is true when \mG is finite (cf. $\#1$).  
In general, let $\chi \in X(\Gamma)$ $-$then to conclude that
\[
i_{\Gamma \ra G} \hsx \chi \in X(G), 
\]
it suffices to show that its restriction to every finite subgroup \mH of \mG is a virtual character (cf. $\S13$, $\#14$).  
So consider
\[
r_{G \ra H}(i_{\Gamma \ra G} \hsx \chi)
\]
or still, take in the above $\Gamma_1 = H$, $\Gamma_2 = \Gamma$, $\phi = \chi$, and consider
\[
\sum\limits_{s \in S} \hsx \chi^\#(U_s) \phi_s(\chi).
\]
Here
\[
\Phi_s(\chi) 
\ = \ 
i_{x_s\Gamma x_s^{-1} \cap H \ra H} \hsx  \chisubs,
\]
where $\chisubs$ is the restriction of $\chi^s$ to $\Gamma(s) \equiv x_s \Gamma x_s^{-1} \cap H$, a finite group.  
But now
\[
\chisubs \in X(\Gamma(s)) \implies i_{\Gamma(s) \ra H} \hsx \chi \in X(H),
\]
which finishes the proof.
\\[-.25cm]
\end{x}

\begin{x}{\small\bf \un{N.B.}} \ 
If \mG is finite, then the arrow
\[
i_{\Gamma \ra G}:CL(\Gamma) \ra CL(G)
\]
sends characters of $\Gamma$ to characters of \mG but this need not be true if $\dim G > 0$
(cf. $\#13$) ($1_\Gamma$ is a character of $\Gamma$ but the induced class function
\[
i_\Gamma \ra G^1 \Gamma
\]
is identically zero, a virtual character, not a character).
\end{x}
\vspace{0.3cm}

\begin{x}{\small\bf REMARK} \ 
Let \mG be a compact connected semisimple Lie group, $\Torus \subset G$ a maximal torus $-$then
\[
\widehat{G} \longleftrightarrow \td\widehat{T} \cap \ov{C} \qquad \text{(cf. $\S10$, $\#23$).}
\]
While the theory developed above gives rise to an arrow
\[
X(T) \ra X(G),
\]
it does not respect the foregoing parameterization which can only be accomplished by a more 
sophisticated version of the preceding process.
\\[-.25cm]
\end{x}

\[
\text{APPENDIX}
\]

There is a different approach to induction which is suggested by A, II, $\S9$, $\#1$.  

So let \mG be a compact Lie group, $\Gamma \subset G$ a closed Lie subgroup.
\\[-.25cm]

{\small\bf CONSTRUCTION} \ 
Let $(\theta, \E)$ be a finite dimensional unitary representation of $\Gamma$ and denote by
$\E_{\Gamma, \theta}^G$ the space of all \mE-valued measurable functions $f$ on \mG such that 
$f(x \hsx \gamma) = \theta(\gamma^{-1}) f(x)$ $(x \in G, \ \gamma \in \Gamma$) subject to  
\[
\int_{G/\Gamma} \  \norm{f}^2 \  \td_{G/\Gamma} 
\ < \ 
\infty.
\]
Then the prescription 
\[
(\Ind_{\Gamma, \theta}^G(x) f) (y)
\ = \ 
f(x^{-1} y)
\]
defines a representation $\Ind_{\Gamma, \theta}^G$ of \mG on $\E_{\Gamma, \theta}^G$, 
the representation of \mG \un{induced} by $\theta$.
\\[-.25cm]

{\small\bf \un{N.B.}} \ 
The inner product
\[
\langle f, g \ranglesubtheta 
\ = \ 
\int_{G/\Gamma} \  \langle f, g \rangle \  \td_{G/\Gamma}
\]
equips  $\E_{\Gamma, \theta}^G$ with the structure of a Hilbert space and $\Ind_{\Gamma, \theta}^G$ is a unitary representation.
\\[-.25cm]

{\small\bf EXAMPLE} \ 
Take $\theta$ to be the trivial representation of $\Gamma$ on $\E = \Cx$ $-$then 
$\E_{\Gamma, \theta}^G = L^2(G/\Gamma)$.
\\[-.25cm]

[Note: \ 
When $\Gamma = \{e\}$, $\E_{\Gamma, \theta}^G = L^2(G)$ and 
\[
\Ind_{\Gamma, \theta}^G 
\ = \
L,
\]
the left translation representation of \mG (cf. $\S1$, $\#5$).]


\chapter{
$\boldsymbol{\S}$\textbf{1}.\quad  ORBITAL INTEGRALS}
\setlength\parindent{2em}
\setcounter{theoremn}{0}
\renewcommand{\thepage}{B II \S1-\arabic{page}}

\qquad Let \mG be a compact group.
\vspace{0.3cm}

\begin{x}{\small\bf DEFINITION} \ 
Given $f \in C(G)$ and $\gamma \in G$, put
\[
\sO(f,\gamma)
\ = \ 
\int_G \ 
f(x \hsy \gamma \hsy x^{-1} \ \td_G(x), 
\]
the \un{orbital integral} of $f$ at $\gamma$.
\\[-.25cm]
\end{x}

\begin{x}{\small\bf LEMMA} \ 
The function $\sO(f)$ defined by the assignment
\[
\gamma \ra \sO(f,\gamma)
\]
is a continuous class function on \mG, i.e., is an element of $CL(G)$.
\\[-.25cm]
\end{x}

\begin{x}{\small\bf RAPPEL} \ 
If $f \in C(G)_\fin$, then 
\[
\langle f, \chisubPi \rangle 
\ = \ 
0
\]
for all but finitely many $\Pi$.
\\[-.25cm]
\end{x}

\begin{x}{\small\bf LEMMA} \ 
Suppose that $f \in C(G)_\fin$ $-$then $\forall \ \gamma \in G$, 
\[
\sO(f,\gamma)
\ = \ 
\sum\limits_{\Pi \in \widehat{G}} \ 
\tr(\Pi^*(f)) \hsx \chisubPi(\gamma)
\qquad \text{(cf. A, III, $\S1$, $\#3$),}
\]
the sum on the right being finite.
\\[-.25cm]

PROOF \ 
Apply I, $\S2$, $\#19$, to get 
\[
\sO(f)
\ = \ 
\sum\limits_{\Pi \in \widehat{G}} \  
\langle \hsy \sO(f), \chisubPi \hsy \rangle_{\chisubPi},
\]
where the series converges in $L^2(G)$.  
But
\begin{align*}
\langle \hsy \sO(f), \chisubPi \hsy \rangle \ 
&=\ 
\int_G \hsx \sO(f,\gamma) \ov{\chisubPi(\gamma)} \hsx \td_G(\gamma) 
\\[11pt]
&=\ 
\int_G \ 
\bigg(\int_G \hsx f(x \hsy \gamma \hsy x^{-1}) \hsx \td_G(x) \bigg) \hsx \ov{\chisubPi(\gamma)} \ \td_G(\gamma) 
\\[11pt]
&=\ 
\int_G \ 
\bigg(\int_G \ f(x \hsy \gamma \hsy x^{-1})  \hsx \ov{\chisubPi(\gamma)}\ \td_G(\gamma) \bigg)\ \td_G(x) 
\\[11pt]
&=\
\int_G \ 
\bigg(\int_G \ f(\gamma) \hsx \ov{\chisubPi(x^{-1} \hsy \gamma \hsy x)} \ \td_G(\gamma) \bigg)\ \td_G(x)  
\\[11pt]
&=\ 
\int_G \ 
\bigg(\int_G \hsx f(\gamma) \hsx \ov{\chisubPi(\gamma)} \ \td_G(\gamma) \bigg)\ \td_G(x)  
\\[11pt]
&=\ 
\int_G \ 
\langle f, \chisubPi \rangle \ \td_G(x)
\\[11pt]
&=\ 
\langle f, \chisubPi \rangle.
\end{align*}
Therefore $\langle \hsy \sO(f), \chisubPi \hsy \rangle  = 0$ for all but finitely many $\Pi$, thus the almost everywhere equality
\[
\sO(f)
\ = \ 
\sum\limits_{\Pi \in \widehat{G}} \ \langle \hsy \sO(f), \chisubPi \hsy \rangle_{\chisubPi}
\]
is that of two continuous functions, thus is valid everywhere.  
Finally, from the definitions, 
\begin{align*}
\langle f, \chisubPi \rangle \ 
&=\ 
\int_G \ f(x) \hsx \ov{\chisubPi(x)}\ \td_G(x) 
\\[11pt]
&=\ 
\int_G \ f(x) \hsy \chisubPiStar(x)\ \td_G(x) 
\\[11pt]
&=\ 
\int_G \ f(x) \hsy \tr(\Pi^*(x)) \ \td_G(x) 
\\[11pt]
&=\ 
\tr(\Pi^*(f)).
\end{align*}
\end{x}


\chapter{
$\boldsymbol{\S}$\textbf{2}.\quad  KERNELS}
\setlength\parindent{2em}
\setcounter{theoremn}{0}
\renewcommand{\thepage}{B II \S2-\arabic{page}}

\qquad Let $(X,\mu)$, $(Y,\nu)$ be $\sigma$-finite measure spaces.
\\

\begin{x}{\small\bf NOTATION} \ 
Given $K \in L^2(X \times Y)$, define $T_K:L^2(Y) \ra L^2(X)$ by
\[
(T_K \phi)(x)
\ = \ 
\int_Y \ K(x,y) \phi(y) \  \td\nu(y).
\]
\\[-1.25cm]
\end{x}

\begin{x}{\small\bf THEOREM} \ 
The map $K \ra T_K$ is a linear isometry of $L^2(X \times Y)$ onto $L_\text{HS}(L^2(Y), L^2 (X))$.
\\[-.25cm]
\end{x}

\begin{x}{\small\bf NOTATION} \ 
Given
\[
\begin{cases}
\ K_1 \in L^2(X \times Y) \\[3pt]
\ K_2 \in L^2(Y \times Z)
\end{cases}
,
\]
define their \un{convolution}
\[ 
K_1 * K_2 \hsx \in \hsx L^2(X \times Z)
\]
by
\[
(K_1 * K_2) (x,z) 
\ = \ 
\int_Y \ K_1(x,y) K_2(y,z) \  \td\nu(y).
\]

[Note: \ 
The underlying measure-theoretic assumption is again $\sigma$-finiteness 
(which is needed infra for Fubini).]
\\[-.25cm]
\end{x}

\begin{x}{\small\bf THEOREM} \ 
\[
T_{K_1} \circ T_{K_2} 
\ = \ 
T_{K_1 \hsy * \hsy K_2}.
\]
\\[-1.25cm]
\end{x}


\begin{x}{\small\bf APPLICATION} \ 
Take $X = Y = Z$ $-$then
\[
\begin{cases}
\ T_{K_1}:L^2(X) \ra L^2(X) \\[3pt]
\ T_{K_2}:L^2(X) \ra L^2(X)
\end{cases}
\]
are Hilbert-Schmidt, hence
\[
T_{K_1 \hsy * \hsy K_2}:L^2(X) \ra L^2(X)
\]
is trace class.
\\[-.25cm]
\end{x}

\begin{x}{\small\bf LEMMA} \ 
Take $X = Y = Z$ and put $K = K_1 * K_2$ $-$then
\[
\tr(T_K) 
\ = \ 
\int_X \  K (x,x) \  \td\mu(x).
\]

PROOF \ 
\allowdisplaybreaks
\begin{align*}
\tr(T_K) \ 
&=\ 
\tr(T_{K_1} \circ T_{K_2})
\\[11pt]
&=\ 
\langle T_{K_2},T_{K_1}^* \rangle_{\text{HS}}
\\[11pt]
&=\ 
\langle K_2,K_1^* \rangle
\\[11pt]
&=\ 
\int_X \  \int_X \  K_2(y,x) \hsx \ov{K_1^*(y,x)} \  \td\mu(y) \  \td\mu(x) \ 
\\[11pt]
&=\ 
\int_X \  \int_X \  K_2(y,x) \hsx \ov{\ov{K_1(x,y)}} \  \td\mu(y) \  \td\mu(x) 
\\[11pt]
&=\ 
\int_X \  \int_X \  K_1(x,y) \hsy K_2(y,x) \  \td\mu(y) \  \td\mu(x) 
\\[11pt]
&=\ 
\int_X \  K_1 *K_2(x,x) \  \td\mu(x)
\\[11pt]
&=\ 
\int_X \  K (x,x) \  \td\mu(x).
\end{align*}
\\[-0.75cm]
\end{x}

\begin{x}{\small\bf REMARK} \ 
It can happen that $K_1 = K_2$ a.e. (so $T_{K_1} = T_{K_2}$), yet
\[
\int_X \  K_1(x,x) \  \td\mu(x)
\ \neq \ 
\int_X \  K_2(x,x) \  \td\mu(x).
\]

[E.g.: \ 
Take $X = Y = [0,1]$, $K_1 \equiv 0$, $K_2 = \chisubDelta$ ($\Delta$ the diagonal).]
\\[-.25cm]
\end{x}

\begin{x}{\small\bf THEOREM} \ 
Let \mX be a locally compact Hausdorff space, $\mu$ a $\sigma$-finite Radon measure on \mX.  
Suppose that $K \in L^2(X \times X)$ is separately continuous and $T_K$ is trace class $-$then 
the function 
\[
x \ra K(x,x)
\]
is integrable on \mX and 
\[
\tr(T_K) 
\ = \ 
\int_X \  K(x,x) \  \td\mu(x).
\]
\\[-1.25cm]
\end{x}

\[
\text{APPENDIX}
\]

{\small\bf LEMMA} \ 
Let \mM be a compact $C^\infty$ manifold, $\mu$ a smooth measure on \mM, 
$T:L^2(M) \ra C^{2k}(M)$ $(k > \frac{1}{4} \dim M)$ $-$then \mT is trace class.
\\[-.25cm]

PROOF \ 
Let $\Delta$ be a Laplacian on \mM and write 
\[
T 
\ = \ 
(1 - \Delta)^{-k} \hsx (1 - \Delta)^k \hsy T.
\]
Then
\[
(1 - \Delta)^k \hsy T (L^2(M)) \ \subset \  C(M) \ \subset \  L^\infty(M)
\]
so $(1 - \Delta)^k \hsy T$ is Hilbert-Schmidt.  
On the other hand, by Sobolev theory, 
\[
(1 - \Delta)^{-k} \hsy L^2(M) 
\ \subset \ 
H^{2k} (M) 
\ \subset \ 
C(M),
\]
thus  $(1 - \Delta)^{-k} \hsy T$ is Hilbert-Schmidt.


\chapter{
$\boldsymbol{\S}$\textbf{3}.\quad  THE LOCAL TRACE FORMULA}
\setlength\parindent{2em}
\setcounter{theoremn}{0}
\renewcommand{\thepage}{B II \S3-\arabic{page}}

\qquad Let \mG be a compact group.
\\[-.25cm]

\begin{x}{\small\bf NOTATION} \ 
Denote by $\pisubLR$ the representation of $G \times G$ on $L^2(G)$ given by
\[
(\pisubLR(x_1, x_2) \hsy f)(x) 
\ = \ 
f(x_1^{-1} \hsy x \hsy x_2) \qquad \text{(cf. A, III, \S2, \#1).}
\]
\\[-1.25cm]
\end{x}

\begin{x}{\small\bf LEMMA} \ 
$\pisubLR$ is unitary.
\\[-.2cm]

PROOF \ 
\allowdisplaybreaks
\begin{align*}
\norm{\pisubLR(x_1, x_2)\hsy f}^2 \ 
&=\ 
\int_G \ 
\abs{(\pisubLR(x_1, x_2)\hsy f)(x)}^2 \ \td_G(x)
\\[11pt]
&=\ 
\int_G \ 
\abs{f(x_1^{-1} \hsy x \hsy x_2)}^2 \ \td_G(x)
\\[11pt]
&=\ 
\int_G \
 \abs{f(x)}^2 \ \td_G(x)
\\[11pt]
&=\ 
\norm{f}^2.
\end{align*}

Given $f_1$, $f_2 \in C(G)$, define $f \in C(G \times G)$ by
\[
f(x_1,x_2) 
\ = \ 
f_1(x_1) f_2(x_2),
\]
and let
\[
\pisubLR(f) 
\ = \ 
\int_G \ \int_G \hsx f_1(x_1) \hsy f_2(x_2) \hsy \pisubLR(x_1,x_2) \ \td_G(x_1) \ \td_G(x_2).
\]
Then $\forall \ \phi \in L^2(G)$, 
\allowdisplaybreaks
\begin{align*}
(\pisubLR(f) \phi)(x) \ 
&=\ 
\int_G \ \int_G \ f_1(x_1) f_2(x_2) \phi(x_1^{-1} \hsy x \hsy x_2) \ \td_G(x_1) \ \td_G(x_2)
\\[11pt]
&=\ 
\int_G \ K_f(x,y) \hsy \phi(y) \ \td_G(y),
\end{align*}
where
\[
K_f(x,y)
\ = \ 
\int_G \ f_1(x z) f_2(z y) \ \td_G(z).
\]
Therefore $\pisubLR(f)$ is an integral operator on $L^2(G)$ with kernel $K_f(x,y)$.
\\[-.25cm]
\end{x}

\begin{x}{\small\bf CONSTRUCTION} \ 
\\[-.2cm]

\qquad \textbullet \quad
Given $f_1 \in C(G)$, put
\[
K_{f_1}(x,y) 
\ = \  
f_1(x y^{-1}) \qquad (x,y \in G).
\]
Then
\[
K_1 \in L^2(G \times G) \qquad (K_1 = K_{f_1})
\]
and 
\allowdisplaybreaks
\begin{align*}
\big(T_{K_1} \phi\big)(x) \ 
&=\ 
\int_G \hsx K_1(x,y) \phi(y) \hsx \td_G(y)
\\[11pt]
&=\ 
\int_G \hsx f_1(xy^{-1}) \phi(y) \hsx \td_G(y) 
\\[11pt]
&=\ 
\int_G \hsx f_1(y)\phi(y^{-1}x) \hsx \td_G(y)
\\[11pt]
&=\ 
\int_G \hsx f_1(y) (L(y) \phi)(x) \hsx \td_G(y)
\\[11pt]
&=\ 
(L(f)\phi) (x).
\end{align*}

\qquad \textbullet \quad
Given $f_2 \in C(G)$, put
\[
K_{f_2}(x,y) 
\ = \ 
f_2(x^{-1} y) \qquad (x,y \in G).
\]
Then
\[
K_2 \in L^2(G \times G) \qquad (K_2 = K_{f_2})
\]
and 
\allowdisplaybreaks
\begin{align*}
\big(T_{K_2} \phi\big)(x) \ 
&=\ 
\int_G \ K_2(x,y) \phi(y) \ \td_G(y)
\\[11pt]
&=\ 
\int_G \ f_2(x^{-1} \hsy y) \phi(y) \ \td_G(y)
\\[11pt]
&=\ 
\int_G \ f_2(y) \hsy \phi(x \hsy y) \ \td_G(y)  
\\[11pt]
&=\ 
\int_G \ f_2(y) \hsy (R(x) \phi)(y) \ \td_G(y)
\\[11pt]
&=\ 
(R(f)\phi) (x).
\end{align*}
\\[-.75cm]
\end{x}

\begin{x}{\small\bf LEMMA} \ 
Let $f_1$, $f_2 \in C(G)$ and let $f = f_1 \hsy f_2$ $-$then
\[
K_f 
\ = \ 
K_{f_1} * K_{f_2}.
\]

PROOF \ 
\allowdisplaybreaks
\begin{align*}
(K_1 * K_2) (x,y) \ 
&=\ 
\int_G \hsx K_1(x,z) K_2(z,y) \hsx \td_G(z)
\\[11pt]
&=\ 
\int_G \hsx f_1(x z^{-1}) f_2(z^{-1} y) \hsx \td_G(z)
\\[11pt]
&=\ 
\int_G \hsx f_1(fx z) f_2(z y) \hsx \td_G(z).
\end{align*}
Since the kernels of 
\[
\begin{cases}
\ T_{K_1} : L^2(G) \ra L^2(G) \\[3pt]
\ T_{K_2} : L^2(G) \ra L^2(G)
\end{cases}
\]
are square integrable, it follows that these operators are Hilbert-Schmidt.  
But
\[T_{K_1} \circ T_{K_2} 
\ = \ 
T_{K_1 * K_2} \qquad \text{(cf. $\S2$, $\#4$).}
\]
Therefore $T_{K_1 * K_2} $ is trace class, i.e., $T_{K_f}$ is trace class, i.e., $\pi_{L,R}(f)$ is trace class.
\\[-.25cm]
\end{x}

\begin{x}{\small\bf LEMMA} \ 
\[
\tr(\pi_{L,R}(f))
\ = \ 
\int_G \hsx K_f(x,x) \hsx \td_G(x) \qquad \text{(cf. $\S2$, $\#6$).}
\]
\\[-1.25cm]
\end{x}

\begin{x}{\small\bf RAPPEL} \ 
Let
\[
f \in \spanx_\Cx (L^2(G) * L^2(G)) \subset C(G).
\]
Then
\[
f(e) 
\ = \ 
\sum\limits_{\Pi \in \widehat{G}} \ \td_\Pi \tr(\Pi(f)) \qquad \text{(cf. I, $\S2$, $\#15$).}
\]

We have
\allowdisplaybreaks
\begin{align*}
K_f(x,y) \ 
&=\ 
\int_G \ f_1(xz)  \hsx f_2(zy) \ \td_G(z) 
\\[11pt]
&=\ 
\int_G \ f_1(u)\hsy  f_2(x^{-1} \hsy u \hsy y) \ \td_G(u)
\\[11pt]
&=\ 
\int_G \ f_1(u)  \hsy f_{2,x,y}(u) \ \td_G(u).
\end{align*}
Put now
\[
F_{x,y}(v)
\ = \ 
\int_G \ f_1(u)  \hsy f_{2,x,y}(v^{-1}\hsy u) \ \td_G(u).
\]
Then
\allowdisplaybreaks
\begin{align*}
F_{x,y}(v) \
&=\ 
\int_G \ f_1(u)  \hsx \widecheck{f}_{2, x, y}(u^{-1} v) \ \td_G(u)
\\[11pt]
&=\ 
f_1 * \widecheck{f}_{2, x, y}(v)
\end{align*}

$\implies$ 
\allowdisplaybreaks
\begin{align*}
K_f(x,y) \ 
&=\ 
F_{x,y}(e) 
\\[11pt]
&=\ 
\sum\limits_{\Pi \in \widehat{G}} \ 
\td_\Pi \hsy \tr(\Pi(F_{x,y}))
\\[11pt]
&=\ 
\sum\limits_{\Pi \in \widehat{G}} \ 
\td_\Pi \hsy \tr(\Pi(f_1 * \widecheck{f}_{2, x, y}))
\\[11pt]
&=\ 
\sum\limits_{\Pi \in \widehat{G}} \ 
\td_\Pi \hsy \tr(\Pi(f_1) \hsy \Pi(\widecheck{f}_{2, x, y}))
\\[11pt]
&=\ 
\sum\limits_{\Pi \in \widehat{G}} \ 
\td_\Pi \hsy \tr(\Pi(f_1) \hsy \Pi(\delta_x * \widecheck{f_2} * \delta_{y^{-1}}))
\\[11pt]
&=\ 
\sum\limits_{\Pi \in \widehat{G}} \ 
\td_\Pi  \hsy \tr(\Pi(f_1) \hsy \Pi(\delta_x) \hsy \Pi(\widecheck{f_2}) \hsy \Pi(\delta_{y^{-1}}))
\\[11pt]
&=\ 
\sum\limits_{\Pi \in \widehat{G}} \ 
\td_\Pi \hsy \tr(\Pi(\delta_{y^{-1}}) \hsy \Pi(f_1) \hsy \Pi(\delta_x) \hsy \Pi(\widecheck{f_2}))
\\[11pt]
&=\ 
\sum\limits_{\Pi \in \widehat{G}} \ 
\td_\Pi \tr(\Pi(y^{-1}) \hsy \Pi(f_1) \hsy \Pi(x)  \hsy \Pi(\widecheck{f_2}))
\end{align*}

$\implies$ 
\allowdisplaybreaks
\begin{align*}
\tr(\pi_{L,R}(f)) \ 
&=\ 
\int_G \hsx K_f(x,x) \hsx \td_G(x) 
\\[11pt]
&=\ 
\sum\limits_{\Pi \in \widehat{G}} \ 
\td_\Pi \hsy \tr\bigg(\bigg(
\int_G \ 
\Pi(x^{-1}) \hsy \Pi(f_1) \hsy \Pi(x) \hsx \td_G(x)\bigg)
\circ \Pi(\widecheck{f_2}) 
\bigg)
\\[11pt]
&=\ 
\sum\limits_{\Pi \in \widehat{G}} \ 
\td_\Pi \hsy \tr\bigg(\bigg(
\int_G \ 
\Pi(x) \hsy \Pi(f_1) \hsy \Pi(x^{-1}) \ \td_G(x)\bigg)
\circ \Pi(\widecheck{f_2}) 
\bigg).
\end{align*}
\\[-.75cm]
\end{x}

\begin{x}{\small\bf SUBLEMMA} \ 
$\forall \ \phi \in C(G)$, the operator
\[
\int_G \ \Pi(x) \Pi(\phi) \Pi(x^{-1}) \ \td_G(x)
\]
intertwines $\Pi$, hence is a scalar multiple of the identity (cf. I, $\S1$, $\#15$), call it $\lambda_\phi$.
\\[-.25cm]
\end{x}

\begin{x}{\small\bf \un{N.B.}} \ 
\[
\lambda_\phi
\ = \ 
\int_G \
\Pi(x) \hsx \Pi(\phi) \hsx \Pi(x^{-1})
\hsx \td_G(x)
\]
\qquad\qquad $\implies$
\[
\lambda_\phi \td_\Pi
\ = \ 
\tr\bigg( \int_G \ \Pi(x) \hsx \Pi(\phi) \hsx \Pi(x^{-1}) \ \td_G(x) \bigg)
\]
\qquad\qquad $\implies$
\[
\lambda_\phi
\ = \ 
\frac{\tr(\Pi(\phi))}{\td_\Pi}.
\]
Therefore
\allowdisplaybreaks
\begin{align*}
\tr(\pi_{L,R} (f)) \ 
&=\ 
\sum\limits_{\Pi \in \widehat{G}} \ \td_\Pi \tr 
\bigg(\bigg(
\int_G \ \Pi(x) \hsx \Pi(f_1)\hsx  \Pi(x^{-1}) \ \td_G(x)\bigg) \circ \Pi(\widecheck{f_2})\bigg)
\\[11pt]
&=\ 
\sum\limits_{\Pi \in \widehat{G}} \ \td_\Pi \lambda_{f_1} \hsx \tr(\Pi(\widecheck{f_2}))
\\[11pt]
&=\ 
\sum\limits_{\Pi \in \widehat{G}} \ \td_\Pi
\frac{\tr(\Pi(f_1))}{\td_\Pi} 
\hsx
\tr(\Pi(\widecheck{f_2}))
\\[11pt]
&=\ 
\sum\limits_{\Pi \in \widehat{G}} \ 
\tr(\Pi(f_1)) \hsx \tr(\Pi(\widecheck{f_2}))
\\[11pt]
&=\ 
\sum\limits_{\Pi \in \widehat{G}} \ \tJ(\Pi,f)
\end{align*}
if 
\[
\tJ(\Pi,f) 
\ = \ 
\tr(\Pi(f_1)) \hsx \tr(\Pi(\widecheck{f_2})).
\]
\\[-1.25cm]
\end{x}

\begin{x}{\small\bf SUBLEMMA} \ 
$\forall \ \phi \in C(G)$, 
\allowdisplaybreaks
\begin{align*}
\tr(\Pi(\widecheck{\phi})) \ 
&=\ 
\tr\bigg(\int_G \ \widecheck{\phi}(x) \hsx \Pi(x) \ \td_G(x)\bigg)
\\[11pt]
&=\ 
\tr\bigg(\int_G \ \phi(x^{-1})\hsx  \Pi(x) \ \td_G(x)\bigg)
\\[11pt]
&=\ 
\tr\bigg(\int_G \ \phi(x) \hsx \Pi(x^{-1}) \ \td_G(x)\bigg)
\\[11pt]
&=\ 
\tr\bigg(\int_G \ \phi(x) \hsx \Pi^*(x) \ \td_G(x)\bigg)
\\[11pt]
&=\ 
\tr(\Pi^*(\phi)).
\end{align*}
\\[-.75cm]
\end{x}

\begin{x}{\small\bf \un{N.B.}} \ 
Consequently, 
\[
\tJ(\Pi,f) 
\ = \ 
\tr(\Pi(f_1))\hsx  \tr(\Pi^*(f_2)).
\]

There is another way to manipulate
\[
\int_G \hsx K_f(x,x) \hsx \td_G(x)
\]
which then leads to a second formula for
\[
\tr(\pi_{L,R} (f)).
\]

To wit:
\[
\int_G \ K_f(x,x) \ \td_G(x)
\ = \ 
\int_G \ \int_G \hsx f_1(x \hsy z \hsy x^{-1}) \hsx f_2(z) \ \td_G(z) \ \td_G(x)
\]
or still, for any $y \in G$, 
\[
\int_G \ \int_G \ f_1(x \hsy z \hsy x^{-1}) \hsx f_2(y \hsy z \hsy y^{-1}) \ \td_G(z) \ \td_G(x).
\]
Now multiply through by $\td_G(y)$ and integrate with respect to $y$:
\allowdisplaybreaks
\begin{align*}
\tr(\pi_{L,R}(f)) \ 
&=\ 
\int_G \ \tr(\pi_{L,R}(f)) \ \td_G(y)
\\[11pt]
&=\  
\int_G \ \int_G \ \int_G \    f_1(x \hsy z \hsy x^{-1}) f_2(y \hsy z \hsy y^{-1})  \ \td_G(z) \ \td_G(x) \ \td_G(y)
\\[11pt]
&=\ 
\int_G \ \bigg(\int_G \ f_1(x \hsy z \hsy x^{-1}) \ \td_G(x) \bigg) \bigg(\int_G \  f_2(y \hsy z \hsy y^{-1}) \ \td_G(y) \ \td_G(z)
\\[11pt]
&=\  
\int_G \ \sO(f_1, z) \hsx \sO(f_2, z) \ \td_G(z).
\end{align*}
\\[-.75cm]
\end{x}

\begin{x}{\small\bf DEFINITION} \ 
Given $f = f_1 f_2$, the \un{local trace formula} is the relation
\allowdisplaybreaks
\begin{align*}
\sum\limits_{\Pi \in \widehat{G}} \ \tJ(\Pi,f) \ 
&=\ 
\sum\limits_{\Pi \in \widehat{G}} \ \tr(\Pi(f_1)) \hsx \tr(\Pi^*(f_2))
\\[11pt]
&=\ 
\int_G \ \sO(f_1, z) \hsx \sO(f_2, z) \ \td_G(z).
\end{align*}
\\[-.75cm]
\end{x}

Let \mG be a compact connected semisimple Lie group, $\Torus \subset G$ a maximal torus.
\\[-.25cm]

\begin{x}{\small\bf RAPPEL} \ 
For any continuous function $f \in C(G)$, 
\[
\int_G \hsx f(x) \hsx \td_G(x)
\ = \ 
\frac{1}{\abs{W}} \hsx 
\int_\Torus \ \abs{\Delta(t)}^2 \ 
\int_G \ f(x \hsy t \hsy x^{-1}) \ \td_G(x) \ \td_\Torus(t) \qquad \text{(cf. I,  $\S9$, $\#12$)}
\]
or still, 
\[
\int_G \ f(x) \ \td_G(x)
\ = \ 
\frac{1}{\abs{W}} \ 
\int_\Torus \ \abs{\Delta(t)}^2 \hsx \sO(f,t) \ \td_\Torus(t).
\]
As above, let $f = f_1 f_2$ $-$then
\allowdisplaybreaks
\begin{align*}
\int_G \ &K_f(x,x) \ \td_G(x) \ 
\\[11pt]
&=\ 
\int_G \ \int_G \ f_1(x \hsy z \hsy x^{-1}) \hsx  f_2(z)  \ \td_G(z) \ \td_G(x)
\\[11pt]
&=\
\frac{1}{\abs{W}} \hsx 
\int_\Torus \ \abs{\Delta(t)}^2 \hsx
\int_G \  f_1(x \hsy y \hsy t \hsy y^{-1} \hsy x^{-1})  f_2(y \hsy t \hsy y^{-1}) 
\ \td_G(y)
\ \td_T(t)
\ \td_G(x)
\\[11pt]
&=\ 
\frac{1}{\abs{W}} \ 
\int_\Torus \ \abs{\Delta(t)}^2 \hsx
\int_G \ 
\bigg(\int_G \ f_1(x \hsy y \hsy t \hsy y^{-1} \hsy x^{-1}) \ \td_G(x) \bigg) 
 f_2(y \hsy t \hsy y^{-1}) \hsx \td_G(y)  
\ \td_\Torus(t)
\\[11pt]
&=\ 
\frac{1}{\abs{W}} \ 
\int_\Torus \ \abs{\Delta(t)}^2 \hsx
\int_G \ 
\bigg(\int_G \ f_1(x \hsy t \hsy x^{-1}) \ \td_G(x) \bigg) 
 f_2(y \hsy t \hsy y^{-1}) \ \td_G(y) 
 \ \td_\Torus(t)
\\[11pt]
&=\ 
\frac{1}{\abs{W}} \ 
\int_\Torus \ \abs{\Delta(t)}^2 \hsx
\bigg(\int_G \ f_1(x \hsy t \hsy x^{-1}) \ \td_G(x) \bigg) 
\bigg(\int_G \ f_2(y \hsy t \hsy y^{-1}) \ \td_G(y)  \bigg) 
\ \td_\Torus(t)
\\[11pt]
&=\  
\frac{1}{\abs{W}} \ 
\int_\Torus \ \abs{\Delta(t)}^2 \hsx \sO(f_1,t) \hsx \sO(f_2,t) \ \td_\Torus(t).
\end{align*}
\end{x}


%% file: _C.tex
\chapter{
$\boldsymbol{\S}$\textbf{1}.\quad  TOPOLOGICAL TERMINOLOGY}
\setlength\parindent{2em}
\setcounter{theoremn}{0}
\renewcommand{\thepage}{C I \S1-\arabic{page}}

\begin{x}{\small\bf DEFINITION} \ 
A topological space \mX is \un{compact} if every open cover of \mX has a finite subcover.
\\[-.25cm]
\end{x}

\begin{x}{\small\bf DEFINITION} \ 
A topological space \mX is \un{locally compact} if every point in \mX has a neighborhood basis consisting of compact sets. 
\\[-.25cm]
\end{x}

\begin{x}{\small\bf LEMMA} \ 
A Hausdorff space \mX is locally compact iff every point in \mX has a compact neighborhood.  
\\[-.25cm]
\end{x}

\begin{x}{\small\bf APPLICATION} \ 
Every compact Hausdorff space is locally compact.
\\[-.25cm]
\end{x}

\begin{x}{\small\bf EXAMPLE} \ 
$\R$ is a locally compact Hausdorff space.
\\[-.25cm]
\end{x}

\begin{x}{\small\bf EXAMPLE} \ 
$\Q$ is a Hausdorff space but it is not locally compact ($\Q$ is first category while a locally compact Hausdorff space is second category).
\\[-.25cm]
\end{x}

\begin{x}{\small\bf LEMMA} \ 
An open subset of a locally compact Hausdorff space is locally compact.
\\[-.25cm]
\end{x}

\begin{x}{\small\bf LEMMA} \ 
A closed subset of a locally compact Hausdorff space is locally compact.
\\[-.25cm]
\end{x}

\begin{x}{\small\bf LEMMA} \ 
In a locally compact Hausdorff space, the intersection of an open set with a closed set is locally compact.
\\[-.25cm]
\end{x}

\begin{x}{\small\bf EXAMPLE} \ 
The semiclosed, semiopen interval $[0,1[$ is locally compact.
\\[-.5cm]

[In fact, 
\[
[0,1[ 
\ \ =\  
[-1,1[ \ \cap \  [0,1].]
\]
\\[-1.25cm]
\end{x}


\begin{x}{\small\bf DEFINITION} \ 
A topological group is a Hausdorff topological space \mG equipped with a group structure such that the function from 
$G \times G$ to \mG defined by $(x,y) \ra x y^{-1}$ is continuous or still, as is equivalent:
\\[-.5cm]

\textbullet \quad The function $G \times G \ra G$ that sends $(x,y) $ to $x y$ is continuous.
\\[-.5cm]

\textbullet \quad The function $G \ra G$ that sends $x$ to $x^{-1}$ is continuous.
\\[-.25cm]
\end{x}

If \mG is a topological group and if $H \subset G$ is a subgroup, then the set $G / H$ is to be given the quotient topology.
\\[-.25cm]

\begin{x}{\small\bf LEMMA} \ 
The space $G / H$ is Hasudorff iff \mH is closed.
\\[-.25cm]
\end{x}

\begin{x}{\small\bf DEFINITION} \ 
A \un{locally compact} (\un{compact}) group is a topological group \mG that is both locally compact (compact) and Hasudorff. 
\\[-.25cm]
\end{x}

\begin{x}{\small\bf LEMMA} \ 
If \mG is a locally compact group and if \mH is a closed subgroup, then $G / H$ is a locally compact Hausdorff space.
\\[-.25cm]
\end{x}

\begin{x}{\small\bf LEMMA} \ 
If \mG is a locally compact group and if \mH is a closed normal subgroup, then $G / H$ is a locally compact group. 
\\[-.25cm]
\end{x}

\begin{x}{\small\bf LEMMA} \ 
If \mG is a locally compact group  and if \mH is a locally compact subgroup, then \mH is closed in \mG.
\\[-.25cm]
\end{x}

\begin{x}{\small\bf LEMMA} \ 
If \mG is a locally compact group, then a subgroup \mH is open iff the quotient $G / H$ is discrete.
\\[-.25cm]
\end{x}

\begin{x}{\small\bf LEMMA} \ 
If \mG is a compact group, then a subgroup \mH is open iff the quotient $G / H$ is finite.
\\[-.25cm]
\end{x}


\begin{x}{\small\bf LEMMA} \ 
If \mG is a locally compact group, then every open subgroup of \mG is closed and every finite index closed subgroup of \mG is open.  
\\[-.25cm]
\end{x}

\begin{x}{\small\bf DEFINITION} \ 
A topologial space \mX is \un{totally disconnected} if the connected components of \mX are singletons.
\\[-.25cm]
\end{x}

\begin{x}{\small\bf EXAMPLE} \ 
$\Q$ is totally disconnected.
\\[-.25cm]
\end{x}

\begin{x}{\small\bf LEMMA} \ 
If \mG is a totally disconnected locally compact group, then $\{e\}$ has  a neighborhood basis consisting of open-compact subgroups.
\\[-.25cm]
\end{x}

\begin{x}{\small\bf LEMMA} \ 
If \mG is a totally disconnected compact group, then $\{e\}$ has  a neighborhood basis consisting of open-compact normal subgroups.
\\[-.25cm]
\end{x}

\begin{x}{\small\bf DEFINITION} \ 
A topological space \mX is \un{0-dimensional} if every point of \mX has a neighborhood basis consisting of open-closed sets.
\\[-.25cm]
\end{x}

\begin{x}{\small\bf EXAMPLE} \ 
$\Q$ is 0-dimensional.
\\[-.25cm]
\end{x}

\begin{x}{\small\bf LEMMA} \ 
A locally compact Hausdorff space is 0-dimensional iff it is totally disconnected.
\\[-.5cm]

[Note: \ 
In such a space, every point has a neighborhood basis consisting of open-compact sets.]
\\[-.25cm]
\end{x}

\begin{x}{\small\bf REMARK} \ 
It is false that the continuous image of a 0-dimensional locally compact Hausdorff space is again 0-dimensional.
\\[-.5cm]

[To see this, recall that every compact metric space is the continuous image of the Cantor set.]
\\[-.25cm]
\end{x}


\begin{x}{\small\bf LEMMA} \ 
If \mG is a locally compact 0-dimensional group and if \mH is a closed subgroup of \mG, then $G / H$ is 0-dimensional.
\\[-.25cm]
\end{x}

\begin{x}{\small\bf LEMMA} \ 
A 0-dimensional $\tT_1$ space is totally disconnected.
\\[-.25cm]
\end{x}

\begin{x}{\small\bf REMARK} \ 
There are totally disconnected metric spaces which are not 0-dimensional.
\\[-.25cm]
\end{x}


\chapter{
$\boldsymbol{\S}$\textbf{2}.\quad  INTEGRATION THEORY}
\setlength\parindent{2em}
\setcounter{theoremn}{0}
\renewcommand{\thepage}{C I \S2-\arabic{page}}

\qquad Let \mX be a locally compact Hausdorff space.
\\[-.25cm]

\begin{x}{\small\bf DEFINITION} \ 
A \un{Radon measure} is a measure $\mu$ defined on the Borel $\sigma$-algebra of \mX subject to the following conditions. 
\\[-.5cm]

1. \ $\mu$ is finite on compacta, i.e., for every compact set $K \subset X$, $\mu(K) < \infty$.
\\[-.5cm]

2. \ $\mu$ is outer regular, i.e., for every Borel set $A \subset X$, 
\[
\mu(A) \ = \ 
\inf\limits_{U \supset A} \mu(U), 
\]
where $U \subset X$ is open.
\\[-.5cm]

3. \ $\mu$ is inner regular, i.e., for every open set $A \subset X$, 
\[
\mu(A) \ = \ 
\sup\limits_{K \subset A} \mu(K), 
\]
where $K \subset X$ is compact.
\\[-.25cm]
\end{x}

\begin{x}{\small\bf RAPPEL} \ 
If \mX is a locally compact Hausdorff space and if \mX is second countable, then for any open subset $U \subset X$, 
there exist compact sets 
$K_1 \subset K_2 \subset \cdots$ 
such that 
\[
U \ = \ 
\bigcup\limits_{n = 1}^\infty \hsx K_n.
\]
\\[-1.25cm]
\end{x}

\begin{x}{\small\bf APPLICATION} \ 
If $(X, \mu)$ is a Radon measure space and if \mX is second countable, then \mX is $\sigma$-finite.
\\[-.25cm]
\end{x}

\begin{x}{\small\bf RIESZ REPRESENTATION THEOREM} \ 
Let \mX be a locally compact Hausdorff space.  
Suppose that 
$\Lambda : C_c(X) \ra \Cx$ 
is a positive linear functional $-$then there exists a unique Radon measure $\mu$ on \mX such that 
$\forall \ f \in C_c(X)$, 
\[
\Lambda \hsy  f 
\ = \ 
\int_X \ f(x) \ \td\mu(x).
\]
\\[-1.25cm]
\end{x}

Let \mG be a locally compact group.
\\[-.25cm]

\begin{x}{\small\bf DEFINITION} \ 
A \un{left Haar measure} on \mG is a Radon measure $\mu_G \neq 0$ which is left invariant, i.e., 
$\forall \ x \in G$ and $\forall$ Borel set $A \subset G$, $\mu_G(x A) = \mu_G(A)$.
\\[-.5cm]

[Note: \ 
Equivalently, a Radon measure $\mu \neq 0$ is a left Haar measure on \mG if 
$\forall \ f \in C_c(G)$ and $\forall \ y \in G$, 

\[
\int_G  f(y \hsy x) \ \td \mu(x) 
\ = \ 
\int_G  f(x) \ \td \mu(x).]
\]
\\[-1.25cm]
\end{x}

\begin{x}{\small\bf THEOREM} \ 
\mG admits a left Haar measure and if $\mu_{G_1}$, $\mu_{G_2}$ are two such, then 
$\mu_{G_1} = c \mu_{G_2}$ $(\exists \ c > 0)$.
\\[-.25cm]
\end{x}

\begin{x}{\small\bf LEMMA} \ 
Every nonempty open subset of \mG has positive left Haar measure.
\\[-.25cm]
\end{x}

\begin{x}{\small\bf LEMMA} \ 
Every compact subset of \mG has finite left Haar measure.
\\[-.25cm]
\end{x}

\begin{x}{\small\bf \un{N.B.}} \ 
The definition of a right Haar measure on \mG is analogous.
\\[-.25cm]
\end{x}

Given $x \in G$ and a Borel set $A \subset X$, let
\[
\mu_{G,x} (A) 
\ = \ 
\mu_G(A x).
\]
Then $\mu_{G,x}$ is a left Haar measure on \mG:
\[
\mu_{G,x} (y A) 
\ = \ 
\mu_G(y \hsy A \hsy x) 
\ = \ 
\mu_G(A \hsy x) 
\ = \ 
\mu_{G,x} (A).
\]
The uniqueness of left Haar measure now implies that there is a unique positive real number 
$\Delta_G(x)$ such that
\[
\mu_{G,x}
\ = \ 
\Delta_G(x) \mu_G.
\]
\\[-1.25cm]

\begin{x}{\small\bf LEMMA} \ 
$\Delta_G : G \ra \R_{> 0}^\times$ \ is independent of the choice of $\mu$.
\\[-.25cm]
\end{x}


\begin{x}{\small\bf LEMMA} \ 
$\Delta_G : G \ra \R_{> 0}^\times$ \  is a continuous homomorphism.
\\[-.25cm]
\end{x}

\begin{x}{\small\bf DEFINITION} \ 
$\Delta_G$ is called the \un{modular function} of \mG.
\\[-.25cm]
\end{x}

So, $\forall \ f \in C_c(G)$ and $\forall \ y \in G$, 
\[
\int_G \ f(x \hsy y^{-1}) \ \td \mu_G (x) 
\ = \ 
\Delta_G(y) \  \int_G \ f(x) \ \td \mu_G (x).
\]
\\[-1.25cm]

\begin{x}{\small\bf LEMMA} \ 
$\forall \ f \in C_c(G)$, 
\[
\int_G \ \widecheck{f}(x) \ \td \mu_G(x)
\ = \ 
\int_G \ \frac{f(x)}{\Delta_G(x)} \ \td \mu_G(x).
\]

[Note: \ 
As usual, $\widecheck{f}(x) = f(x^{-1})$.]
\\[-.25cm]
\end{x}

\begin{x}{\small\bf \un{N.B.}} \ 
The positive linear functional that assigns to each $f \in C_c(G)$ the common value of the two members of this equality is a right Haar integral. 
\\[-.25cm]
\end{x}

\begin{x}{\small\bf LEMMA} \ 
If $\phi:G \ra G$ is a topological automorphism, then there is a unique positive real number $\delta_G(\phi)$ such that
$\forall \ f \in C_c(G)$, 
\[
\int_G \ f(\phi^{-1}(x)) \ \td \mu_G(x) 
\ = \ 
\delta_G(\phi) \ \int_G \hsx f(x) \ \td \mu_G(x).
\]

[The positive linear functional
\[
f \ra \int_G \ f(\phi^{-1}(x)) \ \td \mu_G(x) 
\]
is a left Haar integral.]
\\[-.5cm]

[Note: \ 
The arrow $\phi \ra \delta_G(\phi)$ is a homomorphism: \ 
$\delta_G(\phi_1 \hsy \phi_2) \ = \ \delta_G(\phi_1) \hsy \delta_G(\phi_2)$.]
\\[-.25cm]
\end{x}

\begin{x}{\small\bf EXAMPLE} \ \  
If \mV is a \ real finite \ dimensional \ vector \ space \ and \ if 
$T: V \ra V$ is an invertible linear transformation, then per ``Lebesgue measure'',
\[
\int_V \hsx f(T^{-1}(x)) \ \td x
\ = \ 
\abs{\tdet T} \hsx \int_V \hsx f(x) \  \td x,
\]
so here
\[
\delta_V(T) 
\ = \ 
\abs{\tdet{T}}.
\]
\\[-1.25cm]
\end{x}

\begin{x}{\small\bf EXAMPLE} \ 
Define $I_y:G \ra G$ by $I_y(x) = y \hsy x \hsy y^{-1}$ $-$then 
\allowdisplaybreaks
\begin{align*}
\int_G \ f(I_y^{-1}(x)) \ \td\mu_G(x)
&=\ 
\int_G \ f(y^{-1}x \hsy y) \ \td\mu_G(x)
\\[11pt]
&=\
\int_G \ f(x \hsy y) \ \td\mu_G(x)
\\[11pt]
&=\
\Delta_G(y^{-1}) \ \int_G \ f(x) \ \td\mu_G(x),
\end{align*}
which implies that 
\[
\delta_G(I_y) 
\ = \ 
\Delta_G(y^{-1}).
\]
\\[-1.25cm]
\end{x}

\begin{x}{\small\bf LEMMA} \ 
If $\phi:G \ra G$ is a topological automorphism, then $\forall \ y \in G$,
\[
\Delta_G (\phi(y)) 
\ = \ 
\Delta_G(y).
\]

[On the one hand, 
\allowdisplaybreaks
\begin{align*}
\int_G \ f(\phi (x \hsy y^{-1})) \ \td\mu_G(x) 
&=\ 
\Delta_G(y) \ \int_G \ f(\phi(x)) \ \td\mu_G(x)
\\[11pt]
&=\ 
\Delta_G(y) \ \delta_G(\phi^{-1}) \ \int_G \ f(x) \ \td\mu_G(x)
\end{align*}
and, on the other hand,
\allowdisplaybreaks
\begin{align*}
\int_G \ f(\phi (x \hsy  y^{-1})) \  \td\mu_G(x) 
&=\ 
\int_G \ f(\phi (x) \hsy \phi(y^{-1})) \ \td\mu_G(x) 
\\[11pt]
&=\ 
\delta_G(\phi^{-1}) \ \int_G \ f(x \hsy \phi (y)^{-1}) \ \td\mu_G(x) 
\\[11pt]
&=\ 
\delta_G(\phi^{-1}) \Delta_G(\phi(y)) \ \int_G \ f(x) \ \td\mu_G(x).
\end{align*}
Therefore
\[
\Delta_G(y)
\ = \ 
\Delta_G(\phi(y)).]
\]
\\[-1.25cm]
\end{x}

\begin{x}{\small\bf LEMMA} \ 
If $G_1$, $G_2$ are locally compact groups and if $\mu_{G_1}$, $\mu_{G_2}$ are left Haar measures 
per $G_1$, $G_2$, then $\mu_{G_1} \times \mu_{G_2}$ is a left Haar measure per $G_1 \times G_2$ and
\[
\Delta_{G_1 \times G_2} (x_1, x_2) 
\ = \ 
\Delta_{G_1}(x_1)  \Delta_{G_2}(x_2).
\]
\\[-1.25cm]
\end{x}

Let \mG be a locally compact group, \mX and \mY two closed subgroups of \mG.
\\[-.25cm]

\begin{x}{\small\bf DEFINITION} \ 
The pair $(X,Y)$ is \un{admissible} if the following conditions are satisfied.
\\[-.5cm]

\qquad \textbullet \quad The intersection $X \cap Y$ is compact.
\\[-.5cm]

\qquad \textbullet \quad The multiplication $X \times Y \ra G$ is an open map.
\\[-.5cm]

\qquad \textbullet \quad The set of products $X Y$ exhausts \mG up to a set of Haar measure 0 (left or right).
\\[-.25cm]
\end{x}

\begin{x}{\small\bf EXAMPLE} \ 
Using the notation of $\#19$, work with $G_1 \times G_2$ and take 
$X = G_1 \times \{e_2\}$, $Y = \{e_1\} \times G_2$ $-$then the pair $(X,Y)$ is admissible.
\\[-.25cm]
\end{x}

\begin{x}{\small\bf THEOREM} \ 
Suppose that the pair $(X,Y)$ is admissible.  
Fix left Haar 
measures $\mu_X$, $\mu_Y$, on \mX, \mY $-$then there is a unique left Haar measure $\mu_G$ on \mG such 
that $\forall \ f \in C_c(G)$, 
\[
\int_G \hsx f \hsx \td\mu_G 
\ = \ 
\int_{X \times Y} \hsx f(x \hsy y) \hsx \frac{\Delta_G(y)}{\Delta_Y(y)} \ \td\mu_X(x) \  \td\mu_Y(y).
\] 
\\[-1.25cm]
\end{x}

\begin{x}{\small\bf \un{N.B.}} \ 
Specializing the setup to that of \# 21 leads back to \#19.  
\\[-.5cm]

[Note: \ 
\begin{align*}
\Delta_{G_1 \times G_2} (e_1, x_2) \ 
&= \ 
\Delta_{G_1} (e_1) \Delta_{G_2} (x_2) 
\\[5pt]
&=\ 
\Delta_{G_2} (x_2)
\end{align*}
thereby cancelling the factor in the denominator.]
\\[-.25cm]
\end{x}

\begin{x}{\small\bf LEMMA} \ 
If \mG is a locally compact group and if $H \subset G$ is a closed normal subgroup, then 
$\restr{\Delta_G}{H} = \Delta_H$.
\\[-.25cm]
\end{x}

\begin{x}{\small\bf APPLICATION} \ 
In the setup of \# 22, assume in addition that \mY is normal $-$then $\forall \ f \in C_c(G)$, 
\[
\int_G \ f \ \td \mu_G 
\ = \ 
\int_{X \times Y}  \ f(x y) \ \td \mu_X(x) \ \td \mu_Y(y).
\]

[Note: \ 
Given $x \in X$, the restriction 
\[
\begin{cases}
\ I_{x^{-1}} : Y \ra Y \\[3pt]
\ \hspace{1.1cm} y \ra x^{-1} \hsy y \hsy x
\end{cases}
\]
is an automorphism of \mY and 
\[
\Delta_G(x \hsy y) 
\ = \ 
\Delta_X(x) \hsx \Delta_Y(y) \hsx \delta_Y \big(I_{x^{-1}}\big).]
\]
\\[-1.25cm]
\end{x}

Let \mG be a locally compact group, \mX and \mY two closed subgroups of \mG.
\\[-.25cm]

\begin{x}{\small\bf DEFINITION} \ 
\mG is the \un{topological semidirect product} of \mX and \mY if every element $z \in G$ can be expressed 
in a unique manner as a product $z = x y$ $(x \in X, \ y \in Y)$ and if the multiplication 
$X \times Y \ra G$ is a homeomorphism.
\\[-.25cm]
\end{x}

\begin{x}{\small\bf \un{N.B.}} \ 
A priori, the multiplication 
$X \times Y \ra G$ 
is a continuous bijection, thus the condition is satisfied if the multiplication 
$X \times Y \ra G$ 
is an open map, this being automatic whenever \mG is second countable.
\\[-.5cm]

[Under these circumstances, \mG is the union of a sequence of compact sets (cf. \#2), so the same is true of 
$X \times Y$.  
But \mG is a locally compact Hausdorff space, hence is a Baire space.]
\\[-.5cm]

[Note: \ 
If \mA is a Baire space and if $\{A_n: \ n \in \N\}$ is a closed covering of \mA, then at least one of the $A_n$ must 
contain an open set.]
\\[-.25cm]
\end{x}

If \mG is the topological semidirect product of \mX and \mY, then $X \cap Y = \{e\}$ and the pair 
$(X,Y)$ is admissible.  
Therefore the theory is applicable in this situation.
\\[-.25cm]

\begin{x}{\small\bf \un{N.B.}} \ 
In general, the arrow $(x,y) \ra x y$ is not an isomorphism of groups but this will be the case if every element of 
\mX commutes with every element of \mY or, equivalently, if \mX and \mY are normal subgroups of \mG, i.e., 
if \mG is the topological direct product of \mX and \mY.
\end{x}


\chapter{
$\boldsymbol{\S}$\textbf{3}.\quad  UNIMODULARITY}
\setlength\parindent{2em}
\setcounter{theoremn}{0}
\renewcommand{\thepage}{C I \S3-\arabic{page}}

\qquad Let \mG be a locally compact group, $\mu_G$ a left Haar measure on \mG.
\\[-.25cm]

\begin{x}{\small\bf DEFINITION} \ 
\mG is \un{unimodular} if $\Delta_G \equiv 1$.
\\[-.25cm]
\end{x}

\begin{x}{\small\bf \un{N.B.}} \ 
\mG is unimodular iff $\mu_G$ is a right Haar measure on \mG.
\\[-.25cm]
\end{x}

\begin{x}{\small\bf EXAMPLE} \ 
Take for \mG the group of all real matrices of the form
\[
\begin{pmatrix}
1 &&x\\
0 &&y
\end{pmatrix}
\qquad (y \neq 0)
\]
$-$then
\[
\Delta_G
\begin{pmatrix}
1 &&x\\
0 &&y
\end{pmatrix}
\ = \ 
\abs{y},
\]
thus \mG is not unimodular.
\\[-.25cm]
\end{x}

\begin{x}{\small\bf LEMMA} \ 
\mG is unimodular iff $\forall \ f \in C_c(G)$, 
\[
\int_G \ f(x^{-1}) \ \td\mu_G(x) 
\ = \ 
\int_G \ f(x) \ \td\mu_G(x)
\qquad \text{(cf. $\S2$, $\#13$).}
\]
\\[-1.25cm]
\end{x}

\begin{x}{\small\bf LEMMA} \ 
\\[-.5cm]

\qquad \textbullet \quad Every locally compact abelian group is unimodular.
\\[-.5cm]

\qquad \textbullet \quad Every compact group is unimodular.
\\[-.5cm]

\qquad \textbullet \quad Every discrete group is unimodular.
\\[-.25cm]
\end{x}

\begin{x}{\small\bf LEMMA} \ 
 Every locally compact group that coincides with its closed commutator subgroup is unimodular.
\\[-.25cm]
\end{x}

\begin{x}{\small\bf LEMMA} \ 
Every open subgroup of a unimodular locally compact group is unimodular.
\\[-.25cm]
\end{x}

\begin{x}{\small\bf LEMMA} \ 
Every closed normal subgroup of a unimodular locally compact group is unimodular.
\\[-.5cm]

[Note: \ 
A closed subgroup of a unimodular locally compact group is not necessarily unimodular.]
\\[-.25cm]
\end{x}

\begin{x}{\small\bf LEMMA} \ 
Let \mG be a locally compact group, $Z(G)$ its center $-$then \mG is unimodular iff $G/Z(G)$ is unimodular.
\\[-.25cm]
\end{x}

Let $G$ be a locally compact group, $H \subset \hsy G$ a closed subgroup ($H$ is then a locally compact subgroup) (cf. \S 1, \#8).
\\[-.25cm]

\begin{x}{\small\bf DEFINITION} \ 
\mH is a \un{cocompact} subgroup if the quotient $G/H$ is compact.
\\[-.25cm]
\end{x}

\begin{x}{\small\bf LEMMA} \ 
If \mG admits a unimodular cocompact subgroup $H \subset G$, then \mG is unimodular.
\end{x}


\chapter{
$\boldsymbol{\S}$\textbf{4}.\quad  INTEGRATION ON HOMOGENEOUS SPACES}
\setlength\parindent{2em}
\setcounter{theoremn}{0}
\renewcommand{\thepage}{C I \S4-\arabic{page}}

\qquad Let \mG be a locally compact group, $H \subset G$ a closed subgroup.
\\[-.25cm]

\begin{x}{\small\bf \un{N.B.}}\ 
The quotient $G / H$ is a locally compact Hausdorff space (cf. $\S 1, \ \#14$).
\\[-.25cm]
\end{x}

Fix left Haar measures
\[
\begin{cases}
\ \mu_G \hspace{0.25cm} \text{on}  \hspace{0.25cm} G\\[3pt]
\ \mu_H \hspace{0.25cm} \text{on}  \hspace{0.25cm}  H
\end{cases}
.
\]
\\[-1.25cm]

\begin{x}{\small\bf NOTATION} \ 
Given $f \in C_c(G)$, define $f^H \in C_c(G / H)$ by the rule 
\[
f^H (x H) 
\ = \ 
\int_H \ f(x \hsy y) \ \td \mu_H (y).
\]
\\[-1.25cm]
\end{x}

\begin{x}{\small\bf LEMMA} \ 
The arrow
\[
f \ra f^H
\]
sends $C_c(G)$ onto $C_c(G/H)$.
\\[-.25cm]
\end{x}

\begin{x}{\small\bf DEFINITION} \ 
A Radon measure $\mu \neq 0$ on the Borel $\sigma$-algebra of $G/H$ is said to be an \un{invariant measure} if 
$\forall \ x \in G$ and $\forall$ Borel set $A \subset G/H$, $\mu(x A) = \mu(A)$.
\\[-.5cm]

[Note: \ 
If $H = \{e\}$, then ``invariant measure'' = ``left Haar measure''.]
\\[-.25cm]
\end{x}

\begin{x}{\small\bf THEOREM} \ 
There exists an invariant measure $\mu_{G/H}$ on $G/H$ iff $\restr{\Delta_G}{H} = \Delta_H$ and when this is so, 
$\mu_{G/H}$ is unique up to a positive scalar factor.
\\[-.5cm]

[Note: \ 
Matters are automatic if \mH is compact or if both \mG and \mH are unimodular.]
\\[-.25cm]
\end{x}

\begin{x}{\small\bf \un{N.B.}} \ 
If \mH is a normal closed subgroup of \mG, then $\restr{\Delta_G}{H} = \Delta_H$.
\\[-.5cm]

[For a left Haar measure on $G/H$ is an invariant measure.]
\\[-.25cm]
\end{x}

\begin{x}{\small\bf THEOREM} \ 
There is a unique choice for $\mu_{G/H}$ such that $\forall \ f \in C_c(G)$, 
\[
\int_G \ f(x) \ \td \mu_G(x) 
\ = \ 
\int_{G/H} \ f^H (\xdot) \ \td\mu_{G/H}(\xdot) \qquad (\xdot = x H).
\]

[Note: \ 
Bear in mind that $\mu_G$, $\mu_H$ have been fixed at the beginning.]
\\[-.25cm]
\end{x}

\begin{x}{\small\bf \un{N.B.}} \ 
This formula is valid for all $f \in L^1(G)$.
\\[-.25cm]
\end{x}

\begin{x}{\small\bf LEMMA} \ 
Let $H_1 \subset G$, $H_2 \subset G$ be closed subgroups of \mG with $H_1 \subset H_2$ $-$then 
$G/H_2$ and $H_2/H_1$ admit finite invariant measures iff $G/H_1$ admits a finite invariant measure.
\\[-.25cm]
\end{x}

\begin{x}{\small\bf APPLICATION} \ 
If $G/H$ has a finite invariant measure and if \mH is unimodular, then \mG is unimodular.
\\[-.5cm]

[Let \mK be the kernel of $\Delta_G$ $-$then $\restr{\Delta_G}{H} = \Delta_H \equiv 1$, thus $H \subset K$ and so $G/K$ has a finite invariant measure (as does $(K/H)$.  
But $G/K$ is a locally compact group.  
Therefore $G/K$ is actually a compact group (its Haar measure being finite) and this implies that $\Delta_G(G)$ is a compact subgroup of $\R_{>0}^\times$, 
hence $\Delta_G(G) = \{1\}$, i.e., \mG is unimodular.]
\\[-.25cm]
\end{x}

\begin{x}{\small\bf \un{N.B.}} \ 
Suppose that $H \subset G$ is a unimodular cocompact subgroup $-$then $G/H$ admits a finite invariant measure $\mu_{G/H}$.
\\[-.5cm]

[In fact, $G$ is necessarily unimodular (cf. \S3, \#11), from which the existence 
of $\mu_{G/H}$.  
But $\mu_{G/H}$ is Radon, hence finite on compacta, hence in particular, 
\[
\mu_{G/H} (G/H) \ < \ \infty.]
\]

[Note: \ 
Take $G = \SL(2,\R)$ and let
\[
H 
\ = \ 
\bigg\{X : X = 
\begin{pmatrix}
a &b\\
0 &d
\end{pmatrix}
\in \SL(2,\R) \bigg\}
.
\]
Then \mG is unimodular but \mH is not unimodular.  
Therefore $G/H$ does not admit an invariant measure even though \mH is a cocompact subgroup.]
\\[-.25cm]
\end{x}

\begin{x}{\small\bf LEMMA} \ 
Let $H_1 \subset G$, $H_2 \subset G$ be closed subgroups of \mG with $H_1$ normalizing $H_2$ and $H_1 H_2$ closed in \mG 
$-$then the following are equivalent.
\\[-.5cm]

\qquad \textbullet \quad $H_1 H_2 / H_1$ admits a finite invariant measure.
\\[-.5cm]

\qquad \textbullet \quad $H_2 / H_1 \hsx \cap \hsx H_2$ admits a finite invariant measure.
\\[-.5cm]

[Note: \ 
There is a commutative diagram
\[
\begin{tikzcd}
{H_2} 
\arrow[rr,shift right=0.45,dash] \arrow[rr,shift right=-0.45,dash] 
\ar{d} 
&&{H_2} \ar{d}\\
{H_2/H_1 \cap H_2} \ar{rr}[swap]{\phi}
&&{H_1 H_2 / H_1} 
\end{tikzcd}
,
\]
where

\[
\phi(x_2 (H_1 \cap H_2)) 
\ = \ 
x_2 H_1.
\]
The vertical arrows are continuous and open.  
Therefore the bottom horizontal arrow is a homeomorphism.]
\\[-.25cm]
\end{x}

\begin{x}{\small\bf APPLICATION} \ 
Suppose that \mG is the topological semidirect product
of \mX and \mY (cf. $\S 2, \# 26$) and take \mY normal $-$then $G = X Y$ and $X \cap Y = \{e\}$.  
Therefore $G/X$ has a finite invariant measure iff \mY has a finite invariant measure.
\end{x}


\chapter{
$\boldsymbol{\S}$\textbf{5}.\quad  INTEGRATION ON LIE GROUPS}
\setlength\parindent{2em}
\setcounter{theoremn}{0}
\renewcommand{\thepage}{C I \S5-\arabic{page}}

\qquad Suppose that \mM is an orientable $n$-dimensional $C^\infty$ manifold which we take to be second countable.  
Let $\omega$ be a positive $n$-form on \mM $-$then the theory leads to a positive linear functional
\[
f \ra \int_M \ f \omega \qquad (f \in C_c(M))
\]
from which a Radon measure $\mu_\omega$.
\\[-.25cm]

Assume now that \mG is a Lie group with Lie algebra $\mathfrak{g}$.  
Let $L_x:G \ra G$ be left translation $y \ra x y$ by $x$.
\\[-.25cm]

\begin{x}{\small\bf DEFINITION} \ 
A differential form $\omega$ on \mG is \un{left invariant} if $\forall \ x \in G$, $L_x^* \hsy \omega = \omega$.
\\[-.25cm]
\end{x}

\begin{x}{\small\bf NOTATION} \ 
Given $X \in \fg$, let $\widetilde{X}$ be the corresponding left invariant vecor field on \mG.
\\[-.25cm]
\end{x}

Let $n = \dim G$ (= $\dim \fg$) and fix a basis $X_1, \ldots, X_n$ for $\fg$.  
Define 1-forms $\omega^1, \ldots, \omega^n$ on \mG by the condition 
$\omega^i (\widetilde{X}_j) = \delta_j^i$.
\\[-.25cm]

\begin{x}{\small\bf \un{N.B.}} \ 
The $\omega^i$ are left invariant.
\\[-.25cm]

Put
\[
\omega 
\ = \ 
\omega^1 \wedge \cdots \wedge \omega^n.
\]
Then $\forall \ x \ \in G$, 
\begin{align*}
L_x^* \hsy \omega \ 
&= \ 
L_x^* (\omega^1 \wedge \cdots \wedge \omega^n)
\\[11pt]
&=\ 
L_x^*\omega^1 \wedge \cdots \wedge L_x^*\omega^n
\\[11pt]
&=\ 
\omega^1 \wedge \cdots \wedge \omega^n
\\[11pt]
&=\ 
\omega.
\end{align*}
I.e.: $\omega$ is a left invariant $n$-form on \mG.
\\[-.25cm]
\end{x}

\begin{x}{\small\bf LEMMA} \ 
$\omega$ is nowhere vanishing on \mG.
\\[-.25cm]
\end{x}

\begin{x}{\small\bf LEMMA} \ 
\mG can be oriented so as to render $\omega$ positive. 
\\[-.25cm]

[Note: \ 
The orientation of \mG depends on the choice of a basis for $\fg$.  
If $Y_1, \ldots, Y_n$ is another basis, then the resulting orientation of \mG does not change iff the linear transformation 
$X_i \ra Y_i$ $(1 \leq i \leq n)$ has positive determinant.]
\\[-.5cm]
\end{x}

\begin{x}{\small\bf SCHOLIUM} \ 
The assignment 
\[
f \ra \int_G \hsx f \omega \qquad (f \in C_c(G))
\]
is a positive linear functional.
\\[-.25cm]
\end{x}

\begin{x}{\small\bf LEMMA} \ 
The Radon measure $\mu_\omega$ is a left Haar measure.
\\[-4pt]

PROOF \ 
$\forall \ x \in G$, $L_x:G \ra G$ is an orientation preserving diffeomorphism, so $\forall \ f \in C_c(G)$, 
\begin{align*}
\int_G \ f \hsx \td \mu_\omega \ 
&=\ 
\int_G \ f \omega
\\[11pt]
&=\ 
\int_G \ (f \circ L_x) \hsy L_x^* \hsy \omega
\\[11pt]
&=\ 
\int_G \ (f \circ L_x) \hsy \omega
\\[11pt]
&=\ 
\int_G \ (f \circ L_x)  \ \td \mu_\omega.
\end{align*}
\\[-.75cm]
\end{x}

\begin{x}{\small\bf REMARK} \ 
Any subset \mS of \mG which is contained in an at most countable union of smooth images of $C^\infty$ manifolds of dimension $< \ \dim G$ 
has zero left Haar measure.
\\[-.25cm]
\end{x}

\begin{x}{\small\bf THEOREM} \ 
$\forall \ x \in G$, 
\[
\Delta_G(x) 
\ = \ 
\frac{1}{\abs{\det \Ad(x)}}.
\]
\\[-1.25cm]
\end{x}


\begin{x}{\small\bf EXAMPLE} \ 
Every connected nilpotent Lie group \mN is unimodular.
\\[-.5cm]

[If $X \in \fn$ (the Lie algebra of \mN), then ad(\mX) is nilpotent, thus tr(ad$(X)) = 0$ and so
\begin{align*}
\det \Ad(\exp X) \ 
&=\ 
\det e^{\ad(X)} 
\\[8pt]
&=\ 
e^{\tr(\ad(X))}
\\[8pt]
&=\
1.]
\end{align*}
\\[-1.25cm]
\end{x}

\begin{x}{\small\bf LEMMA}\ 
A 1-dimensional representation of a connected semisimple Lie group is trivial.
\\[-.25cm]
\end{x}

\begin{x}{\small\bf APPLICATION} \ 
The restriction of $\Delta_G$ to any semisimple analytic subgroup of \mG is $\equiv 1$.
\\[-.25cm]
\end{x}

\begin{x}{\small\bf THEOREM} \ 
Suppose that \mG is a reductive Lie group in the Harish-Chandra class $-$then \mG is unimodular.
\\[-0.25cm]

PROOF \ 
First decompose \mG as the product $^0G \times V$, where \mV is a central vector group (possibly trivial) and 
\[
^0G 
\ = \ 
\bigcap\limits_\chi \ \Ker \chi ,
\]
the $\chi$ running through the set of continuous homomorphisms $G \ra \R_{>0}^\times$.  
This done, take for a left Haar measure on \mG the product of the left Haar measures on $^0G$ and \mV.  
Since \mV is unimodular, it will be enough to deal with $^0G$ (cf. $\S 2, \ \# 19$).  
Fix a maximal compact subgroup \mK of \mG $-$then \mK is a maximal compact subgroup of $^0G$ and 
$^0G = K G_\tss$, thus $\forall \ k \in K$, $\forall \ x \in G_\tss$, 
\[
\Delta_{\hsy ^0G} (k x) 
\ = \ 
\Delta_{\hsy ^0G} (k) \hsx  \Delta_{\hsy ^0G}  (x) 
\ = \ 
1 \cdot 1 
\ = \ 
1.
\]

[Note: \ 
$G_\tss$ is the analytic subgroup of \mG corresponding to $\fg_\tss$ (the ideal in $\fg$ spanned by $[\fg.\fg]$).  
It is closed and normal.]
\\[-.25cm]
\end{x}

Maintaining the supposition that \mG is a reductive Lie group in the Harish-Chandra class, consider an Iwasawa decomposition $G = KAN$.
\\[-.25cm]

\begin{x}{\small\bf \un{N.B.}} \ 
\mN is a normal subgroup of $AN$ and $AN$ is the topological semidirect product of \mA and \mN.
\\[-.5cm]

[Note: \ 
$A N$ is second countable so there are no technical issues.]
\\[-.25cm]
\end{x}

\begin{x}{\small\bf LEMMA} \ 
\[
\Delta_{AN} (an) 
\ = \ 
\frac{1}{\abs{\det \Ad(an)}} 
\ = \ 
\frac{1}{e^{2 \rho (\log a)}}.
\]

[Note: \ 
Here $2 \rho$ is the sum of the positive roots of $(\fg, \fa)$ counted with multiplicities.]
\\[-.25cm]
\end{x}

Since the pair $(K, AN)$ is admissible and since $\Delta_G \equiv 1$, it follows from $\S 2, \ \#22$ that $\forall \ f \in C_c(G)$, 
\begin{align*}
\int_G \hsx f \hsx \td \mu_G \ 
&=\ 
\int_{K \times A N} \hsx f(k \hsy a \hsy n) \frac{\Delta_G(an)}{\Delta_{AN}(an)} \hsx \td \mu_K (k)\hsx \td \mu_{AN}(an)
\\[8pt]
&=\ 
\int_{K \times A N} \hsx f(k \hsy a \hsy n) \frac{1}{\Delta_{AN}(an)} \hsx \td \mu_K (k) \hsx \td \mu_{AN}(an)
\\[8pt]
&=\ 
\int_{K \times A N} \hsx f(k \hsy a \hsy n) e^{2 \rho (\log a)}  \hsx \td \mu_K (k) \hsx \td \mu_{AN}(an)
\\[8pt]
&=\ 
\int_{K \times A \times N} \hsx f(k \hsy a \hsy n) e^{2 \rho (\log a)}  \hsx \td \mu_K (k) \hsx \td \mu_{A}(a)  \hsx \td \mu_{N}(n).
\end{align*}

[Note: \ 
To be completely precise, fix left Haar measures $\mu_K$, $\mu_A$, $\mu_N$ on \mK, \mA, \mN $-$then there is a unique determination 
of the left Haar measure $\mu_G$ on \mG such that for any $f \in C_c(G)$ the function 
\[
(k, a, n) \ra f(k \hsy a \hsy n)
\]
lies in 
\[
C_c(K \times A \times N)
\]
and 
\[
\int_G \hsx f \hsx \td \mu_G \ 
\ = \ 
\int_{K \times A \times N} \ f(k \hsy a \hsy n) e^{2 \rho (\log a)} \ \td \mu_K (k) \ \td \mu_A (a) \ \td \mu_N (n).]
\]
\\[-0.25cm]

\begin{x}{\small\bf LEMMA} \ 
\begin{align*}
\Delta_{AN}(an) \ 
&=\ 
\Delta_A(a) \hsx \Delta_N(n) \hsx \delta_N(I_{a^{-1}}) \qquad \text{(cf. $\S 2, \ \# 25$)}
\\[11pt]
&=\ 
\delta_N(I_{a^{-1}}).
\end{align*}

[Note: \ 
\mA is abelian and \mN is nilpotent \ldots \hsy .]
\\[-.25cm]

So, $\forall \ f \in C_c(G)$, 
\allowdisplaybreaks
\begin{align*}
\int_{K \times N \times A} \ 
&f(k \hsy a \hsy n) \hsx \td \mu_K(k) \ \td \mu_N(n) \ \td \mu_A(a) 
\\[11pt]
&=\ 
\int_{K \times N \times A} \hsx f(k \hsy a \hsy a^{-1} n \hsy a) \ \td \mu_K(k) \ \td \mu_N(n) \ \td \mu_A(a) 
\\[11pt]
&=\ 
\int_{K \times N \times A} \ f(k a I_{a^{-1}} (n)) \ \td \mu_K(k) \ \td \mu_N(n) \ \td \mu_A(a) 
\\[11pt]
&=\ 
\int_{K \times N \times A} \ f(k \hsy a \hsy n) \hsx \delta_N(I_a)  \ \td \mu_K(k) \ \td \mu_A(a) \ \td \mu_N(n) 
\qquad \text{(cf. $\S 2, \ \#15$)}
\\[11pt]
&=\ 
\int_{K \times A \times N} \ f(k \hsy a \hsy n) \hsx \Delta_{AN}(a^{-1} n) \ \td \mu_K(k) \ \td \mu_A(a) \ \td \mu_N(n)
\\[11pt]
&=\ 
\int_{K \times A \times N} \ f(k \hsy a \hsy n) \hsx e^{2 \rho (\log a)} \ \td \mu_K(k) \ \td \mu_A(a) \ \td \mu_N(n)  
\\[11pt]
&=\ 
\int_G \hsx f \hsx \td \mu_G.
\end{align*}

[Note: \ 
As a corollary, 
\begin{align*}
\int_{A \times N \times K} \
&f(a \hsy n \hsy k) \ \td \mu_A(a) \ \td \mu_N(n) \ \td \mu_K(k)  
\\[11pt]
&=\ 
\int_{K \times N \times A} \ \widecheck{f}(k^{-1} n^{-1} a^{-1}) \ \td \mu_K(k) \ \td \mu_N(n) \ \td \mu_A(a) 
\\[11pt]
&=\ 
\int_{K \times N \times A} \ \widecheck{f}(k \hsy n \hsy a) \ \td \mu_K(k) \ \td \mu_N(n) \ \td \mu_A(a) 
\qquad \text{(cf. $\S 3, \ \#4$)} 
\\[11pt]
&=\ 
\int_G \ \widecheck{f} \ \td \mu_G
\\[11pt]
&=\ 
\int_G \ f  \ \td\mu_G ,
\end{align*}
\mG being unimodular (cf. $\# 13$).]
\\[-.25cm]
\end{x}

Let \mM be the centralizer of $\fa$ in \mK and put $\ov{N} = \Theta N$ $-$then the map
\[
(\ov{n}, m, a, n) \ra \ov{n} \hsy m \hsy a \hsy n
\]
is an open bijection of 
$\ov{N} \times M \times A \times N$ onto an open submanifold $\ov{N} M A N \subset G$.
\\[-.25cm]

\begin{x}{\small\bf LEMMA} \ 
The complement of $\ov{N} M A N$ in \mG is a set of Haar measure 0.
\\[-.25cm]

[Using the Bruhat decomposition, the said complement is seen to be a finite union of smooth images of $C^\infty$ manifolds of 
dimension $< \dim G$ so one can quote $\# 8$.]
\\[-.25cm]
\end{x}

The pair $(\ov{N}, MAN)$ is therefore admissible, hence $\forall \ f \in C_c(G)$ (cf. $\S2, \ \#22$), 
\begin{align*}
\int_G \ f  \ \td \mu_G \ 
&=\ 
\int_{\ov{N} \times MAN} \ 
f(\ov{n}\hsy m \hsy a \hsy n) \hsx \frac{\Delta_G(m \hsy a \hsy n)}{\Delta_{MAN}(m \hsy a \hsy n)} 
\ \td \mu_{\ov{N}}(\ov{n}) \ \td \mu_{MAN}(m \hsy a \hsy n)
\\[11pt]
&=\ 
\int_{\ov{N} \times MAN} \ f(\ov{n} \hsy m \hsy a \hsy n) \hsx \frac{1}{\Delta_{MAN}(m \hsy a \hsy n)} 
\ \td \mu_{\ov{N}}(\ov{n}) \ \td \mu_{MAN}(m \hsy a \hsy n)
\\[11pt]
&=\ 
\int_{\ov{N} \times MAN} \ 
f(\ov{n} \hsy m \hsy a \hsy n) \hsx e^{2 \rho (\log a)}\hsx \td \mu_{\ov{N}}(\ov{n}) 
\ \td \mu_{MAN}(m \hsy a \hsy n)
\\[11pt]
&=\ 
\int_{\ov{N} \times M \times A \times N} \ 
f(\ov{n} \hsy m \hsy a \hsy n) \hsx e^{2 \rho (\log a)}\hsx \td \mu_{\ov{N}}(\ov{n}) 
\ \td \mu_{M}(m)
\ \td \mu_{A}(a)
\ \td \mu_{N}(n).
\end{align*}
\\[-.75cm]

\begin{x}{\small\bf RAPPEL} \ 
Let \mV be a finite dimensional real Hilbert space $-$then the canonical Haar measure $\td V$ on \mV is that in which the parallelpiped 
determined by an orthonormal basis has unit measure.
\\[-.5cm]

[Spelled out, if $\{X_1, \ldots, X_n\}$ is an orthonormal basis for \mV and if \mQ is the set of all points 
$X = \ds \sum\limits_{i = 1}^d \hsx c_i X_i \ (c_i \in \R)$ with $0 \leq c_i \leq 1$, then $\ds\int_Q \hsx \td V = 1$.]
\\[-.25cm]

[Note: \ 
Matters are independent of the particular choice of an orthonormal basis since the transition matrix between any two such is orthogonal, 
hence the absolute value of its determinant is 1.]
\\[-.25cm]
\end{x}

\begin{x}{\small\bf SUBLEMMA} \ 
Let \mV be a finite dimensional real Hilbert space; let 
$V_1 \subset V$, $V_2 \subset V$ be subspaces.  
Suppose that 
$T:V_1 \ra V_2$ is a bijective linear transformation $-$then $\forall \ \phi \in C_c(V_2)$, 
\[
\int_{V_2} \ \phi \ \td V_2 
\ = \ 
\abs{\det \tT} \ \int_{V_1} \hsx \phi \circ \tT \ \td V_1,
\]
where the determinant is computed relative to an orthonormal basis in $V_1$ and an orthonormal basis in $V_2$.
\\[-.25cm]
\end{x}

\begin{x}{\small\bf \un{N.B.}} \ 
Symbolically, 
\[
\td V_2 
\ = \ 
\abs{\det \tT} \td V_1.
\]
\\[-1.25cm]
\end{x}

\begin{x}{\small\bf CONVENTION} \ 
Extend the Killing form on $\fg_{ss} \times \fg_{ss}$ to a nondegenerate symmetric bilinear form 
$\tB:\fg \times \fg \ra \R$ with the following properties:
\\[-8pt]

\qquad \textbullet \quad B is Ad \mG invariant.
\\[-.5cm]

\qquad \textbullet \quad B is $\theta$-invariant.
\\[-.5cm]

\qquad \textbullet \quad $\fk$ and $\fp$ are orthogonal under B.
\\[-.5cm]

\qquad \textbullet \quad B is positive definite on $\fp$ and negative definite on $\fk$.
\\[-.25cm]
\end{x}

\begin{x}{\small\bf \un{N.B.}} \ 
The bilinear form 
\[
(X,Y)_\theta 
\ = \ 
-\tB(X,\theta Y) \qquad (X, \ Y \in \fg)
\]
equips $\fg$ with the structure of a real Hilbert space.
\\[-.25cm]

Relative to this data, any subspace $\fl$ of $\fg$ carries a canonical Haar measure $\td \fl$, 
an instance being the Lie algebra $\fl$ of a closed Lie subgroup \mL of \mG.
\\[-.25cm]
\end{x}

\begin{x}{\small\bf EXAMPLE} \ 
$\fk$ and $\fp$ are orthogonal and $\td \fg = \td \fk \td \fp$.
\\[-.5cm]

[Note: \ 
The orthogonal projections $\tE_\fk$, $\tE_\fp$, of $\fg$ onto $\fk$, $\fp$ are given by
\[
\begin{cases}
\ \ds\tE_\fk = \frac{1 + \theta}{2} \\[8pt]
\ \ds\tE_\fp = \frac{1 - \theta}{2}
\end{cases}
\]
respectively.]
\\[-.25cm]
\end{x}

\begin{x}{\small\bf CONSTRUCTION} \ 
Choose an open neighborhood $N_0$ of 0 in $\fl$ and an open 
neighborhood $N_e$ of $e$ in \mL such that exp is an analytic diffeomorphism of $N_0$ 
onto $N_e$.  
Normalize the left Haar measure $\mu_L$ on \mL in such a way that $\forall \ f \in C_c(N_e)$, 
\[
\int_{N_e} \ f \ \td \mu_L 
\ = \ 
\int_{N_0} \ F \ \td \fl,
\]
where 
\[
F(X) 
\ = \ 
f(\exp X) \abs{\hsy\det\bigg(\frac{1 - e^{-\ad(x)}}{\ad(x)}\bigg)}.
\]
This fixes $\mu_L$ uniquely, call it $\td L$, and its definition is independent of the choice of $N_0$.
\\[-.25cm]
\end{x}

\begin{x}{\small\bf \un{N.B.}} \ 
If $L$ is compact, put
\[
\vol(L) 
\ = \ 
\int_L \ \td L
\]
and term $\ds\frac{1}{\vol(L)} \hsx \td L$ the \un{normalized Haar measure} of $L$.
\end{x}
\vspx

Now write after Iwasawa \ $G = KAN$, thus $\forall \ f \in C_c(G)$,
\[
\int_G \ f \ \td \mu_G 
\ = \ 
\int_{K \times A \times N} \ f(k \hsy a \hsy n) \hsx e^{2 \rho (\log a)} 
\ \td \mu_K(k) \ \td \mu_A(a) \ \td \mu_N(n).
\]
On the right hand side, take
\[
\td \mu_K(k) 
\ = \ 
\frac{1}{\vol(K)} \hsx \td K, \ 
\td \mu_A(a)  \ = \ \td A, \ 
\td \mu_N(n)  \ = \ \td N.
\]
Then these choices determine $\td \mu_G$ uniquely, denote it by the symbol $\td_{st} G$ and refer to it as the 
\un{standard Haar measure} of \mG.
\\[-.25cm]

\begin{x}{\small\bf LEMMA} \ 
\[
\td G 
\ = \ 
2^{-\frac{1}{2} \dim N} \ e^{2 \rho (\log a)} \ \td K \hsx \td A \hsx \td N.
\]

PROOF \ 
It suffices to show that 
\[
\td \fg 
\ = \ 
2^{-\frac{1}{2} \dim N} \hsx \td \fk \hsx \td \fa \hsx \td \fn.
\]
To establish this, write
\[
\fp 
\ = \ 
\fa + \tE_\fp \fn,
\]
the sum being orthogonal, hence
\begin{align*}
\td \fg \ 
&= \ 
\td \fk \hsx \td \fp 
\\[6pt]
&= \ 
\td \fk \hsx \td \fa \hsx \td \tE_\fp \fn
\\[6pt]
&= \ 
\abs{\det \tE_\fp} \hsx \abs{\fn} \hsx \td \fk \hsx \td \fa \hsx \td \fn.
\end{align*}
Choose an orthonormal basis $\tZ_i$ for $\fn$ $-$then 
\[
(\tE_\fp \tZ_i, \tE_\fp \tZ_j)_\theta 
\ = \ 
\delta_{i \hsx j} / 2
\]
which implies that $\sqrt{2} \  \tE_\fp \tZ_i$ is an orthonormal basis for $\tE_\fp \fn$, so
\[
\abs{\det \tE_\fp \fn} 
\ = \ 
\abs{
\begin{matrix}
\hsx \frac{1}{\sqrt{2}}\\[-7pt]
&\cdot\\[-7pt]
&&\cdot\\[-7pt]
&&&\cdot\\[-7pt]
&&&&{\frac{1}{\sqrt{2}} \hsx }
\end{matrix}
}
\ = \ 
2^{-\frac{1}{2} \dim N
.
}
\]

[Note: \ 
\[
\dim N 
\ = \ 
\dim G/K - \rank \hsy G / K.]
\]

Therefore
\begin{align*}
\td_{st} G \ 
&=\ 
e^{2 \rho (\log a)} \bigg( \frac{\td K}{\vol(K)} \bigg) \td A \hsx \td N
\\[11pt]
&=\ 
\frac{1}{\vol(K)} 
\ \frac{2^{-\frac{1}{2} \dim N}}{2^{-\frac{1}{2} \dim N}} 
\ e^{2 \rho (\log a)} \hsx \td K \hsx \td A \hsx \td N
\\[11pt]
&=\ 
\frac{1}{\vol(K)} \ 2^{\frac{1}{2} \dim N} \ \td G.
\end{align*}
\end{x}


\chapter{
$\boldsymbol{\S}$\textbf{1}.\quad  TRANSVERSALS}
\setlength\parindent{2em}
\setcounter{theoremn}{0}
\renewcommand{\thepage}{C II \S1-\arabic{page}}

\qquad Let \mG be a locally compact group.
\\[-.25cm]

\begin{x}{\small\bf SUBLEMMA} \ 
Fix $x \in G$ $-$then for any open neighborhood \mU of $e$ there exists an open neighborhood \mV of $x$ such that $V^{-1} V \subset U$.
\\[-.25cm]
\end{x}

\begin{x}{\small\bf DEFINITION} \ 
A subgroup $\Gamma \subset G$ is a \un{discrete subgroup} if the relative topology on $\Gamma$ is the discrete topology.
\\[-.25cm]
\end{x}

\begin{x}{\small\bf LEMMA} \ 
A subgroup $\Gamma \subset G$ is discrete iff there exists an open neighborhood \mU of $e$ (in \mG) such that $\Gamma \cap U = \{e\}$.
\\[-.25cm]
\end{x}

\begin{x}{\small\bf THEOREM} \ 
Suppose that $\Gamma \subset G$ is a discrete subgroup $-$then $\Gamma$ is closed in \mG, hence $G/\Gamma$ is a 
locally compact Hausdorff space (cf. I, $\S1$, $\#14$).
\\[-.25cm]
\end{x}

\begin{x}{\small\bf EXAMPLE} \ 
\\[-.5cm]

\qquad \textbullet \quad $G = \R$, $\Gamma = \Z$.
\\[-.2cm]

\qquad \textbullet \quad $G = \A$, $\Gamma = \Q$. 
\\[-.2cm]

\qquad \textbullet \quad $G = \I$, $\Gamma = \Q^\times$.
\\[-.25cm]
\end{x}

\begin{x}{\small\bf LEMMA} \ 
Let $\Gamma$ be a discrete subgroup  of \mG $-$then there exists an open neighborhood $U_0$ of $e$ such that 
$U_0 \hsy \gamma \cap U_0 \neq \emptyset$ for all $\gamma \neq e$ in $\Gamma$.
\\[-.5cm]

PROOF \ 
First choose \mU per $\#3$.  
This done, choose \mV per $\#1$ (with $x = e$) and put $U_0 = V$.  
Assume now that $u_0^\prime \in U_0 \hsy \gamma \hsx \cap \hsx U_0$, 
thus $u_0^\prime = u_0 \gamma$ ($\exists \ u_0 \in U_0$), so
\[
\gamma 
\ = \ 
u_0^{-1} u_0^\prime \in U_0^{-1} U_0 = V^{-1} V \hsx \subset \hsx U
\]
\qquad $\implies$ 
\[
\gamma = e.
\]
\\[-1.25cm]
\end{x}

\begin{x}{\small\bf SUBLEMMA} \ 
Let \mH be a closed subgroup of \mG and give $G/H$ the quotient topology $-$then the projection 
$\pi:G \ra G/H$ is an open map.
\\[-.25cm]

[Let $U \subset G$ be a nonempty open set, the claim being that $\pi(U) \subset G/H$ is a nonempty open set.  
But $\pi(U)$ is open iff $\pi^{-1}(\pi(U))$ is open.  
And
\[
\pi^{-1}(\pi(U)) 
\ = \ 
U H 
\ = \ 
\bigcup\limits_{h \in H} \hsx Uh
\]
which is a union of open sets.]
\end{x}
\vspace{0.3cm}

\begin{x}{\small\bf THEOREM} \ 
Suppose that $\Gamma \subset G$ is a discrete subgroup $-$then the projection $\pi:G \ra G/\Gamma$ is a local homeomorphism.
\\[-.5cm]

PROOF \ 
Fix $x \in G$ and choose $U_0$ per $\#6$ to get an open neighborhood $x \hsy U_0$ of $x$ with the property that $\forall \ \gamma \neq e$ in $\Gamma$, 
\[
x \hsy U_0 \hsy \gamma \cap x \hsy U_0 
\ = \ 
x(U_0 \hsy \gamma \cap U_0) 
\ = \ 
\emptyset.
\]
Therefore the arrow $x \hsy U_0 \ra \pi(x \hsy U_0)$ is a continuous bijection, hence is a homeomorphism (cf. $\#7$).
\\[-.25cm]
\end{x}

\begin{x}{\small\bf DEFINITION} \ 
Let $\Gamma$ be a discrete subgroup of \mG $-$then a Borel subset $\gT \subset G$ is a \un{transversal} for $G/\Gamma$ if the restriction of 
$\pi$ to $\gT$ is bijective.
\\[-.25cm]
\end{x}

\begin{x}{\small\bf \un{N.B.}} \ 
In other words, a transversal $\gT$ for $G/\Gamma$ is a Borel subset of \mG which meets each coset exactly once.
\\[-.25cm]
\end{x}

\begin{x}{\small\bf THEOREM} \ 
Suppose that $\Gamma \subset G$ is a discrete subgroup.  
Assume: \mG is second countable $-$then $G/\Gamma$ admits a transversal $\gT$.
\\[-.25cm]
\end{x}

\begin{x}{\small\bf REMARK} \ 
A transversal $\gT$ for $G/\Gamma$ gives rise to a unique section 
$\tau:G/\Gamma \ra \gT \subset G$ $(\pi \circ \tau = \id)$ which is Borel measurable if \mG is second countable.
\\[-.25cm]
\end{x}

\begin{x}{\small\bf \un{N.B.}} \ 
Tacitly, Lie groups are assumed to be second countable (cf. I, $\S5$), hence $\sigma$-compact (cf. I, $\S2$, $\#2$).
\\[-.5cm]

[Note: \ Still, in this situation it is not claimed (nor is it true in general) that smooth sections exist.]
\\[-.25cm]
\end{x}

\begin{x}{\small\bf EXAMPLE} \ 
Take $G = \R$, $\Gamma = \Z$ $-$then $[0,1[$ is a transversal for $\R/\Z$.
\\[-.25cm]
\end{x}

\begin{x}{\small\bf EXAMPLE} \ 
Take $G =\A$, $\Gamma = \Q$ $-$then $\prod\limits_p \hsx \Z_p \times [0,1[$ is a transversal for $\A/\Q$.
\end{x}
\vspace{0.3cm}

\begin{x}{\small\bf EXAMPLE} \ 
Take $G = \I$, $\Gamma = \Q^\times$ 
$-$then $\prod\limits_p \hsx \Z_p^\times \times \R_{>0}^\times$ is a transversal for $\I/\Q^\times$.
\\[-.25cm]
\end{x}

\begin{x}{\small\bf CONVENTION} \ 
The Haar measure on a discrete group $\Gamma$ is the counting measure:
\[
\int_\Gamma \ f(\gamma) \ \td_\Gamma(\gamma) 
\ = \ 
\sum\limits_{\gamma \in \Gamma} \ f(\gamma).
\]

[Note: \ 
$\Gamma$ is unimodular (being discrete).]
\\[-.25cm]
\end{x}

\begin{x}{\small\bf LEMMA} \ 
If $\Gamma \subset G$ is a discrete subgroup and if \mG is second countable, then $\Gamma$ is at most countable.
\\[-.25cm]
\end{x}

\begin{x}{\small\bf LEMMA} \ 
If $\Gamma \subset G$ is a discrete subgroup, if \mG is second countable and if $\gT$ is a transversal for $G/\Gamma$, then 
\allowdisplaybreaks
\begin{align*}
G \ 
&=\ 
\bigcup\limits_{\gamma \in \Gamma} \ \gT \gamma \qquad \text{(disjoint union)},
\\[11pt]
\int_G \ f \ \td\mu_G
&=\ 
\sum\limits_{\gamma \in \Gamma} \ \int_{\gT \gamma} \  f \ \td\mu_G 
\\[11pt]
&=\ 
\int_{\gT} \  f^\Gamma \circ \pi  \ \td\mu_G.
\end{align*}


[Note: \ 
$\forall x \in \gT$, 
\allowdisplaybreaks
\begin{align*}
(f^\Gamma \circ \pi) (x) \ 
&=\ 
f^\Gamma (x\Gamma)
\\[11pt]
&=\ 
\int_\Gamma \ f(x \gamma) \ \td \mu_\Gamma(\gamma)
\\[11pt]
&=\ 
\sum\limits_{\gamma \in \Gamma} \ f(x \gamma).]
\end{align*}
\\[-.75cm]
\end{x}

\begin{x}{\small\bf RAPPEL} \ 
If \mG is unimodular and if $\mu_G$ is fixed, then $G/\Gamma$ admits an invariant measure $\mu_{G/\Gamma}$ 
characterized by the condition that for all $f \in C_c(G)$, 
\[
\int_G \hsx f(x) \hsx \td\mu_G(x) 
\ = \ 
\int_{G/\Gamma} \hsx \ f^\Gamma (\dot{x}) \hsx \td \mu_{G/\Gamma} (\dot{x}) \qquad (\dot{x} = x\Gamma) 
\qquad \text{(cf. I, $\S4$, $\#7$).}
\]
\\[-1.25cm]
\end{x}

\begin{x}{\small\bf THEOREM} \ 
If $\Gamma \subset G$ is a discrete subgroup, if \mG is second countable, 
if $\gT$ is a transversal for $G/\Gamma$, 
if \mG is unimodular and if $\mu_G$ is fixed, then $\forall \ f \in C_c(G)$, 
\[
\int_{G/\Gamma} \hsx \ f^\Gamma \hsx \td \mu_{G/\Gamma}
\ = \ 
\int_{\gT} \hsx \ f^\Gamma \circ \pi \td\mu_G.
\]

[Simply assemble the foregoing data.]
\\[-.2cm]

[Note: \ 
Since the $f^\Gamma$ $(f \in C_c(G))$ exhaust $C_c(G/\Gamma)$ (cf. I, $\S4$, $\#3$), 
it follows that $\forall \ \phi \in C_c(G/\Gamma)$, 
\[
\int_{G/\Gamma} \ \ \phi \ \td \mu_{G/\Gamma}
\ = \ 
\int_{\gT} \ \ \phi\circ \pi \ \td\mu_G.
\]
In particular, this holds for all $\phi$ if $G/\Gamma$ is compact.]
\\[-.25cm]
\end{x}

\begin{x}{\small\bf DEFINITION} \ 
Let $\Gamma$ be a discrete subgroup of \mG $-$then a Borel subset $\fF \subset G$ is a 
\un{fundamental domain} for $G/\Gamma$ if it differs from a transversal by a set of Haar measure 0 (left or right).
\\[-.25cm]
\end{x}


\begin{x}{\small\bf EXAMPLE} \ 
Take $G = \R$, $\Gamma = \Z$ $-$then $[0,1]$ is a fundamental domain for $\R/\Z$.
\\[-.25cm]
\end{x}

\begin{x}{\small\bf \un{N.B.}} \ 
What was said in $\#21$ goes through verbatim if ``transversal'' is replaced by ``fundamental domain''.
\end{x}


\chapter{
$\boldsymbol{\S}$\textbf{2}.\quad  LATTICES}
\setlength\parindent{2em}
\setcounter{theoremn}{0}
\renewcommand{\thepage}{C II \S2-\arabic{page}}


\qquad Let \mG be a second countable locally compact group, $\Gamma \subset G$ a discrete subgroup.
\\[-.25cm]

\begin{x}{\small\bf NOTATION} \ 
Given a finite subset $\Delta \subset \Gamma$, let $G_\Delta$ denote the centralizer of $\Delta$ in \mG.
\\[-.5cm]
\end{x}

\begin{x}{\small\bf \un{N.B.}} \ 
$G_\Delta$ is closed in \mG.
\\[-.25cm]
\end{x}

\begin{x}{\small\bf LEMMA} \ 
$G_\Delta \Gamma$ is closed in \mG.
\\[-.5cm]

PROOF \ 
Let $x_n \in G_\Delta$ and $\gamma_n \in \Gamma$ be sequences such that $x_n \hsy \gamma_n$ converges to a limit $x$ 
$-$then 
the claim is that $x \in G_\Delta \Gamma$.  
To begin with, $\forall \ \gamma \in \Delta$, 
\begin{align*}
x^{-1} \gamma x \ 
&=\ 
\lim\limits_{n \ra \infty} \big(\gamma_n^{-1} \hsy x_n^{-1} \hsy \gamma \hsy x_n \hsy \gamma_n\big)
\\[8pt]
&=\ 
\lim\limits_{n \ra \infty} \gamma_n^{-1} \hsy \gamma \hsy \gamma_n.
\end{align*}
Since $\Gamma$ is discrete, $\exists \ n_0(\gamma)$:
\begin{align*}
n \geq n_0(\gamma) 
&\implies 
\gamma_n^{-1} \hsy \gamma \hsy \gamma_n = \gamma_{n+1}^{-1} \hsy \gamma \hsy \gamma_{n+1}
\\[8pt]
&\implies
\gamma_{n+1} \gamma_n^{-1} \in G_\Delta.
\end{align*}
But $\Delta$ is finite, thus $\exists \ n_0$ independent of the choice of $\gamma$ such that 
\[
n \geq n_0 
\implies
\gamma_n 
\ = \ 
y_n \gamma_{n_0} \qquad (y_n \in G_\Delta)
\]
\hspace{2cm} $\implies$
\[
x_n \hsy \gamma_n 
\ = \ 
x_n \hsy y_n \hsy \gamma_{n_o} 
\ = \ 
z_n \hsy \gamma_{n_0} \qquad (z_n \in G_\Delta)
\]
\hspace{2cm} $\implies$
\[
z_n 
\ = \ 
x_n \hsy \gamma_n \hsy \gamma_{n_0}^{-1} \lra x \hsy \gamma_{n_0}^{-1} \qquad (n \ra \infty)
\]
\hspace{2cm} $\implies$
\[
x \hsy \gamma_{n_0}^{-1} \in G_\Delta 
\implies 
x \in G_\Delta \Gamma.
\]
\\[-1.25cm]
\end{x}

\begin{x}{\small\bf NOTATION} \ 
Given $\gamma \in \Gamma$, $G_\gamma$ is its centralizer in \mG and $\Gamma_\gamma$ $(= G_\gamma \hsx \cap \hsx \Gamma)$ 
is its centralizer in $\Gamma$.
\\[-.25cm]
\end{x}

\begin{x}{\small\bf \un{N.B.}} \ 
$G_\gamma$ is a closed subgroup of \mG, as is $\Gamma_\gamma$ (cf. $\S1$, $\#4$).
\\[-.25cm]
\end{x}

\begin{x}{\small\bf LEMMA} \ 
$G_\gamma \hsy \Gamma$ is closed in \mG (cf. $\#3$ (take $\Delta = \{\gamma\}$)).
\\[-.25cm]
\end{x}

\begin{x}{\small\bf SUBLEMMA} \ 
If \mH is a closed subgroup of \mG, if $\pi:G \ra G/H$ is the projection and if \mF is a closed subset of \mG that is the union of cosets $x H$, 
then $\pi(F)$ is closed in $G/H$.
\\[-.25cm]
\end{x}

\begin{x}{\small\bf APPLICATION} \ 
The image of
\[
G_\gamma \hsy \Gamma 
\ = \ 
\bigcup\limits_{x \in G_\gamma} \ x \hsy \Gamma
\]
in $G / \Gamma$ is closed, hence is a locally compact Hausdorff space.
\\[-.25cm]
\end{x}

\begin{x}{\small\bf REMARK} \ 
The projection $\pi:G \ra G /  \Gamma$ is an open map (cf. $\S 1$, $\# 7$) but, in general, it is not a closed map.
\\[-.5cm]

[Take $G = \R$, $\Gamma = \Z$ and view $\R / \Z$ as $[0,1[$ equipped with the topology in which an open basis 
consists of all sets $]a,b \hsy[$ $(0 < a < b < 1)$ and of all sets
$[0,a \hsy[ \hsx \cup \hsx ]b, 1[$ $(0 < a < b < 1)$ $-$then
\[
A 
\ = \ 
\bigg\{ \frac{3}{2}, \hsx \frac{9}{4}, \ldots, n + 2^{-n}, \ldots\bigg\}
\]
is closed in $\R$ but
\[
\pi(A) 
\ = \ 
\bigg\{ \frac{1}{2}, \hsx \frac{1}{4}, \ldots, \frac{1}{2^n}, \ldots\bigg\}
\]
is not closed in $[0,1[$\hsx.]
\\[-.25cm]
\end{x}

Considered as families of subsets of \mG, $G_\gamma \hsy \Gamma / \Gamma$ and $\pi(G_\gamma)$ are identical:  
The elements of $G_\gamma \hsy \Gamma / \Gamma$ are the cosets $x\Gamma$ with $x \in G_\gamma \hsy \Gamma$ and the elements of 
$\pi(G_\gamma)$ are the cosets of $x \hsy \Gamma$ with $x \in G_\gamma$.

\begin{x}{\small\bf LEMMA} \ 
The identity map
\[
\{x \hsy \Gamma : x \in G_\gamma \hsy \Gamma\} \ra \{x \hsy \Gamma : x \in G_\gamma\}
\]
is a homeomorphism.
\\[-.5cm]

[Note: \ 
That is to say, the two topologies are the same.]
\\[-.25cm]
\end{x}

\begin{x}{\small\bf \un{N.B.}} \ 
One may then identify $\pi(G_\gamma)$ with $G_\gamma \hsy \Gamma /\Gamma$ which is therefore closed in $G/\Gamma$ (cf. $\# 8$).
\end{x}
\vspace{0.3cm}

\begin{x}{\small\bf NOTATION} \ 
Let 
\[
r: G_\gamma \hsy \Gamma /\Gamma \ra G_\gamma / G_\gamma \hsx \cap \hsx \Gamma
\]
be the arrow defined by 
\[
r (x \hsy \Gamma) 
\ = \ 
x(G_\gamma \hsx \cap \hsx \Gamma).
\]
\\[-1.25cm]
\end{x}

\begin{x}{\small\bf \un{N.B.}} \ 
$r$ is bijective. 
\\[-.25cm]
\end{x}

\begin{x}{\small\bf THEOREM} \ 
$r$ is a homeomorphism.
\\[-.25cm]

This is not completely obvious and it will be best to break the proof into two parts.
\\[-.25cm]
\end{x}

\begin{x}{\small\bf LEMMA} \ 
$r$ carries open subsets of $G_\gamma \hsy \Gamma /\Gamma$ onto open subsets of $G_\gamma / G_\gamma \hsx \cap \hsx \Gamma$. 
\\[-.5cm]

PROOF \ 
An open subset of $G_\gamma \hsy \Gamma /\Gamma$ is a subset $\{x\Gamma : x \in X\}$, where $X \subset G_\gamma$, such that 
$X \Gamma$ is open in $G_\gamma \hsy \Gamma$ viewed as a subspace of \mG.   
Since
\[
X(G_\gamma \hsx \cap \hsx \Gamma) 
\ = \ 
X \Gamma \hsx \cap \hsx G_\gamma,
\]
it follows that $X(G_\gamma \hsx \cap \hsx \Gamma)$ is an open subset of $G_\gamma$ in its relative topology as a subspace of \mG, 
thus by the very definition of the topology on $G_\gamma / G_\gamma \hsx \cap \hsx \Gamma$, 
\[
r\{x \hsy \Gamma : x \in X\} 
\ = \ 
\{x (G_\gamma \hsx \cap \hsx \Gamma) : x \in X\}
\]
is an open subset of $G_\gamma / G_\gamma \hsx \cap \hsx \Gamma$.  
\\[-.25cm]
\end{x}

\begin{x}{\small\bf LEMMA} \ 
$r^{-1}$ carries open subsets of $G_\gamma / G_\gamma \hsx \cap \hsx \Gamma$ onto open subsets of $G_\gamma \hsy \Gamma / \Gamma$.
\\[-.5cm]

PROOF \ 
Let $\{y (G_\gamma \hsx \cap \hsx \Gamma) : y \in Y\}$ $(Y \subset G_\gamma)$ be an open subset of $G_\gamma /  G_\gamma \hsx \cap \hsx \Gamma$ 
$-$then 
$Y(G_\gamma \hsx \cap \hsx \Gamma)$ is an open subset of $G_\gamma$, so
\[
\pi(Y(G_\gamma \hsx \cap \hsx \Gamma))
\ = \ 
\{y \Gamma: y \in Y\}
\]
is open in $G_\gamma \hsy \Gamma / \Gamma$ (see the Appendix infra) or still, 
\[
\{y \Gamma : y \in Y\} 
\ = \ 
r^{-1}\{y(G_\gamma \hsx \cap \hsx \Gamma) : y \in Y\}
\]
is open in $G_\gamma \hsy \Gamma / \Gamma$.
\\[-.25cm]
\end{x}

\begin{x}{\small\bf EXAMPLE} \ 
Take $G = \R$, $\Gamma = \Z$ and $H = \sqrt{2} \hsx \Z$ $-$then the argument used in
$\# 15$ is applicable if the $G_\gamma$ there is replaced by \mH, thus the map
\[
H / H \cap \Gamma \ra H + \Gamma / \Gamma
\]
is continuous.  
Nevertheless, it is not a homeomorphism.
\\[-.5cm]

[$H \cap \Gamma$ is trivial so $H / H \cap \Gamma$ is isomorphic to $\Z$ and carries the discrete topology.  
Meanwhile, $H + \Gamma = \sqrt{2} \hsx \Z + \Z$ is dense in $\R$, hence
\[
H + \Gamma / \Gamma 
\ = \ 
\sqrt{2} \hsx \Z + \Z /  \Z
\]
is dense in $\R/\Z \approx \T$.  
It is isomorphic to $\Z$ as a group but it is not discrete since every nonempty open subset of $\T$ intersects it in an infinite set implying thereby 
that none of its finite subsets are open.]
\\[-.5cm]

[Note: \ 
The difference here is this: $G_\gamma \hsy \Gamma / \Gamma$ is locally compact but $H + \Gamma / \Gamma$ is not locally compact.]
\\[-.25cm]
\end{x}

\begin{x}{\small\bf DEFINITION} \ 
$\Gamma$ is said to be a \un{lattice} if $G / \Gamma$ admits a finite invariant measure (cf. I, $\S 4, \ \# 4$), 
$\Gamma$ being termed \un{uniform} or \un{nonuniform} according to whether $G / \Gamma$ is compact or not.
\\[-.25cm]
\end{x}

\begin{x}{\small\bf \un{N.B.}} \ 
If there is a lattice in \mG, then \mG is necessarily unimodular (cf. I, $\S 3, \ \# 11$ and I, $\S 4, \ \# 10$).
\\[-.5cm]

[Note: \ 
A discrete cocompact subgroup is necessarily a uniform lattice \ldots \hsy .]
\\[-.25cm]
\end{x}

\begin{x}{\small\bf EXAMPLE} \ 
$\Z$ is a uniform lattice in $\R$.
\\[-.25cm]
\end{x}

\begin{x}{\small\bf EXAMPLE} \ 
$\SL(2,\Z)$ is a nonuniform lattice in $\SL(2,\R)$ .
\\[-.25cm]
\end{x}

\begin{x}{\small\bf THEOREM} \ 
Suppose that $\Gamma \subset G$ is a uniform lattice $-$then $\forall \ \gamma \in \Gamma$, $G_\gamma / \Gamma_\gamma$ is compact.
\\[-.5cm]

PROOF \ 
$G_\gamma \hsy \Gamma / \Gamma$ is closed in $G / \Gamma$, hence is compact (this being the case of $G/\Gamma$).  
On the other hand,
\[
r: G_\gamma \hsy \Gamma / \Gamma \ra G_\gamma / G_\gamma \hsx \cap \hsx \Gamma \qquad (= G_\gamma / \Gamma_\gamma)
\]
is a homeomorphism (cf. $\#14$).
\\[-.5cm]

[Note: \ 
Consequently, $\Gamma_\gamma \subset G_\gamma$ is a uniform lattice and $G_\gamma$ is unimodular.]
\\[-.25cm]
\end{x}

\begin{x}{\small\bf NOTATION} \ 
$[\Gamma]$ is a set of representatives for the $\Gamma$-conjugacy classes in $\Gamma$.
\\[-.25cm]

Put
\[
\ggS 
\ = \ 
\coprod\limits_{\gamma \in [\Gamma]} \ G / \Gamma_\gamma \times \{\gamma\}
\]
and define $\psi: \ggS \ra G$ by the rule
\[
\psi(x \hsy \Gamma_\gamma, \gamma) 
\ = \ 
x \hsy \gamma \hsy x^{-1}.
\]
\\[-1.25cm]
\end{x}

\begin{x}{\small\bf \un{N.B.}} \ 
$\Gamma_\gamma$ is a discrete subgroup of \mG, thus $\Gamma_\gamma$ is closed in \mG 
(cf. $\S 1$, $\# 4$) and therefore the quotient $G / \Gamma_\gamma$ is a locally compact Hausdorff space from which it follows that 
$\ggS$ is a locally compact Hausdorff space. 
\\[-.25cm]
\end{x}

\begin{x}{\small\bf DEFINITION} \ 
Let \mX and \mY be locally compact Hausdorff spaces, $f:X \ra Y$ a continuous function $-$then $f$ is \un{proper} if for every 
compact subset \mK of \mY, the inverse image $f^{-1}(K)$ is a compact subset of \mX.
\\[-.25cm]
\end{x}

\begin{x}{\small\bf THEOREM} \ 
Suppose that $\Gamma \subset G$ is a uniform lattice $-$then $\psi$ is a proper map.  
\\[-.25cm]
\end{x}

\begin{x}{\small\bf NOTATION} \ 
Given $\gamma \in \Gamma$, let
\[
[\gamma]_G 
\ = \ 
\{x \hsy \gamma \hsy x^{-1} : x \in G\}.
\]
\\[-1.25cm]
\end{x}

\begin{x}{\small\bf APPLICATION} \ 
In the uniform situation, for any compact subset $K \subset G$, 
\[
\{\gamma \in [\Gamma] : [\gamma]_G \cap K \neq \emptyset\}
\]
is finite.
\\[-.25cm]
\end{x}

\begin{x}{\small\bf LEMMA} \ 
A proper map $f:X \ra Y$ is closed:
\[
S \subset X \ \text{closed} \ 
\implies 
f(S) \subset Y \ \text{closed}.
\]
\\[-1.25cm]
\end{x}

\begin{x}{\small\bf APPLICATION} \ 
In the uniform situation, $\forall \ \gamma \in \Gamma$, $[\gamma]_G$ is closed.

[In fact, 
\[
[\gamma]_G 
\ = \ 
\psi\bigg( \bigcup\limits_{\gamma_0 \hsy \in \hsy [\gamma]_G \hsy \cap \hsy [\Gamma]} \ G / \gamma_0 \times \{\gamma_0\}\bigg).]
\]
\\[-1.25cm]
\end{x}

\begin{x}{\small\bf \un{N.B.}} \ 
Accordingly, $[\gamma]_G$ is a locally compact Hausdorff space and the canonical arrow
\[
G / G_\gamma \ra [\gamma]_G
\]
is a homeomorphism.
\end{x}

\[
\text{APPENDIX}
\]

Denote by $\restr{\pi}{G_\gamma}$ the restriction of $\pi:G \ra G/\Gamma$ to $G_\gamma$.
\\

\qquad {\small\bf CRITERION} \ 
Suppose that there exist nonempty open sets
\[
U 
\hsx \subset \hsx G_\gamma, 
\quad 
V \hsx \subset \hsx
G_\gamma \hsy \Gamma / \Gamma
\]
such that the restriction of $\restr{\pi}{G_\gamma}$ to \mU is an open continuous map of \mU onto \mV $-$then $\restr{\pi}{G_\gamma}$ is open.
\\[-.5cm]

PROOF \ 
Given $x \in G_\gamma$ and an open neighborhood \mW of $x$ in $G_\gamma$, it suffices to show that $(\restr{\pi}{G_\gamma})(W)$ 
contains an open neighborhood $N_x$ of $(\restr{\pi}{G_\gamma})(x)$.  
So fix a point $y \in U$ and put
\[
\widetilde{U} 
\ = \ 
U \cap y x^{-1} W,
\]
an open neighborhood of $y$ in \mU, thus the image $(\restr{\pi}{G_\gamma})(\widetilde{U})$ is an open subset of 
$G_\gamma \hsy \Gamma /\Gamma$ or still, $\widetilde{U} \Gamma$ is an open subset of $G_\gamma \hsy \Gamma$, hence
\[
x \hsy y^{-1} \widetilde{U} \Gamma 
\ = \ 
(x \hsy y^{-1} U \hsy \cap \hsy W ) \Gamma
\]
is an open subset of $G_\gamma \hsy \Gamma$ and 
\[
x \in x \hsy y^{-1} U \hsy \cap \hsy W.
\]
Put now
\[
N_x 
\ = \ 
\{z \Gamma : z \in x \hsy y^{-1} U \hsy \cap \hsy W\}.
\]
Then $N_x$ is an open subset of $G_\gamma \hsy \Gamma /\Gamma$ contained in 
\[
(\restr{\pi}{G_\gamma})(W) 
\ = \ 
\{w \Gamma : w \in W\}
\]
to which $(\restr{\pi}{G_\gamma})(x)$ belongs.
\\[-.25cm]

There is a commutative diagram
\[
\begin{tikzcd}[sep=huge]
{G_\gamma } 
\ar{rr}{\restr{\pi}{G_\gamma}}
\arrow[d,shift right=0.5,dash] \arrow[d,shift right=-0.5,dash]
&&{G_\gamma \hsy \Gamma / \Gamma} \\
{G_\gamma } \ar{rr}[swap]{\pi_\gamma}
&&{G_\gamma / G_\gamma \hsx \cap \hsx \Gamma} \ar{u}[swap]{r^{-1}}
\end{tikzcd}
\]

and $G_\gamma \hsy \Gamma /\Gamma$ is a locally compact Hausdorff space (cf. $\# 8$), thus is a Baire space.
\\[-.25cm]

\qquad {\small\bf LEMMA} \ 
$\restr{\pi}{G_\gamma}$ is an open map.
\\[-.5cm]

PROOF \ 
The quotient $G_\gamma / G_\gamma \hsx \cap \hsx \Gamma$ is second countable, hence $\sigma$-compact, hence
\[
G_\gamma / G_\gamma \hsx \cap \hsx \Gamma
\ = \ 
\bigcup\limits_{n = 1}^\infty \ K_n,
\]
where $K_1, K_2, \ldots$ are compact.  
In view of $\# 15$,
\[
r^{-1} : G_\gamma / G_\gamma \hsx \cap \hsx \Gamma  \ra G_\gamma \hsy \Gamma /\Gamma
\]
is continuous and one-to-one, so $\forall \ n$ the restriction of $r^{-1}$ to $K_n$ is a homeomorphism of $K_n$ onto 
$L_n \equiv r^{-1} (K_n)$:
\[
G_\gamma \hsy \Gamma / \Gamma
\ = \ 
\bigcup\limits_{n = 1}^\infty \ L_n,
\]
a countable union of compacta.  
Being Baire, it therefore follows that $\exists \ n \in \N$ and a nonempty open subset \mV of $G_\gamma \hsy \Gamma /\Gamma$ such that 
$V \subset r^{-1}(L_n)$.  
Put
\[
U 
\ = \ 
(\restr{\pi}{G_\gamma})^{-1} (V).
\]
Then $U \subset G_\gamma$ is nonempty and open and the restriction of $\pi_\gamma$ to \mU is an open continuous map of \mU 
onto $r(V)$ or still, the restriction of $\restr{\pi}{G_\gamma}$ to \mU is an open continuous map of \mU onto \mV.

\chapter{
$\boldsymbol{\S}$\textbf{3}.\quad  UNIFORMLY INTEGRABLE FUNCTIONS}
\setlength\parindent{2em}
\setcounter{theoremn}{0}
\renewcommand{\thepage}{C II  \S3-\arabic{page}}

\qquad Let \mG be a unimodular locally compact group and, generically, let $\fU$ be a compact symmetric neighborhood of the identity in \mG.
\\

\begin{x}{\small\bf NOTATION} \ 
Given a continuous function $f$ on \mG, put
\[
f_\fU(y) 
\ = \ 
\sup\limits_{x, z \in \fU} \hsx \abs{f(x \hsy y \hsy z)} \qquad (y \in G).
\]
\\[-1.25cm]
\end{x}

\begin{x}{\small\bf LEMMA} \ 
$f_\fU \in C(G)$, i.e., is a continuous function on \mG.
\\[-.25cm]
\end{x}

\begin{x}{\small\bf DEFINITION} \ 
A continuous function $f$ on \mG is said to be \un{uniformly integrable} if there exists a $\fU$ such that $f_\fU \in \hsx L^1(G)$.
\\[-.25cm]
\end{x}

\begin{x}{\small\bf \un{N.B.}}  \ 
Since $\abs{f} \leq f_\fU$, it is clear that if $f$ is uniformly integrable, then $f$ is integrable: $f \in L^1(G)$.
\\[-.25cm]
\end{x}

\begin{x}{\small\bf NOTATION} \ 
Write $C_\UN(G)$ for the set of continuous functions on \mG that are uniformly integrable.
\\[-.25cm]
\end{x}

\begin{x}{\small\bf LEMMA} \ 
\[
C_c(G) \ \subset \ C_\UN(G) \ \subset \ C_0(G).
\]

[Note: \ 
As usual, $C_c(G)$ is the set of continuous functions on \mG that are compactly supported and $C_0(G)$ is the set of continuous functions 
on \mG that vanish at infinity.]
\\[-.25cm]
\end{x}

\begin{x}{\small\bf LEMMA} \ 
\[
C_\UN(G) \ \subset \ L^2(G).
\]

[Integrable functions in $C_0(G)$ are square integrable.]
\\[-.25cm]
\end{x}

\begin{x}{\small\bf EXAMPLE} \ 
Take $G = \R$ $-$then $\ds f(x) = e^{-x^2}$ is uniformly integrable.
\\[-.25cm]
\end{x}

\begin{x}{\small\bf LEMMA} \ 
If $f$, $g, \in C_\UN(G)$, then $f * g \in C_\UN(G)$.
\\[-.5cm]

[Working with a common $\fU$, 
\begin{align*}
(f * g)_\fU(y) \ 
&=\ 
\sup\limits_{x, z \in \fU} \ \abs{\int_G \ f(u)  \hsy g(u^{-1}  \hsy x  \hsy y  \hsy z) \ \td \mu_G(u)}
\\[11pt]
&=\ 
\sup\limits_{x, z \in \fU} \ \abs{\int_G \ f(x \hsy u) g(u^{-1}  \hsy y  \hsy z) \ \td \mu_G(u)}
\\[11pt]
&\leq\ 
\sup\limits_{x, z \in \fU} \ \int_G \ \abs{f(x \hsy u) g(u^{-1}  \hsy y  \hsy z)} \ \td \mu_G(u)
\\[11pt]
&\leq\ 
\int_G \ f_\fU(u) g_\fU(u^{-1} y)\ \td \mu_G(u)
\\[11pt]
&=\ 
(f_\fU * g_\fU) (y),
\end{align*}
which suffices.]
\\[-.5cm]

[Note: \ The convolution $f * g$ is continuous.]
\\[-.25cm]
\end{x}

Let $H \subset G$ be a closed subgroup and assume that \mH is unimodular and cocompact.
\\[-.2cm]

\begin{x}{\small\bf NOTATION} \ 
$L^2(G/H)$ is the Hilbert space associated with $\mu_{G/H}$ (the invariant measure on $G/H$ per I, $\S4$, $\#5$).
\\[-.25cm]
\end{x}

\begin{x}{\small\bf NOTATION} \ 
$L(G/H)$ is the left translation representation of \mG on $L^2(G/H)$.
\\[-.25cm]
\end{x}

\begin{x}{\small\bf THEOREM} \ 
Let $f \in C_\UN(G)$ $-$then
\[
L_{G/H} (f) 
\ = \ 
\int_G \ f(x) \hsx L_{G/H} (x) \ \td\mu_G(x)
\]
is an integral operator on $L^2(G/H)$ with continuous kernel
\[
K_f(x, y) 
\ = \ 
\int_H \ f(x \hsy h \hsy y^{-1}) \ \td\mu_H(h).
\]

Since
\[
C(G/H \times G/H) \ \subset \ L^2(G/H \times G/H),
\]
it follows that $\forall \ f \in C_\UN(G)$, $L_{G/H}(f)$ is Hilbert-Schmidt, hence is compact.
\\[-.25cm]
\end{x}

\begin{x}{\small\bf SUBLEMMA} \ 
Let \mU be a unitary representation of \mG on a Hilbert space $\sH$ with the property that $\forall \ f \in C_c(G)$, the operator
\[
U(f) 
\ = \ 
\int_G \ f(x) \hsx U(x) \ \td\mu_G(x)
\]
is compact $-$then \mU is discretely decomposable, a given irreducible unitary representation of \mG occurring at most a finite number of times in 
the orthogonal decomposition of \mU.
\\[-.5cm]

[Note: \ If \mG is a Lie group, then one can replace $C_c(G)$ by $C_c^\infty(G)$.]
\\[-.25cm]
\end{x}

\begin{x}{\small\bf \un{N.B.}} \ 
If \mG is second countable, then $\sH$ is separable.
\\[-.25cm]
\end{x}

\begin{x}{\small\bf APPLICATION} \ 
Take $\sH = L^2(G/H)$, $U = L_{G/H}$ $-$then there exist nonnegative integers $m(\Pi,L_{G/H})$ $(\Pi \in \widehat{G})$ such that 
\[
L_{G/H} 
\ = \ 
\widehat{\bigoplus\limits_{\Pi \in \widehat{G}}} \ m(\Pi,L_{G/H}) \hsx \Pi.
\]

[Note: \ Per usual, $\widehat{G}$ is the set of unitary equivalence classes of irreducible unitary representations of $G$.]
\end{x}


\chapter{
$\boldsymbol{\S}$\textbf{4}.\quad  THE SELBERG TRACE FORMULA}
\setlength\parindent{2em}
\setcounter{theoremn}{0}
\renewcommand{\thepage}{C II \S4-\arabic{page}}

\qquad Let \mG be a second countable locally compact group, $\Gamma \subset G$ a discrete subgroup.  
Assume: $\Gamma$ is a uniform lattice $-$then $G/\Gamma$ is cocompact and \mG is necessarily unimodular 
(cf. \S2, \#19).

Working with $L^2(G/\Gamma)$, there is an orthogonal decomposition
\[
L_{G/\Gamma} \ = \ \widehat{\bigoplus\limits_{\Pi \in \widehat{G}}} \  m(\Pi,L_{G/\Gamma}) \Pi
\qquad (cf. \ \S3, \ \#15),
\]
the multiplicities $m(\Pi,L_{G/\Gamma})$ being certain nonnegative integers.

\begin{x}{\small\bf RAPPEL} \ 
$\forall \ f \in C_\UN(G)$, $L_{G/\Gamma} (f)$ is an integral operator on $L_{G/\Gamma}^2$ with continuous kernel
\[
K_f(x,y) 
\ = \ 
\sum\limits_{\gamma \in \Gamma} \ f(x \hsy \gamma \hsy y^{-1}) \qquad \text{(cf. $\S3$, $\#12$).}
\]

[Note: \ 
This implies that $L_{G/\Gamma} (f)$ is Hilbert-Schmidt.]
\\[-.25cm]
\end{x}

\begin{x}{\small\bf CONVENTION} \ 
Fix a Haar measure $\mu_G$ on \mG, take the counting measure on $\Gamma$, and normalize the invariant measure $\mu_{G/\Gamma}$ on $G/\Gamma$ by the stipulation 
\[
\int_G 
\ = \ 
\int_{G/\Gamma} \ \int_\Gamma \qquad \bigg(= \int_{G/\Gamma} \ \sum\limits_\Gamma\bigg).
\]
\\[-.5cm]

If $f = g*g^*$ $(g \in C_\UN(G))$, then $f \in C_\UN(G)$ (cf. $\S3$, $\#9$),
\[
L_{G/\Gamma} (f) 
\ = \ 
L_{G/\Gamma} (g) L_{G/\Gamma} (g)^*
\]
is trace class and (cf. B, II, $\S2$, $\#8$)
\[
\tr(L_{G/\Gamma} (f)) 
\ = \ 
\int_{G/\Gamma} \ K_f(\dot{x}, \dot{x}) \ \td \mu_{G/\Gamma}(\dot{x}) \qquad (\dot{x} = x \Gamma).
\]
\\[-1.25cm]
\end{x}

\begin{x}{\small\bf REMARK} \ 
The assumption that $f = g*g^*$ $(g \in C_\UN(G))$ is not restrictive.  
For if $f = g*h^*$ $(g, \ h \in C_\UN(G))$, put
\[
T(g, h) 
\ = \ 
g* h^*
\]
and using the same letter for the diagonal, note that
\[
T(g, h) 
\ = \ 
\frac{1}{4} \big(T(g + h) - T(g - h) - \sqrt{-1} \hsx T(g - \sqrt{-1} \hsx h) + \sqrt{-1} \hsx T(g + \sqrt{-1} \hsx h) \big).
\]

Let $\chisubGmodGamma$ be the characteristic function of $G/\Gamma$, i.e., the function $\equiv 1$. 
Choose $\alpha \in C_c(G) : \alpha^\Gamma = \chisubGmodGamma$ (cf. I, $\S4$, $\#3$), thus $\forall \ x \in G$, 
\[
\alpha^\Gamma (x \Gamma) 
\ = \ 
\sum\limits_{\gamma \in \Gamma} \ \alpha(x \hsy \gamma) 
\ = \ 
1.
\]

One can then write
\allowdisplaybreaks
\begin{align*}
\sum\limits_{\Pi \in \widehat{G}} \ m\big(\Pi, L_{G/\Gamma}\big) \tr(\Pi(f)) \ 
&=\ 
\tr( L_{G/\Gamma}(f))
\\[11pt]
&=\ 
\int_{G/\Gamma} \ K_f(\dot{x}, \dot{x}) \ \td\mu_{G/\Gamma} \dot{x}
\\[11pt]
&=\ 
\int_\gT \ K_f(x \hsy \Gamma, x \hsy \Gamma) \ \td\mu_G(x) \qquad \text{(cf. $\S1$, $\#21$)}
\\[11pt]
&=\ 
\int_\gT \ \bigg( \sum\limits_{\gamma \in \Gamma} \alpha(x \hsy \gamma) \bigg) \hsy K_f(x \hsy \Gamma, x \hsy \Gamma) 
\ \td\mu_G(x)
\\[11pt]
&=\ 
\int_\gT \  \sum\limits_{\gamma \in \Gamma} \ \alpha(x \hsy \gamma) K_f(x \hsy \gamma, x \hsy \gamma) \ \td\mu_G(x)
\\[11pt]
&=\ 
\int_\gT \ (\alpha K_f)^\Gamma \circ \pi(x) \ \td\mu_G(x)
\\[11pt]
&=\ 
\int_G \ \alpha(x) K_f(x, x) \ \td\mu_G(x)
\\[11pt]
&=\ 
\int_G \ \alpha(x) \ \sum\limits_{\gamma \in \Gamma} \  f(x \hsy \gamma \hsy x^{-1})  \ \td\mu_G(x)
\\[11pt]
&=\ 
\sum\limits_{\gamma \in \Gamma} \ \int_G \ \alpha(x)  f(x \hsy \gamma \hsy x^{-1})  \ \td\mu_G(x).
\end{align*}
\\[-.75cm]
\end{x}

\begin{x}{\small\bf NOTATION} \ 
For any $\gamma \in \Gamma$,
\[
\begin{cases}
\ G_\gamma = \ \text{centralizer of $\gamma$ in \mG}\\[3pt]
\ \Gamma_\gamma = \ \text{centralizer of $\gamma$ in $\Gamma$}
\end{cases}
.
\]
\\[-1.25cm]
\end{x}

\begin{x}{\small\bf RAPPEL} \ 
$\Gamma_\gamma$ is a uniform lattice in $G_\gamma$ (cf. $\S2$, $\#22$).
\\[-.5cm]

[Note: \ 
Consequently, $G_\gamma$ is unimodular.]
\end{x}
\vspace{0.3cm}

\begin{x}{\small\bf NOTATION} \ 
For any $\gamma \in \Gamma$,
\[
\begin{cases}
\ [\gamma]_\Gamma = \ \text{conjugacy class of $\gamma$ in $\Gamma$}\\[3pt]
\ [\gamma]_G = \ \text{conjugacy class of $\gamma$ in \mG}
\end{cases}
.
\]
\\[-1.25cm]
\end{x}

\begin{x}{\small\bf RAPPEL} \ 
There are canonical bijections
\[
\begin{cases}
\ \Gamma/\Gamma_\gamma \ra [\gamma]_\Gamma \\[3pt]
\ G/G_\gamma \ra [\gamma]_G.
\end{cases}
\]
\\[-1.25cm]
\end{x}

Returning to the computation, break the sum over $\Gamma$ into conjugacy classes in $\Gamma$, the contribution from 
\[
[\gamma]_\Gamma \ = \ \{\delta \hsy \gamma \hsy \delta^{-1} : \delta \in \Gamma / \Gamma_\gamma\}
\]
being
\begin{align*}
\sum\limits_{\delta \in \Gamma / \Gamma_\gamma} \
\int_G \ \alpha(x) f(x \delta \hsy \gamma \hsy \delta^{-1} x^{-1}) \ d\mu_G(x) \ 
&=\ \sum\limits_{\delta \in \Gamma / \Gamma_\gamma} \ 
\int_G \ \alpha(x \delta^{-1}) f(x \hsy \delta \hsy  x^{-1}) \ d\mu_G(x)
\\[8pt]
&=\ \int_G \ 
\bigl 
( \sum\limits_{\delta \in \Gamma / \Gamma_\gamma} 
\alpha(x \delta^{-1}) \bigr ) f(x \hsy \gamma x^{-1}) \ d\mu_G(x).
\end{align*}

\begin{x}{\small\bf CONVENTION} \ 
Supplementing the agreements in $\#2$, fix a Haar measure $\mu_{G_\gamma}$ on $G_\gamma$, take the counting measure on $\Gamma_\gamma$, and 
normalize the invariant measure $\mu_{G_\gamma/\Gamma_\gamma}$ on $G_\gamma/\Gamma_\gamma$ by the stipulation 
\[
\int_{G_\gamma} 
\ = \ 
\int_{G_\gamma/\Gamma_\gamma} \ \int_{\Gamma_\gamma} 
\qquad \bigg(= \int_{G_\gamma/\Gamma_\gamma} \ \sum\limits_{\Gamma_\gamma} \ \bigg).
\]
Next, fix $\mu_{G/G_\gamma}$ via
\[
\int_G 
\ = \ 
\int_{G/G_\gamma} \ \int_{G_\gamma}.
\]
Finally, make the identification
\[
G/\Gamma_\gamma 
\ \approx \ 
\big(G/G_\gamma\big) / \big(G_\gamma/\Gamma_\gamma \big)
\]
and put
\[
\int_{G/\Gamma_\gamma}
\ = \ 
\int_{G/G_\gamma} \ \int_{G_\gamma /\Gamma_\gamma}.
\]

Moving on
\begin{align*}
\int_G \ \bigg(\sum\limits_{\delta \in \Gamma/\Gamma_\gamma} 
\ \alpha(x \hsy \delta^{-1})\bigg) f(x \hsy \gamma \hsy x^{-1}) 
\ \td \mu_G(x)\ 
&=\ 
\int_{G/G_\gamma} 
\ \int_{G_\gamma} 
\ \int_{\Gamma/\Gamma_\gamma} \cdots
\\[11pt]
&=\ 
\int_{G/G_\gamma} 
\ \int_{G_\gamma/\Gamma_\gamma} 
\ \int_{\Gamma_\gamma} 
\ \int_{\Gamma/\Gamma_\gamma} \cdots \hsx .
\end{align*}
But
\[
 \int_{\Gamma_\gamma} \ \int_{\Gamma/\Gamma_\gamma} \ \alpha(x \hsy \eta \hsy \delta^{-1})
\]
is $\equiv 1$, leaving
\[
\int_{G/G_\gamma} \ \int_{G_\gamma/\Gamma_\gamma} \cdots 
\ = \ 
\int_{G /\Gamma_\gamma} \cdots .
\]

Summary:
\[
\tr(L_{G/\Gamma} (f)) 
\ = \ 
\sum\limits_{\gamma \in [\Gamma]} \ \int_{G /\Gamma_\gamma} \ f(x \hsy \gamma \hsy x^{-1}) 
\ \td \mu_{G/\Gamma_\gamma} (\dot{x}),
\]
the sum being taken over a set of representatives for the $\Gamma$-conjugacy classes in $\Gamma$ (cf. $\S2$, $\#23$).
\\[-.25cm]
\end{x}

\begin{x}{\small\bf \un{N.B.}} \ 
\allowdisplaybreaks
\begin{align*}
\int_{G /\Gamma_\gamma}  
f(x \hsy \gamma x^{-1}) \ \td \mu_{G/\Gamma_\gamma} (\dot{x}) \ 
&=\ 
\int_{G/G_\gamma} 
\ \bigg(\int_{G_\gamma/\Gamma_\gamma} 
f(x \hsy \eta \hsy \gamma \hsy \eta^{-1} \hsy x^{-1}) 
\ \td\mu_{G_\gamma/\Gamma_\gamma} (\dot{\eta}) \bigg)
\ \td \mu_{G/G_\gamma} (\dot{x})
\\[11pt]
&=\ 
\int_{G/G_\gamma} \ \bigg(\int_{G_\gamma/\Gamma_\gamma} 
f(x \hsy \gamma \hsy x^{-1}) 
\ \td\mu_{G_\gamma/\Gamma_\gamma} (\dot{\eta}) \bigg)
\ \td \mu_{G/G_\gamma} (\dot{x})
\\[11pt]
&=\ 
\int_{G/G_\gamma} \ f(x \hsy \gamma \hsy x^{-1}) \hsx \bigg(\int_{G_\gamma/\Gamma_\gamma} \ 
\td\mu_{G_\gamma/\Gamma_\gamma} \bigg) 
\ \td \mu_{G/G_\gamma} (\dot{x})
\\[11pt]
&=\ 
\vol(G_\gamma/\Gamma_\gamma) \int_{G/G_\gamma} \ f(x \hsy \gamma \hsy x^{-1}) \ \td \mu_{G/G_\gamma} (\dot{x}).
\end{align*}
\\[-.75cm]
\end{x}

\begin{x}{\small\bf DEFINITION} \ 
Given $f \in C_\UN(G) * C_\UN(G)$, the \un{Selberg trace formula} is 
the relation 
\[
\sum\limits_{\Pi \in \widehat{G}} \ m\big(\Pi, L_{G/\Gamma}\big) \tr(\Pi(f)) 
\ = \ 
\sum\limits_{\gamma \in [\Gamma]} \ \vol(G_\gamma/\Gamma_\gamma) \
\int_{G/G_\gamma} 
\ f(x \hsy \gamma \hsy x^{-1}) \ \td \mu_{G/G_\gamma} (\dot{x}),
\]
their common value being
\[
\tr(L_{G/\Gamma}(f)).
\]
\\[-1.25cm]
\end{x}

\begin{x}{\small\bf REMARK} \ 
Suppose that \mG is a Lie group $-$then
\[
C_c^\infty(G) * C_c^\infty(G) 
\ = \ 
C_c^\infty(G) 
\qquad \text{(Dixmier-Malliavin)}, 
\]
Since
\[
C_c^\infty(G) \subset C_\UN(G), 
\]
it follows that the Selberg trace formula is valid for all $f \in C_c^\infty(G)$.
\\[-.25cm]
\end{x}

Let \mG be a second countable locally compact group, $\Gamma \subset G$ a uniform lattice.
\\[-.25cm]

\begin{x}{\small\bf LEMMA} \ 
Let 
$\chi:G \ra \bT$ 
be a unitary character $-$then the multiplicity of 
$\chi$ in $L^2(G/\Gamma)$ is 1 if $\chi(\Gamma) = \{1\}$ and 0 otherwise.
\\[-.25cm]
\end{x}

Now \ take \ \mG \ abelian \ and identify $\widehat{G}$ with the unitary character group of $G:\Pi \longleftrightarrow \chi$, 
the Fourier transform being defined by 
\[
\tr(\Pi(f)) 
\ = \ 
\widehat{f}(x) 
\ = \ 
\int_G \ f(x) \chi(x) \ \td \mu G(x).
\]
\\[-1.25cm]

\begin{x}{\small\bf NOTATION} \ 
\[
\Gamma^\perp 
\ = \ 
\{\chi \in \widehat{G}: \chi(\gamma) = 1 \ \forall \ \gamma \in \Gamma\}.
\]
\\[-1.25cm]
\end{x}

\begin{x}{\small\bf \un{N.B.}} \ 
Therefore
\[
\begin{cases} 
\ \chi \in \Gamma^\perp \implies m(\chi,L_{G/\Gamma}) = 1\\[3pt]
\ \chi \notin \Gamma^\perp \implies m(\chi,L_{G/\Gamma}) = 0
\end{cases}
.
\]

The Selberg trace formula thus simplifies:
\\[-.2cm]

\qquad \textbullet \quad Matters on the ``spectral side'' reduce to 
\[
\sum\limits_{\gamma \in \Gamma^\perp} \ \widehat{f}(\chi).
\]

\qquad \textbullet \quad Matters on the ``geometric side'' reduce to 
\[
\vol(G/\Gamma) \ \sum\limits_{\gamma \in \Gamma} \ f(\gamma).
\]
\\[-.125cm]
\end{x}

\begin{x}{\small\bf DEFINITION} \ 
The relation
\[
\sum\limits_{\chi \in \Gamma^\perp} \widehat{f}(\chi) 
\ = \ 
\vol(G / \Gamma) \ \sum\limits_{\gamma \in \Gamma} f(\gamma)
\]
is the 
\un{Poisson summation formula}
\index{Poisson summation formula} 
(cf. A, III, \S4, \#7) 
(in that situation 
\[
\vol(G / \Gamma) \ = \ \frac{\abs{G}}{\abs{\Gamma}} \bigg).
\]
\end{x}


\chapter{
$\boldsymbol{\S}$\textbf{5}.\quad  FUNCTIONS OF REGULAR GROWTH}
\setlength\parindent{2em}
\setcounter{theoremn}{0}
\renewcommand{\thepage}{C II \S5-\arabic{page}}

\qquad Let \mG be a second countable locally compact group, $\Gamma \subset G$ a uniform lattice.  
While $C_{\UN}(G)$ is theoretically convenient, there is a larger class of functions that can be fed into the Selberg trace formula.
\\[-.25cm]

\begin{x}{\small\bf DEFINITION} \ 
Let $\phi \in C(G) \hsx \cap \hsx L^1(G)$ be nonnegative $-$then $\phi$ is said to be of \un{regular growth} if there is a compact symmetric 
neighborhood $\fU$ of the identity in \mG and a positive constant \mC (depending on $\phi$ and $\fU$) such that $\forall \ y \in G$,
\[
\phi(y) 
\ \leq \
C \hsx \int_U \hsx \phi(x \hsy y) \hsx \td \mu_G(x).
\]
\\[-1.25cm]
\end{x}

\begin{x}{\small\bf \un{N.B.}} \ 
In terms of the characteristic function $\chisubfU$ of $\fU$, $\forall \ y \in G$,
\begin{align*}
(\chi_U * \phi)(y) \ 
&=\
\int_G \ \chisubfU(x) \phi(x^{-1} y) \ \td \mu_G(x)
\\[8pt]
&=\
\int_\fU  \ \phi(x^{-1} y) \ \td \mu_G(x)
\\[8pt]
&=\
\int_{\fU^{-1}}  \ \phi(x \hsy y) \ \td \mu_G(x)
\\[8pt]
&=\
\int_\fU  \ \phi(x \hsy y) \ \td \mu_G(x).
\end{align*}
\\[-.75cm]
\end{x}

\begin{x}{\small\bf EXAMPLE} \ 
Take $G = \R^n$ and fix a real number $r > 0$ such that
\[
\int_{\R^n} \ \frac{1}{(1 + \norm{Y})^r} \hsx \td Y 
< 
\infty.
\]
Given $\fU$, fix a real number $N > 0$ such that $\forall \ X \in \fU$, 
\[
\big(1 + \norm{Y}\big)^{-r} 
\leq 
N\big(1+ \norm{X + Y}\big)^{-r}.
\]
Then
\begin{align*}
\frac{1}{\vol(\fU)} \ \int_\fU \hsx \frct{\td X}{\big(1+ \norm{X + Y}\big)^r} \ 
&\geq \ 
\frac{1}{\vol(\fU)} \ \int_\fU \ \frct{\td X}{\big(1+ \norm{Y}\big)^r} 
\\[8pt]
&=\ 
\big(1 + \norm{Y}\big)^{-r}.
\end{align*}
Therefore
\[
\phi (Y)
\ = \ 
\big(1 + \norm{Y}\big)^{-r} 
\]
is of regular growth.
\\[-.25cm]
\end{x}

\begin{x}{\small\bf EXAMPLE} \ 
Let \mG be a connected semisimple Lie group with finite center and fix a real number $r > 0$ such that 
\[
\int_G \hsx \Gsym (y) (1 + \sigma(y))^{-r} \ \td_G(y) \ < \ \infty.
\]
Given $\fU$, fix a real number $M > 0$ such that $\forall \ x \in \fU$, 
\[
\Gsym (y) 
\ \leq \ 
M \Gsym (x \hsy y)
\]
and fix a real number $N > 0$ such that $\forall \ x \in \fU$,  
\[
(1 + \sigma(y))^{-r} 
\ \leq \ 
N (1 + \sigma(x \hsy y))^{-r}.  
\]
Then
\begin{align*}
\frac{MN}{\vol(\fU)} \ \int_\fU \ \Gsym (x \hsy y) (1 &+ \sigma(x \hsy y))^{-r} \ \td \mu_G(x) \ 
\\[11pt]
&\geq \ 
\frac{1}{\vol(\fU)} \ \int_\fU \hsx \Gsym (y) (1 + \sigma(y))^{-r} \ \td \mu_G(x) \ 
\\[11pt]
&=\ 
\Gsym (y) (1 + \sigma(y))^{-r}.
\end{align*}
Therefore
\[
\phi(y) 
\ = \ 
\Gsym (y) (1 + \sigma(y))^{-r} 
\]
is of regular growth.
\\[-.25cm]
\end{x}

\begin{x}{\small\bf DEFINITION} \ 
Let $f$ be a continuous function on \mG $-$then $f$ is \un{admissible} if there exists a function $\phi$ of 
regular growth such that $\forall \ y \in G$,  
\[
\abs{f(y)}
\ \leq \ 
\phi(y) \quad (\leq C ( \chisubfU * \phi)(y)).
\]

[Note: \ 
Admissible functions are integrable.]
\\[-.25cm]
\end{x}

\begin{x}{\small\bf EXAMPLE} \ 
The rapidly decreasing functions on $\R^n$ are admissible (cf. $\#3$).
\\[-.25cm]
\end{x}

\begin{x}{\small\bf LEMMA} \ 
If $f \in C_\text{UN}(G)$, then $f$ is admissible.
\\[-.5cm]

PROOF \ 
$\forall \ y \in G$, $\abs{f(y)} \leq f_\fU(y)$.  
And
\begin{align*}
f_\fU(y) \ 
&=\ 
\sup\limits_{u, z \in \fU} \hsx \abs{f(u \hsy y \hsy  z)}
\\[11pt]
&\leq \ 
\sup\limits_{u, z \in \fU} \hsx \abs{f(u \hsy x \hsy y  \hsy z)} \qquad (x \in \fU)
\end{align*}
\hspace{2cm} $\implies$
\begin{align*}
f_\fU(y) \ 
&=\ 
\frac{\vol(\fU)}{\vol(\fU)} \hsx f_\fU(y)
\\[11pt]
&\leq \ 
\frac{1}{\vol(\fU)} \ \int_\fU \hsx f_\fU(y) \  \td \mu_G(x) 
\\[11pt]
&\leq \ 
\frac{1}{\vol(\fU)} \ \int_\fU \ \sup\limits_{x, z \in \fU} \hsx \abs{f(u \hsy x \hsy y \hsy z)} \ \td \mu_G(x)
\\[11pt]
&= \ 
\frac{1}{\vol(\fU)} \ \int_\fU \ f_\fU (x \hsy y) \ \td \mu_G(x).
\end{align*}
Therefore $f_\fU$ is of regular growth, hence $f$ is admissible.
\\[-.25cm]
\end{x}

\begin{x}{\small\bf LEMMA} \ 
Suppose that $\abs{f} \leq \abs{g}$, where $g$ is admissible, say $\abs{g} \leq \psi$ $-$then 
$f$ is admissible (clear) as is $f * f$.

[For
\begin{align*}
\abs{f * f} \ 
&\leq \ 
\abs{f} * \abs{f}
\\[11pt]
&\leq \ 
\abs{g} * \abs{g} 
\\[11pt]
&\leq \ 
\psi * \abs{g}.
\end{align*}
And $\psi * \abs{g}$ is of regular growth
\begin{align*}
\psi * \abs{g} \ 
&\leq \ 
(C \chisubfU * \psi) * \abs{g}
\\[11pt]
&= \ 
C (\chisubfU * (\psi * \abs{g})).
\end{align*}
The condition of admissible is then met by 
\[
\phi \ = \ \psi  * \abs{g}.]
\]

[Note: \ 
If $f_1$, $f_2 \in C(G) \hsx \cap \hsx L^1(G)$ and if $f_1$ is admissible, then $f_1 * f_2$ is admissible.  
Proof: 
\begin{align*}
\abs{f_1 * f_2} \ 
&\leq \ 
\abs{f_1} * \abs{f_2} 
\\[11pt]
&\leq \ 
\phi_1 * \abs{f_2}
\\[11pt]
&\leq \ 
C(\chisubfU * (\phi_1 * \abs{f_2})).]
\end{align*}
\\[-.75cm]
\end{x}

\begin{x}{\small\bf DEFINITION} \ 
A series of functions $f_1, f_2, \ldots$ on a locally compact Hausdorff space \mX is 
\un{locally dominantly absolutely convergent} (ldac) if for every compact set $K \subset X$ there exists a positive constant 
$M_K$ such that $\forall \ k \in K$,
\[
\sum\limits_n \hsx \abs{f_n(k)} \ < \ M_K.
\]
\\[-1.25cm]
\end{x}

\begin{x}{\small\bf CRITERION} \ 
Let $f \in C(G) \hsx \cap \hsx L^1(G)$.  
Assume: \ The operator $L_{G/\Gamma}(f)$ is trace class and the series
\[
\sum\limits_{\gamma \in \Gamma} \ f(x \gamma y^{-1})
\]
is ldac on $G \times G$ to a separately continuous function $-$then the Selberg trace formula obtains:
\[
\tr\big(L_{G/\Gamma}(f)\big) 
\ = \ 
\sum\limits_{\gamma \in [\Gamma]} \ \vol(G_\gamma/\Gamma_\gamma) \ 
\int_{G/G_\gamma} \hsx f(x \hsy \gamma \hsy  x^{-1}) \ \td \mu_{G/G_\gamma} (\dot{x}),
\]
the sum on the right hand side being absolutely convergent.
\\[-.5cm]

[First of all, 
\[
\tr\big(L_{G/\Gamma}(f)\big) 
\ = \ 
\int_{G/\Gamma} \ \sum\limits_{\gamma \in \Gamma} \ f(x \hsy \gamma \hsy x^{-1}) \ \td \mu_{G/\Gamma} (\dot{x}) 
\qquad \text{(cf. B, II, $\S2$, $\#8$).} 
\]
Proceeding, fix a compact set $K \subset G$: $K \Gamma = G$ (cf. $\#11$ infra) and choose $M_K > 0$:
\[
k, \ \ell \in K 
\implies 
\sum\limits_{\gamma \in \Gamma} \ \abs{f(k \hsy \gamma \hsy \ell^{-1})} 
\ < \ M_K.
\]
Here, of course, the ldac condition is per $K \times K^{-1} \subset G \times G$.  
Given $x$, $y \in G$, $\exists \ \gamma_x, \ \gamma_y \in \Gamma$: $x \hsy \gamma_x, \ y \hsy \gamma_y \in K$, so
\[
\sum\limits_{\gamma \in \Gamma} \ \abs{f(x \hsy  \gamma \hsy  y^{-1})} 
\ = \ 
\sum\limits_{\gamma \in \Gamma} \ \abs{f(x \hsy \gamma_x \hsy \gamma \hsy \gamma_y^{-1} \hsy y^{-1})} 
\ < \ M_K,
\]
from which 
\[
M_K \vol(G/\Gamma) 
\ \geq \ 
\int_{G/\Gamma} \hsx 
\sum\limits_{\gamma \in \Gamma} \ f(x \hsy \gamma \hsy y^{-1}) \ \td \mu_{G/\Gamma} (\dot{x}).
\]
Now interchange sum and integral, the ensuing formal manipulation being justified by Fubini.]
\\[-.25cm]
\end{x}

\begin{x}{\small\bf SUBLEMMA} \ 
There exists a compact set $K \subset G$ such that $K \Gamma = G$.
\\[-.5cm]

[Let \mU be an open neighborhood of $e$ such that $\ov{U}$ is compact $-$then the collection 
$\{\pi(xU): x \in G\}$ is an open covering of $G/\Gamma$, thus there is a finite subcollection
\[
\pi(x_1 U), \pi(x_2 U), \ldots, \pi(x_n U)
\]
that covers $G/\Gamma$ and one may take
\[
K 
\ = \ 
x_1 \ov{U} \hsx \cup \hsx x_2 \ov{U} \hsx \cup \cdots \hsx \cup \hsx x_n \ov{U}.
\]
Indeed, 
\[
G/\Gamma 
\ = \ 
\{k \Gamma : k \in K\},
\]
so given $x \in G$, 
\[
x \hsy \Gamma = k \hsy \Gamma \ (\exists \ k) \ 
\implies \ 
x = k \hsy \gamma  \ (\exists \ \gamma) \ 
\implies \ 
x \in K \Gamma.]
\]

[Note: \ 
It can be shown that \mK contains a transversal $\gT$ which is therefore relatively compact.]
\\[-.25cm]
\end{x}

Suppose that $f$ is admissible $-$then $\forall \ x, \ y \in G$, 
\begin{align*}
\abs{f(x \hsy \gamma \hsy y^{-1})} \ 
&\leq \ 
\phi(x \hsy \gamma \hsy y^{-1})
\\[11pt]
&\leq \ 
C \hsx \int_\fU \ \phi(u \hsy x \hsy \gamma \hsy y^{-1}) \ \td \mu_G(u).
\end{align*}
\\[-.25cm]

\begin{x}{\small\bf REMARK} \ 
Fix $x, \ y \in G$ $-$then $\forall \ \gamma_1, \ \gamma_2 \in \Gamma$, 
\[
\fU \hsy x \hsy \gamma_1 \hsy y^{-1} \ \cap \ \fU \hsy x \hsy \gamma_2 \hsy y^{-1} \ \neq \ \emptyset
\]
iff
\[
\gamma_2 \hsy \gamma_1^{-1} \in x^{-1} \hsy \fU^{-1} \hsy \fU \hsy x.
\]

[In one direction, 
\[
u_1 \hsy x \hsy \gamma_1 \hsy y^{-1} 
\ = \ 
u_2 \hsy x \hsy \gamma_2 \hsy y^{-1}
\]
\hspace{2cm} $\implies$
\[
u_1 \hsy x \hsy \gamma_1 
\ = \ 
u_2 \hsy x \hsy \gamma_2
\]
\hspace{2cm} $\implies$
\[
u_1 \hsy x 
\ = \ 
u_2 \hsy x \hsy \gamma_2 \hsy \gamma_1^{-1}
\]
\hspace{2cm} $\implies$
\[
u_2^{-1} \hsy u_1 \hsy  x 
\ = \ 
x \hsy \gamma_2 \hsy \gamma_1^{-1}
\]
\hspace{2cm} $\implies$
\[
x_1^{-1} \hsy u_2^{-1} \hsy u_1 \hsy x 
\ = \ 
\gamma_2 \hsy \gamma_1^{-1}.]
\]

Since $\Gamma$ is discrete, the compact set $x^{-1} \hsy \fU^{-1} \hsy \fU \hsy x$ contains a finite number \mN of elements of $\Gamma$.  
So, for fixed $x$, $y$, not more than $\mN$ of the 
$\fU \hsy x \hsy \gamma_2 \hsy y^{-1}$ can intersect 
$\fU \hsy x \hsy \gamma_1 \hsy y^{-1}$.
\\[-.25cm]
\end{x}

\begin{x}{\small\bf \un{N.B.}} \ 
Consider the case when $N = 1$.  
Since it is always true that $e \in x^{-1} \fU^{-1} \fU x$, in this situation the $\fU \hsy  x \hsy \gamma \hsy y^{-1}$ are disjoint, hence
\begin{align*}
\sum\limits_{\gamma \in \Gamma} \ \int_\fU \ \phi(u \hsy x \hsy \gamma \hsy y^{-1}) \ \td \mu_G (u) \ 
&\leq \  
\int_G \ \phi  \ \td\mu_G 
\\[11pt]
&< \ 
\infty.
\end{align*}
\\[-.75cm]
\end{x}

\begin{x}{\small\bf RAPPEL} \ 
If $\mu$ is a measure, then 
\begin{align*}
\sum\limits_{i = 1}^n \hsx \mu(X_i) \ 
&=\ 
\mu\bigg(\bigcup\limits_{i = 1}^n \ X_i \bigg) 
+ 
\mu\bigg(
\bigcup\limits_{\substack{i = 1\\ i < j}}^n \ 
\bigcup\limits_{j = 1}^n \ X_i \hsx \cap \hsx X_j \bigg)
\\[11pt]
&\hspace{1cm} 
+ \mu\bigg(
\bigcup\limits_{\substack{i = 1\\ i < j}}^n \ 
\bigcup\limits_{\substack{j = 1\\ j < k}}^n \ 
\bigcup\limits_{k = 1}^n \ X_i \hsx \cap \hsx X_j \hsx \cap \hsx X_k \bigg)
+ \cdots + 
\mu\bigg(\bigcap\limits_{i = 1}^n \ X_i \bigg).
\end{align*}
\\[-.75cm]
\end{x}

\begin{x}{\small\bf LEMMA} \ 
Fix $x$, $y \in G$ $-$then
\[
\sum\limits_{\gamma \in \Gamma} \ \abs{f(x \hsy \gamma \hsy y^{-1})}
\ \leq \ 
N C \int_G \ \phi \ \td\mu_G 
\ < \ \infty.
\]
\\[-1.25cm]
\end{x}

\begin{x}{\small\bf \un{N.B.}} \ 
More is true: The series
\[
\sum\limits_{\gamma \in \Gamma} \ f(x \hsy \gamma \hsy y^{-1})
\]
is ldac on $G \times G$ to a continuous function .
\\[-.5cm]

[The point is that the preceding estimate is uniform in $x$ and $y$ if these variables are confined to compacta $K_x$ and $K_y$.]
\\[-.5cm]

[Note: \ 
Consequently, 
\[
\text{$f$ admissible $\implies$ $L_{G /\Gamma}(f)$ Hilbert-Schmidt.]}
\]
\\[-1.25cm]
\end{x}

\begin{x}{\small\bf THEOREM} \ 
If $f$ is admissible and if $L_{G /\Gamma}(f)$ is trace class, then the Selberg trace formula obtains (cf. $\#10$).
\\[-.25cm]
\end{x}

\begin{x}{\small\bf \un{N.B.}} \ 
\[
\text{$f$ admissible $\implies$ $f * f$ admissible (cf. $\#8$).}
\]
Therefore
\[L_{G /\Gamma}(f * f) 
\ = \ 
L_{G /\Gamma}(f) L_{G /\Gamma}(f)
\]
is trace class and the foregoing is applicable.
\\[-.25cm]
\end{x}

Specialize now to the case when \mG is a connected semisimple Lie group with finite center.

\begin{x}{\small\bf RAPPEL} \ 
$C^1(G)$ is the $L^1$-Schwartz space of \mG.  
It is closed under convolution and contains $C_c^\infty(G)$ as a dense subspace.
\end{x}
\vspace{0.3cm}

Let $f \in C^1(G)$ and take $r > 0$ per $\#4$ $-$then there exists a constant $C > 0$ such that
\[
\abs{f(y)} 
\ \leq \ 
C \hsx 
\abs{-\hspace{-.12cm}\circ \hspace{-.12cm}-}^2
(y) \hsy (1 + \sigma(y))^{-r}
\qquad (y \in G).
\]
Therefore $f$ is admissible.
\\

\begin{x}{\small\bf LEMMA} \ 
$L_{G /\Gamma}(f)$ is trace class.  
\\[-.5cm]

[Using the theory of the parametrix, write
\[
f 
\ = \ 
g * \mu + f * \nu,
\]
where $g \in C^1(G)$ (a certain derivative of $f$), 
$\mu \in C_c^p(G)$, $\nu \in C_c^\infty(G)$, so
\[
L_{G /\Gamma}(f)
\ = \ 
L_{G /\Gamma}(g) \hsy L_{G /\Gamma}(\mu) + L_{G /\Gamma}(f) \hsy L_{G /\Gamma}(\nu).
\]
The functions
\[
f, \ g, \ \mu, \ \nu
\]
are admissible, hence the operators
\[
L_{G /\Gamma}(f), \ 
L_{G /\Gamma}(g), \ 
L_{G /\Gamma}(\mu), \ 
L_{G /\Gamma}(\nu)
\]
are Hilbert-Schmidt.]
\\[-.25cm]
\end{x}

\begin{x}{\small\bf SCHOLIUM} \ 
$\forall \ f \in C^1(G)$, the Selberg trace formula obtains.
\\[-.25cm]
\end{x}

\begin{x}{\small\bf \un{N.B.}} \ 
The assignment
\[
f \ra \tr\big(L_{G/\Gamma}(f)\big)
\]
is continuous in the topology of $C^1(G)$.
\\[-.5cm]

[Note: \ Analogously, the assignment
\[
f \ra \tr\big(L_{G/\Gamma}(f)\big)
\]
is continuous in the topology of $C_c^\infty(G)$, i.e., is a distribution on \mG.]
\\[-.25cm]
\end{x}


\[
\text{APPENDIX}
\]

By way of reconciliation, consider the case when \mG is finite and use the notation of A, III, $\S3$ and $\S4$ $-$then 
given $f \in C(G)$, $\phi \in C(G/\Gamma)$, we have
\[
\big(L_{G/\Gamma}(f) \hsy \phi\big) (x) 
\ = \ 
\sum\limits_{y \in G} \ K_f (x,y) \hsy \phi(y),
\]
where in this context
\[
K_f (x,y) 
\ = \ 
\frac{1}{\abs{\Gamma}} \ \sum\limits_{\gamma \in \Gamma} \ f(x \hsy \gamma \hsy y^{-1}).
\]
Here
\[
\begin{cases}
\ \mu_G = \ \text{counting measure on \mG}\\[3pt]
\ \mu_\Gamma = \ \text{counting measure on $\Gamma$}
\end{cases}
.
\]

Write
\[
G 
\ = \ 
\coprod\limits_{k = 1}^n \ x_k \hsy \Gamma.
\]
Then for any $f \in C(G)$, 
\begin{align*}
\int_G \ f \hsx \td \mu_G \ 
&=\ 
\sum\limits_{x \in G} \ f(x)
\\[11pt]
&=\ 
\int_{G/\Gamma} \ f^\Gamma \ \td \mu_{G/\Gamma}
\\[11pt]
&=\ 
\sum\limits_{k = 1}^n \ \sum\limits_{\gamma \in \Gamma} \ f(x_k \gamma),
\end{align*}
so $\mu_{G/\Gamma}$ is counting measure on $G/\Gamma$.
\\

Now explicate matters:
\allowdisplaybreaks
\begin{align*}
(L_{G/\Gamma}(f)\phi) \hsy (x) \ 
&=\ 
\sum\limits_{y \in G} \ K_f (x,y) \hsy \phi(y)
\\[11pt]
&=\ 
\sum\limits_{k = 1}^n \ \sum\limits_{\gamma \in \Gamma} \  K_f(x,x_k \hsy \gamma) \hsy \phi(x_k \hsy \gamma)
\\[11pt]
&=\ 
\sum\limits_{k = 1}^n \ \sum\limits_{\gamma \in \Gamma} \  K_f(x,x_k) \hsy \phi(x_k)
\\[11pt]
&=\ 
\sum\limits_{k = 1}^n \ \abs{\Gamma} \cdot K_f(x,x_k) \hsy \phi(x_k)
\\[11pt]
&=\ 
\sum\limits_{k = 1}^n \ \abs{\Gamma} \cdot \frac{1}{\abs{\Gamma}} 
\ \sum\limits_{\gamma \in \Gamma} \ f(x_k \hsy \gamma \hsy x_k^{-1}) \hsy  \phi(x_k)
\\[11pt]
&=\ 
\sum\limits_{k = 1}^n \ \sum\limits_{\gamma \in \Gamma} \ f(x_k \hsy \gamma \hsy x_k^{-1}) \hsy  \phi(x_k)
\end{align*}
which establishes that $L_{G/\Gamma}(f)$ is an integral operator on $C(G/\Gamma)$ with kernel
\[
\sum\limits_{\gamma \in \Gamma} \ f(x \hsy \gamma \hsy y^{-1}),
\]
this being the ``$K_f$'' of $\S4$, $\#1$.
\\[-.5cm]

There is more to be said.  
Thus given $f \in C(G)$, we have
\begin{align*}
\tr\big(L_{G/\Gamma}(f)\big) 
&=\ 
\sum\limits_{x \in G} \hsx \frac{1}{\abs{\Gamma}} \ \sum\limits_{\gamma \in \Gamma} \ f(x \hsy \gamma \hsy x^{-1})
\qquad \text{(cf. A, III, $\S3$, $\#8$)}
\\[11pt]
&=\ 
\sum\limits_{i = 1}^n \ \frac{1}{\abs{\Gamma_{\gamma_i}}} \ \sO(f,\gamma_i) 
\hspace{1.65cm} \text{(cf. A, III, $\S4$, $\#2$)}.
\end{align*}
Here
\[
\abs{\Gamma} 
\ = \ 
\{\gamma_1, \ldots, \gamma_n\}
\]
while
\begin{align*}
\sO(f,\gamma_i) \ 
&=\ 
\sum\limits_{x \in G} \hsx f(x \gamma_i x^{-1})
\\[11pt]
&=\ 
\abs{G_{\gamma_i}} \hsx \sum\limits_{x \in G /  G_{\gamma_i}} \hsx f(x \gamma_i x^{-1}).
\end{align*}
Therefore
\begin{align*}
\vol(G_{\gamma_i}/\Gamma_{\gamma_i})
&=\ 
\abs{\frac{G_{\gamma_i}}{\Gamma_{\gamma_i}}}
\\[11pt]
&=\ 
[G_{\gamma_i} : \Gamma_{\gamma_i}].
\end{align*}
\\[-.75cm]

\qquad {\small\bf \un{N.B.}}
The Haar measure $\mu_G$ (or $\mu_{G_\gamma}$) and $\mu_\Gamma$ (or $\mu_{\Gamma_\gamma}$) 
are counting measures, hence the invariant measure $\mu_{G/\Gamma}$ (or $\mu_{G_\gamma / \Gamma_\gamma}$) 
is counting measure, hence the invariant measure $\mu_{G/G_\gamma}$ per
\[
\int_G 
\ = \ 
\int_{G/G_\gamma} \hsx \int_{G_\gamma}
\]
is counting measure, its total volume being
\[
[G : G_\gamma] 
\ = \ 
\frac{\abs{G}}{\abs{G_\gamma}}.
\]
Finally, the invariant measure $\mu_{G/\Gamma_\gamma}$ per
\[
\int_{G/\Gamma_\gamma} 
\ = \ 
\int_{G/G_\gamma} \hsx \int_{G_\gamma / \Gamma_\gamma}
\]
is counting measure and 
\[
\vol(G/\Gamma_\gamma)
\ = \ 
\vol(G/G_\gamma) \hsx \vol(G_\gamma / \Gamma_\gamma),
\]
i.e., 
\[
[G : \Gamma_\gamma]
\ = \ 
[G : G_\gamma] \hsx [G_\gamma : \Gamma_\gamma],
\]
i.e., 
\[
\frac{\abs{G}}{\abs{\Gamma_\gamma}} 
\ = \ 
\frac{\abs{G}}{\abs{G_\gamma}} \hsx [G_\gamma : \Gamma_\gamma]
\]
\hspace{2cm} $\implies$
\[
\frac{\abs{G_\gamma}}{\abs{\Gamma_\gamma}} 
\ = \ 
[G_\gamma : \Gamma_\gamma].
\]

Matters are thus consistent, so the bottom line is that the global trace formula of A, III, $\S4$, $\#6$ is in this context 
the Selberg trace formula.


\chapter{
$\boldsymbol{\S}$\textbf{6}.\quad  DISCRETE SERIES}
\setlength\parindent{2em}
\setcounter{theoremn}{0}
\renewcommand{\thepage}{C II \S6-\arabic{page}}

\qquad Let \mG be a unimodular locally compact group.
\\[-.25cm]

\begin{x}{\small\bf DEFINITION} \ 
Let $\Pi$ be an irreducible unitary representation of \mG on a Hilbert space $V(\Pi)$ $-$then $\Pi$ is 
\un{square integrable} if $\exists \ v \neq 0$ in $V(\Pi)$ such that the coefficient 
\[
x \ra \langle \Pi(x)v, v \rangle
\]
is square integrable on \mG.
\\[-.25cm]
\end{x}

\begin{x}{\small\bf THEOREM} \ 
If $\Pi$ is square integrable, then for all $v_1$, $v_2 \in V(\Pi)$, the coefficient
\[
x \ra \langle \Pi(x)v_1, v_2 \rangle
\]
lies in $L^2(G)$ and there exists a unique positive real number $\td_\Pi$ (depending on the normalization of the Haar measure 
on \mG but independent of $v_1$, $v_2$) such that 
\[
\int_G \hsx \abs{\langle \Pi(x) v_1, v_2 \rangle}^2 \hsx \td \mu_G(x) 
\ = \ 
\frac{1}{\td_\Pi} \hsx \norm{v_1}^2 \norm{v_2}^2.
\]
\\[-1.25cm]
\end{x}

\begin{x}{\small\bf DEFINITION} \ 
$\td_\Pi$ is called the \un{formal dimension} of $\Pi$.
\\[-.5cm]

[Note: \ If \mG is compact, then every irreducible unitary representation of \mG is square integrable and 
$\td_\Pi$ is the dimension of $\Pi$ in the usual sense provided $\ds\int_G \hsx \td \mu_G = 1.$]  
\\[-.25cm]
\end{x}

\begin{x}{\small\bf NOTATION} \ 
$\widehat{G}_d$ is the subset of $\widehat{G}$ comprised of the square integrable representations and is called the 
\un{discrete series} for \mG.
\\[-.5cm]

[Note: \ 
$\widehat{G}_d$ may very well be empty (e.g., take $G = \R$).]
\\[-.25cm]
\end{x}

\begin{x}{\small\bf REMARK} \ 
If $\widehat{G}_d$ is nonempty, then the center of \mG is compact (the converse being false).
\\[-.25cm]
\end{x}

\begin{x}{\small\bf \un{N.B.}} \ 
The elements of $\widehat{G}_d$ are precisely those irreducible unitary representations of \mG which occur as irreducible 
subrepresentations of the left translation representation of \mG on $L^2(G)$.
\\[-.25cm]
\end{x}

\begin{x}{\small\bf NOTATION} \ 
Given a $\Pi \in \widehat{G}_d$, let
\[
\phi_{.,.}(x) 
\ = \ 
\langle \Pi(x) \cdot, \cdot \rangle \qquad (x \in G)
\]
stand for a generic coefficient.
\\[-.25cm]
\end{x}

\begin{x}{\small\bf THEOREM} \ 
Suppose that $\Pi$ is square integrable $-$then $\forall \ v_1$, $v_2$, $\forall \ w_1, \ w_2$ in $V(\Pi)$,
\[
\int_G \hsx \phi_{v_1, v_2} (x) \ov{\phi_{w_1, w_2} (x)} \td\mu_G(x) 
\ = \ 
\frac{1}{\td_\Pi} \langle v_1, w_1\rangle \hsx \ov{\langle v_2, w_2 \rangle}.
\]
\\[-1.25cm]
\end{x}

\begin{x}{\small\bf APPLICATION} \ 
\[
\phi_{v_1, v_2} * \phi_{w_1, w_2} 
\ = \ 
\frac{1}{\td_\Pi} \langle v_1, w_2\rangle \hsx \phi_{w_1, v_2}.
\]

[Note: \ 
If $v_1 = v_2 = w_1 = w_2$ is a unit vector, call it $v$ and abbreviate $\phi_{v,v}$ to $\phi$, then
\[
\norm{\phi}_2^2 
\ = \ 
\frac{1}{\td_\Pi} 
\quad \text{and} \quad 
\phi * \phi 
\ =\ 
\frac{1}{\td_\Pi} \phi.]
\]
\\[-1.25cm]
\end{x}

\begin{x}{\small\bf DEFINITION} \ 
Let $\Pi$ be an irreducible unitary representation of \mG on a Hilbert space $V(\Pi)$ $-$then 
$\Pi$ is \un{integrable} if $\exists  \ v \neq 0$ in $V(\Pi)$ such that the coefficient
\[
x \ra \langle \Pi(x)v,v \rangle
\]
is integrable on \mG.
\\[-.25cm]
\end{x}

\begin{x}{\small\bf \un{N.B.}} \ 
The coefficient
\[
x \ra \langle \Pi(x) v, v\rangle
\]
is bounded and $L^1$, hence is $L^2$.  Therefore
\[
\text{``$\Pi$ integrable'' $\implies$ ``$\Pi$ square integrable''}
\]
but the converse is false.
\\[-.25cm]
\end{x}

\begin{x}{\small\bf THEOREM} \ 
If $\Pi$ is integrable, then there exists a dense subspace $V(\Pi)^\sim$ of $V(\Pi)$ such that for all $v_1$, $v_2 \in V(\Pi)^\sim$ 
the coefficient 
\[
x \ra \langle \Pi(x)v_1,v_2 \rangle
\]
lies in $L^1(G)$.
\\[-.5cm]

[Note: \ 
If $\phi_{v,v} \in L^1(G)$, then one can take
\[
V(\Pi)^\sim 
\ = \ 
\Pi(C_c(G))v.] 
\]
\\[-1.25cm]
\end{x}

Take \mG second countable and assume that $\Pi \in \widehat{G}$ is integrable, say 
$\phi_{v,v} \in L^1(G)$ $-$then $\forall \ f \in C_c(G)$,
\[
\phi_{\Pi(f)v,\Pi(f)v} \in L^1(G) \qquad \text{(cf. $\#12$).} 
\]
Put $v_0 = \Pi(f)v$, normalized by $\norm{v_0} = 1$, and let
\[
\phi_0 
\ = \ 
\td_\Pi \hsy \phi_{v_0,v_0}.
\]
\\[-1.25cm]


\begin{x}{\small\bf \un{N.B.}} \ 
\begin{align*}
\phi_0 *\phi_0 \ 
&=\ 
\td_\Pi \hsx  \phi_{v_0,v_0} * \td_\Pi \hsx \phi_{v_0,v_0} 
\\[11pt]
&=\ 
\td_\Pi^2 \hsx  \phi_{v_0,v_0} * \phi_{v_0,v_0}
\\[11pt]
&=\ 
\td_\Pi^2 \hsx  \frac{1}{\td_\Pi} \hsx \phi_{v_0,v_0}  
\qquad \text{(cf. $\#9$)}
\\[11pt]
&=\ 
\td_\Pi \hsx \phi_{v_0,v_0}
\\[11pt]
&=\ 
\phi_0.
\end{align*}

It is also clear that $\phi_0^* = \phi_0$ and $\Pi(\phi_0)v_0 = v_0$.
\\[-.25cm]
\end{x}

\begin{x}{\small\bf NOTATION} \ 
If $\Pi$ is an irreducible unitary representation of \mG and if $\pi$ is a unitary representation of \mG , then 
\[
\tI_G(\Pi,\pi)
\]
is the set of intertwining operators between $\Pi$ and $\pi$.
\\[-.25cm]
\end{x}

\begin{x}{\small\bf LEMMA} \ 
For any unitary representation $\pi$ of \mG, $\pi(\ov{\phi}_0)$ is the orthogonal projection onto
\[
\{T v_0 : T \in \tI_G(\Pi,\pi)\}.
\]

[Note: \ 
It's $\pi(\ov{\phi}_0)$ , not $\pi(\phi_0) \ldots$ .]
\\[-.25cm]
\end{x}

Suppose that $\Gamma \hsx \subset \hsx G$ is a uniform lattice and take $\pi = L_{G/\Gamma}$.
\\[-.25cm]

\begin{x}{\small\bf APPLICATION} \ 
\[
L_{G/\Gamma} (\ov{\phi}_0)
\]
is trace class and 
\[
\tr (L_{G/\Gamma} \ov{\phi}_0))
\ = \ 
\dim \tI_G(\Pi,L_{G/\Gamma} ) 
\ = \ 
m(\Pi,L_{G/\Gamma}).
\]
\\[-1.25cm]
\end{x}
\begin{x}{\small\bf THEOREM} \ 
The series
\[
\sum\limits_{\gamma \in \Gamma} \hsx \phi_0(x \hsy \gamma \hsy y^{-1})
\]
is ldac on $G \times G$ to a separately continuous function. 
\\[-.5cm]

PROOF \ 
Let $K \subset G$ be compact and let
\[
n(K) 
\ = \ 
\abs{\Gamma \cap K^{-1} \spt(f) \hsy \spt(f)^{-1}K}.
\]
Then the arrow 
\[
\spt(f)^{-1} K \ra G/\Gamma
\]
is at most $n(K)$-to-1 and $\forall \ x \in K$, 

\allowdisplaybreaks
\begin{align*}
\sum\limits_{\gamma \in \Gamma} \ \abs{\phi_0(x \hsy \gamma \hsy y^{-1})} \ 
&=\ 
\td_\Pi \ 
\sum\limits_{\gamma \in \Gamma} \ 
\abs{\langle \Pi(x \hsy \gamma \hsy y^{-1}) \hsy v_0, v_0 \rangle}
\\[11pt]
&=\ 
\td_\Pi \ 
\sum\limits_{\gamma \in \Gamma} \ 
\abs{\langle \Pi(x \hsy \gamma \hsy y^{-1}) \hsy v_0, \Pi(f) \hsy v \rangle}
\\[11pt]
&=\ 
\td_\Pi \ 
\sum\limits_{\gamma \in \Gamma} \ 
\abs{\langle \Pi(f) \hsy v, \Pi(x \hsy \gamma \hsy y^{-1}) \hsy v_0 \rangle}
\\[11pt]
&=\ 
\td_\Pi \ 
\sum\limits_{\gamma \in \Gamma} \ 
\abs{\int_G f(z) \ 
\langle \Pi(z) \hsy v, \Pi(x \hsy \gamma \hsy y^{-1}) \hsy v_0 \rangle 
\ \td\mu_G(z)}
\\[11pt]
&=\ 
\td_\Pi \ 
\sum\limits_{\gamma \in \Gamma} \ 
\abs{\int_G \ f(z) \ 
\langle v, \Pi(z^{-1} \hsy x \hsy \gamma \hsy y^{-1}) \hsy v_0 \rangle 
\ \td\mu_G(z)}
\\[11pt]
&=\ 
\td_\Pi \ 
\sum\limits_{\gamma \in \Gamma} \ 
\abs{\int_G \  f(z^{-1}) \ 
\langle v, \Pi(z \hsy x \hsy \gamma \hsy y^{-1}) \hsy v_0 \rangle 
\ \td\mu_G(z)}
\\[11pt]
&=\ 
\td_\Pi \ 
\sum\limits_{\gamma \in \Gamma} \ 
\abs{\int_G \  f(z^{-1}) \ 
\langle \ov{\Pi(z \hsy x \hsy \gamma \hsy y^{-1}) \hsy v_0,v} \rangle 
\ \td\mu_G(z)}
\\[11pt]
&\leq\ 
\td_\Pi \ 
\sum\limits_{\gamma \in \Gamma} \ 
\int_G \ \abs{f(z^{-1})} \ 
\abs{\langle \ov{\Pi(z \hsy x \hsy \gamma \hsy y^{-1}) \hsy v_0,v} \rangle} 
\ \td\mu_G(z)
\\[11pt]
&=\ 
\td_\Pi \ 
\sum\limits_{\gamma \in \Gamma} \ 
\int_G \ \abs{f(z^{-1})} \ \abs{\langle \Pi(z \hsy x \hsy \gamma \hsy y^{-1}) \hsy v_0,v \rangle} 
\ \td\mu_G(z)
\\[11pt]
&\leq\ 
\td_\Pi \ 
\norm{f}_\infty  \ \sum\limits_{\gamma \in \Gamma} \ 
\int_{\spt(f)^{-1_K}} \hsx \abs{\langle \Pi(z \hsy \gamma \hsy y^{-1}) \hsy v_0,v \rangle} 
\ \td\mu_G(z) 
\\[11pt]
&\leq\ 
\td_\Pi \ 
\norm{f}_\infty  \hsx n(K) \ \int_{G/\Gamma} \ 
\sum\limits_{\gamma \in \Gamma} \hsx \abs{\langle \Pi(\dot{z} \gamma \hsy y^{-1}) \hsy v_0,v \rangle} 
\ \td\mu_{G/\Gamma}(\dot{z})
\\[11pt]
&=\ 
\td_\Pi \norm{f}_\infty  \hsx n(K) \ \int_{G} \ 
\abs{\langle \Pi(z y^{-1}) v_0,v \rangle} 
\ \td\mu_G(z)
\\[11pt]
&=\ 
\td_\Pi \ 
\norm{f}_\infty  \hsx n(K) \ \int_{G} \ 
\abs{\langle \Pi(z) v_0,v \rangle} \ \td\mu_G(z)
\\[11pt]
&=\ 
\td_\Pi  \ 
\norm{f}_\infty n(K) \norm{\phi_{v_0, v}}_1.
\end{align*}
And
\begin{align*}
\norm{\phi_{v_0, v}}_1\ 
&\leq \ 
\int_G \ \int_G \ \abs{f(y) \langle \Pi(x \hsy y) v, v \rangle } 
\ \td\mu_G(x) \ \td\mu_G(y)
\\[11pt]
&\leq \ 
\norm{f}_1 \norm{\phi_{v,v}}_1 
\\[11pt]
&< \ \infty, 
\end{align*}
thereby settling the ldac condition (and then some (no restriction on ``$y$'')), 
leaving the claim of separate continuity which can be left to the reader.
\\[-.25cm]
\end{x}

The operator $L_{G/\Gamma} (\ov{\phi}_0)$ is trace class (cf. $\#16$).  
So, in view of what has been said above, the criterion of $\S5$, $\#10$ is applicable.
\\[-.25cm]

\begin{x}{\small\bf SCHOLIUM} \ 
\[
m(\Pi, L_{G/\Gamma})
\ = \ 
\sum\limits_{\gamma \in [\Gamma]} \ \vol(G_\gamma/\Gamma_\gamma) 
\ \int_{G/G_\gamma} \ \ov{\phi_0(x \hsy \gamma \hsy x^{-1})} 
\ \td \mu_{G/G_\gamma}(\dot{x}),
\]
the sum on the right hand side being absolutely convergent.
\\[-.25cm]
\end{x}

\begin{x}{\small\bf REMARK} \ 
There are circumstances in which the integral
\[
\int_{G/G_\gamma} \ \ov{\phi_0(x \hsy \gamma \hsy x^{-1})} 
\ \td \mu_{G/G_\gamma}(\dot{x})
\]
vanishes for all $\gamma$ except $\gamma = e$, hence then
\allowdisplaybreaks
\begin{align*}
m(\Pi, L_{G/\Gamma}) \ 
&=\ 
\vol(G/\Gamma) \hsx \ov{\phi_0(e)}
\\[11pt]
&=\ 
\vol(G/\Gamma) \hsx \td_\Pi \hsx \langle \Pi(e) v_0, v_0 \rangle
\\[11pt]
&=\ 
\vol(G/\Gamma) \hsx \td_\Pi \hsx \langle v_0, v_0 \rangle
\\[11pt]
&=\ 
\vol(G/\Gamma) \hsx \td_\Pi.
\end{align*}
Therefore $m(\Pi, L_{G/\Gamma})$ is positive, so $\Pi$ definitely occurs in $L_{G/\Gamma}$.
\\[-.5cm]

[Note: \ 
To run a reality check, take \mG finite, $\Gamma = \{e\}$ 
$-$then $\vol(G/\Gamma) = \vol(G) = 1$ and $\forall \ \Pi \in \widehat{G}$, 
\[
m(\Pi, L_{G/\Gamma}) 
\ = \ 
\td_\Pi \qquad \text{(cf. A, II, $\S5$, $\#8$ and A, III, $\S3$, $\#15$).]}
\]
\\[-1.25cm]
\end{x}

\begin{x}{\small\bf \un{N.B.}} \ 
The situation envisioned in $\#19$ is realized if \mG is a connected semisimple Lie group with finite center and if $\Gamma$ has no elements 
of finite order other than the identity.
\\[-.25cm]
\end{x}

\begin{x}{\small\bf LEMMA} \ 
If \mG is a Lie group and if $f \in C_c^\infty(G)$, then the series
\[
\sum\limits_{\gamma \in \Gamma} \ \phi_0(x \hsy \gamma \hsy y^{-1})
\]
is a $C^\infty$ function of $x, y$.
\end{x}


%% file: __refs.tex
\centerline{\textbf{\large REFERENCES}}
\setcounter{page}{1}
\setcounter{theoremn}{0}
\renewcommand{\thepage}{References-\arabic{page}}
\vspace{0.75cm}

\[
\text{ARTICLES}
\]

\begin{rf}
A. Borel and Harish-Chandra, Arithmetic Subgroups of algebraic groups, Ann. of Math. 75 (1962), 485-535.
\end{rf}

\begin{rf}
A. Selberg, Discontinuous groups and harmonic analysis, Proc. Internat. Congr. of Math., Stockholm (1962) 177-189
\end{rf}

\begin{rf}
T. Tamagawa, On Selberg’s trace formula, J. Fac. Sci. Univ. Tokyo Sect. IA Math. 8 (1960) 363-386.
\end{rf}

\begin{rf}
G. Warner, Selberg’s trace formula for non-uniform lattices:  The R-rank I case, Adv. in Math. Studies 6 (1979), 1-142.
\end{rf}
\vspace{0.25cm}

\[
\text{BOOKS}
\]

\begin{rf}
Christian Berg, Jens Peter Reus Christiansen, Harmonic Analysis on Semigroups, Springer-Verlag, 1984.
\end{rf}

\begin{rf}
A. Borel, Introduction aux groupes arithmetiques, Publlcations de l’Institut Mathematique de l’Universite de Strasbourg XV, Hermann, Paris, 1969.
\end{rf}

\begin{rf}
Z. I. Borevich and I. R. Shafarevich, Number Theory, Academic Press, New York, 1966.
\end{rf}

\begin{rf}
Glen E. Bredon, Topology and Geometry, Springer-Verlag, New York, 1993.
\end{rf}

\begin{rf}
K. Chandrasekharan, Arithmetical Functions, Springer-Verlag, New York, 1970.
\end{rf}

\begin{rf}
George Csordas, Linear Operators, Fourier Transforms and the Riemann (capital xi)-Function, In:  Some Topics on Value Distribution and Differentiability in Complex and p-Adic Analysis, Beijing Science Press (2008), pp. 188-218.
\end{rf}

\begin{rf}
J. Dieudonne, Foundations of Modern Analysis, Academic Press, New York, 1960.
\end{rf}

\begin{rf}
Harish-Chandra, Automorphic Forms on Semi-Simple Lie Groups, Lecture Notes in Mathematics Vol. 62, Springer-Verlag, Berlin and New York, 1976.
\end{rf}

\begin{rf}
Edwin Hewitt and Kenneth A, Ross, Abstract Harmonic Analysis I, Structure of Topological Groups Integration Theory Group Representations, Springer-Verlag, New York, 1963.
\end{rf}

\begin{rf}
J. Humphreys, Arithmetic Groups, Lecture Notes in Mathematics Vol. 789 Springer-Verlag, Berlin and New York, 1980.
\end{rf}

\begin{rf}
R. P. Langlands, On the Functional Equations Satisfied by Eisenstein Series, Lecture Notes in Mathematics Vol. 544, Springer-Verlag, Berlin and New York, 1976.
\end{rf}

\begin{rf}
Lynn H, Loomis, An Introduction to Abstract Harmonic Analysis, D. Van Nostrand, Princeton, New Jersey, 1953.
\end{rf}

\begin{rf}
Patrick Morandi, Field and Galois Theory, Springer-Verlag, 1996.
\end{rf}

\begin{rf}
Leopoldo Nachbin, The Haar Integral, D. Van Nostrand, Princeton, New Jersey, 1965.
\end{rf}

\begin{rf}
. Scott Osborne and Garth Warner, The Theory of Eisenstein Systems, Academic Press, 1981.
\end{rf}

\begin{rf}
M. S. Raghunathan, Discrete Subgroups of Lie Groups, Springer-Verlag, Berlin and New York, 1972.
\end{rf}

\begin{rf}
V. S. Varadarajan, Harmonic Analysis on Real Reductive Groups, Lecture Notes in Mathematics Vol. 576, Springer-Verlag, Berlin and New York, 1977.
\end{rf}

\begin{rf}
Garth Warner, Harmonic Analysis on Semi-Simple Lie Groups I, Springer-Verlag, New York, 1972.
\end{rf}

\begin{rf}
Garth Warner, Harmonic analysis on Semi-Simple Lie Groups II, Springer-Verlag, New York, 1972.
\end{rf}

\begin{rf}
Andre Weil, Basic Number Theory, Springer-Verlag, New York, 1967.
\end{rf}